\documentclass[twoside,a4paper,titlepage]{report}
\usepackage{amsmath}
\usepackage{amssymb}
\usepackage{euscript}
\usepackage{xspace}
\usepackage{enumerate}
\usepackage{amsthm}
\usepackage[all]{xy}
\usepackage[german,american]{babel}
\newcommand{\blank}{{\llcorner\!\lrcorner}}
\newcommand{\Mp}[2][60]{${#2}$\hbox{-}\penalty#1\hspace*{0pt}}
\newcommand{\Mpn}[2][10000]{${#2}$\hbox{-}\penalty#1\hspace*{0pt}}

\newcommand{\Frechet}{Fr\'e\-chet\xspace}
\newcommand{\Grothendieck}{Gro\-then\-dieck\xspace}
\newcommand{\Hochschild}{Hoch\-schild\xspace}

\newcommand{\Cstar}{C\Mp{{}^\ast}}

\newcommand{\C}{{\mathbb{C}}}
\newcommand{\R}{{\mathbb{R}}}

\newcommand{\Z}{{\mathbb{Z}}}
\newcommand{\Ztwo}{{\mathbb{Z}_2}}
\newcommand{\N}{{\mathbb{N}}}
\newcommand{\Q}{{\mathbb{Q}}}

\newcommand{\CCINF}{\mathop{\mathrm{C_c^\infty}}}

\newcommand{\NBC}[1][]{\mathrm{C}^{#1}}
\newcommand{\CINF}{\mathrm{C}^\infty}
\newcommand{\ABC}{\mathrm{AC}}
\newcommand{\HO}{\mathcal{O}}

\newcommand{\Sch}[1][1]{\ell^{#1}}

\newcommand{\Mat}[1][n]{\mathbb{M}_{\,#1}}        
\newcommand{\Hils}[1][H]{\EuScript{#1}}         

\newcommand{\VS}[1][V]{{\mathsf{#1}}}
\newcommand{\vs}[1][v]{{\mathsf{#1}}}
\newcommand{\MA}[1][A]{{\mathsf{#1}}}
\newcommand{\ma}[1][a]{{\mathsf{#1}}}

\newcommand{\GR}{{\Gamma}}            

\newcommand{\Unse}[1]{#1^{+}}                   

\newcommand{\ket}[1]{\mathopen|{#1}\mathclose\rangle}
\newcommand{\bra}[1]{\mathopen\langle{#1}\mathclose|}
\newcommand{\5}[2]{\langle{{#1},{#2}}\rangle}   

\newcommand{\cl}[1]{{\overline{#1}}}            
\newcommand{\bipol}[1]{#1^{\diamondsuit}}       
\newcommand{\coco}[1]{{#1}^{\scriptscriptstyle \heartsuit}}
\newcommand{\lin}{\mathop\mathrm{lin}}          
\newcommand{\cllin}{\mathop{\cl{\mathrm{lin}}}} 

\newcommand{\Ker}{\mathop\mathrm{Ker}}          
\newcommand{\Coker}{\mathop\mathrm{Coker}}      

\newcommand{\tr}{\mathop{\mathrm{tr}}\nolimits}

\newcommand{\ID}[1][{}]{\mathrm{id}_{{#1}}}     

\newcommand{\const}{\mathrm{const}}

\newcommand{\indlim}{\mathop{\lim\limits_{\longrightarrow}}\nolimits}
\newcommand{\prolim}{\mathop{\lim\limits_{\longleftarrow}}\nolimits}
\newcommand{\sep}{\mathop{\mathrm{sep}}}

\theoremstyle{plain}
\newtheorem{theorem}{Theorem}[chapter]
\newtheorem{proposition}[theorem]{Proposition}
\newtheorem{corollary}[theorem]{Corollary}
\newtheorem{lemma}[theorem]{Lemma}
\newtheorem{deflemma}[theorem]{Definition and Lemma}
\newtheorem{axiom}[theorem]{Axiom}
\newtheorem{claim}[theorem]{Claim}

\theoremstyle{definition}
\newtheorem{definition}[theorem]{Definition}

\theoremstyle{remark}
\newtheorem{example}[theorem]{Example}
\newtheorem{warning}[theorem]{Warning}

\theoremstyle{plain}

\newcommand{\hot}{\mathbin{\hat{\otimes}}}

\newcommand{\prot}{\mathbin{\hat{\otimes}_\pi}}
\newcommand{\protu}{\mathbin{\otimes_\pi}}
\newcommand{\indt}{\mathbin{\hat{\otimes}_\iota}}

\newcommand{\HH}{\operatorname{HH}}
\newcommand{\HC}{\operatorname{HC}}
\newcommand{\HD}{\operatorname{HD}}
\newcommand{\HE}{\operatorname{HE}}

\newcommand{\HP}{\operatorname{HP}}
\newcommand{\HA}{\operatorname{HA}}

\newcommand{\KK}{\operatorname{KK}}

\newcommand{\opp}{{\mathrm{opp}}}

\newcommand{\an}{{\mathrm{an}}}

\newcommand{\Null}{\operatorname{null}}

\newcommand{\chern}{\operatorname{ch}}

\newcommand{\rank}{\operatorname{rank}}

\newcommand{\Cpl}[1]{#1^{\mathrm{c}}}

\newcommand{\CBS}{\mathfrak{S}}

\newcommand{\BOUND}{\mathfrak{Bound}}
\newcommand{\COMP}{\mathfrak{Comp}}
\newcommand{\FINE}{\mathfrak{Fine}}
\newcommand{\EQUI}{\mathfrak{Equi}}

\newcommand{\Lin}{\mathcal{B}}                  
\newcommand{\Fin}{\mathcal{F}}                  
\newcommand{\Comp}{\mathcal{K}}                 
\newcommand{\Endo}{\mathcal{B}}                 

\newcommand{\Ad}{\mathop\mathrm{Ad}}            

\newcommand{\even}{\mathrm{even}}
\newcommand{\odd}{\mathrm{odd}}

\newcommand{\Cat}{\mathfrak{C}}                 
\newcommand{\Mor}{\operatorname{Mor}}           
\newcommand{\PRO}{\operatorname{pro}}           

\newcommand{\dissect}{\mathop{\mathfrak{dis}}}

\newcommand{\defeq}{\mathrel{:=}}

\newcommand{\Tanil}{\mathcal{T}}
\newcommand{\Tpnil}{\mathcal{T}_{\mathrm{p}}}
\newcommand{\Janil}{\mathcal{J}}
\newcommand{\Jpnil}{\mathcal{J}_{\mathrm{p}}}
\newcommand{\Tens}{\mathrm{T}}
\newcommand{\Jens}{\mathrm{J}}

\newcommand{\LL}{\mathfrak{L}}
\newcommand{\II}{\mathfrak{I}}

\newcommand{\opt}[1]{\langle{#1}\rangle}

\newcommand{\congto}{\mathrel{\overset{\cong}{\longrightarrow}}}
\newcommand{\prto}{\twoheadrightarrow}
\newcommand{\injto}{\rightarrowtail}

\newcommand{\ev}{\mathrm{ev}}

\newcommand{\Barn}{\mathrm{Bar}}

\newcommand{\lanilcur}{la\-nil\-cur\xspace}
\newcommand{\lanilcurs}{la\-nil\-curs\xspace}
\newcommand{\lonilcur}{lo\-nil\-cur\xspace}
\newcommand{\lonilcurs}{lo\-nil\-curs\xspace}

\newcommand{\LLH}[1]{{[\![ #1 ]\!]}}
\newenvironment{digression}{\begin{quotation}\small\sffamily}{\end{quotation}}

\newcommand{\Xtower}{\EuScript{X}}
\newcommand{\Xmap}{\xi}

\begin{document}

\begin{titlepage}
  \begin{center}
    \Huge\hspace*{1pt}\\[6cm]
    Ralf Meyer\\[1cm]
    Analytic cyclic cohomology\\[1cm]
    1999
  \end{center}
\end{titlepage}

\tableofcontents

\chapter{Introduction}
\label{sec:introduction}

The main results of this thesis are a proof of excision for entire cyclic
cohomology and the construction of a Chern-Connes character for Fredholm
modules without any summability conditions having values in a variant of
entire cyclic cohomology.

The excision theorem in entire cyclic cohomology asserts that if $\MA[K]
\injto \MA[E] \prto \MA[Q]$ is an algebra extension with a bounded linear
section $s \colon \MA[Q]\to \MA[E]$, then there is a six-term exact sequence
$$
\xymatrix{
  {\HE^0(\MA[Q])} \ar[r]^{i^\ast} &
    {\HE^0(\MA[E])} \ar[r]^{p^\ast} &
      {\HE^0(\MA[K])} \ar[d] \\
  {\HE^1(\MA[K])} \ar[u] &
    {\HE^1(\MA[E])} \ar[l]^{p^\ast} &
      {\HE^1(\MA[Q])} \ar[l]^{i^\ast}
  }
$$
of the associated entire cyclic cohomology groups $\HE^\ast(\blank)$.
Recently, this has been proved also by Puschnigg~\cite{puschnigg98:excision}
using quite different methods.

Entire cyclic cohomology is a relative of periodic cyclic cohomology that
accommodates also certain ``infinite dimensional'' cohomology contributions.
This can be made more precise in terms of the Chern-Connes character in
\Mpn{K}homology.  In periodic cyclic cohomology, the Chern-Connes character
can be defined only for \emph{finitely summable} Fredholm modules.  That is,
the commutators $[x,F]$ that are required to be compact for a Fredholm module
are even contained in a Schatten ideal $\Sch[p](\Hils)$ for some~$p$.
Unfortunately, there are many interesting Fredholm modules that are not
finitely summable.  In entire cyclic cohomology, the summability condition can
be weakened to \emph{\Mpn{\theta}summability} \cite{connes88:entire}.  This
allows certain ``infinite'' Fredholm modules, but still imposes a serious
restriction.  In particular, Fredholm modules over \Cstar{}algebras usually
are not \Mpn{\theta}summable.

However, a small modification of the definition of entire cyclic cohomology
allows to define a Chern-Connes character for Fredholm modules without any
summability restriction.  For entire cyclic cohomology we consider families of
multi-linear maps $\phi_n \colon \MA^{\otimes n} \to \C$, $n \in \N$, such
that for all \emph{bounded} sets $S \subset \MA$, there is a constant $C(S)$
such that $|\phi_n(\ma_1, \dots, \ma_n)| \le C$ for all $n \in \N$, $\ma_1,
\dots, \ma_n \in S$.  If we instead require this ``entire growth'' condition
only for \emph{compact} sets~$S$, we get more linear maps.  It turns out that
we get sufficiently many cochains to write down a Chern-Connes character for
Fredholm modules without imposing any summability condition
(Section~\ref{sec:Chern_character}).  The resulting character coincides with
the usual character for finitely summable Fredholm modules and probably also
for \Mpn{\theta}summable Fredholm modules.  However, I have not checked the
latter.

If we impose the entire growth condition only on \emph{finite} subsets instead
of bounded or compact subsets, we obtain a theory that is considered already
by Connes~\cite{connes94:ncg}.

So far we have not specified for which algebras entire cyclic cohomology is
defined.  Usually, entire cyclic cohomology is defined on locally convex
topological algebras.  However, only the bounded subsets are used in its
definition to formulate the growth constraint mentioned above.  Moreover, we
may use the compact or the finite subsets instead of the bounded subsets.
Therefore, a more natural domain of definition is the category of
\emph{complete bornological algebras}.  A \emph{bornological vector space} is
a vector space~$\VS$ together with a collection~$\CBS$ of subsets, satisfying
certain axioms already formulated by Bourbaki~\cite{bourbaki:TVS}.  A standard
example is the collection of all bounded subsets of a locally convex
topological vector space.  We call sets $S \in \CBS$ \emph{small}.  It appears
that bornological vector spaces have not yet been used by people studying
cyclic cohomology theories.  Hence we outline some basic analysis in
Chapter~\ref{cha:bornologies}.  We do not need much analysis, because we
really want to do \emph{algebra}.  Analysis mainly has to set the stage for
the algebra to go through smoothly.

It is impossible to understand entire cyclic cohomology using only topological
vector spaces.  Unless~$A$ is a Banach algebra, it appears to be impossible to
describe the semi-norms of a locally convex topology on $\sum A^{\hot n}$ such
that the continuous linear functionals are precisely the families of
multi-linear maps satisfying the entire growth condition.  In adittion, for
the purposes of the Chern-Connes character it is reasonable to endow a
\Cstar{}algebra with the bornology generated by the compact subsets instead of
the bornology of all bounded subsets.  This change of bornology cannot be
described as choosing a different topology.  The only reasonable alternative
to bornological vector spaces are inductive systems.  These are used quite
successfully by Puschnigg in~\cite{puschnigg98:cyclic} and are necessary to
handle \emph{local} theories.  However, for the purposes of entire cyclic
cohomology it suffices to work with bornological vector spaces.  We will see
in Appendix~\ref{app:bornologies_inductive} that inductive systems and
bornological vector spaces are closely related.

The original definition of periodic cyclic cohomology by Connes suggests the
following route to compute the periodic cyclic cohomology of an algebra~$\MA$.
Firstly, compute the \emph{\Hochschild cohomology} $\HH^\ast(\MA)$.  Secondly,
use Connes's long exact ``SBI'' sequence to obtain the \emph{cyclic
  cohomology} $\HC^\ast(\MA)$ from $\HH^\ast (\MA)$.  Finally, the
\emph{periodic cyclic cohomology} $\HP^\ast (\MA)$ is the inductive limit
\begin{equation}  \label{eq:HP_indlim}
  \HP^\ast(\MA) = \indlim \HC^{\ast + 2k}(\MA).
\end{equation}

This powerful strategy works, for example, for smooth manifolds and group
algebras.  However, it fails miserably if we cannot compute the \Hochschild
cohomology of~$\MA$.  This happens, for instance, for the Schatten ideals
$\Sch[p](\Hils)$ with $p \gg 1$.  Nevertheless, Cuntz~\cite{cuntz97:bivariant}
was able to compute $\HP^0(\Sch[p]) = \C$ and $\HP^1(\Sch[p]) = 0$ for all~$p$
without knowing anything about the \Hochschild or cyclic cohomology of
$\Sch[p](\Hils)$.  This is possible because periodic cyclic cohomology has
better homological properties than \Hochschild or cyclic cohomology.  It is
invariant under smooth homotopies, stable, and satisfies excision.  For many
algebras, we can compute periodic cyclic cohomology using long exact sequences
that follow from these properties.  For example, there are six term exact
sequences computing the periodic cyclic cohomology of crossed products by $\Z$
and~$\R$.  This is why the result of Cuntz and Quillen~\cite{cuntz97:excision}
that periodic cyclic cohomology satisfies excision is so important.

Connes's original definition of entire cyclic cohomology is elementary but
gives us no hint how to compute it.  Connes takes the complex computing
periodic cyclic cohomology and enlarges it by allowing cochains with infinite
support that satisfy the entire growth condition.  However, entire cyclic
cohomology is no longer related to \Hochschild cohomology.  As a result, we do
not know how to cut the huge complex defining entire cyclic cohomology to a
manageable size.  It is already a hard problem to compute the entire cyclic
cohomology of the algebra $\CINF(S^1)$ of smooth functions on the circle.  The
obvious conjecture that $\HE^\ast\bigl( \CINF(S^1) \bigr)$ is equal to the de
Rham homology of the circle, $\C$ in even and odd degree, has been verified
only recently by Puschnigg using excision.  At the moment, excision is the
most powerful method to compute entire cyclic cohomology groups.

We use a more complicated but also more explanatory definition of entire
cyclic cohomology and call it \emph{analytic cyclic cohomology} instead.  This
avoids confusion in connection with the Chern-Connes character, because its
range is not the entire cyclic cohomology that is usually considered.  The
definition of analytic cyclic cohomology is analogous to the description of
periodic cyclic cohomology given by Cuntz and Quillen in \cite{cuntz95:cyclic}
and~\cite{cuntz97:excision}.

The \emph{tensor algebra} $\Tens\MA$ of~$\MA$ can be realized as the even part
of the algebra of non-commutative differential forms $\Omega\MA$, endowed with
the \emph{Fedosov product}.  The latter is a deformation of the ordinary
product of differential forms (see Appendix~\ref{app:Omega_Tens}).  We endow
$\Omega\MA$ with the convex bornology $\CBS_\an$ generated by sets of the form
$\{ \opt{\ma_0} d\ma_1 \dots d\ma_n \mid n \in \N, \ma_0, \dots, \ma_n \in
S\}$ with small $S \subset \MA$.  The notation $\opt{\ma_0} d\ma_1 \dots
d\ma_n$ stands for $\ma_0 d\ma_1 \dots d\ma_n$ or $d\ma_1 \dots d\ma_n$.  The
bornology $\CBS_\an$ is quite obvious from the entire growth condition.  The
completion of $\Tens \MA \subset \Omega \MA$ with respect to the bornology
$\CBS_\an$ is the \emph{analytic tensor algebra} $\Tanil\MA$.

The \emph{X-complex} $X(\MA)$ of an algebra~$\MA$ is a \Mpn{\Ztwo}graded
complex with
$$
X_0(\MA) \defeq \MA, \qquad
X_1(\MA) \defeq \Omega^1\MA / [\Omega^1\MA, \MA] \defeq
\Omega^1\MA / b(\Omega^2\MA).
$$
The boundary maps in the X-complex are induced by $d \colon \MA \to
\Omega^1\MA$ and $b \colon \Omega^1\MA \to \MA$.  Thus we obtain a
low-dimensional quotient of the \Mpn{(B,b)}bicomplex.  For a general
algebra~$\MA$, this complex is not particularly interesting.  However, it is
good enough if~$\MA$ is \emph{quasi-free}, that is, homologically
\Mpn{1}dimensional.  If~$\MA$ is quasi-free, then $X(\MA)$ computes the
periodic cyclic cohomology of~$\MA$.

The \emph{analytic cyclic cohomology} $\HA^\ast(\MA)$, $\ast = 0,1$, of a
complete bornological algebra~$\MA$ is defined as the cohomology of the
complex $X(\Tanil\MA)$.  Furthermore, we define a corresponding homology
theory $\HA_\ast(\MA)$ and a bivariant theory $\HA^\ast(\MA; \MA[B])$ by
considering the complex of bounded linear maps $X(\Tanil\MA) \to
X(\Tanil\MA[B])$.  If~$\MA$ is a locally convex algebra endowed with the
bornology of bounded subsets, then $\HA^\ast(\MA)$ is naturally isomorphic to
the entire cyclic cohomology $\HE^\ast(\MA)$.

This definition of analytic cyclic cohomology already suggests the right tools
for its study.  The purely algebraic tensor algebra $\Tens\MA$ is universal
for bounded linear maps $\MA \to \MA[B]$ into complete bornological algebras
in the sense that any such map can be extended uniquely to a bounded
homomorphism $\Tens\MA \to \MA[B]$.  In Cuntz's treatment of bivariant
\Mpn{K}theory~\cite{cuntz97:bivariant}, tensor algebras play a prominent role.
The identity map induces a natural projection $\Tens\MA \to \MA$.  Let
$\Jens\MA$ be the kernel of this map.  The extension $\Jens\MA \injto \Tens\MA
\prto \MA$ is universal among all extensions of~$\MA$ with a bounded linear
section.  This universal property characterizes $\Tens\MA$ uniquely up to
smooth homotopy.

The analytic tensor algebra $\Tanil\MA$ is no longer universal for arbitrary
linear maps, but we can characterize the class of linear maps that can be
extended to bounded homomorphisms $\Tanil\MA \to \MA[B]$.  We call these maps
\emph{\lanilcurs{}}.  The extension $\Janil\MA \injto \Tanil\MA \prto \MA$ is
no longer universal among all extensions of~$\MA$, but it enjoys an analogous
universal property among those extensions with \emph{analytically nilpotent
  kernel}.  A complete bornological algebra~$\MA[N]$ is called
\emph{analytically nilpotent} iff for all small subsets $S \subset \MA[N]$,
the set $\bigcup_{n \in \N} S^n$ is small.  The universal property of the
extension $\Janil\MA \injto \Tanil\MA \prto \MA$ determines $\Tanil\MA$
uniquely up to smooth homotopy.  This is an important ingredient in the proof
of the Excision Theorem.  The algebra $\Tanil\MA$ has the property that all
(allowable) extensions $\MA[N] \injto \MA[E] \prto \Tanil\MA$ with
analytically nilpotent kernel split by a bounded homomorphism.  This can be
viewed as a condition of non-singularity because of its similarity to a
condition for smooth algebraic varieties.  Since the algebra $\Janil\MA$ is
analytically nilpotent, we may consider $\Tanil\MA$ as an ``non-singular
infinitesimal thickening'' of~$\MA$.  Consequently, the process of
replacing~$\MA$ by its analytic tensor algebra $\Tanil\MA$ can be viewed as a
``non-commutative resolution of singularities by an infinitesimal
thickening''.

These are the tools to handle the analytic tensor algebra.  The X-complex is
susceptible to homological algebra.  It is closely related to the complex
$\Omega^1\MA / b(\Omega^2\MA) \overset{b}{\longrightarrow} \MA$ that computes
the \Hochschild cohomology groups $\HH^0(\MA)$ and $\HH^1(\MA)$.  We write
$X_\beta(\MA)$ for this complex.

The even cohomology of $X(\MA)$ can be described as the space of smooth
homotopy classes of traces on~$\MA$.  That is, $H^0 \bigl( X(\MA) \bigr)$ is
obtained by making the functor $\HH^0(\MA)$ that associates to~$\MA$ the space
of traces on~$\MA$ invariant under smooth homotopies.  Thus we can view
elements of $\HA^0(\MA)$ as homotopy classes of traces on $\Tanil\MA$.  The
universality of the extension $\Janil\MA \injto \Tanil\MA \prto \MA$ implies
that any trace on an analytically nilpotent extension of~$\MA$ can be pulled
back to a trace on $\Tanil\MA$ whose homotopy class is uniquely determined.
Thus we can view $\HA^0(\MA)$ as the space of homotopy classes of traces on
analytically nilpotent extensions of~$\MA$.

In Section~\ref{sec:X_Tanil_HA}, we derive some basic results about analytic
cyclic cohomology.  We prove homotopy invariance for absolutely continuous
homotopies and stability with respect to generalized algebras of trace class
operators.  We construct the Chern-Connes character from \Mpn{K}theory to
analytic cyclic homology.  These results are quite straightforward using the
techniques developed by Cuntz and Quillen~\cite{cuntz95:cyclic}.

The Excision Theorem is proved using the tools mentioned above.  See
Section~\ref{sec:excision_outline} for an outline of the proof.  One part of
the proof is homological algebra: We have to verify that certain bimodules
form a projective resolution and compute the associated commutator quotient
complex.  The other part is universal algebra: We have to show that a certain
algebra~$\LL$ is \emph{analytically quasi-free}, that is, the natural
projection $\Tanil\LL \to \LL$ has a splitting homomorphism.  Such a splitting
homomorphism is constructed explicitly using the universal property of
analytic tensor algebras.

In addition, we apply the techniques used to study analytic cyclic cohomology
to periodic cyclic cohomology.  The formalism of non-commutative resolution of
singularities is available as well.  We only have to use different notions of
nilpotence and a different version of the tensor algebra.  Actually, the right
tensor algebra for periodic cyclic cohomology is not an algebra but a
projective system of algebras.  This is observed in~\cite{cuntz97:excision}.
Therefore, to handle periodic cyclic cohomology properly, we study it on the
category of ``pro-algebras'' (Section~\ref{sec:pro_algebras}).  All the
desirable homological properties, in particular excision, continue to hold in
that generality.

The connecting map in $\HP^\ast$ for an extension $\MA[K] \injto \MA[E] \prto
\MA[Q]$ maps $\HP^{\ast}(\MA[K]) \to \HP^{\ast+1}(\MA[Q])$.
By~\eqref{eq:HP_indlim} we expect it to map $\HC^k(\MA[K])$ into
$\HC^{f(k)}(\MA[Q])$ for a suitable function~$f$.
Puschnigg~\cite{puschnigg99:excision} shows that we can achieve $f(k) = 3k+3$
and that this estimate is optimal for a large class of extensions, including
the universal extensions $\Jens\MA \injto \Tens\MA \prto \MA$.  Our proof of
excision yields the same estimate and, in addition, more precise results
concerning relative cyclic cohomology.  An equivalent formulation of the
excision theorem asserts that $\HP^\ast (\MA[K]) \cong \HP^\ast (\MA[E] :
\MA[Q])$.  We obtain estimates for the dimension shift that occurs when
switching between $\HC^\ast (\MA[K])$ and $\HC^\ast (\MA[E] : \MA[Q])$
(Theorem~\ref{the:excision_cyclic}).

Appendix~\ref{app:bornologies} contains some proofs omitted in
Chapter~\ref{cha:bornologies}.  The remaining Appendices mostly contain
well-known algebra, carried over to complete bornological vector spaces.  In
particular, we prove the homotopy invariance of the X-complex for quasi-free
algebras.

\chapter{Bornologies}
\label{cha:bornologies}

This preparatory section contains a brief review of some important concepts of
bornological functional analysis.  We define complete bornological vector
spaces, bounded linear maps, convergence of sequences in bornological vector
spaces, completions, and completed bornological tensor products and illustrate
these concepts by examples.  Of great interest are \Frechet spaces endowed
with the precompact or the bounded bornology and vector spaces without
additional structure endowed with the finest possible bornology.  Often
difficult theorems of analysis are required to prove the assertions in the
examples.  Since the examples mainly serve to illustrate the definitions, we
usually omit proofs.

There is a close relationship between complete bornological vector spaces and
inductive systems of Banach spaces.  Put in a nutshell, the category of
complete bornological vector spaces is equivalent to a full subcategory of the
category of inductive systems of Banach spaces, namely the subcategory of
inductive systems with injective structure maps.  This relationship is quite
fundamental and used heavily throughout this thesis because it is needed to
prove the existence of completions.  In addition, it is necessary to work with
inductive systems to understand the \emph{local cyclic cohomology} of
Puschnigg~\cite{puschnigg98:cyclic}.  However, as long as we deal only with
non-local theories like entire and periodic cyclic cohomology it suffices to
work with bornological vector spaces.  Therefore, we avoid direct use of
inductive systems as much as possible and relegate this topic to
Appendix~\ref{app:bornologies}.

I am not sufficiently familiar with bornological analysis to comment
accurately on its history.  Let me just point out the following.  The basic
definitions can be found in the present form in Bourbaki~\cite{bourbaki:TVS}.
Completions and tensor products of bornological vector spaces are studied by
Henri Hogbe-Nlend \cite{hogbe-nlend70:completions}.  Important contributions
are due to Lucien Waelbroeck.  Hogbe-Nlend has written good books about the
general theory of bornological vector spaces, namely
\cite{hogbe-nlend71:theorie}, \cite{hogbe-nlend73:techniques}, and
\cite{Hogbe-Nlend:Born}.  The textbook~\cite{Hogbe-Nlend:Born} is written at
an elementary level.  We will usually adopt the terminology
of~\cite{Hogbe-Nlend:Born} with a few exceptions.

It is important to avoid a prejudice in favor of the bounded bornology on a
topological vector space.  The bounded bornology is usually not so
well-behaved.  This is good for analysts.  They can single out the good spaces
from the bad ones by requiring the bounded bornology to have some desirable
property.  Algebraists, however, want all objects to be well-behaved.  By
choosing a smaller bornology, we can often improve the behavior of a space.
For instance, the bounded bornology on a projective tensor product $V \prot W$
of \Frechet spaces may have nothing to do with the bounded bornologies on the
factors $V$ and~$W$.  This follows from the negative answer to Grothendieck's
\emph{Probl\`eme des Topologies} \cite{grothendieck55:produits}.  However,
this problem completely disappears if we work with the precompact bornology.
Moreover, the analytic cyclic cohomology of \Cstar{}algebras gets interesting
once we endow them with the precompact bornology.

\section{Basic Definitions}
\label{sec:basic_def}

All vector spaces shall be over the complex numbers~$\C$.  Let~$\VS$ be a
vector space.  A subset $S\subset\VS$ is a \emph{disk} iff it is circled and
convex.  The \emph{disked hull}~$\bipol{S}$ is the circled convex hull of~$S$.
If $S\subset\VS$ is a disk, define the semi-normed space~$\VS_S$ to be the
linear span of~$S$ endowed with the semi-norm $\|\blank\|_S$ whose unit ball
is~$S$.  The disk~$S$ is called \emph{norming} iff~$\VS_S$ is a normed space
and \emph{completant} iff~$\VS_S$ is a Banach space.

\subsection{Bornologies}
\label{sec:def_bornologies}

Let~$\CBS$ be a collection of subsets of~$\VS$.  We call~$\CBS$ a \emph{convex
(vector) bornology} and the pair $(\VS,\CBS)$ a \emph{convex bornological
vector space} iff the following conditions hold:
\begin{enumerate}[(1)]
  \item $\{\vs\} \in \CBS$ for all $\vs \in \VS$;
  \item if $S \in \CBS$ and $T \subset S$, then $T \in \CBS$;
  \item if $S_1, S_2 \in\CBS$, then $S_1 + S_2 \in \CBS$;
  \item if $S \in \CBS$, then $\bipol{S} \in \CBS$.
\end{enumerate}
We call $S\subset \VS$ \emph{small} iff $S \in \CBS$, avoiding the more common
term ``bounded''.  The conditions above imply that $S_1 \cup S_2 \in \CBS$ if
$S_1, S_2 \in \CBS$ and $c \cdot S \in \CBS$ for all $c \in \C$, $S \in \CBS$.

The bornology~$\CBS$ is \emph{separating} iff all small disks $S \in \CBS$ are
norming.  If~$\CBS$ is separating, then we call $(\VS,\CBS)$ a \emph{separated
convex bornological vector space}.  The bornology~$\CBS$ is \emph{completant}
iff each $S\in\CBS$ is contained in a \emph{completant} small disk
$S'\in\CBS$.  If~$\CBS$ is completant, then we call~$(\VS,\CBS)$ a
\emph{complete bornological vector space}.  A completant bornology is
necessarily separating.

In the following, we usually write just~$\VS$ for a bornological vector space
and~$\CBS(\VS)$ for the bornology on~$\VS$.  We let $\CBS_d(\VS)$ be the
family of all small \emph{disks} in~$\VS$.  If~$\VS$ is complete,
$\CBS_c(\VS)$ denotes the family of all \emph{completant} small disks
in~$\VS$.  If $S_1,S_2$ are (completant) small disks, then $S_1 + S_2$ is a
(completant) small disk as well.  Thus $\CBS_d(\VS)$ and $\CBS_c(\VS)$ are
directed sets with inclusion as order relation.  Moreover, $\CBS_d(\VS)$ and
$\CBS_c(\VS)$ are closed under arbitrary intersections.  Hence we can
well-define the \emph{completant disked hull}~$\coco{S}$ of $S \subset \VS$ as
the smallest completant small disk containing~$S$.  Of course, the
set~$\coco{S}$ may fail to exist if~$\VS$ is incomplete.

\begin{example}
  Let~$\VS$ be a locally convex topological vector space.  Let $\BOUND =
  \BOUND(\VS)$ be the collection of all \emph{bounded} subsets of~$\VS$.  That
  is, $S \in \BOUND(\VS)$ iff~$S$ is absorbed by all neighbourhoods of the
  origin.  The bounded sets always form a convex vector bornology on~$\VS$
  called the \emph{bounded bornology} (it is called ``von Neumann bornology''
  in~\cite{Hogbe-Nlend:Born}).  The bornology $\BOUND(\VS)$ is separating if
  the topology of~$\VS$ is Hausdorff and completant if~$\VS$ is a complete
  topological vector space.  Actually, it suffices to assume that~$\VS$ is
  ``sequentially complete''.
  
  The case of (semi)normed spaces is particularly important.  Let~$\VS$ be a
  (semi)normed space with unit ball~$B$.  Then $S \subset \VS$ is bounded iff
  $S \subset c \cdot B$ for some constant $c>0$.  We call a (semi)normed space
  with the bounded bornology \emph{primitive}.  It is explained in
  Appendix~\ref{app:bornologies} how general convex bornological vector spaces
  are pieced together out of primitive spaces.
\end{example}

\begin{example}
  Let~$\VS$ be a locally convex topological vector space.  Let $\COMP =
  \COMP(\VS)$ be the collection of all \emph{precompact} subsets of~$\VS$.
  That is, $S \in \COMP(\VS)$ iff for all neighbourhoods of the origin~$U$
  there is a \emph{finite} set $F\subset \VS$ such that $S \subset F + U$.
  If~$\VS$ is complete, then $S\subset \VS$ is precompact iff its closure is
  compact.  In general, $S\subset \VS$ is precompact iff its closure taken in
  the completion of~$\VS$ is compact.  The precompact subsets $\COMP(\VS)$
  always form a convex vector bornology on~$\VS$, called the \emph{precompact
  bornology}.  The bornology $\COMP(\VS)$ is separating if the topology
  of~$\VS$ is Hausdorff and completant if~$\VS$ is a complete topological
  vector space.
\end{example}

\begin{example}
  Let~$\VS$ be just a vector space.  The \emph{fine bornology} $\FINE(\VS)$
  on~$\VS$ is the smallest possible vector bornology on~$\VS$.  That is,
  $S\in\FINE(\VS)$ if and only if there are finitely many points
  $\vs_1,\dots,\vs_n \in \VS$ such that~$S$ is contained in the disked hull of
  $\{\vs_1,\dots,\vs_n\}$.  The fine bornology is always completant because
  finite dimensional topological vector spaces are complete.
\end{example}

\begin{example}
  The analytic tensor algebra and many other universal objects are constructed
  as follows.  We start with a vector space~$\VS$ and a collection~$\CBS$ of
  subsets not satisfying the axioms for a convex bornology.  The \emph{convex
  bornology generated by~$\CBS$} is defined as the smallest collection~$\CBS'$
  of subsets of~$\VS$ that verifies the axioms for a convex vector bornology
  and contains~$\CBS$.

  This bornology always exists.  It may fail to be separating and will usually
  not be complete.  This can be helped in a second step by taking the
  completion of $(\VS, \CBS')$.
\end{example}

\subsection{Bounded maps}
\label{sec:bounded_maps}

A linear map $l\colon \VS_1 \to \VS_2$ between bornological vector spaces is
\emph{bounded} iff~$l$ maps small sets to small sets, that is, $l(S) \in
\CBS(\VS_2)$ for all $S \in \CBS(\VS_1)$.  A set~$L$ of linear maps $\VS_1 \to
\VS_2$ is \emph{equibounded} iff $L(S) \defeq \{ l(x) \mid l \in L, x \in S\}$
is small in~$\VS_2$ for all $S \in \CBS(\VS_1)$.  A \emph{bornological
isomorphism} is a bounded linear map $l \colon \VS_1 \to \VS_2$ with a bounded
(two-sided) inverse.

An \Mpn{n}linear map $l\colon \VS_1 \times \dots \times \VS_n \to \VS_{n+1}$
is \emph{bounded} iff $l(S_1 \times \dots \times S_n) \in \CBS(\VS_{n+1})$ for
all $S_j \in \CBS(\VS_j)$.  A set~$L$ of \Mpn{n}linear maps $\VS_1 \times
\dots \times \VS_n \to \VS_{n+1}$ is \emph{equibounded} iff $L(S_1 \times
\dots \times S_n) \in \CBS(\VS_{n+1})$ for all $S_j \in \CBS(\VS_j)$.

\begin{definition}
  A \emph{complete bornological algebra} is a complete bornological vector
  space~$\MA$ with a bounded, associative multiplication $\MA \times \MA \to
  \MA$.
\end{definition}

\begin{example}
  Let~$\VS_j$ be \Frechet spaces endowed with the bounded or the precompact
  bornology.  Then a continuous linear map $\VS_1 \to \VS_2$ is necessarily
  bounded and vice versa.  The reason is that for metrizable spaces,
  continuity is equivalent to ``sequential continuity''.  In other words, the
  functors $\BOUND$ and $\COMP$ from the category of \Frechet spaces to the
  category of complete bornological vector spaces are fully faithful.
  
  This carries over to multi-linear maps.  An \Mpn{n}linear map $\VS_1 \times
  \dots \times\VS_n \to \VS_{n+1}$ is bounded if and only if it is jointly
  continuous.  This is related to the fact that for \Frechet spaces,
  separately continuous multi-linear maps are automatically jointly
  continuous.  Hence a \Frechet algebra, endowed with either the bounded or
  the precompact bornology, is a complete bornological algebra.
\end{example}

\begin{example}  \label{exa:R_convolution}
  A locally convex topological vector space~$V$ is an \emph{LF-space} iff
  there is an increasing sequence $(V_n)_{n\in \N}$ of subspaces of~$V$ such
  that each~$V_n$ is a \Frechet space with respect to the subspace topology,
  $\bigcup V_n = V$, and~$V$ carries the finest locally convex topology making
  the inclusions $V_n \to V$ continuous for all $n\in\N$.  Thus a linear map
  $V \to W$ with locally convex range~$W$ is continuous iff its restriction
  to~$V_n$ is continuous for all $n \in \N$.  A bounded subset of an
  LF-space~$V$ is necessarily contained in a \Frechet subspace~$V_n$ for some
  $n \in \N$ (see~\cite{Treves:67}).
  
  Let~$\VS_j$ be LF-spaces.  Endow all~$\VS_j$ with the bounded or the
  precompact bornology.  Then an \Mpn{n}linear map $l\colon \VS_1 \times \dots
  \times \VS_n \to \VS_{n+1}$ is bounded if and only if it is separately
  continuous.  Indeed, boundedness easily implies separate continuity.
  Conversely, separate continuity implies joint continuity for \Frechet
  spaces.  Hence the restriction of~$l$ to all \Frechet subspaces of $\VS_1
  \times \dots \times \VS_n$ is jointly continuous and therefore bounded.
  Since a bounded subset of an LF-space is already contained in a \Frechet
  subspace, the map~$l$ is bounded.

  Thus bounded bilinear maps need not be jointly continuous.  We call
  LF-spaces with a separately continuous associative multiplication
  \emph{LF-algebras}.  Thus an LF-algebra is a complete bornological algebra
  with respect to either the bounded or the precompact bornology.
  
  Interesting LF-algebras are the convolution algebras of compactly supported
  smooth functions on smooth groupoids.  For instance, endow the LF-space
  $\CCINF(\R)$ of smooth compactly supported functions on the real line with
  the convolution product
  $$
  (f \ast g)(t) = \int_{\R} f(s) g(t-s) \,ds.
  $$
  It is easy to verify that this multiplication is separately continuous but
  not jointly continuous.  To prove the latter assertion, consider the
  continuous seminorm
  $$
  p(f) \defeq \sum_{n\in\N} |f^{(n)}(n)| \qquad \forall f \in \CCINF(\R).
  $$
  Assume that there are continuous seminorms $q_1,q_2$ on $\CCINF(\R)$ with
  $p(f\ast g) \le q_1(f)\cdot q_2(g)$.  Restrict~$q_1$ to functions with
  support in $[-1,1]$.  By continuity, this restricted seminorm $q_1(f)$ can
  only take into account suprema of a finite number of derivatives of~$f$.
  Hence, for sufficiently big~$n$, we have no control on the $n$th derivative
  of~$f$ near~$0$.  Since $(f\ast g)^{(n)} = f^{(n)} \ast g$, the modulus of
  $(f\ast g)^{(n)}(n)$ can become arbitrarily big while $q_1(f)\cdot q_2(g)$
  remains bounded.
\end{example}

\begin{example}
  Let~$\VS_j$ be fine spaces and let~$\VS[W]$ be an arbitrary convex
  bornological vector space.  Then any \Mpn{n}linear map $\VS_1 \times \dots
  \times \VS_n \to \VS[W]$ is bounded.  Consequently, the functor $\FINE$ that
  associates to a vector space~$\VS$ the bornological vector space
  $(\VS,\FINE)$ is a fully faithful functor from the category of vector spaces
  to the category of complete bornological vector spaces.  Any algebra becomes
  a complete bornological algebra when endowed with the fine bornology.
\end{example}

\begin{example}
  Assume that the bornology of~$\VS$ is the convex bornology generated by a
  collection~$\CBS$ of subsets of~$\VS$.  Then a linear map $l\colon \VS \to
  \VS[W]$ is bounded iff $l(S) \in \CBS(\VS[W])$ for all $S \in \CBS$.
\end{example}

\subsection{Bornological convergence}
\label{sec:born_convergence}

A sequence $(\vs_n)_{n \in \N}$ in a convex bornological vector space~$\VS$ is
\emph{(bornologically) convergent} towards $\vs_\infty \in \VS$ iff there is a
null-sequence of scalars $(\epsilon_n)_{n \in \N}$ and a small set $S \in
\CBS(\VS)$ such that $\vs_n - \vs_\infty \in \epsilon_n S$ for all $n \in \N$.
Equivalently, there is a small disk $S \in \CBS(\VS)$ such that $\vs_n \in
\VS_S$ for all $n \in \N \cup \{\infty\}$ and $\lim_{n \to \infty} \| \vs_n -
\vs_\infty\|_S = 0$.

A subset~$S$ of a bornological vector space~$\VS$ is \emph{(bornologically)
closed} iff it is sequentially closed for bornologically convergent sequences.
The bornologically closed sets satisfy the axioms for the closed sets of a
topology (though in general not a vector space topology).  Hence the
\emph{(bornological) closure} of a set is well-defined.  The closure of a
vector subspace is again a vector subspace.  Bounded maps respect convergence
of sequences and therefore are continuous with respect to the topology defined
by the bornologically closed sets.

\emph{Bornological Cauchy sequences} and \emph{absolutely convergent series}
are defined as Cauchy sequences or absolutely convergent series in~$\VS_S$ for
a suitable small disk~$S$.

\begin{example}
  In a \Frechet space (or an LF-space), a sequence~$(\vs_n)$ converges with
  respect to the topology if and only if it is bornologically convergent with
  respect to the bounded or precompact bornology.  Hence a subset of a
  \Frechet space is bornologically closed iff it is topologically closed.
\end{example}

\begin{example}
  In a fine space~$\VS$, a sequence~$(\vs_n)$ converges iff all~$\vs_n$ are in
  a fixed finite dimensional subspace and converge in that finite dimensional
  subspace.  Thus a subset $S \subset \VS$ is bornologically closed iff its
  intersection with all finite dimensional subspaces is closed in the usual
  sense.
\end{example}

\begin{warning}
  Let~$\VS$ be a complete bornological vector space.  A completant subset
  of~$\VS$ need not be closed.  Hence the \emph{closed disked hull} of a
  subset~$S$, defined as the smallest closed disked subset of~$\VS$
  containing~$S$, may be bigger than the \emph{completant} disked
  hull~$\coco{S}$.  We will not use closed disked hulls at all.
\end{warning}

\begin{warning}
  In most applications coming from topological vector spaces, the bornological
  closure of a small set is again small.  For example, the closure of a
  bounded or precompact set is bounded or precompact, respectively.  However,
  in general the closure of a small set need not be small.
\end{warning}

\section{Constructions with bornological vector spaces}
\label{sec:bornology_construct}

Most constructions that can be done with ordinary vector spaces can be carried
over to the category of (complete) bornological vector spaces.  There are
\emph{direct products}, \emph{direct sums}, \emph{projective limits}, and
\emph{inductive limits}.  These are defined by the well-known universal
properties.  We should not bother about these rather abstract matters before
we really need them.  Therefore, the discussion of these constructions is
relegated to Appendix~\ref{app:bornology_construct}.  See
also~\cite{Hogbe-Nlend:Born}.

\subsection{Subspaces, quotients, extensions}
\label{sec:extensions}

There are natural bornologies on \emph{subspaces} and \emph{quotients}.  Thus
a bounded linear map has a kernel and a cokernel.  Let~$\VS$ be a complete
bornological vector space and $\VS[W] \subset \VS$ a vector subspace.  A
subset $S \subset \VS[W]$ is small in the \emph{subspace bornology} iff $S \in
\CBS(\VS)$.  A subset $S \subset \VS/\VS[W]$ is small in the \emph{quotient
  bornology} iff $S = T \bmod \VS[W]$ for some $T \in \CBS(\VS)$.  The
subspace bornology on~$\VS[W]$ is completant iff~$\VS[W]$ is bornologically
closed in~$\VS$.  The quotient bornology on $\VS / \VS[W]$ is completant iff
it is separating iff~$\VS[W]$ is bornologically closed.

\begin{example}
  Let~$V$ be a topological vector space and $W\subset V$ a subspace, endowed
  with the subspace topology.  The subspace bornology induced by $\BOUND(V)$
  is equal to $\BOUND(W)$.  The subspace bornology induced $\COMP(V)$ is equal
  to $\COMP(W)$.  That is, a subset of~$W$ is bounded (precompact) in the
  subspace topology iff it is bounded (precompact) as a subset of~$V$.
\end{example}

\begin{example}
  Let~$V$ be a \Frechet space and let $W \subset V$ be a closed subspace.
  Then the quotient bornology on $V/W$ induced by $\COMP(V)$ is equal to
  $\COMP(V/W)$.  That is, a precompact subset of the \Frechet space~$V/W$ is
  the image of a precompact subset of~$V$.  This is an important and
  non-trivial theorem of functional analysis.  It fails if~$V$ is not
  metrizable.

  Moreover, it fails for the bounded bornology.  There may be a bounded subset
  of~$V/W$ that cannot be lifted to a bounded subset of~$V$.
\end{example}

\begin{example}
  The quotient and subspace bornologies induced by a fine bornology are always
  fine.
\end{example}

An important use of subspace and quotient bornologies is to formulate the
general concept of an \emph{extension}.  A diagram
$$
\VS'  \overset{i}{\longrightarrow}
\VS   \overset{p}{\longrightarrow}
\VS''
$$
with complete bornological vector spaces $\VS',\VS,\VS''$ and bounded linear
maps $i,p$ is called an \emph{extension} (of complete bornological vector
spaces) iff $i\colon \VS' \to \VS$ is a bornological isomorphism onto its
image, endowed with the subspace bornology from~$\VS$; $p\circ i=0$; and the
map $\VS/i(\VS') \to \VS''$ induced by~$p$ is a bornological isomorphism.

A \emph{bounded linear section} for an extension is a bounded linear map
$s\colon \VS'' \to \VS$ such that $p\circ s = \ID$.  The existence of a
bounded linear section implies that $\VS \cong \VS' \oplus \VS''$ as
bornological vector spaces.  The vector spaces in an extension often carry
additional structure: they are algebras or complexes.  Unfortunately, most
tools of universal algebra and homological algebra only apply to extensions
with a bounded linear section.  The bounded linear section is used as raw
material to construct homomorphisms or chain maps.

\begin{definition}
  We briefly write $(i,p) \colon \VS' \injto \VS \prto \VS''$ or $\VS' \injto
  \VS \prto \VS''$ iff $\VS',\VS,\VS'',i,p$ form the data of an extension.  An
  extension is called \emph{allowable} iff it has a bounded linear section.
\end{definition}

\subsection{Completions}
\label{sec:completion}

We define the \emph{completion}~$\Cpl{\VS}$ of a convex bornological vector
space~$\VS$ as the universal \emph{complete} target for bounded linear maps
out of~$\VS$.  That is, $\Cpl{\VS}$ is a complete bornological vector space
with a bounded linear map $\natural \colon \VS \to \Cpl{\VS}$ such that any
bounded linear map $l \colon \VS \to \VS[W]$ with complete range~$\VS[W]$ can
be factored uniquely as $\Cpl{l} \circ \natural$ for a bounded linear map
$\Cpl{l} \colon \Cpl{\VS} \to \VS[W]$.  This universal property
determines~$\Cpl{\VS}$ and the map $\natural \colon \VS \to \Cpl{\VS}$
uniquely up to bornological isomorphism.  The completion is introduced by
Hogbe-Nlend in~\cite{hogbe-nlend70:completions}.  To construct the completion
explicitly, we have to use the relationship between bornological vector spaces
and inductive systems.  Therefore, the proof of existence of completions is
relegated to Appendix~\ref{app:bornologies}.  There we also prove that we can
``extend'' bounded multi-linear maps to completions:

\begin{lemma}  \label{lem:Cpl_multi}
  Let $\VS_1, \dots, \VS_n$ be convex bornological vector spaces and
  let~$\VS[W]$ be a complete bornological vector space.  
  Let $l \colon \VS_1 \times \dots \times \VS_n \to \VS[W]$ be a bounded
  \Mpn{n}linear map.  Then there is a unique bounded \Mpn{n}linear map
  $\Cpl{l}\colon \Cpl{\VS_1} \times \dots \times \Cpl{\VS_n} \to \VS[W]$ such
  that $\Cpl{l} \circ (\natural_1 \times \dots \times \natural_n) = l$.
  
  Furthermore, if~$L$ is an equibounded set of \Mpn{n}linear maps $\VS_1
  \times \dots \times \VS_n \to \VS[W]$, then the set of extensions $\Cpl{L}
  \defeq \{ \Cpl{l} \mid l \in L\}$ is an equibounded set of \Mpn{n}linear
  maps $\Cpl{l} \colon \Cpl{\VS_1} \times \dots \times \Cpl{\VS_n} \to
  \VS[W]$.
\end{lemma}

In particular, the completion of a bornological algebra is a complete
bornological algebra.  The completion in the bornological framework is not as
well-behaved as the completion of a locally convex topological vector space.
Therefore we will mainly use its universal property and
Lemma~\ref{lem:Cpl_multi}.  In some places, we need the following
concrete information about the bornology of the completion:

\begin{lemma}  \label{lem:completion_bornology}
  Let~$\VS$ be a convex bornological vector space and let $\natural \colon \VS
  \to \Cpl{\VS}$ be the natural map.  Each small subset of~$\Cpl{\VS}$ is
  contained in $\coco{\bigl( \natural(S) \bigr)}$ for suitable $S \in
  \CBS(\VS)$.  In particular, each point of~$\Cpl{\VS}$ is contained in a set
  of the form $\coco{\bigl( \natural(S) \bigr)}$.
\end{lemma}

\begin{proof}
  This follows immediately from the explicit construction of the completion in
  Appendix~\ref{app:bornologies}.  Alternatively, it can be deduced from the
  universal property as follows.  Consider the subspace~$\VS'$ of $\Cpl{\VS}$
  spanned by the sets of the form $\coco{\bigl( \natural(S) \bigr)}$.  Endow
  it with the bornology~$\CBS'$ of all subsets of sets of the form
  $\coco{\bigl( \natural(S) \bigr)}$.  Verify that $(\VS', \CBS')$ is a
  complete bornological vector space and show that the universal property
  of~$\Cpl{\VS}$ implies that $(\VS', \CBS')$ has the universal property as
  well.
\end{proof}

\begin{warning}
  The most serious problem with completions is that the natural map $\natural
  \colon \VS \to \Cpl{\VS}$ may fail to be injective.  Therefore, the phrase
  that~``$\Cpl{l}$ extends~$l$'' should not be taken too literally.
\end{warning}

\begin{example}
  There is a separated convex bornological vector space~$\VS$ such that the
  only bounded linear map $l \colon \VS \to \VS[W]$ with complete
  range~$\VS[W]$ is the zero map.  Thus $\Cpl{\VS} = \{0\}$ is the zero space.
  The following example is taken from~\cite{Hogbe-Nlend:Born}.
  
  Let $\VS = t \cdot \C[t]$ be the vector space of polynomials without
  constant coefficient.  For all $\epsilon > 0$, let $\| f \|_\epsilon$ be the
  maximum of the polynomial $f \in t \cdot \C[t]$ on the interval $[0,
  \epsilon]$.  This is a \emph{norm} for each $\epsilon > 0$.  A subset $S
  \subset \VS$ is declared small iff $\|S\|_\epsilon < \infty$ for some
  $\epsilon > 0$, that is, $S$ is bounded with respect to the norm
  $\|\blank\|_\epsilon$ for some $\epsilon>0$.  This defines a separated
  convex bornology on~$\VS$.  Let $l \colon \VS \to \VS[W]$ be a bounded
  linear map with complete range~$\VS[W]$.  We claim that $l = 0$.
  
  Let~$\VS_\epsilon$ be the completion of~$\VS$ with respect to the norm
  $\|\blank\|_\epsilon$.  Thus~$\VS_\epsilon$ is isomorphic to $\NBC_0( \left]
  0, 1 \right] )$.  Since~$\VS[W]$ is complete, we can extend~$l$ to a bounded
  linear map $l_\epsilon \colon \VS_\epsilon \to \VS[W]$.  Pick a continuous
  function $f \in \VS_\epsilon$ vanishing in a neighborhood of~$0$.  Then~$f$
  is annihilated by the restriction map $ r_{\epsilon, \epsilon'} \colon
  \VS_\epsilon \to \VS_{\epsilon'}$ for suitably small $\epsilon' \in \left]
  0, \epsilon \right[$.  However, $l_\epsilon = l_{\epsilon'} \circ
  r_{\epsilon, \epsilon'}$ factors through this restriction map.  Hence
  $l_\epsilon(f) = 0$.  Since $f(0) = 0$ for all $f \in \VS_\epsilon$, the
  space of functions vanishing in a neighborhood of~$0$ is dense
  in~$\VS_\epsilon$.  Thus $l_\epsilon = 0$ and hence $l = 0$.
\end{example}

\begin{example}
  Let~$\VS$ be a primitive space.  Then~$\Cpl{\VS}$ is naturally isomorphic to
  the Hausdorff completion of~$\VS$.  It is elementary to verify the universal
  property.
\end{example}

\begin{example}
  Let~$V$ be a \Frechet space and let $W\subset V$ be a dense subspace,
  endowed with the precompact bornology.  Then the completion of~$W$ is equal
  to~$V$, endowed with the precompact bornology.  Since $(V, \COMP)$ is
  complete, the inclusion $W \to V$ induces a natural bounded linear map $l
  \colon \Cpl{(W,\COMP)} \to (V,\COMP)$.  The assertion is that this map is a
  bornological isomorphism.
  
  All precompact subsets of~$V$ are contained in a set of the form~$\coco{S}$
  with a precompact set $S \subset W$.  In fact, we can choose for~$S$ the
  points of a null-sequence in~$W$.  This implies at once that~$l$ is a
  quotient map.  That is, $l$ is surjective and each small subset of $(V,
  \COMP)$ is the image of a small subset of $\Cpl{(W, \COMP)}$.  It remains to
  show that~$l$ is injective.  This follows easily from
  Lemma~\ref{lem:completion_injective} and the observation that topological
  and bornological null-sequences are the same both in $(V, \COMP)$ and
  $(W, \COMP)$.
\end{example}

If we have a subset $S \subset \VS$ of a bornological vector space, we define
its \emph{linear span} $\lin S \subset \VS$ to be the vector subspace of~$\VS$
generated by~$S$ endowed with the subspace bornology from~$\VS$.  The
\emph{completant linear hull} $\cllin S$ is defined as the completion of $\lin
S$.  By Lemma~\ref{lem:completion_injective}, if~$\VS$ is complete then we
have natural injective bounded linear maps $\lin S \to \cllin S \to \VS$.

\begin{warning}
  There is no reason to believe that the ``completion bornology'' on $\cllin
  S$ is the subspace bornology from~$\VS$.  This would imply that $\cllin S$
  is a bornologically closed subspace of~$\VS$.  However, each point in
  $\cllin S$ is the limit of a sequence in $\lin S$.  In general, there may be
  points in the bornological closure of $\lin S$ that are not limits of
  sequences.
  
  Thus the completant linear hull of~$S$ has to be distinguished from the
  closed linear span of~$S$ in~$\VS$.  We will not use closed linear spans
  unless they are equal to completant linear hulls.
\end{warning}

\subsection{Completed bornological tensor products}
\label{sec:hot}

The \emph{completed (bornological) tensor product} $\VS_1 \hot \VS_2$ of two
convex bornological vector spaces is defined as the universal complete target
for bounded bilinear maps $\VS_1 \times \VS_2 \to \blank$.  That is, part of
the structure of $\VS_1 \hot \VS_2$ is a bounded bilinear map $b \colon \VS_1
\times \VS_2 \to \VS_1 \hot \VS_2$; the bornological vector space $\VS_1 \hot
\VS_2$ is complete; and composition with~$b$ gives rise to a bijection between
bounded linear maps $\VS_1 \hot \VS_2 \to \VS[W]$ and bounded bilinear maps
$\VS_1 \times \VS_2 \to \VS[W]$ for all complete bornological vector
spaces~$\VS[W]$.  This universal property determines $\VS_1 \hot \VS_2$
and~$b$ uniquely up to bornological isomorphism.

To construct $\VS_1 \hot \VS_2$, first endow the algebraic tensor product
$\VS_1 \otimes \VS_2$ with the convex bornology generated by the \emph{bismall
sets} $S_1 \otimes S_2 = \{\vs_1 \otimes \vs_2 \mid \vs_1 \in S_1,\ \vs_2 \in
S_2\}$.  Thus a subset $S \subset \VS_1 \otimes \VS_2$ is small if and only if
it is contained in the disked hull of a bismall set.  We always endow $\VS_1
\otimes \VS_2$ with this bornology.  It is usually not complete but has the
right universal property for a tensor product in the category of convex
bornological vector spaces.  Hence the completion $\VS_1 \hot \VS_2 \defeq
\Cpl{(\VS_1 \otimes \VS_2)}$ is a model for the completed bornological tensor
product.  The universal property follows immediately from the universal
property of completions and of $\VS_1 \otimes \VS_2$.  It is clear that $\VS_1
\hot \VS_2$ is functorial for bounded linear maps in both variables.

\begin{proposition}  \label{pro:tensor_associative}
  The completed tensor product is associative.  The \Mpn{n}fold completed
  tensor product $\VS_1 \hot \dots \hot \VS_n$ is universal for bounded
  \Mpn{n}linear maps $\VS_1 \times \dots \times \VS_n \to \VS[W]$ with
  complete range~$\VS[W]$.
\end{proposition}

\begin{digression}
  Associativity appears to be a trivial property but one should not take it
  for granted.  It fails for \Grothendieck's completed inductive tensor
  product~$\indt$ that is universal for separately continuous bilinear
  maps~\cite{grothendieck55:produits}.  The reason is that separately
  continuous bilinear maps may fail to extend to the completions.  Hence $(V_1
  \indt \C) \indt V_2$ and $V_1 \indt (\C \indt V_2)$ need not be isomorphic.
\end{digression}

\begin{proof}
  The corresponding assertions for the uncompleted tensor product are easy.
  It is clear that the bornologies on $(\VS_1 \otimes \VS_2) \otimes \VS_3$
  and $\VS_1 \otimes (\VS_2 \otimes \VS_3)$ coincide with the bornology
  generated by the trismall sets $S_1 \otimes S_2 \otimes S_3$.  Thus the
  uncompleted bornological tensor product is associative and the \Mpn{n}fold
  uncompleted tensor product is universal for bounded \Mpn{n}linear maps.
  
  To carry this over to the completed tensor product, we have to extend
  bounded multi-linear maps to the completions.
  Lemma~\ref{lem:Cpl_multi} implies that $\Cpl{\VS_1} \hot \Cpl{\VS_2}
  \cong \Cpl{(\VS_1 \otimes \VS_2)}$ for all convex bornological vector spaces
  because both sides are complete bornological vector spaces universal for
  bounded bilinear maps $\VS_1 \times \VS_2 \to \VS[W]$ with complete
  range~$\VS[W]$.  Especially, $\VS_1 \hot (\VS_2 \hot \VS_3) \cong
  \Cpl{(\VS_1 \otimes (\VS_2 \otimes \VS_3))}$.
  Proposition~\ref{pro:tensor_associative} follows.
\end{proof}

The functor~$\hot$ is additive, even for infinite sums.  It commutes with
inductive limits, that is, $(\indlim \VS_i) \hot \VS[W] \cong \indlim {(\VS_i
  \hot \VS[W])}$ for all inductive systems of complete bornological vector
spaces $(\VS_i)_{i\in I}$ and complete bornological vector spaces~$\VS[W]$.
This follows from~\eqref{eq:indlim_universal_multi}.  However, $\hot$ does not
commute with infinite direct products.

\begin{lemma}
  Let $\MA$ and~$\MA[B]$ be complete bornological algebras.  Then the natural
  multiplication makes $\MA \hot \MA[B]$ a complete bornological algebra.
\end{lemma}

\begin{proof}
  View the multiplications in $\MA$ and~$\MA[B]$ as bounded linear maps $\MA
  \hot \MA \to \MA$, $\MA[B] \hot \MA[B] \to \MA[B]$.  Then we can tensor them
  to get a bounded linear map $\MA \hot \MA \hot \MA[B] \hot \MA[B] \to \MA
  \hot \MA[B]$.  Up to a flip of the tensor factors, this is the linearized
  version of the natural multiplication in the tensor product.  Hence the
  natural multiplication is bounded.
\end{proof}

\begin{example}
  Let $\VS_1$ and~$\VS_2$ be fine spaces.  The completed bornological tensor
  product $\VS_1\hot \VS_2$ is naturally isomorphic to the algebraic tensor
  product $\VS_1\otimes \VS_2$ endowed with the fine bornology.
\end{example}

\begin{example}
  Let $\VS_1$ and~$\VS_2$ be primitive spaces.  Then $\VS_1 \hot \VS_2$ is
  naturally isomorphic to the \emph{completed projective tensor product}
  $\VS_1 \prot \VS_2$ endowed with the primitive bornology.  In fact, this is
  the definition of the projective tensor
  product~\cite{grothendieck55:produits}.
\end{example}

The example of \Frechet spaces deserves to be called a theorem.

\begin{theorem}  \label{the:tensor_Frechet}
  Let\/ $\VS_1$ and~$\VS_2$ be \Frechet spaces and let\/ $\VS_1 \prot \VS_2$ be
  their completed projective tensor product~\cite{grothendieck55:produits}.
  The natural bilinear map $\natural\colon \VS_1 \times \VS_2 \to \VS_1 \prot
  \VS_2$ induces bornological isomorphisms
  \begin{gather*}
    (\VS_1, \COMP ) \hot (\VS_2, \COMP ) \cong (\VS_1 \prot \VS_2, \COMP);
    \\
    (\VS_1, \BOUND) \hot (\VS_2, \BOUND) \cong (\VS_1 \prot \VS_2, \CBS ).
  \end{gather*}
  Here~$\CBS$ denotes the bornology of all $S \subset \VS_1 \prot \VS_2$ that
  are contained in a set of the form $\coco{(B_1 \otimes B_2)}$ with bounded
  sets $B_1 \in \BOUND(\VS_1)$ and $B_2 \in \BOUND(\VS_2)$.
  
  The bornology $\CBS$ is equal to $\BOUND( \VS_1 \prot \VS_2 )$ in the
  following cases.  If both\/ $\VS_1$ and~$\VS_2$ are Banach spaces;
  if\/~$\VS_2$ is arbitrary and\/ $\VS_1 = L^1(M, \mu)$ is the space of
  integrable functions on a locally compact space with respect to some Borel
  measure; or if\/~$\VS_2$ is arbitrary and\/~$\VS_1$ is nuclear.
\end{theorem}

This theorem is proved in Appendix~\ref{app:tensor_Frechet}.  Although I have
not found Theorem~\ref{the:tensor_Frechet} in the literature, it may be known
to the experts.  The proof is an application of \Grothendieck's fundamental
theorems about compact subsets of projective tensor products of \Frechet
spaces.

\begin{corollary}  \label{cor:tensor_LF}
  Let\/ $\VS_1$ and\/~$\VS_2$ be nuclear LF-spaces.  Then $(\VS_1, \BOUND) \hot
  (\VS_2, \BOUND)$ is isomorphic to \Grothendieck's inductive tensor product\/
  $\VS_1 \indt \VS_2$, endowed with the bounded bornology.  The inductive
  tensor product\/ $\VS_1 \indt \VS_2$ is again a nuclear LF-space.
\end{corollary}

This corollary is proved in Appendix~\ref{app:tensor_Frechet}.

\begin{example}  \label{exa:tensor_smooth_manifold}
  If~$M$ is a smooth manifold, denote by $\CCINF(M)$ the LF-space of smooth
  compactly supported functions on~$M$.  We endow $\CCINF(M)$ with the bounded
  bornology, which equals the precompact bornology.  Thus $S \in
  \CBS\bigl(\CCINF(M) \bigr)$ iff there is a compact subset $K \subset M$ such
  that all functions $f\in S$ have support in~$K$ and all total derivatives
  $D^{(n)}f$, $f\in S$, are uniformly bounded on~$K$.  The space $\CCINF(M)$
  is nuclear.  Therefore, Corollary~\ref{cor:tensor_LF} implies that
  $$
  \CCINF(M) \hot \CCINF(N) \cong \CCINF(M \times N)
  $$
  for all smooth manifolds $M,N$.
\end{example}

\subsection{Spaces of bounded linear maps}
\label{sec:Lin}

It is clear that the bounded linear maps $\VS_1 \to \VS_2$ form a vector space
$\Lin(\VS_1; \VS_2)$.  The equibounded families of linear maps $\VS_1 \to
\VS_2$ form a convex bornology $\EQUI(\VS_1; \VS_2)$ on $\Lin(\VS_1; \VS_2)$.
We always endow $\Lin(\VS_1; \VS_2)$ with this bornology.  It is separating
if~$\VS_2$ is separated and completant if~$\VS_2$ is complete.  To prove the
completantness, take an equibounded subset $S \subset \Lin(\VS_1; \VS_2)$.  If
$l_j \in S$, $\lambda_j \in \C$, $\sum |\lambda_j| \le 1$, then the series
$\sum \lambda_j l_j$ converges pointwise on~$\VS_1$.  The pointwise limit of
$\sum \lambda_j l_j$ is a bounded operator $\VS_1 \to \VS_2$, and the set of
all limits of such series is an equibounded completant small disk in $\Lin(
\VS_1; \VS_2)$.

Composition of linear maps gives rise to a bounded bilinear map
$$
\Lin(\VS_2; \VS_3) \times \Lin(\VS_1; \VS_2) \to \Lin(\VS_1; \VS_3)
$$
because $S_1 \circ S_2$ is equibounded if both $S_1$ and $S_2$ are
equibounded.  In particular,

\begin{proposition}  \label{pro:Endo_CBA}
  Let~$\VS$ be a complete bornological vector space.  Then $\Endo(\VS) \defeq
  \Lin(\VS; \VS)$ is a complete bornological algebra.
\end{proposition}

There is a natural bornological isomorphism $\Lin(\C; \VS) \cong \VS$.  The
evaluation of linear maps is a bounded bilinear map $\Lin(\VS_1; \VS_2) \times
\VS_1 \to \VS_2$.

\begin{proposition}  \label{pro:adjoint_associativity}
  The functors $\hot$ and $\Lin$ are related by \emph{adjoint associativity}.
  That is, there is a natural bornological isomorphism
  \begin{equation}  \label{eq:adjoint_associativity}
    \Lin(\VS_1 \hot \VS_2; \VS_3) \cong \Lin(\VS_1; \Lin(\VS_2; \VS_3)).
  \end{equation}
\end{proposition}

\begin{proof}
   An equibounded family of bounded linear maps $\VS_1 \hot \VS_2 \to \VS_3$
  corresponds to an equibounded family of bounded bilinear maps $\VS_1 \times
  \VS_2 \to \VS_3$.  By definition, boundedness for bilinear maps $b\colon
  \VS_1 \times \VS_2 \to \VS_3$ means that if we fix $S_1 \in \CBS(\VS_1)$ and
  let~$\vs_1$ vary in~$\VS_1$, then the family of linear maps $\VS_2 \to
  \VS_3$, $\vs_2 \mapsto b(\vs_1, \vs_2)$ maps any small set $S_2 \in
  \CBS(\VS_2)$ into a small set in~$\VS_3$.  This is precisely the condition
  for a subset of $\Lin(\VS_1; \Lin(\VS_2; \VS_3))$ to be equibounded.
\end{proof}

\subsection{Smooth and absolutely continuous homotopies}
\label{sec:homotopy}

We need tensor products mainly to handle homotopies.  A smooth homotopy is a
bounded homomorphism $\MA \to \CINF([0,1]) \hot \MA[B]$.  Although we do not
need any concrete information about $\CINF([0,1]) \hot \MA[B]$ for the
purposes of the general theory, it is helpful to describe how this space looks
like in the familiar examples.  Besides smooth homotopies, we consider
\emph{absolutely continuous homotopies}.  These are not needed for the general
theory but are more general and easier to handle than smooth homotopies.

The \Frechet algebra $\CINF([0,1])$ is endowed with the bounded
bornology, which equals the precompact bornology.  If~$\MA$ is a complete
bornological algebra, we define
$$
\CINF([0,1]; \MA) \defeq \CINF([0,1]) \hot \MA.
$$
This is again a complete bornological algebra.  For $t \in [0,1]$, let $\ev_t
\colon \CINF([0,1]) \to \C$ be evaluation at~$t$.  We get induced bounded
homomorphisms $\ev_t \defeq \ev_t \hot \ID[\MA] \colon \CINF([0,1]; \MA) \to
\MA$.

A \emph{smooth homotopy} is a bounded homomorphism $\Phi \colon \MA \to
\CINF([0,1]; \MA[B])$.  For $t \in [0,1]$, let $\Phi_t \defeq \ev_t \circ
\Phi$.  One can check that~$\Phi$ is determined uniquely by the map $[0,1]
\mapsto \Lin(\MA; \MA[B])$, $t \mapsto \Phi_t$.  Two bounded homomorphisms
$\phi_0, \phi_1 \colon \MA \to \MA[B]$ are \emph{smoothly homotopic} iff there
is a smooth homotopy $\Phi \colon \MA \to \CINF([0,1]; \MA[B])$ such that
$\Phi_t = \phi_t$ for $t = 0,1$.

Smooth homotopy of homomorphisms is an equivalence relation.  It is not quite
obvious that it is transitive.  This is proved by Cuntz
in~\cite{cuntz97:bivariant}.  The idea is to reparametrize the interval $[0,
1]$ by an increasing smooth bijection $h \colon [0,1] \to [0,1]$ with the
property that $h^{(n)}(t) = 0$ for all $n \in \N$, $t = 0,1$.  In this way,
any smooth homotopy $\MA \to \CINF([0,1]; \MA[B])$ can be replaced by a smooth
homotopy $H \colon \MA \to \CINF([0,1]; \MA[B])$ with the additional property
that $H^{(n)}(t) = 0$ for all $n \in \N$, $t = 0,1$.  Homotopies $H,H'$ with
this additional property and $H_1 = H'_0$ can be glued together without
producing a jump in the derivatives at the gluing point.

\begin{example}
  Let $\VS$ and~$\VS[W]$ be \Frechet spaces.  Let~$\CBS$ be the precompact or
  the bounded bornology.  Theorem~\ref{the:tensor_Frechet} implies that
  $\CINF([0,1]) \hot (\VS[W], \CBS) \cong (\CINF([0,1]) \prot \VS[W], \CBS)$
  because $\CINF([0,1])$ is nuclear.  Furthermore, the projective tensor
  product $\CINF([0,1]) \prot \VS[W]$ is isomorphic, as a topological vector
  space, to the space of smooth functions $[0,1] \to \VS[W]$ with the natural
  topology (see~\cite{Treves:67}).
  
  The functor $\VS \mapsto (\VS,\CBS)$ from \Frechet spaces to complete
  bornological vector spaces is fully faithful.  Thus there is a bijection
  betwen bounded linear maps $(\VS, \CBS) \to \CINF\bigl([0,1]; (\VS[W], \CBS)
  \bigr)$ and continuous linear maps $\VS \to \CINF([0,1]) \prot \VS[W]$.
  Hence we get the usual notion of smooth homotopy.
\end{example}

Unfortunately, analytic and periodic cyclic cohomology are not invariant under
arbitrary continuous homotopies.  We need a certain amount of
differentiability.  Actually, we do not need more than just one derivative,
and we only need the derivative to be Lebesgue integrable.

Let $\ABC([0,1])$ be the completion of $\CINF([0,1])$ with respect to the norm
$$
\| f \|_{ac} \defeq |f(0)| + \int_0^1 |f'(t)| \,dt,
$$
where~$dt$ is the usual Lebesgue measure.  Let $L^1([0,1]) \defeq L^1([0,1];
dt)$ be the Banach space of functions on $[0,1]$ that are integrable with
respect to Lebesgue measure, with the standard norm $\| g \|_1 \defeq \int_0^1
|g(t)| \,dt$.  Let $\partial / \partial t \colon \CINF([0,1]) \to \CINF([0,1])
\subset L^1([0,1])$ be the differentiation map $f \mapsto f'$.  The norm $\|
\blank \|_{ac}$ is defined so as to make the linear map $(\ev_0, \partial /
\partial t) \colon \CINF([0,1]) \to \C \oplus L^1([0,1])$ isometric.  Its
range is dense.  Hence we get an isomorphism of Banach spaces $\ABC([0,1])
\cong \C \oplus L^1([0,1])$.  We endow $\ABC([0,1])$ with the bounded
bornology.  The evaluation maps $\ev_t \colon \CINF([0,1]) \to \C$ extend to
$\ABC([0,1])$ in a bounded way because $|\ev_t(f)| = |f(0) + \int_0^t
f'(s)\,ds| \le \|f\|_{ac}$.  It follows easily that $\|f g\| \le 2\|f\| \cdot
\|g\|$.  Hence $\ABC([0,1])$ is a Banach algebra.  Let $f \colon [0,1] \to
\C$.  Then $f \in \ABC([0,1])$ iff~$f'$ exists almost everywhere and $f' \in
L^1([0,1])$.  This is equivalent to~$f$ being \emph{absolutely continuous}.

Let~$\MA$ be a complete bornological algebra.  Define $\ABC([0,1]; \MA) \defeq
\ABC([0,1]) \hot \MA$.  This is again a complete bornological algebra.  An
\emph{absolutely continuous homotopy} is a bounded homomorphism $\Phi \colon
\MA \to \ABC([0,1]; \MA[B])$.  Two bounded homomorphisms $\phi_0, \phi_1
\colon \MA \to \MA[B]$ are \emph{AC-homotopic} iff there is an absolutely
continuous homotopy $\Phi \colon \MA \to \ABC([0,1]; \MA[B])$ such that
$\Phi_t = \phi_t$ for $t = 0,1$.  It is trivial that AC-homotopy is an
equivalence relation.

\begin{example}
  Let~$\VS$ be a \Frechet space endowed with the bounded bornology.
  Theorem~\ref{the:tensor_Frechet} yields a bornological isomorphism
  $$
  \ABC([0,1]) \hot \VS \cong
  \VS \oplus L^1([0,1]) \hot \VS \cong
  \VS \oplus L^1([0,1]) \prot \VS \cong
  \VS \oplus L^1([0,1]; \VS).
  $$
  All spaces above carry the bounded bornology.  The space $L^1([0,1];
  \VS)$ is the space of Lebesgue integrable functions $[0,1] \to \VS$,
  see~\cite{grothendieck55:produits}.  Thus $\ABC([0,1]; \VS)$ can be computed
  quite explicitly.  This is an advantage of $\ABC([0,1])$ over the algebra of
  continuously differentiable functions $\NBC[1]([0,1])$.
\end{example}

\chapter{Analytic cyclic cohomology}
\label{cha:HA}

\section[Analytic tensor algebras and a-nilpotent algebras]{Analytic tensor
  algebras and analytically nilpotent algebras}
\label{sec:Tanil}

We define the \emph{analytic differential envelope $\Omega_\an\MA$} of a
complete bornological algebra~$\MA$ as the completion of the usual
differential envelope $\Omega\MA$ with respect to a certain bornology.  We
refer the reader to~\cite{connes94:ncg} for a discussion of the algebra
$\Omega\MA$ over a \emph{non-unital} algebra~$\MA$.  Cuntz and Quillen
\cite{cuntz95:algebra}, \cite{cuntz95:cyclic} work with unital algebras and
therefore use a slightly different definition of $\Omega\MA$.  The
\emph{analytic tensor algebra} $\Tanil\MA$ is defined as the even part of
$\Omega_\an\MA$ endowed with the Fedosov product as multiplication.  Thus
$\Tanil\MA$ is a completion of the purely algebraic \emph{tensor algebra}
$\Tens\MA$ with respect to a certain bornology.  The idea that the Fedosov
product on $\Omega\MA$ can be used to obtain the tensor algebra of~$\MA$ is
due to Cuntz and Quillen~\cite{cuntz95:algebra}.

Let $\tau_{\MA} \colon \Tanil\MA \to \MA$ be the natural projection that sends
a differential form to its degree~$0$ component.  The kernel $\Janil\MA \defeq
\Ker \tau_{\MA}$ is the ideal of even forms of degree at least~$2$.  The
bounded linear map $\sigma_{\MA} \colon \MA \congto \Omega^0\MA \subset
\Tanil\MA$ is a natural section for~$\tau_{\MA}$.  Thus we get a natural
allowable extension of complete bornological algebras
\begin{equation}  \label{eq:Tanil_extension}
  \xymatrix{
    {\Janil\MA\;} \ar@{>->}[r] &
      {\Tanil\MA} \ar@{->>}[r]^{\tau_{\MA}} &
        {\MA.} \ar@/^/@{.>}[l]^{\sigma_{\MA}}
    }
\end{equation}

Analytic cyclic cohomology is defined as the cohomology of the X-complex of
$\Tanil\MA$.  Thus we divide the definition of analytic cyclic cohomology into
two steps, the functor $\Tanil$ and the X-complex functor.  The X-complex is a
rather small object and can be handled by techniques of homological algebra.
The main problem is to understand the analytic tensor algebra $\Tanil\MA$.
Following Cuntz and Quillen~\cite{cuntz95:cyclic}, we interpret $\Tanil\MA$ as
a ``resolution of singularities'' of~$\MA$.

To fill this philosophy with meaning, we have to explain the notion of
``non-singularity'' that is relevant here.  In algebraic geometry, there is
the following criterion for an affine variety to be smooth.  The commutative
algebra~$A$ describes a smooth affine variety iff each extension $N \injto E
\prto A$ of commutative algebras with nilpotent~$N$ can be split by a
homomorphism $A \to E$.  We can carry over this definition to non-commutative
spaces by simply dropping the commutativity assumption.  An algebra that has
the property that all (non-commutative) nilpotent extensions of it split by a
homomorphism is called \emph{quasi-free} by Cuntz and
Quillen~\cite{cuntz95:algebra}.  Quasi-free algebras can be characterized by
several equivalent conditions (\ref{deflem:quasi_free}).  Most notably, an
algebra is quasi-free iff it has \Hochschild homological dimension at
most~$1$.  In particular, a commutative algebra is quasi-free iff it is smooth
and at most \emph{\Mpn{1}dimensional}.

An algebra~$\MA[N]$ is \emph{nilpotent} iff $\MA[N]^k = \{0\}$ for some $k \in
\N$.  This notion of nilpotence and the corresponding notion of quasi-freeness
is purely algebraic and not related to the analytic tensor algebra.
Therefore, we have to consider a larger class of ``nilpotent'' algebras.

\begin{definition}  \label{def:a_nilpotent}
  A convex bornological algebra~$\MA[N]$ is \emph{analytically nilpotent}
  iff
  $$
  S^\infty \defeq \bigcup_{n \in \N} S^n
  $$
  is a small subset of~$\MA[N]$ for all $S \in \CBS(\MA[N])$.
\end{definition}

We usually abbreviate the adverb ``analytically'' by the prefix ``a-'',
writing \emph{a-nilpotent} instead of analytically nilpotent.

Let~$\MA[N]$ be a complete a-nilpotent algebra.  If $S \in \CBS(\MA[N])$, then
$c \cdot S \in \CBS(\MA[N])$ for all $c \in \R$ and hence $(c \cdot S)^\infty
= \bigcup c^n S^n$ is small.  Thus all elements $x \in \MA[N]$ have ``spectral
radius zero'' in the following sense.  Let $\HO(\{0\} \subset \C^{\,n})_0$ be
the algebra of germs of analytic functions~$f$ in a neighborhood of the origin
in~$\C^{\,n}$ satisfying $f(0) = 0$.  If $x_1, \dots, x_n \in \MA[N]$ are
commuting elements, then there is a well-defined bounded functional calculus
$\HO(\{0\} \subset \C^{\,n})_0 \to \MA[N]$ for the elements $x_1, \dots, x_n$.
We may interpret the definition of analytic nilpotence as the condition that
we have a functional calculus for ``germs of non-commutative analytic
functions in infinitely many variables'', as long as these variables all
remain in a fixed small set $S \subset \MA[N]$.

\begin{definition}  \label{def:a_quasi_free}
  A complete bornological algebra~$\MA[R]$ is \emph{analytically quasi-free}
  iff there is a bounded homomorphism $\upsilon \colon \MA[R] \to
  \Tanil\MA[R]$ that is a section for $\tau_{\MA[R]} \colon \Tanil\MA[R] \to
  \MA[R]$.  That is, $\tau_{\MA[R]} \circ \upsilon = \ID[{\MA[R]}]$.
\end{definition}

\begin{definition}  \label{def:a_nilpotent_ext}
  An \emph{analytically nilpotent extension} (of~$\MA$) is an
  \textbf{allowable} extension of complete bornological algebras $\MA[N]
  \injto \MA[E] \prto \MA$ with a-nilpotent kernel~$\MA[N]$.
  
  A \emph{universal analytically nilpotent extension} (of~$\MA$) is an
  \textbf{allowable} extension of complete bornological algebras $\MA[N]
  \injto \MA[R] \prto \MA$ with a-nilpotent~$\MA[N]$ and
  a-quasi-free~$\MA[R]$.
\end{definition}

We will see below that a complete bornological algebra~$\MA[R]$ is
a-quasi-free iff each a-nilpotent extension $\MA[N] \injto \MA[E] \prto
\MA[R]$ has a bounded splitting homomorphism.  Thus analytic quasi-freeness is
the notion of ``non-commutative non-singularity'' corresponding to analytic
nilpotence.  The advantage of Definition~\ref{def:a_quasi_free} is that it is
easier to verify.

The \emph{geometric picture} behind Definition~\ref{def:a_nilpotent_ext} is
that of an \emph{infinitesimal thickening} and a \emph{resolution of
  singularities}.  An extension $\MA[N] \injto \MA[R] \prto \MA$ can be viewed
as an embedding of the geometric space~$\widehat{\MA}$ corresponding to~$\MA$
as a closed subspace of the space~$\widehat{\MA[R]}$ corresponding
to~$\MA[R]$.  The complement $\widehat{\MA[R]} \setminus \widehat{\MA}$ is the
space~$\widehat{\MA[N]}$ described by the algebra~$\MA[N]$.  The analytic
nilpotence of~$\MA[N]$ can be interpreted as saying that the
space~$\widehat{\MA[N]}$ is infinitesimally small (think of the commutative
examples $\HO(\{0\} \subset \C^{\,n})_0$).  The algebra~$\MA[E]$ in an
a-nilpotent extension $\MA[N] \injto \MA[E] \prto \MA$ can therefore be viewed
as an ``infinitesimal thickening'' of~$\widehat{\MA}$.  A universal
a-nilpotent extension may be thought of as a non-commutative resolution of
singularities by an infinitesimal thickening.

These definitions help to improve our understand of the analytic tensor
algebra $\Tanil\MA$ because the extension $\Janil\MA \injto \Tanil\MA \prto
\MA$ in~\eqref{eq:Tanil_extension} is a universal a-nilpotent extension
of~$\MA$.  In particular, such universal a-nilpotent extensions exist.
Furthermore, as the name \emph{uni\/}versal suggests, there is essentially
only one universal a-nilpotent extension of~$\MA$.  More precisely, $\MA$ may
have many universal a-nilpotent extensions, but they are all smoothly homotopy
equivalent.  This is the assertion of the Uniqueness
Theorem~\ref{the:uniqueness}.  The notion of smooth homotopy equivalence of
extensions is explained below and implies that the entries of the extension
are smoothly homotopy equivalent.  Thus if $\MA[N] \injto \MA[R] \prto \MA$ is
any universal analytically nilpotent extension of~$\MA$, then~$\MA[R]$ is
smoothly homotopy equivalent to $\Tanil\MA$.  This consequence will be used in
the proof of the Excision Theorem~\ref{the:excision_analytic}.

The Uniqueness Theorem allows us to extend homology theories from
``non-singular'' spaces to arbitrary spaces by a resolution of singularities.
A \emph{homology theory for non-singular spaces} is simply a smooth homotopy
functor~$F$ from the category of analytically quasi-free algebras to, say,
Abelian groups.  Such a homology theory can be extended to arbitrary complete
bornological algebras by putting $\bar{F} (\MA) \defeq F(\MA[R])$,
where~$\MA[R]$ is any non-singular infinitesimal thickening of~$\MA$.  In
fact, we can take the standard choice $\MA[R] = \Tanil\MA$ and put $\bar{F}
(\MA) \defeq F(\Tanil\MA)$.  In Section~\ref{sec:X_Tanil_HA} we choose for~$F$
the cohomology of the X-complex to obtain analytic cyclic cohomology.

A similar method can be applied to define the de Rham homology of singular
algebraic varieties.  For smooth varieties, we can define de Rham homology by
the de Rham complex of differential forms.  For singular varieties, this is
not the right definition.  Instead, the singular variety~$S$ is embedded in a
smooth variety~$R$.  The de Rham complex of~$R$ is ``localized at the
subvariety~$S$''.  Roughly speaking, we consider germs of differential forms
in a neighborhood of $S \subset R$.  The cohomology of this localized de Rham
complex is taken as the de Rham homology of~$S$.

The results of this section depend on the close relationship between analytic
tensor algebras and a-nilpotent algebras.  The mediator between them is the
notion of a \emph{\lanilcur{}}.  This name is an abbreviation for
\emph{bounded \textbf{l}inear map with \textbf{an}alytically
\textbf{nil}potent \textbf{cur}vature}.

The \emph{curvature} of a bounded linear map $l\colon \MA \to \MA[B]$ between
complete bornological algebras is the bounded bilinear map $\omega_l \colon
\MA \times \MA \to \MA[B]$ defined by
\begin{equation}  \label{eq:def_curvature}
  \omega_l(\ma_1, \ma_2) \defeq
  l(\ma_1 \cdot \ma_2) - l(\ma_1) \cdot l(\ma_2).
\end{equation}
We often view~$\omega_l$ as a bounded linear map $\omega_l \colon \MA \hot \MA
\to \MA[B]$.  The curvature measures the deviation of~$l$ from being a
homomorphism.

\begin{definition}  \label{def:lanilcur}
  A bounded linear map $l \colon \MA \to \MA[B]$ is a \emph{\lanilcur{}} iff
  its curvature is a-nilpotent in the sense that $\omega_l(S, S)^\infty
  \subset \MA[B]$ is small for all $S \in \CBS(\MA)$.
\end{definition}

The machinery of universal nilpotent extensions applies to other situations as
well, where we have analogues of a-nilpotent algebras, \lanilcurs, and
analytic tensor algebras.  We will examine the case of ``pro-algebras'' in
Section~\ref{sec:pro_systems}.  To minimize the number of results that have to
be checked again, we proceed axiomatically and base everything on a few
results.  These are called ``axioms'' to distinguish them from those results
that follow from the axioms.

The first axiom~\ref{axiom:Tanil_universal} is the universal property of
$\Tanil\MA$.  If a bounded linear map $l \colon \MA \to \MA[B]$ between
complete bornological algebras is a \lanilcur, then there is a unique bounded
homomorphism $\LLH{l} \colon \Tanil\MA \to \MA[B]$ extending~$l$ in the sense
that $\LLH{l} \circ \sigma_{\MA} = l$.  In addition, linear maps of the form
$f \circ \sigma_{\MA}$ for a bounded homomorphism $f \colon \Tanil\MA \to
\MA[B]$ are automatically \lanilcurs.  Thus the analytic tensor algebra can be
defined as the universal target for \lanilcurs.  This determines the algebra
$\Tanil\MA$ and the linear map $\sigma_{\MA} \colon \MA \to \Tanil\MA$
uniquely up to isomorphism.  Conversely, we can define \lanilcurs as those
linear maps that can be extended to bounded homomorphisms $\Tanil\MA \to
\MA[B]$.  Thus \lanilcurs and analytic tensor algebras determine each other.

The second axiom~\ref{axiom:Janil_nilpotent} is that the kernel $\Janil\MA$ of
the extension~\eqref{eq:Tanil_extension} is always a-nilpotent.  The third and
fourth axioms \ref{axiom:nilpotent_lanilcur}
and~\ref{axiom:lanilcur_nilpotent} relate a-nilpotent algebras and \lanilcurs.
They imply that nilpotent algebras are a-nilpotent and that an
algebra~$\MA[N]$ is a-nilpotent iff each bounded linear map into it is a
\lanilcur.  Thus a-nilpotent algebras are determined if we know the definition
of \lanilcurs.  Conversely, \lanilcurs are determined if we know the
definition of a-nilpotent algebras.  A bounded linear map $l \colon \MA \to
\MA[B]$ is a \lanilcur iff its curvature can be factored as $\omega_l = f
\circ \omega'$ with a bounded homomorphism $f \colon \MA[N] \to \MA[B]$, a
bounded linear map $\omega' \colon \MA \hot \MA \to \MA[N]$, and
a-nilpotent~$\MA[N]$.

Finally, we need the Extension Axiom~\ref{axiom:extension} and the Homotopy
Axiom~\ref{axiom:homotopy}, asserting that the class of a-nilpotent algebras
is closed under forming extensions and tensor products with $\CINF([0,1])$.
The Extension Axiom is equivalent to the assertion that the composition of two
\lanilcurs is again a \lanilcur and implies that analytic tensor algebras are
a-quasi-free.  We prove these axioms and reduce all other results to them.
The first main consequence is the \emph{Universal Extension
Theorem}~\ref{the:universal_ext}.  It asserts that universal a-nilpotent
extensions satisfy a universal property among arbitrary a-nilpotent
extensions, with uniqueness up to smooth homotopy.  The Uniqueness
Theorem~\ref{the:uniqueness} asserts that all universal a-nilpotent extensions
are smoothly homotopy equivalent.  This is a straightforward consequence of
the universal property.  Moreover, we get some equivalent characterizations of
a-quasi-free algebras (Theorem~\ref{the:aqf_definitions}).

We use the universal property of $\Tanil\MA$ to determine the
\Mp{\Tanil\MA}bimodule $\Omega^1(\Tanil\MA)$.  This will be needed to
determine the X-complex of $\Tanil\MA$.  The idea of this proof goes back to
Cuntz and Quillen \cite[Proposition~2.6]{cuntz95:algebra}.  A more concrete
proof is relegated to the Appendices.  The advantage of the universal algebra
proof is that it easily generalizes to other tensor algebras that are related
to different notions of nilpotence.

Since the composition of two \lanilcurs is again a \lanilcur, they form the
morphisms of a category.  The functor~$\Tanil$ is naturally defined on the
\lanilcur category and is a smooth homotopy functor for \lanilcurs.  A nice
consequence is the ``Goodwillie Theorem'' that if $\MA[N] \injto \MA[E] \prto
\MA$ is an a-nilpotent extension, then $\Tanil\MA[N]$ is smoothly contractible
and $\Tanil\MA$ is a smooth deformation retract of $\Tanil\MA[E]$.  This
implies that excision in analytic cyclic cohomology holds for a-nilpotent
extensions.  However, we will prove excision in complete generality and will
not use the Goodwillie Theorem for that purpose.

At least the definitions of analytic tensor algebras, analytic nilpotence, and
\lanilcurs are given by Puschnigg~\cite{puschnigg98:cyclic} in the setting of
inductive systems of Banach algebras.  The analytic tensor algebra $\Tanil\MA$
is essentially the same as Puschnigg's ``universal linear deformation''
\cite{puschnigg98:cyclic} in the bornological framework (however, Puschnigg's
definition is much more complicated).  The analogue of analytic nilpotence is
called ``topological nilpotence'' in~\cite{puschnigg98:cyclic}, the analogue
of a \lanilcur is called an ``almost multiplicative map''
\cite{puschnigg98:cyclic}.  However, Puschnigg does not study topological
nilpotence at all and does not state results like the Universal Extension
Theorem or the Uniqueness Theorem.  Instead the almost multiplicative maps
occupy a prominent place in his theory.  In this thesis, \lanilcurs are
regarded only as a technical tool.  They do not show up in the statements of
theorems.

\subsection{Definition of the analytic tensor algebra}
\label{sec:def_Tanil}

Let~$\MA$ be a complete bornological algebra.  For $S \in \CBS(\MA)$, define
\begin{align*}
   \opt{S}(dS)^\infty
 &\defeq
   S(dS)^\infty \cup (dS)^\infty \cup S \subset
   \Omega\MA,
 \\
   \opt{S}(dS)^\even
  &\defeq
   \opt{S}(dS)^\infty \cap \Omega^\even\MA,
 \\
   \opt{S}(dS)^\odd
  &\defeq
   \opt{S}(dS)^\infty \cap \Omega^\odd\MA.
\end{align*}
The notation~$\opt{S}$ above will be used frequently in the following.  If~$S$
is a subset or an element of a bornological algebra~$\MA$, then~$\opt{S}$
denotes $S\cup \{1\}$, where~$1$ is the unit in the
\emph{unitarization}~$\Unse{\MA}$.  For example,
$$
\opt{S}\cdot T\cdot \opt{S} = S\cdot T\cdot S\cup S\cdot T\cup T\cdot S.
$$
Sometimes we consider $\opt{S} \subset \MA$, leaving out~$1$.  The expression
$\opt{\ma_0} d\ma_1 \dots d\ma_n$ denotes $\ma_0 d\ma_1 \dots d\ma_n$ or
$d\ma_1 \dots d\ma_n$.  This is a convenient notation for a generic
differential form.

\begin{warning}
  If~$\MA$ is unital with unit~$1$, then the differential forms $1\cdot d\ma_1
  \dots d\ma_n$ and $d\ma_1 \dots d\ma_n$ are \emph{different}.  That is,
  $1\in \MA$ is not a unit in $\Omega\MA$.
\end{warning}

\begin{definition}  \label{def:Omega_nil}
  Let~$\CBS_\an$ be the convex bornology on $\Omega\MA$ generated by the sets
  $\opt{S}(dS)^\infty$ with $S\in\CBS(\MA)$.  The \emph{analytic differential
    envelope} $\Omega_\an\MA$ of~$\MA$ is the completion of the convex
  bornological vector space $(\Omega\MA, \CBS_\an)$.  Let
  $\Omega^\even_\an\MA$ and $\Omega^\odd_\an\MA$ be the even and odd part of
  $\Omega_\an\MA$.
\end{definition}

The universal property of the completion asserts that bounded linear maps
$\Omega_\an \MA \to \VS[W]$ with complete range~$\VS[W]$ are in bijection with
bounded linear maps $(\Omega\MA, \CBS_\an) \to \VS[W]$.  Since the
bornology~$\CBS_\an$ is generated by the sets $\opt{S} (dS)^\infty$, a linear
map $l \colon \Omega\MA \to \VS[W]$ is bounded with respect to~$\CBS_\an$ iff
$l(\opt{S} (dS)^\infty) \in \CBS(\VS[W])$ for all $S\in \CBS_c(\MA)$.  Thus
bounded linear maps $\Omega_\an \MA \to \VS[W]$ are in bijection with linear
maps $l\colon \Omega\MA \to \VS[W]$ satisfying $l(\opt{S} (dS)^\infty) \in
\CBS(\VS[W])$ for all $S \in \CBS(\MA)$.

By Lemma~\ref{lem:completion_bornology}, a subset of $\Omega_\an\MA$ is small
iff it is contained in $\coco{ (\opt{S} (dS)^\infty)}$ for some $S \in
\CBS_c(\MA)$.  Similarly, a subset of $\Omega_\an^\even\MA$ is small iff it is
contained in $\coco{ (\opt{S} (dS)^\even)}$ for some $S \in \CBS_c(\MA)$.

If $C \in \R$, $S \in \CBS(\MA)$, then $C\cdot S \in \CBS(\MA)$ and hence $S
\cup \bigcup_{n\ge 1} C^n \opt{S} (dS)^n \in \CBS_\an$.  Consequently, the
operator $f(N) \colon \Omega\MA \to \Omega\MA$ that multiplies a form of
degree~$n$ by~$f(n)$ is bounded with respect to the bornology~$\CBS_\an$ if
$f\colon \Z_+ \to \C$ is of at most exponential growth.  Conversely, $f(N)$ is
unbounded if~$f$ does not have exponential growth.  Thus~$f(N)$ can be
extended to a bounded linear map $f(N) \colon \Omega_\an\MA \to \Omega_\an\MA$
iff~$f$ has at most exponential growth.  If we let $f \colon \Z_+ \to \{0,1\}$
be the characteristic function of the set of even integers, we find that the
even and odd parts $\Omega_\an^\even\MA$ and $\Omega_\an^\odd\MA$ are
complementary direct summands of $\Omega_\an\MA$.  That is,
$\Omega_\an^\even\MA \oplus \Omega_\an^\odd\MA \cong \Omega_\an\MA$.  Thus
$\Omega_\an\MA$ is a \Mpn{\Ztwo}graded complete bornological vector space.  If
we let $f \colon \Z_+ \to \{0,1\}$ be the characteristic function of $\{j\}$,
we find that $\Omega^j\MA \subset \Omega_\an\MA$ is a direct summand.

The familiar operators $d$, $b$ and~$b'$ on $\Omega\MA$ (see \eqref{eq:def_b}
and~\eqref{eq:def_bPrime} for definitions of $b$ and~$b'$) are bounded with
respect to the bornology~$\CBS_\an$.  This is quite trivial for the
differential~$d$.  The equations \eqref{eq:def_bPrime} and~\eqref{eq:def_b}
imply that $b$ and~$b'$ map $\opt{S}(dS)^n$ into $\sum_{j=0}^n (-1)^j
\opt{T}(dT)^{n-1}$ for $T\defeq S\cup S^2\cup \{0\}$.  Thus if $f(n)\defeq
n+1$, then $f(N)^{-1}\circ b$ and $f(N)^{-1}\circ b'$ map $\opt{S}(dS)^\infty$
into the disked hull of $\opt{T}(dT)^\infty$ and therefore are bounded.
Since~$f$ has only linear growth, the operator~$f(N)$ is bounded, so that $b$
and~$b'$ are bounded.  Hence $d$, $b$, and~$b'$ extend to bounded linear maps
$\Omega_\an\MA \to \Omega_\an\MA$.

\begin{lemma}  \label{lem:Omega_born_alg}
  Let~$\MA$ be a complete bornological algebra.  Then $\Omega_\an\MA$ is a
  complete bornological DG-algebra in the sense that it is a differential
  \Mpn{\Ztwo}graded algebra with a completant bornology for which the grading
  operator, the differential, and the multiplication are bounded.
\end{lemma}

\begin{proof}
  It is clear that the differential and the grading are bounded operators on
  $\Omega_\an \MA$.  It suffices to verify that the multiplication is bounded
  on $(\Omega\MA, \CBS_\an)$ because bounded bilinear maps extend to the
  completion by Lemma~\ref{lem:Cpl_multi}.  Of course, the
  multiplication is the usual multiplication of non-commutative differential
  forms~\cite{connes94:ncg}.

  Since~$b'$ is bounded, $\opt{S} (dS)^\infty \cdot S = b'(\opt{S}
  (dS)^\infty) \in \CBS_\an$ for all $S \in \CBS_c(\MA)$.  Thus
  $\opt{S} (dS)^\infty \cdot \opt{S} (dS)^\infty \in \CBS_\an$ for all $S \in
  \CBS_c(\MA)$.  That is, the multiplication on $(\Omega\MA, \CBS_\an)$ is
  bounded.
\end{proof}

On any differential \Mpn{\Ztwo}graded algebra, we can deform the product to
the \emph{Fedosov product}
$$
x \odot y \defeq x \cdot y - (-1)^{\deg x} dx \cdot dy.
$$
It is bounded if the grading operator, the differential, and the ordinary
multiplication~$\cdot$ are bounded.  A routine computation shows that the
Fedosov product is associative.

\begin{definition}  \label{def:Tanil}
  Let~$\Tanil\MA$ be the even part of $\Omega_\an\MA$ endowed with the
  Fedosov product~$\odot$.  This is called the \emph{analytic tensor algebra
    of~$\MA$}.  Let~$\Janil\MA$ be the kernel of the natural map
  $\tau_{\MA}\colon \Tanil\MA \to \MA$ that maps a form to its degree zero
  part.  That is, $\Janil\MA$ is the closed linear span of $\sum_{j=1}^\infty
  \Omega^{2j} \MA$.
\end{definition}

Lemma~\ref{lem:Omega_born_alg} implies that $\Tanil\MA$ and~$\Janil\MA$
are complete bornological algebras.  Clearly, $\tau_{\MA} \colon \Tanil\MA \to
\MA$ is a bounded homomorphism.  The inclusion of the degree zero part
$\sigma_{\MA} \colon \MA \to \Tanil\MA$ is a natural bounded linear section
for~$\tau_{\MA}$.  Thus $\Janil\MA \injto \Tanil\MA \prto \MA$ is an allowable
extension of~$\MA$.  The maps~$\sigma_{\MA}$ and~$\tau_{\MA}$ will occur
repeatedly.  If no confusion can arise, the subscript~$\MA$ is omitted
and~$\MA$ is considered as a subspace of~$\Tanil\MA$ via~$\sigma_{\MA}$.

\subsection{Properties of analytically nilpotent algebras}
\label{sec:analytically_nilpotent}

In~\ref{def:a_nilpotent} we have defined analytic nilpotence for not
necessarily complete bornological algebras to make the following observation:

\begin{lemma}  \label{lem:nilpotent_complete}
  If\/~$\MA[N]$ is a-nilpotent, so is the completion\/~$\Cpl{\MA[N]}$.
\end{lemma}

\begin{proof}
  By Lemma~\ref{lem:completion_bornology}, any small set $S \subset
  \Cpl{\MA[N]}$ is contained in a set~$\coco{T}$ with $T \in \CBS(\MA[N])$.
  The boundedness of the multiplication in~$\Cpl{\MA[N]}$ implies that
  $(\coco{T})^n \subset \coco{(T^n)}$.  Thus $S^\infty \subset
  (\coco{T})^\infty \subset \coco{(T^\infty)}$ is small in~$\Cpl{\MA[N]}$.
\end{proof}

\begin{lemma}  \label{lem:functional_calculus}
  Let~$\MA[N]$ be a complete a-nilpotent bornological algebra and $x \in
  \MA[N]$.  Let $f = \sum_{j=1}^\infty f_j z^j$ be a power series without
  constant coefficient and positive radius of convergence.  Then
  $\sum_{j=1}^\infty f_j x^j$ converges bornologically in~$\MA[N]$.  We write
  $f(x)$ for the limit of this series.
\end{lemma}

\begin{proof}
  Since~$f$ has positive radius of convergence, we have $|f_j| \le C^j$ for
  some $C \in \left] 0, \infty\right[$.  Let $S \subset \MA[N]$ be the disked
  hull of $\{ 2C\cdot x \}$.  Since~$\MA[N]$ is a-nilpotent, there is a
  completant small disk $T \subset \MA[N]$ containing $S^\infty$.  Thus $2^j
  f_j x^j \in T$ for all~$j$ and hence
  $$
  \biggl\|\sum_{j=N}^\infty f_j x^j \biggr\|_T \le
  \sum_{j=N}^\infty 2^{-j} \|2^j f_j x^j\|_T \le
  \sum_{j=N}^\infty 2^{-j} \le
  2^{-N+1}.
  $$
  Hence $\sum_{j=1}^\infty f_j x^j$ is an absolutely convergent series in
  the Banach space~$\MA[N]_T$.
\end{proof}

Thus in a complete a-nilpotent algebra we have a functional calculus with
germs of holomorphic functions in a neighborhood of the origin.  Roughly
speaking, all elements of an a-nilpotent algebra have spectral radius~$0$.  It
is left to the interested reader to formulate an analogous functional calculus
in infinitely many (non-commuting) variables.  By the way, in a general
complete bornological algebra, for example in $\C[t]$, we only have a
polynomial functional calculus.

We end this section by proving that the class of a-nilpotent algebras is
closed under some constructions.  Most notably, we prove the Extension Axiom
and the Homotopy Axiom.

\begin{lemma}  \label{lem:nilpotent_inherit_triv}
  Analytic nilpotence is inherited by subalgebras, quotients by closed ideals,
  direct products, direct sums, projective limits, and inductive limits.
\end{lemma}

\begin{proof}
  This follows immediately from the definitions.
\end{proof}

\begin{digression}
  We have not made Lemma~\ref{lem:nilpotent_inherit_triv} an axiom because we
  will not need it in this section.  Moreover, some assertions of it fail for
  more general notions of nilpotence, others follow from the axioms.  The
  class of locally nilpotent pro-algebras considered in
  Section~\ref{sec:pro_algebras} is not closed under taking direct sums and
  inductive limits.  The axioms we will meet below imply that products,
  projective limits, closed subalgebras, and quotients by allowable ideals are
  again a-nilpotent.

  The proof for products and projective limits is essentially the same.  Let
  $(\MA[N]_i)_{i \in I}$ be a projective system of a-nilpotent algebras.  The
  identity maps $\Null( \prolim \MA[N]_i) \to \Null(\MA[N]_i) \to \MA[N]_i$
  form a compatible family of \lanilcurs because all~$\MA[N]_i$ are
  a-nilpotent.  Thus we get a compatible family of bounded homomorphisms
  $\Tanil\bigl( \Null (\prolim \MA[N]_i) \bigr) \to \MA[N]_i$ and thus a
  bounded homomorphism into $\prolim \MA[N]_i$.  Its restriction to $\Null
  (\prolim \MA[N]_i)$ is equal to the ``identity map'' and is a \lanilcur
  because it can be extended to a bounded homomorphism on~$\Tanil$.  It
  follows that $\prolim \MA[N]_i$ is a-nilpotent.

  Let~$\MA[N]$ be an a-nilpotent algebra and let $\MA[N]' \subset \MA[N]$ be a
  closed subalgebra.  The identity map $\Null(\MA[N]') \to \MA[N]' \subset
  \MA[N]$ is a \lanilcur because any bounded linear map into the a-nilpotent
  algebra~$\MA[N]$ is a \lanilcur.  Hence it extends to a bounded homomorphism
  $f \colon \Tanil \Null(\MA[N]') \to \MA[N]$.  Since~$\MA[N]'$ is a closed
  subalgebra, the range of this map is contained in~$\MA[N]'$ (when working in
  a more general category, this step requires some care).  Hence the identity
  map $\Null(\MA[N]') \to \MA[N]'$ is a \lanilcur because it can be extended
  to the bounded homomorphism~$f$.

  Let $\pi \colon \MA[N] \prto \MA[N]''$ be a quotient map with a bounded
  linear section~$s$.  If~$\MA[N]$ is a-nilpotent, then so is~$\MA[N]''$.  We
  can factor the ``identity map'' $\Null(\MA[N]'') \to \MA[N]''$ as
  $$
  \Null(\MA[N]'') \overset{s}{\longrightarrow}
  \Null(\MA[N]) \overset{\ID}{\longrightarrow}
  \MA[N] \overset{\pi}{\longrightarrow}
  \MA[N]''.
  $$
  Thus it is a \lanilcur, so that~$\MA[N]''$ is a-nilpotent.
\end{digression}

\begin{axiom}[Extension Axiom]  \label{axiom:extension}
  Let $(i,p) \colon \MA[N]'' \injto \MA[N] \prto \MA[N]'$ be an (allowable)
  extension of complete bornological algebras.  If\/ $\MA[N]''$
  and\/~$\MA[N]'$ are a-nilpotent, then~$\MA[N]$ is a-nilpotent as well.
\end{axiom}

We can prove this statement without using a bounded linear section.  However,
we only need it for allowable extensions.

\begin{proof}
  Let $S \in \CBS(\MA[N])$.  We have to show that~$S^\infty$ is small.  Let
  $S_1 \defeq \bigl( p(4 \cdot S) \bigr)^\infty$.  This set is small
  in~$\MA[N]'$ by a-nilpotence.  Since $p \colon \MA[N] \prto \MA[N]'$ is a
  quotient map, there is $S_2 \in \CBS(\MA[N])$ with $p(S_2) =
  \frac{1}{2}(S_1)$.  Let
  $$
  S_3 \defeq (2S - S_2) \cap\MA[N]'' \in \CBS(\MA[N]'').
  $$
  Then
  $$
  S \subset \tfrac{1}{2} (S_2 + S_3) \subset \bipol{(S_2 \cup S_3)}
  $$
  because for all $x \in S$, there is $\tilde{x} \in S_2$ with
  $p(\tilde{x}) = 2p(x)$.  That is, $2x - \tilde{x} \in \MA[N]''$ and thus $2x
  - \tilde{x} \in S_3$ by definition of~$S_3$.  Let
  $$
  S_4 \defeq S_3 \cup (2S_2^2 - S_2) \cap \MA[N]'' \in \CBS(\MA[N]'').
  $$
  We claim that
  \begin{equation}
    \label{eq:extension_temp}
    S_2^2 \subset \tfrac{1}{2}(S_2+S_4)\subset\bipol{(S_2\cup S_4)}.
  \end{equation}
  To see this, observe that $S_1 = S_1^\infty$ implies
  $$
  p(S_2 \cdot S_2) = p(S_2) \cdot p(S_2) = \tfrac{1}{4}(S_1\cdot S_1)\subset
  \tfrac{1}{4} S_1 = p\bigl(\tfrac{1}{2} S_2\bigr)
  $$
  Thus for all $x \in S_2^2$, there is $\tilde{x} \in S_2$ with
  $p(2x) = p(\tilde{x})$ and therefore $2x - \tilde{x} \in \MA[N]''$ and
  $2x - \tilde{x} \in S_4$ by definition of~$S_4$.  Let
  $$
  S_5 \defeq S_2 \cup (\opt{S_2} \cdot S_4)^\infty \cdot \opt{S_2}, \qquad
  S_6 \defeq \bipol{S_5}.
  $$
  Since~$\MA[N]''$ is an ideal, $\opt{S_2} \cdot S_4 \in \CBS(\MA[N]'')$.
  The sets $S_5, S_6 \subset \MA[N]$ are therefore small because~$\MA[N]''$ is
  a-nilpotent.  Equation~\eqref{eq:extension_temp} implies $S_2^2 \subset
  \bipol{S_5}$, $(\opt{S_2} \cdot S_4)^\infty \opt{S_2} \cdot S_2 \subset
  \bipol{S_5}$, and $S_2 \cdot (\opt{S_2} \cdot S_4)^\infty \opt{S_2} \subset
  \bipol{S_5}$.  Thus $S_6^2 \subset S_6$.  This implies $S_6 = S_6^\infty$.
  Moreover, $S \subset S_6$ because $S \subset \bipol{(S_2\cup S_3)}$.  Thus
  $S^\infty \subset S_6^\infty$ is small.  Consequently, $\MA[N]$ is
  a-nilpotent.
\end{proof}

\begin{definition}  \label{def:tensoring}
  A complete bornological algebra~$\MA[C]$ is called \emph{tensoring} iff
  $\MA[C] \hot \MA[N]$ is a-nilpotent for all a-nilpotent complete
  bornological algebras~$\MA[N]$.
\end{definition}

\begin{lemma}  \label{lem:tensor_a_nilpotent}
  Let~$\MA[C]$ be a complete bornological algebra.  Assume that for all small
  sets $S \in \CBS(\MA[C])$ there is a constant $\lambda > 0$ such that
  $(\lambda \cdot S)^\infty$ is small.  Then~$\MA[C]$ is tensoring.
\end{lemma}

\begin{proof}
  Let~$\MA[N]$ be a-nilpotent.  By Lemma~\ref{lem:nilpotent_complete}, it
  suffices to show that the uncompleted bornological tensor product $\MA[N]
  \otimes \MA[C]$ is a-nilpotent.  Any small subset of $\MA[N] \otimes \MA[C]$
  is contained in a set of the form $\bipol{(S_{\MA[N]} \otimes S_{\MA[C]})}$
  with $S_{\MA[N]} \in \CBS(\MA[N])$, $S_{\MA[C]} \in \CBS(\MA[C])$.  Choose
  $\lambda > 0$ such that $(\lambda \cdot S_{\MA[C]})^\infty$ is small.  Since
  $$
  (S_{\MA[N]} \otimes S_{\MA[C]})^\infty =
  (\lambda^{-1} \cdot S_{\MA[N]})^\infty \otimes
  (\lambda \cdot S_{\MA[C]})^\infty,
  $$
  the set $(S_{\MA[N]} \otimes S_{\MA[C]})^\infty$ is small, and hence so
  is $\bigl( \bipol{(S_{\MA[N]} \otimes S_{\MA[C]})} \bigr)^\infty$.  Thus
  $\MA[N] \otimes \MA[C]$ is a-nilpotent.
\end{proof}

Lemma~\ref{lem:tensor_a_nilpotent} yields immediately that a-nilpotent
algebras and Banach algebras are tensoring.  We will study tensoring algebras
in Section~\ref{sec:tensoring}.  There we prove the converse of
Lemma~\ref{lem:tensor_a_nilpotent} and that a complete bornological algebra is
tensoring iff it is an inductive limit of Banach algebras.  Moreover, a
\Frechet algebra with the precompact bornology is tensoring iff it is an
admissible \Frechet algebra in the sense of~\cite{puschnigg96:asymptotic}.
The algebras $\C[t]$ and $\HO(\C)$ of polynomials and holomorphic functions
on~$\C$ are not tensoring.  It is easy to see that $\C[t] \hot \MA[N]$ is
a-nilpotent iff~$\MA[N]$ is nilpotent (that is, $\MA[N]^k = 0$ for some $k \in
\N$).

\begin{lemma}[Homotopy Axiom]  \label{axiom:homotopy}
  The algebras $\CINF([0,1])$ and $\ABC([0,1])$ are tensoring.
  
  Thus $\CINF([0,1]; \MA[N])$ and $\ABC([0,1]; \MA[N])$ are a-nilpotent
  if\/~$\MA[N]$ is a-nilpotent.
\end{lemma}

\begin{proof}
  Let $S \subset \CINF([0,1])$ be a bounded subset.  Rescale~$S$ such that
  $|f(t)| \le 1/2$ for all $f \in S$, $t \in [0,1]$.  We claim that~$S^\infty$
  is bounded in $\CINF([0,1])$.  This implies that $\CINF([0,1])$ is
  tensoring by Lemma~\ref{lem:tensor_a_nilpotent}.  We have to prove
  that~$S^\infty$ is bounded in the algebras $\NBC[k]([0,1])$ of $k$~times
  continuously differentiable functions for all $k \in \Z_+$.  The Gelfand
  transform for these commutative Banach algebras is the inclusion into
  $\NBC([0,1])$.  Therefore, the spectral radius in $\NBC[k]([0,1])$ is equal
  to the norm in $\NBC([0,1])$ and thus bounded above on~$S$ by~$1/2$.  By
  definition of the spectral radius, the set~$S^\infty$ is norm-bounded in
  $\NBC[k]([0,1])$.
  
  The condition of Lemma~\ref{lem:tensor_a_nilpotent} is trivially verified
  for the Banach algebra $\ABC([0,1])$.  Thus $\ABC([0,1])$ is tensoring.
\end{proof}

\subsection{The interrelations between analytic nilpotence, \lanilcurs, and
  analytic tensor algebras}
\label{sec:triality}

We have now three concepts: The analytic tensor algebra functor~$\Tanil$; the
class of \lanilcurs; and the class of analytically nilpotent algebras.  We
show that they mutually determine each other.

\begin{axiom}  \label{axiom:Tanil_universal}
  Let~$\MA$ be a complete bornological algebra.  The complete bornological
  algebra $\Tanil\MA$ and the bounded linear map $\sigma_{\MA} \colon \MA \to
  \Tanil\MA$  have the following universal property.  If $l \colon \MA \to
  \MA[B]$ is a \lanilcur into a  complete bornological algebra~$\MA[B]$, then
  there is a unique bounded homomorphism $\LLH{l} \colon \Tanil\MA \to \MA[B]$
  extending~$l$ in the sense that $\LLH{l} \circ \sigma_{\MA} = l$.
  Furthermore, any map of the form $f \circ \sigma_{\MA}$ with a bounded
  homomorphism $f \colon \Tanil\MA \to \MA[B]$ is a \lanilcur.
\end{axiom}

\begin{proof}
  By the universal property of completions, bounded homomorphisms $\Tanil\MA
  \to \MA[B]$ are in bijection with homomorphisms $\Tens \MA \to \MA[B]$ that
  map $\opt{S} (dS)^\even$ to a small set in~$\MA[B]$ for all $S\in
  \CBS(\MA)$.  The universal property~\ref{pro:Tens_universal} of $\Tens\MA$
  asserts that bounded homomorphisms $\Tens \MA \to \MA[B]$ are in bijection
  with bounded linear maps $\MA \to \MA[B]$.  In addition, the homomorphism
  $\LLH{l} \colon \Tens\MA \to \MA[B]$ corresponding to a linear map $l \colon
  \MA \to \MA[B]$ satisfies
  \begin{equation}
    \label{eq:Tanil_universal}
    \LLH{l} (\opt{\ma_0} d\ma_1 \dots d\ma_{2n}) =
    l\opt{\ma_0} \cdot
    \omega_l(\ma_1,\ma_2) \cdots \omega_l(\ma_{2n-1},\ma_{2n}).
  \end{equation}
  The sets $\LLH{l} (\opt{S} (dS)^\even) = l\opt{S} \cdot
  \omega_l(S,S)^\infty \subset \MA[B]$ are small for all $S \in \CBS(\MA)$
  iff~$l$ is a \lanilcur.  Consequently, we can extend~$\LLH{l}$ to a bounded
  homomorphism $\Tanil\MA \to \MA[B]$ iff~$l$ is a \lanilcur, and this
  extension is necessarily unique.
\end{proof}

\begin{example}
  Any bounded homomorphism $f \colon \MA \to \MA[B]$ is a \lanilcur.  Since $f
  \circ \tau_{\MA} \circ \sigma_{\MA} = f$, the uniquely determined
  homomorphism $\LLH{f} \colon \Tanil\MA \to \MA[B]$ in
  Axiom~\ref{axiom:Tanil_universal} is equal to $f \circ \tau_{\MA}$.
\end{example}

\begin{example}
  The natural linear section $\sigma_{\MA} \colon \MA \to \Tanil\MA$ is a
  \lanilcur by Axiom~\ref{axiom:Tanil_universal} because we can write it as
  $\ID[{\Tanil\MA}] \circ \sigma_{\MA}$ and the identity map $\ID \colon
  \Tanil\MA \to \Tanil\MA$ is certainly a bounded homomorphism.  The
  maps~$\sigma_{\MA}$ are the most fundamental examples of \lanilcurs.
  Axiom~\ref{axiom:Tanil_universal} asserts that any \lanilcur $l \colon \MA
  \to \MA[B]$ is the composition of a bounded homomorphism $f \colon \Tanil\MA
  \to \MA[B]$ with the \lanilcur~$\sigma_{\MA}$.

  We can easily check directly that~$\sigma_{\MA}$ is a \lanilcur.  Its
  curvature is
  $$
  \omega_\sigma(\ma_1, \ma_2) =
  \ma_1 \cdot \ma_2 - \ma_1 \odot \ma_2 =
  d\ma_1 d\ma_2.
  $$
  The subset $(dS)^\even \subset \Tanil\MA$ is small for all $S \in
  \CBS(\MA)$ by the definition of the bornology $\CBS_\an$.
\end{example}

\begin{axiom}  \label{axiom:Janil_nilpotent}
  The algebra $\Janil\MA$ is a-nilpotent for all complete bornological
  algebras~$\MA$.
\end{axiom}

\begin{proof}
  By Lemma~\ref{lem:nilpotent_complete}, it suffices to show that $(\Jens\MA,
  \odot, \CBS_\an)$ is a-nilpotent because $\Janil\MA$ is the completion of
  this convex bornological algebra.  Thus we have to verify that $T^\infty \in
  \CBS_\an$ for all $T \in \CBS_\an \cap \Jens\MA$.  There is $S \in
  \CBS(\MA)$ such that~$T$ is contained in the disked hull of $\opt{S}(dS)^2
  (dS)^\even \subset (\opt{S} dS)^\even$.  Hence it suffices to prove that
  $(\opt{S} dS)^\infty \in \CBS_\an$.

  We claim that there is $I \in \CBS_\an$ such that $\opt{S} dS \odot I
  \subset \bipol{I}$ and $\opt{S} dS \subset \bipol{I}$ (the letter~$I$ stands
  for ``invariant'').  Since multiplication is bilinear, it follows by
  induction that $(\opt{S} dS)^n \odot \bipol{I} \subset \bipol{I}$ for all $n
  \in \N$ and hence $(\opt{S} dS)^\infty \subset \bipol{I}$.  Thus $(\opt{S}
  dS)^\infty \in \CBS_\an$ as desired.

  Let $S^{(n)} \defeq S \cup S^2 \cup \dots \cup S^n$ and
  $$
  I \defeq \opt{S^{(2)}} \bigcup_{n=1}^\infty (2dS^{(3)})^{n}.
  $$
  Clearly, $\opt{S} dS \subset \bipol{I}$.  Pick $\opt{\ma_0} d\ma_1 \in
  \opt{S} dS$ and $\opt{\ma[b]_0} d\ma[b]_1 \dots d\ma[b]_n \in I$ and
  compute
  \begin{multline*}
    \opt{\ma_0} d\ma_1 \odot
    \opt{\ma[b]_0} d\ma[b]_1 \dots d\ma[b]_n
    \\ =
    \opt{\ma_0} d(\ma_1 \opt{\ma[b]_0}) d\ma[b]_1 \dots d\ma[b]_n
    - \opt{\ma_0} \ma_1 d\opt{\ma[b]_0} d\ma[b]_1 \dots d\ma[b]_n
    - d\opt{\ma_0} d\ma_1 d\opt{\ma[b]_0} d\ma[b]_1 \dots d\ma[b]_n.
  \end{multline*}
  This is a sum of (at most) $3$ terms in $\opt{S^{(2)}} dS^{(2)} dS^{(3)}
  (d2S^{(3)})^{n}$ and thus lies in the disked hull~$\bipol{I}$ of~$I$ as
  desired.
\end{proof}

If~$\VS$ is a complete bornological vector space, then we let $\Null(\VS)$ be
the algebra with underlying vector space~$\VS$ and the zero multiplication.

\begin{axiom}  \label{axiom:nilpotent_lanilcur}
  Let~$\MA[N]$ be a complete bornological algebra.  If the identity map
  $\Null(\MA[N]) \to \MA[N]$ is a \lanilcur, then~$\MA[N]$ is a-nilpotent.
\end{axiom}

\begin{proof}
  The curvature of the identity map $l \colon \Null (\MA[N]) \to \MA[N]$ is
  $\omega_l(x, y) = l(xy) - l(x) l(y) = -xy$ for all $x,y \in \Null (\MA[N])$.
  If this map is a \lanilcur, then $S^\even = (-S^2)^\infty \in \CBS (\MA[N])$
  for all $S \in \CBS_d (\MA[N])$.  Hence $\opt{S} S^\even = S^\infty \in
  \CBS(\MA[N])$.  That is, $\MA[N]$ is a-nilpotent.
\end{proof}

\begin{axiom}  \label{axiom:lanilcur_nilpotent}
  Let $l \colon \MA \to \MA[B]$ be a bounded linear map between complete
  bornological algebras whose curvature factors through an a-nilpotent algebra
  in the following sense.  There are a complete a-nilpotent algebra~$\MA[N]$,
  a bounded linear map $\omega' \colon \MA \hot \MA \to \MA[N]$, and a bounded
  homomorphism $f \colon \MA[N] \to \MA[B]$ such that $\omega_l = f \circ
  \omega'$.  That is, the diagram
  \begin{equation}  \label{eq:nilpotent_lanilcur}
    \begin{gathered}
      \xymatrix{
        {\MA \hot \MA} \ar[r]^-{\omega_l} \ar@{.>}[d]_{\omega'} &
          {\MA[B]} \\
        {\MA[N]} \ar@{-->}[ru]_{f}
        }
    \end{gathered}
  \end{equation}
  commutes.  Then~$l$ is a \lanilcur.
\end{axiom}

\begin{proof}  
  If the curvature of~$l$ factors as in~\eqref{eq:nilpotent_lanilcur}, then
  $\omega_l(S, S)^\infty = f\bigl( \omega'(S \otimes S)^\infty \bigr) \in
  \CBS(\MA[B])$ for all $S \in \CBS(\MA)$ because~$f$ is a homomorphism
  and~$\MA[N]$ is a-nilpotent.  That is, $l$ is a \lanilcur.
\end{proof}

\begin{lemma}  \label{lem:a_nilpotent_test}
  A complete bornological algebra is a-nilpotent iff each bounded linear map
  into it is a \lanilcur.
\end{lemma}

\begin{proof}
  If each map into~$\MA[N]$ is a \lanilcur, then the identity map
  $\Null(\MA[N]) \to \MA[N]$ is a \lanilcur, so that~$\MA[N]$ is a-nilpotent
  by Axiom~\ref{axiom:nilpotent_lanilcur}.  Conversely, if~$\MA[N]$ is
  a-nilpotent, then we can factor the curvature of any bounded linear map $l
  \colon \MA \to \MA[N]$ through the a-nilpotent algebra~$\MA[N]$ by writing
  simply $\omega_l = \ID[{\MA[N]}] \circ \omega_l$.  Thus any bounded linear
  map into~$\MA[N]$ is a \lanilcur by Axiom~\ref{axiom:lanilcur_nilpotent}.
\end{proof}

Therefore, the class of \lanilcurs determines the class of a-nilpotent
algebras.  Conversely, the class of a-nilpotent algebras determines the class
of \lanilcurs:

\begin{lemma}  \label{lem:lanilcur_nilpotent}
  A bounded linear map $l \colon \MA \to \MA[B]$ is a \lanilcur iff its
  curvature can be factored through an a-nilpotent algebra as
  in~\eqref{eq:nilpotent_lanilcur}.
\end{lemma}

\begin{proof}
  By Axiom~\ref{axiom:lanilcur_nilpotent}, $l$ is a \lanilcur if its curvature
  can be factored as in~\eqref{eq:nilpotent_lanilcur}.  Conversely, if~$l$ is
  a \lanilcur, then we have $l = f \circ \sigma_{\MA}$ with a bounded
  homomorphism $f \defeq \LLH{l} \colon \Tanil\MA \to \MA[B]$ by the universal
  property~\ref{axiom:Tanil_universal} of $\Tanil\MA$.  Hence the curvature
  of~$l$ can be factored as $f \circ \omega_\sigma$.  The range
  of~$\omega_\sigma$ is contained in $\Janil\MA$ because $\tau_{\MA} \circ
  \sigma_{\MA} = \ID[{\MA}]$ is a homomorphism.  Thus we can
  view~$\omega_\sigma$ as a bounded linear map $\MA \hot \MA \to \Janil\MA$.
  By Axiom~\ref{axiom:Janil_nilpotent}, $\Janil\MA$ is a-nilpotent.
  Hence~$\omega_l$ factors through an a-nilpotent algebra as asserted.
\end{proof}

\begin{lemma}  \label{lem:lanilcur_compose_hom}
  Let $l \colon \MA \to \MA[B]$ be a \lanilcur, and let $f \colon \MA' \to
  \MA$ and $g \colon \MA[B] \to \MA[B]'$ be bounded homomorphisms.  Then $g
  \circ l$ and $l \circ f$ are \lanilcurs.
  
  If $f \colon \MA \to \MA[B]$ is a bounded homomorphism, there is a unique
  bounded homomorphism $\Tanil f \colon \Tanil\MA \to \Tanil\MA[B]$ such that
  $\Tanil f \circ \sigma_{\MA} = \sigma_{\MA[B]} \circ f$.  That is, the
  diagram
  \begin{equation}  \label{eq:Tanil_functorial_abstract}
    \begin{gathered}
      \xymatrix{
        {\Tanil\MA} \ar@{-->}[r]^{\Tanil f} &
          {\Tanil\MA[B]} \\
        {\MA} \ar[r]_f \ar@{.>}[u]^{\sigma_{\MA}} &
          {\MA[B]} \ar@{.>}[u]^{\sigma_{\MA[B]}}
        }
    \end{gathered}
  \end{equation}
  commutes.  Thus~$\Tanil$ is functorial for bounded homomorphisms.
\end{lemma}

\begin{proof}
  We have $\omega_{g \circ l} = g \circ \omega_l$ and $\omega_{l \circ f} =
  \omega_l \circ (f \hot f)$ because $g$ and~$f$ are multiplicative.  Use
  Lemma~\ref{lem:lanilcur_nilpotent} to factor the curvature of~$l$ through an
  a-nilpotent algebra~$\MA[N]$ as $\omega_l = \phi \circ \omega'$ with a
  bounded homomorphism $\phi \colon \MA[N] \to \MA[B]$ and a bounded linear
  map $\omega' \colon \MA \hot \MA \to \MA[N]$.  Then $\omega_{g \circ l} = (g
  \circ \phi) \circ \omega'$ and $\omega_{l \circ f} = \phi \circ
  \bigl(\omega' \circ (f \hot f) \bigr)$ factor through~$\MA[N]$ as well.
  Thus $g \circ l$ and $l \circ f$ are \lanilcurs by
  Axiom~\ref{axiom:lanilcur_nilpotent}.  In particular, if $f \colon \MA \to
  \MA[B]$ is a bounded homomorphism, then $\sigma_{\MA[B]} \circ f$ is a
  \lanilcur.  By the universal property~\ref{axiom:Tanil_universal}, there is
  a unique bounded homomorphism $\Tanil f \defeq \LLH{\sigma_{\MA[B]} \circ
    f}$ with $\Tanil f \circ \sigma_{\MA} = \sigma_{\MA[B]} \circ f$.  The
  symmetry of the diagram~\eqref{eq:Tanil_functorial_abstract} implies that $f
  \mapsto \Tanil f$ is functorial.
\end{proof}

It is quite easy to obtain the homomorphism $\Tanil f$ from the concrete
description of $\Tanil\MA$:
\begin{equation}
  \label{eq:Tanil_functorial}
  \Tanil f( \opt{\ma_0} d\ma_1 \dots d\ma_{2n}) \defeq
  f\opt{\ma_0} df(\ma_1) \dots df(\ma_{2n}).
\end{equation}

\subsection{Reformulations of the Extension and the Homotopy Axiom}
\label{sec:extension_homotopy_axiom}

We translate the Extension Axiom to equivalent assertions about \lanilcurs and
analytic tensor algebras.

\begin{lemma}[Extension Lemma]  \label{lem:extension}
  The statements (i)---(iii) follow from each other without using the
  Extension or the Homotopy Axiom.  Hence they are all equivalent to the
  Extension Axiom.

  \begin{enumerate}[(i)]
  \item Let $(i,p) \colon \MA[N]'' \injto \MA[N] \prto \MA[N]'$ be an
    allowable extension of complete bornological algebras with a-nilpotent
    kernel and quotient\/ $\MA[N]''$ and\/~$\MA[N]'$.  Then~$\MA[N]$ is
    a-nilpotent.
    
  \item The bounded linear map $\sigma_{\MA}^2 \defeq \sigma_{\Tanil\MA} \circ
    \sigma_{\MA} \colon \MA \to \Tanil\MA \to \Tanil\Tanil\MA$ is a \lanilcur
    for all complete bornological algebras~$\MA$.  The composition of two
    \lanilcurs is again a \lanilcur.  That is, \lanilcurs form a category.
    
  \item The natural projection $\tau_{\Tanil\MA} \colon \Tanil\Tanil\MA \to
    \Tanil\MA$ has a natural bounded splitting homomorphism $\upsilon_{\MA}
    \colon \Tanil\MA \to \Tanil\Tanil\MA$ satisfying $\upsilon_{\MA} \circ
    \sigma_{\MA} = \sigma_{\MA}^2$, namely, $\upsilon_{\MA} \defeq
    \LLH{\sigma_{\MA}^2}$.  Hence the analytic tensor algebra $\Tanil\MA$ is
    a-quasi-free.

  \end{enumerate}
\end{lemma}

Of course, statement~(i) is the Extension Axiom~\ref{axiom:extension}.

\begin{proof}
  \textbf{(i) implies~(ii).}  Let $\tau_{\MA}^2 \defeq \tau_{\MA} \circ
  \tau_{\Tanil\MA} \colon \Tanil\Tanil\MA \to \MA$.  Since~$\tau_{\Tanil\MA}$
  is split surjective with bounded linear section~$\sigma_{\MA}$, we have an
  allowable extension $\Ker \tau_{\Tanil\MA} \injto \Ker \tau_{\MA}^2 \prto
  \Ker \tau_{\MA}$.  By Axiom~\ref{axiom:Janil_nilpotent}, $\Ker
  \tau_{\Tanil\MA} = \Janil(\Tanil\MA)$ and $\Ker \tau_{\MA} = \Janil\MA$ are
  a-nilpotent.  Thus $\MA[N] \defeq \Ker (\tau_{\MA} \circ \tau_{\Tanil\MA})$
  is a-nilpotent by~(i).  Since $\tau_{\MA}^2 \circ \sigma_{\MA}^2 = \ID[\MA]$
  is a homomorphism, the curvature of $\sigma_{\MA}^2$ takes values in
  $\MA[N]$.  That is, it can be factored through the inclusion $\MA[N] \injto
  \Tanil\Tanil\MA$.  Since~$\MA[N]$ is a-nilpotent, $\sigma_{\MA}^2$ is a
  \lanilcur by Axiom~\ref{axiom:lanilcur_nilpotent}.
  
  Let $l \colon \MA \to \MA[B]$ and $l' \colon \MA[B] \to \MA[C]$ be
  \lanilcurs.  By the universal property~\ref{axiom:Tanil_universal} of
  $\Tanil\MA$, we can factor $l = f \circ \sigma_{\MA}$ and $l' = f' \circ
  \sigma_{\MA[B]}$ with bounded homomorphisms $f \defeq \LLH{l} \colon
  \Tanil\MA \to \MA[B]$ and $f' \defeq \LLH{l'} \colon \Tanil\MA[B] \to
  \MA[C]$.  By Lemma~\ref{lem:lanilcur_compose_hom}, we can replace
  $\sigma_{\MA[B]} \circ f$ by $\Tanil f \circ \sigma_{\MA}$.  Thus
  $$
  l' \circ l =
  f' \circ (\sigma_{\MA[B]} \circ f) \circ \sigma_{\MA} =
  f' \circ (\Tanil f \circ \sigma_{\Tanil\MA}) \circ \sigma_{\MA} =
  f' \circ \Tanil f \circ \sigma_{\MA}^2.
  $$
  Since $f' \circ \Tanil f$ is a bounded homomorphism and $\sigma_{\MA}^2$
  is a \lanilcur, $l' \circ l$ is a \lanilcur by
  Lemma~\ref{lem:lanilcur_compose_hom}.  Consequently, \lanilcurs form a
  category.
  
  \textbf{(ii) implies~(iii).}  Since~$\sigma_{\MA}^2$ is a \lanilcur, the
  universal property~\ref{axiom:Tanil_universal} of $\Tanil\MA$ yields a
  bounded homomorphism $\upsilon_{\MA} \defeq \LLH{\sigma_{\MA}^2} \colon
  \Tanil\MA \to \Tanil\Tanil\MA$.  By construction,
  $$
  (\tau_{\Tanil\MA} \circ \upsilon_{\MA}) \circ \sigma_{\MA} =
  \tau_{\Tanil\MA} \circ \sigma_{\MA}^2 =
  \sigma_{\MA} =
  \ID[\Tanil\MA] \circ \sigma_{\MA}.
  $$
  The uniqueness part of Axiom~\ref{axiom:Tanil_universal} implies
  $\tau_{\Tanil\MA} \circ \upsilon_{\MA} = \ID[\Tanil\MA]$.
  Hence~$\upsilon_{\MA}$ is a section for~$\tau_{\Tanil\MA}$.  Thus
  $\Tanil\MA$ is a-quasi-free in the sense of
  Definition~\ref{def:a_quasi_free}.

  \textbf{(iii) implies~(ii).}  Assume that there is a bounded
  homomorphism $\upsilon \colon \Tanil\MA \to \Tanil\Tanil\MA$ satisfying
  $\upsilon \circ \sigma_{\MA} = \sigma_{\MA}^2$.  The universal
  property~\ref{axiom:Tanil_universal} of $\Tanil\MA$ implies that
  $\sigma_{\MA}^2$ is a \lanilcur.

  \begin{digression}
    It is unclear whether the quasi-freeness of $\Tanil\MA$ alone, that is,
    the existence of \textsl{any} bounded splitting homomorphism $\upsilon
    \colon \Tanil\MA \to \Tanil\Tanil\MA$ for~$\tau_{\Tanil\MA}$, suffices to
    prove~(ii).
  \end{digression}
  
  \textbf{(ii) implies~(i).}  Let $(i, p) \colon \MA[N]' \injto \MA[N] \prto
  \MA[N]''$ be an \emph{allowable} extension of complete bornological algebras
  in which $\MA[N]'$ and~$\MA[N]''$ are a-nilpotent.  We choose a bounded
  linear section $s \colon \MA[N]'' \to \MA[N]$ of our allowable extension.
  To verify that~$\MA[N]$ is a-nilpotent, we have to verify that any bounded
  linear map $l \colon \MA[B] \to \MA[N]$ is a \lanilcur by
  Axiom~\ref{axiom:lanilcur_nilpotent}.  Since~$\MA[N]''$ is a-nilpotent, $p
  \circ l$ is a \lanilcur by Lemma~\ref{lem:lanilcur_nilpotent}.  By the
  universal property of $\Tanil\MA$ we can extend it to a bounded homomorphism
  $f \defeq \LLH{p \circ l} \colon \Tanil\MA[B] \to \MA[N]''$.  We consider
  the bounded linear map $l' \colon \Tanil\MA[B] \to \MA[N]$ defined by
  $$
  l' = s \circ f + (\ID - s \circ p) \circ l \circ \tau_{\MA[B]}.
  $$
  By adding $(\ID - s \circ p) \circ l \circ \tau_{\MA[B]}$, we achieve
  that $l' \circ \sigma_{\MA[B]} = l$.  Since $p \circ s = \ID$, the
  composition $p \circ l' = f$ is a homomorphism.  Hence the curvature of~$l'$
  factors through the inclusion $i \colon \MA[N]' \to \MA[N]$.  By assumption,
  $\MA[N]'$ is a-nilpotent and hence~$l'$ is a \lanilcur by
  Axiom~\ref{axiom:lanilcur_nilpotent}.  Since~$\sigma_{\MA[B]}$ is a
  \lanilcur as well, the composition $l = l' \circ \sigma_{\MA[B]}$ is a
  \lanilcur by~(ii).  Thus~(i) follows from~(ii).
\end{proof}

\begin{lemma}[Tensor Lemma]  \label{lem:tensor}
  Let~$\MA[C]$ be a complete bornological algebra.  The assertions (i)---(iii)
  follow from each other without using the Extension or the Homotopy Axiom and
  provide alternative definitions for tensoring algebras.

  \begin{enumerate}[(i)]
  \item If\/~$\MA[N]$ is an a-nilpotent complete bornological algebra, so is
    $\MA[N] \hot \MA[C]$.

  \item If $l \colon \MA \to \MA[B]$ is a \lanilcur, then so is $l \hot
    \ID[{\MA[C]}] \colon \MA \hot \MA[C] \to \MA[B] \hot \MA[C]$.
    
  \item If $f \colon \MA \to \MA[B] \hot \MA[C]$ is a bounded homomorphism,
    then there is a unique bounded homomorphism $g \colon \Tanil\MA \to
    (\Tanil\MA[B]) \hot \MA[C]$ that makes the following diagram commute:
    \begin{equation}  \label{eq:tensor_sigma}
      \begin{gathered}
        \xymatrix{
          {\Tanil\MA} \ar@{-->}[r]^-{g} &
            {(\Tanil\MA[B]) \hot \MA[C]} \\
          {\MA} \ar@{.>}[u]^{\sigma_{\MA}} \ar[r]_-{f} &
            {\MA[B] \hot \MA[C]}
            \ar@{.>}[u]_{\sigma_{\MA[B]} \hot \ID[{\MA[C]}]}
          }
      \end{gathered}
    \end{equation}
  \end{enumerate}
\end{lemma}

\begin{proof}
  \textbf{(i) implies~(ii).}  Let $l \colon \MA \to \MA[B]$ be a \lanilcur.
  By Lemma~\ref{lem:lanilcur_nilpotent}, we can factor the curvature of~$l$
  through an a-nilpotent algebra~$\MA[N]$.  That is, $\omega_l = f \circ
  \omega'$ with a bounded linear map $\omega' \colon \MA \hot \MA \to \MA[N]$
  and a bounded homomorphism $f \colon \MA[N] \to \MA[B]$.  Let $m \colon
  \MA[B] \hot \MA[B] \to \MA[B]$ be the multiplication map.  We have
  $\omega_{l \otimes \ID}( \ma_1 \otimes \ma[c]_1, \ma_2 \otimes \ma[c]_2) =
  \omega_l(\ma_1, \ma_2) \otimes (\ma[c]_1 \cdot \ma[c]_2)$ and thus
  $\omega_{l \hot \ID} = (f \hot \ID[{\MA[C]}]) \circ (\omega' \hot m)$.
  Hence the curvature of $l \hot \ID[{\MA[C]}]$ factors through $\MA[N] \hot
  \MA[C]$.  Since $\MA[N] \hot \MA[C]$ is a-nilpotent by~(i), $l \hot
  \ID[{\MA[C]}]$ is a \lanilcur by Axiom~\ref{axiom:lanilcur_nilpotent}.
  
  \textbf{(ii) implies~(iii).}  The linear map $\sigma_{\MA[B]} \colon \MA[B]
  \to \Tanil\MA[B]$ is a \lanilcur.  Hence $\sigma_{\MA[B]} \hot \ID[{\MA[C]}]
  \colon \MA[B] \hot \MA[C] \to (\Tanil\MA[B]) \hot \MA[C]$ is a \lanilcur
  by~(ii).  By Lemma~\ref{lem:lanilcur_compose_hom}, the composition
  $(\sigma_{\MA[B]} \hot \ID[{\MA[C]}]) \circ f \colon \MA \to (\Tanil\MA[B])
  \hot \MA[C]$ is a \lanilcur.  The universal
  property~\ref{axiom:Tanil_universal} of $\Tanil\MA$ yields that there is a
  unique bounded homomorphism~$g$ making the diagram~\eqref{eq:tensor_sigma}
  commute.
  
  \textbf{(iii) implies~(ii).}  If we apply~(iii) to the bounded homomorphism
  $\ID \colon \MA \hot \MA[C] \to \MA \hot \MA[C]$, we obtain that
  $\sigma_{\MA} \hot \ID[{\MA[C]}] \colon \MA \hot \MA[C] \to (\Tanil\MA) \hot
  \MA[C]$ is a \lanilcur because it can be extended to a bounded homomorphism
  $\Tanil(\MA \hot \MA[C]) \to (\Tanil\MA) \hot \MA[C]$.  A \lanilcur $l
  \colon \MA \to \MA[B]$ can be extended to a bounded homomorphism $\LLH{l}
  \colon \Tanil\MA \to \MA[B]$ by the universal
  property~\ref{axiom:Tanil_universal} of $\Tanil\MA$.  Tensoring with
  $\ID[{\MA[C]}]$, we get a bounded homomorphism $f \colon (\Tanil\MA) \hot
  \MA[C] \to \MA[B] \hot \MA[C]$ satisfying $f \circ (\sigma_{\MA} \hot
  \ID[{\MA[C]}]) = l \hot \ID[{\MA[C]}]$.  Since $\sigma_{\MA} \hot
  \ID[{\MA[C]}]$ is a \lanilcur, so is its composition with the bounded
  homomorphism~$f$ by Lemma~\ref{lem:lanilcur_compose_hom}.  Thus $l \hot
  \ID[{\MA[C]}]$ is a \lanilcur as desired.

  \textbf{(ii) implies~(i).}  Let~$\MA[N]$ be an a-nilpotent complete
  bornological algebra.  We have to show that $\MA[N] \hot \MA[C]$ is
  a-nilpotent.  By Lemma~\ref{lem:a_nilpotent_test}, the identity map
  $\Null(\MA[N]) \to \MA[N]$ is a \lanilcur.  Thus the identity map
  $\Null(\MA[N] \hot \MA[C]) \cong \Null(\MA[N]) \hot \MA[C] \to \MA[N] \hot
  \MA[C]$ is a \lanilcur as well by~(ii).  By
  Axiom~\ref{axiom:nilpotent_lanilcur}, this implies that $\MA[N] \hot \MA[C]$
  is a-nilpotent.
\end{proof}

Smooth and absolutely continuous homotopies are defined as bounded
homomorphisms $\MA \to \MA[B] \hot \MA[C]$ with~$\MA[C]$ either $\CINF([0,1])$
or $\ABC([0,1])$.  Thus~\eqref{eq:tensor_sigma} describes the lifting of
homotopies to analytic tensor algebras.  The Homotopy
Axiom~\ref{axiom:homotopy} and the Tensor Lemma~\ref{lem:tensor} imply:

\begin{lemma}  \label{lem:Tanil_homotopy}
  Let $H \colon \MA \to \ABC([0,1]; \MA[B])$ be an absolutely continuous
  homotopy.  Then there is a unique absolutely continuous homotopy $\hat{H}
  \colon \Tanil\MA \to \ABC([0,1]; \Tanil\MA[B])$ such that $\ev_t \circ
  \hat{H} = \Tanil (\ev_t \circ H)$ for all $t \in [0,1]$.  Thus the
  functor~$\Tanil$ preserves absolutely continuous homotopies.  Similarly,
  $\Tanil$ preserves smooth homotopies.
\end{lemma}

\subsection{Universal analytically nilpotent extensions}
\label{sec:universal_ext}

We introduce some natural terminology to handle extensions.  Let $\MA[K]_i
\injto \MA[E]_i \prto \MA[Q]_i$, $i=0,1$, be extensions of complete
bornological algebras.  A \emph{morphism of extensions} is a triple of bounded
homomorphisms $(\xi, \psi, \phi)$ making the diagram
$$
\xymatrix{
  {\MA[K]_0\;} \ar@{>->}[r] \ar[d]^{\xi} &
    {\MA[E]_0} \ar@{->>}[r] \ar[d]^{\psi} &
      {\MA[Q]_0} \ar[d]^{\phi} \\
  {\MA[K]_1\;} \ar@{>->}[r] &
    {\MA[E]_1} \ar@{->>}[r] &
      {\MA[Q]_1} 
  }
$$
commute.  Of course, a morphism of extensions is determined by the
map~$\psi$ in the middle.  A \emph{smooth homotopy} is a morphism of
extensions $(\Xi, \Psi, \Phi)$ from $\MA[K]_0 \injto \MA[E]_0 \prto \MA[Q]_0$
to $\CINF([0,1]; \MA[K]_1) \injto \CINF([0,1]; \MA[E]_1) \prto \CINF([0,1];
\MA[Q]_1)$.  A smooth homotopy $(\Xi, \Psi, \Phi)$ is called a \emph{smooth
  homotopy relative to~$\phi$} iff $\Phi_t = \phi$ for all $t\in[0,1]$.  If
$\MA[Q]_0 = \MA[Q]_1 = \MA[Q]$, a smooth homotopy relative to~$\ID[{\MA[Q]}]$
is called a \emph{smooth homotopy relative to~$\MA[Q]$}.  If we replace
$\CINF([0,1])$ by $\ABC([0,1])$, we obtain \emph{absolutely continuous
  homotopies} and \emph{relative absolutely continuous homotopies}.

\begin{theorem}[Universal Extension Theorem]  \label{the:universal_ext}
  Let $(\iota, \pi) \colon \MA[N] \injto \MA[R] \prto \MA$ be a universal
  a-nilpotent extension.  Let $(j, p) \colon \MA[K] \injto \MA[E] \prto
  \MA[Q]$ be an a-nilpotent extension.  Let $\phi \colon \MA \to \MA[Q]$ be a
  bounded homomorphism.  Then there is a morphism of extensions
  $$
  \xymatrix{
    {\MA[N]\;} \ar@{>->}[r]^{\iota} \ar@{-->}[d]^{\xi} &
      {\MA[R]} \ar@{->>}[r]^{\pi} \ar@{-->}[d]^{\psi} &
        {\MA} \ar[d]^{\phi} \\
    {\MA[K]\;} \ar@{>->}[r]^{j} &
      {\MA[E]} \ar@{->>}[r]^{p} &
        {\MA[Q]} 
    }
  $$
  lifting~$\phi$.  This lifting is unique up to smooth homotopy relative
  to~$\phi$.
  
  More generally, let $(\xi_t, \psi_t, \phi_t)$, $t = 0,1$, be morphisms of
  extensions as above and let $\Phi \colon \MA \to \CINF([0,1]; \MA[Q])$ be a
  smooth homotopy between $\phi_0$ and~$\phi_1$.  Then~$\Phi$ can be lifted to
  a smooth homotopy $(\Xi, \Psi, \Phi)$ between $(\xi_0, \psi_0, \phi_0)$ and
  $(\xi_1, \psi_1, \phi_1)$.  An analogous statement holds for absolutely
  continuous homotopies instead of smooth homotopies.
\end{theorem}

\begin{proof}
  Let $\upsilon \colon \MA[R] \to \Tanil\MA[R]$ be a splitting homomorphism
  for~$\tau_{\MA[R]}$.  Let $s \colon \MA[Q] \to \MA[E]$ be a bounded linear
  section for~$p$.  Since $p \circ (s \circ \phi \circ \pi) = \phi \circ \pi$
  is a homomorphism, the curvature of the bounded linear map $s \circ \phi
  \circ \pi \colon \MA[R] \to \MA[E]$ has values in~$\MA[K]$.  Since~$\MA[K]$
  is a-nilpotent by assumption, $s \circ \phi \circ \pi$ is a \lanilcur by
  Axiom~\ref{axiom:lanilcur_nilpotent}.  The universal
  property~\ref{axiom:Tanil_universal} of $\Tanil\MA[R]$ yields an associated
  bounded homomorphism $\hat{\psi} \defeq \LLH{s \circ \phi \circ \pi} \colon
  \Tanil\MA[R] \to \MA[E]$.  Define $\psi \defeq \hat{\psi} \circ \upsilon
  \colon \MA[R] \to \MA[E]$.  By construction,
  $$
  (p \circ \hat{\psi}) \circ \sigma_{\MA[R]} =
  p \circ s \circ \phi \circ \pi =
  \phi \circ \pi =
  (\phi \circ \pi \circ \tau_{\MA[R]}) \circ \sigma_{\MA[R]}.
  $$
  The uniqueness part of Axiom~\ref{axiom:Tanil_universal} yields $p \circ
  \hat{\psi} = \phi \circ \pi \circ \tau_{\MA[R]}$ and hence $p \circ \psi =
  \phi \circ \pi \circ \tau_{\MA[R]} \circ \upsilon = \phi \circ \pi$ as
  desired.  Thus $\psi(\MA[N]) \subset \MA[K]$, so that there is a
  unique~$\xi$ making $(\xi, \psi, \phi)$ a morphism of extensions.
  
  The uniqueness of the lifting $(\xi, \psi, \phi)$ up to smooth homotopy
  relative to~$\phi$ follows from the more general statement about the lifting
  of homotopies.  We use the bounded homomorphism $p_\ast \defeq \CINF([0,1];
  p) \colon \CINF([0,1]; \MA[E]) \to \CINF([0,1]; \MA[Q])$ and the bounded
  linear map $s_\ast$ defined similarly.  Since $\Phi_t \circ \pi = p \circ
  \psi_t$ for $t = 0,1$, there is a bounded linear map $l \colon \MA[R] \to
  \CINF([0,1]; \MA[E])$ satisfying $\ev_t \circ l = \psi_t$ for $t = 0,1$ and
  $p_\ast \circ l = \Phi \circ \pi$.  For instance, we can take
  $$
  l \defeq
  s_\ast \circ \Phi \circ \pi
  + (1-t) \otimes (\psi_0 - s \circ \Phi_0 \circ \pi)
  + t \otimes (\psi_1 - s \circ \Phi_1 \circ \pi)
  $$
  
  Since $p_\ast \circ l = \Phi \circ \pi$ is a homomorphism, the curvature
  of~$l$ has values in $\CINF([0,1]; \MA[K])$.  This algebra is a-nilpotent by
  the Homotopy Axiom~\ref{axiom:homotopy}.  Thus~$l$ is a \lanilcur by
  Axiom~\ref{axiom:lanilcur_nilpotent}.  Let $\hat{\Psi} \defeq \LLH{l} \colon
  \Tanil\MA[R] \to \CINF([0,1]; \MA[E])$ be the associated homomorphism by the
  universal property~\ref{axiom:Tanil_universal}.  Define $\Psi \defeq
  \hat{\Psi} \circ \upsilon$.  By construction, $p_\ast \circ \hat{\Psi} \circ
  \sigma_{\MA[R]} = p_\ast \circ l = \Phi \circ \pi = (\Phi \circ \pi\circ
  \tau_{\MA[R]}) \circ \sigma_{\MA[R]}$ and thus $p_\ast \circ \hat{\Psi} =
  \Phi \circ \pi\circ \tau_{\MA[R]}$ by the uniqueness assertion of
  Axiom~\ref{axiom:Tanil_universal}.  Therefore, $p_\ast \circ \Psi = \Phi
  \circ \pi \circ \tau_{\MA[R]} \circ \upsilon = \Phi \circ \pi$.  A similar
  computation shows that $\Psi_t = \ev_t \circ \Psi = \psi_t$ for $t = 0,1$.
  
  The above argument used no other information about smooth homotopies than
  the Homotopy Axiom, stating that $\CINF([0,1])$ is tensoring.  Since
  $\ABC([0,1])$ is tensoring as well, the same argument shows that we can lift
  absolutely continuous homotopies to absolutely continuous homotopies between
  morphisms of extensions.
\end{proof}

\begin{theorem}[Uniqueness Theorem]  \label{the:uniqueness}
  Let~$\MA$ be a complete bornological algebra.  Then $\Janil\MA \injto
  \Tanil\MA \prto \MA$ is a universal a-nilpotent extension.  Up to smooth
  homotopy equivalence of extensions relative to~$\MA$, this is the only
  universal a-nilpotent extension of~$\MA$.  That is, any universal
  a-nilpotent extension $\MA[N] \injto \MA[R] \prto \MA$ is smoothly homotopy
  equivalent relative to~$\MA$ to $\Janil\MA \injto \Tanil\MA \prto \MA$.
  Especially, $\MA[R]$ is smoothly homotopy equivalent to $\Tanil\MA$
  and~$\MA[N]$ is smoothly homotopy equivalent to $\Janil\MA$.

  In addition, any morphism of extensions $(\xi, \psi, \ID)$ between two
  universal a-nilpotent extensions is a smooth homotopy equivalence of
  extensions relative to~$\MA$.  Any two morphisms of extensions $(\xi, \psi,
  \ID)$ are smoothly homotopic.
\end{theorem}

\begin{proof}
  The algebra $\Janil\MA$ is a-nilpotent by Axiom~\ref{axiom:Janil_nilpotent}
  and $\Tanil\MA$ is a-quasi-free by the Extension
  Lemma~\ref{lem:extension}.(iii).  Thus $\Janil\MA \injto \Tanil\MA \prto
  \MA$ is a universal a-nilpotent extension.

  Let $E_i \colon \MA[N]_i \injto \MA[R]_i \prto \MA$, $i = 0,1$, be two
  universal a-nilpotent extensions of~$\MA$ (we write~$E_i$ as an abbreviation
  for these extensions).  By the existence part of the Universal Extension
  Theorem~\ref{the:universal_ext}, the identity map $\ID[{\MA}] \colon \MA \to
  \MA$ can be lifted to morphisms of extensions $(\xi_0, \psi_0, \ID) \colon
  E_0 \to E_1$ and $(\xi_1, \psi_1, \ID) \colon E_1 \to E_0$.  The composition
  $(\xi_1, \psi_1, \ID) \circ (\xi_0, \psi_0, \ID) = (\xi_1 \circ \xi_0,
  \psi_1 \circ \psi_0, \ID)$ is a lifting of the identity map $\ID[{\MA}]$ to
  an automorphism of~$E_0$.  Another such lifting is the automorphism
  $(\ID[{\MA[N]_0}], \ID[{\MA[R]_0}], \ID[{\MA}])$.  By the uniqueness part of
  the Universal Extension Theorem~\ref{the:universal_ext}, these two liftings
  are smoothly homotopic relative to~$\MA$.  An analogous argument shows that
  the composition $(\xi_0, \psi_0, \ID) \circ (\xi_1, \psi_1, \ID) \colon E_1
  \to E_1$ is smoothly homotopic to the identity automorphism relative
  to~$\MA$.  This means that $E_0$ and~$E_1$ are smoothly homotopy equivalent
  relative to~$\MA$.  Especially, $\MA[N]_0$ and~$\MA[N]_1$ are smoothly
  homotopy equivalent algebras and $\MA[R]_0$ and~$\MA[R]_1$ are smoothly
  homotopy equivalent algebras.  In addition, any morphism of extensions
  between $E_0$ and $E_1$ of the form $(\xi, \psi, \ID)$ is a smooth homotopy
  equivalence.  In fact, by the uniqueness part of the Universal Extension
  Theorem~\ref{the:universal_ext}, all morphisms of this form are smoothly
  homotopic relative to~$\MA$.
\end{proof}

\begin{theorem}  \label{the:aqf_definitions}
  Let~$\MA[R]$ be a complete bornological algebra.  The following conditions
  are equivalent to~$\MA[R]$ being analytically quasi-free:
  \begin{enumerate}[(i)]%
    
  \item There is a bounded homomorphism $\upsilon \colon \MA[R] \to
    \Tanil\MA[R]$ that is a section for $\tau_{\MA[R]} \colon \Tanil\MA[R] \to
    \MA[R]$.
    
  \item All allowable a-nilpotent extensions $\MA[N] \injto \MA[E] \prto
    \MA[R]$ have a bounded splitting homomorphism.
    
  \item Let $(\iota, \pi) \colon \MA[N] \injto \MA[E] \prto \MA[Q]$ be an
    a-nilpotent extension and let $\phi \colon \MA[R] \to \MA[Q]$ be a bounded
    homomorphism.  Then there is a bounded homomorphism $\psi \colon \MA[R]
    \to \MA[E]$ lifting~$\phi$, that is, making the diagram
    $$
    \xymatrix{
      &
        {\MA[R]} \ar@{=}[r] \ar@{-->}[d]^{\psi} &
          {\MA[R]} \ar[d]^{\phi} \\
      {\MA[N]\:} \ar@{>->}[r]^{\iota} &
        {\MA[E]} \ar@{->>}[r]^{\pi} &
          {\MA[Q]}
      }
    $$
    commute.
  \end{enumerate}

  Furthermore, if\/~$\MA[R]$ is a-quasi-free, then we get some additional
  informations in the situations of (i) and~(iii).  In~(i), we get that
  $\Janil\MA[R]$ is smoothly contractible and that~$\MA[R]$ is a deformation
  retract of $\Tanil\MA[R]$.  Any two bounded splitting homomorphisms $\MA[R]
  \to \Tanil\MA[R]$ are smoothly homotopic.  In~(iii), we get that the
  lifting~$\psi$ is unique up to smooth homotopy.
\end{theorem}

\begin{proof}
  Assertion~(i) is the definition of analytic quasi-freeness.  If~(i) holds,
  then $0 \injto \MA[R] \prto \MA[R]$ is a universal a-nilpotent extension
  of~$\MA[R]$ because the zero algebra is certainly a-nilpotent.  If we apply
  the Universal Extension Theorem~\ref{the:universal_ext} to this universal
  extension, we get~(iii) with the additional information that the
  lifting~$\psi$ is unique up to smooth homotopy.
  
  If we apply~(iii) to the extension $\MA[N] \injto \MA[E] \prto \MA[R]$, we
  get~(ii).  If we apply~(ii) to the a-nilpotent extension $\Janil\MA[R]
  \injto \Tanil\MA[R] \prto \MA[R]$, then we get~(i).  Thus (i)--(iii) are
  equivalent.  Finally, if~$\MA[R]$ is a-quasi-free, then the Uniqueness
  Theorem~\ref{the:uniqueness} implies that the two universal a-nilpotent
  extensions $0 \injto \MA[R] \prto \MA[R]$ and $\Janil\MA[R] \injto
  \Tanil\MA[R] \prto \MA[R]$ are smoothly homotopy equivalent relative
  to~$\MA[R]$.  Thus~$\MA[R]$ is a smooth deformation retract of
  $\Tanil\MA[R]$ and $\Janil\MA[R]$ is smoothly contractible.  In addition,
  the Uniqueness Theorem asserts that all bounded splitting homomorphisms
  $\MA[R] \to \Tanil\MA[R]$ are smoothly homotopic.
\end{proof}

\begin{corollary}  \label{cor:aqf_qf}
  Analytically quasi-free algebras are quasi-free.
\end{corollary}

\begin{proof}
  Let~$\MA[R]$ be a-quasi-free.  By Theorem~\ref{the:aqf_definitions}, all
  a-nilpotent extensions of~$\MA[R]$ split by a bounded homomorphism.  Since
  \Mpn{2}nilpotent algebras are a-nilpotent by
  Axiom~\ref{axiom:nilpotent_lanilcur}, $\MA[R]$ satisfies condition (vii) of
  Definition~\ref{deflem:quasi_free}.  That is, $\MA[R]$ is quasi-free.
\end{proof}

\subsection{The bimodule $\Omega^1(\Tanil\MA)$}
\label{sec:Omega_I_Tanil}

By Corollary~\ref{cor:aqf_qf}, a-quasi-free algebras are always quasi-free.
By Lemma~\ref{deflem:quasi_free}, it follows that $\Omega^1(\MA[R])$ is
a projective bimodule.  In particular, $\Omega^1(\Tanil\MA)$ is a projective
\Mp{\Tanil\MA}bimodule.  We show even that $\Omega^1(\Tanil\MA)$ is isomorphic
to the free \Mp{\Tanil\MA}bimodule on~$\MA$.  A concrete proof of this
isomorphism is given in Appendix~\ref{app:X_Tanil}.  Here we give a rather
abstract proof based on an idea of Cuntz and Quillen
\cite[Proposition~2.6]{cuntz95:algebra}.

Let~$\MA[B]$ be a complete bornological algebra and let~$\VS$ be a
\Mpn{\MA[B]}bimodule.  The \emph{semi-direct product} $\MA[B] \ltimes \VS$ is
the following complete bornological algebra.  As a bornological vector space,
$\MA[B] \ltimes \VS \cong \MA[B] \oplus \VS$.  The multiplication is
$$
(\ma[b], \vs) \cdot (\ma[b]', \vs') \defeq
(\ma[b] \cdot \ma[b]', \vs \cdot \ma[b]' + \ma[b] \cdot \vs').
$$
In this definition, we used the \Mpn{\MA[B]}bimodule structure to define
products $\VS \cdot \MA[B]$ and $\MA[B] \cdot \VS$, and the multiplication in
the algebra~$\MA[B]$.  It is straightforward to verify that $\MA[B] \ltimes
\VS$ is a complete bornological algebra that fits naturally into a split
extension of complete bornological algebras
$$
\Null(\VS) \injto \MA[B] \ltimes \VS \prto \MA[B].
$$
That is, the summand $\VS \subset \MA[B] \ltimes \VS$ is an ideal with $\VS
\cdot \VS = 0$.  The summand $\MA[B] \subset \MA[B] \ltimes \VS$ is a
subalgebra and the natural projection $\MA[B] \ltimes \VS \to \MA[B]$ is a
bounded homomorphism $\MA[B] \ltimes \VS \to \MA[B]$.  The relevance of
$\MA[B] \ltimes \VS$ is that bimodule homomorphisms $\Omega^1(\MA[B]) \to \VS$
can be converted into splitting homomorphisms $\MA[B] \to \MA[B] \ltimes \VS$
and vice versa.

\begin{lemma}  \label{lem:Omega_I_ltimes}
  Let~$\MA[B]$ be a complete bornological algebra and let~$\VS$ be a
  \Mpn{\MA[B]}bimodule.  The bimodule homomorphisms $\Omega^1(\MA[B]) \to \VS$
  are those linear maps of the form $\opt{x} dy \mapsto \opt{x} \cdot
  \delta(y)$ for all $x, y \in \MA[B]$ with a bounded derivation $\delta
  \colon \MA[B] \to \VS$.
  
  The bounded splitting homomorphisms for the natural projection $\MA[B]
  \ltimes \VS \prto \MA[B]$ are the maps $\MA[B] \to \MA[B] \ltimes \VS \cong
  \MA[B] \oplus \VS$ of the form $(\ID[{\MA[B]}], \delta)$ with a bounded
  derivation $\delta \colon \MA[B] \to \VS$.
\end{lemma}

\begin{proof}
  Since $\Omega^1(\MA[B]) \cong \Unse{\MA[B]} \hot \MA[B]$ is a free left
  \Mpn{\MA[B]}module, the left \Mpn{\MA[B]}module homomorphisms
  $\Omega^1(\MA[B]) \to \VS$ are those linear maps of the form $\hat{\delta}
  \colon \opt{x} dy \mapsto \opt{x} \cdot \delta(y)$ with arbitrary bounded
  linear maps $\delta \colon \MA[B] \to \VS$.  The map~$\hat{\delta}$ is even
  a bimodule homomorphism iff $\delta(xy) = x\delta(y) + \delta(x)y$ for all
  $x,y \in \MA[B]$, that is, $\delta$ is a derivation.  The bounded
  \emph{linear} sections $s \colon \MA[B] \to \MA[B] \ltimes \VS$ are those
  maps of the form $(\ID[{\MA[B]}], \delta)$ with an arbitrary bounded linear
  map $\delta \colon \MA[B] \to \VS$.  The multiplication in $\MA[B] \ltimes
  \VS$ is defined so that~$s$ is multiplicative iff~$\delta$ is a derivation.
\end{proof}

Before we can talk about $\Omega^1(\Tanil\MA)$, we have to solve a notational
problem: The letter~$d$ is already used for the differential in $\Omega\MA
\supset \Tanil\MA$.  Therefore, we use a capital~$D$ for the differential in
$\Omega(\Tanil\MA)$.

\begin{proposition}  \label{pro:Omega_I_Tanil}
  Let~$\MA$ be a complete bornological algebra.  Then the following bounded
  linear maps are bornological isomorphisms:
  \begin{alignat}{2}
  \label{eq:mu3}
    \mu_{\ref{eq:mu3}} & \colon
    \Unse{(\Tanil\MA)} \hot \MA \hot \Unse{(\Tanil\MA)} \congto
    \Omega^1(\Tanil\MA), & \qquad
    x \otimes \ma \otimes y & \mapsto
    x \bigl( D\sigma_{\MA}(\ma) \bigr) y,
  \\
  \label{eq:mu1}
    \mu_{\ref{eq:mu1}} & \colon
    \Unse{(\Tanil\MA)} \hot \MA \congto
    \Tanil\MA, & \qquad
    x \otimes \ma & \mapsto x \odot \sigma_{\MA}(\ma),
  \\
  \label{eq:mu2}
    \mu_{\ref{eq:mu2}} & \colon
    \MA \hot \Unse{(\Tanil\MA)} \congto
    \Tanil\MA, & \qquad
    \ma \otimes x & \mapsto \sigma_{\MA}(\ma) \odot x.
  \end{alignat}
  Thus $\Tanil\MA$ is free as a left and right \Mpn{\Tanil\MA}module on the
  subspace $\MA \cong \sigma_{\MA}(\MA)$ and $\Omega^1(\Tanil\MA)$ is free as
  a \Mp{\Tanil\MA}bimodule on the subspace $D\sigma_{\MA}(\MA) \cong \MA$.
\end{proposition}

\begin{proof}
  Let~$\VS$ be a \Mp{\Tanil\MA}bimodule.  We claim that composition with $D
  \circ \sigma_{\MA} \colon \MA \to \Omega^1 (\Tanil\MA)$ gives rise to a
  bijection between \Mp{\Tanil\MA}bimodule homomorphisms $\Omega^1(\Tanil\MA)
  \to \VS$ and bounded linear maps $\MA \to \VS$.  This is precisely the
  universal property of a free module.  Since universal objects are uniquely
  determined, it follows that~$\mu_{\ref{eq:mu3}}$ is a bornological
  isomorphism.
  
  Lemma~\ref{lem:Omega_I_ltimes} yields a bijection between bimodule
  homomorphisms $\Omega^1(\Tanil\MA) \to \VS$ and bounded splitting
  homomorphisms $\Tanil\MA \to \Tanil\MA \ltimes \VS$ for the natural
  projection $p \colon \Tanil\MA \ltimes \VS \to \Tanil\MA$.  The bimodule
  homomorphism $f \colon \Omega^1(\Tanil\MA) \to \VS$ corresponds to the
  homomorphism $(\ID, f \circ D) \colon \Tanil\MA \to \Tanil\MA \ltimes \VS$.
  By the universal property~\ref{axiom:Tanil_universal} of $\Tanil\MA$,
  composition with~$\sigma_{\MA}$ gives rise to a bijection between bounded
  homomorphisms $\Tanil\MA \to \Tanil\MA \ltimes \VS$ and \lanilcurs $\MA \to
  \Tanil\MA \ltimes \VS$.  Furthermore, the homomorphism associated to a
  \lanilcur $l \colon \MA \to \Tanil\MA \ltimes \VS$ is a section for~$p$ iff
  $p \circ l = \sigma_{\MA}$.  Thus we get a bijection between bimodule
  homomorphisms $\Omega^1(\Tanil\MA) \to \VS$ and bounded linear maps $l
  \colon \MA \to \VS$ with the additional property that $(\sigma_{\MA}, l)
  \colon \MA \to \Tanil\MA \ltimes \VS$ is a \lanilcur.  To a bimodule
  homomorphism $f \colon \Omega^1(\Tanil\MA) \to \VS$, this bijection
  associates the bounded linear map $f \circ D \circ \sigma_{\MA} \colon \MA
  \to \VS$.
  
  To finish the proof of~\eqref{eq:mu3}, it remains to show that each map $\MA
  \to \Tanil\MA \ltimes \VS$ of the form $(\sigma_{\MA}, l')$ with a bounded
  linear map $l' \colon \MA \to \VS$ is a \lanilcur.  We have $\tau_{\MA}
  \circ p \circ (\sigma_{\MA}, l') = \tau_{\MA} \circ \sigma_{\MA} =
  \ID[{\MA}]$.  Since this is a homomorphism, the curvature of $(\sigma_{\MA},
  l')$ factors through $\MA[N] \defeq \Ker (\tau_{\MA} \circ p) \colon
  \Tanil\MA \ltimes \VS \to \MA$.  Since~$p$ is split surjective, we have an
  allowable extension $\Ker p \injto \MA[N] \prto \Ker \tau_{\MA}$.  Since
  $\Ker p = \Null(\VS)$ and $\Ker \tau_{\MA} = \Janil\MA$ are a-nilpotent by
  Axiom~\ref{axiom:nilpotent_lanilcur} and Axiom~\ref{axiom:Janil_nilpotent},
  the Extension Axiom~\ref{axiom:extension} implies that~$\MA[N]$ is
  a-nilpotent.  Thus $(\sigma_{\MA}, l')$ is a \lanilcur by
  Axiom~\ref{axiom:lanilcur_nilpotent}.  The proof of~\eqref{eq:mu3} is
  complete.

  To obtain the isomorphisms \eqref{eq:mu1} and~\eqref{eq:mu2}, we consider
  the natural allowable extension
  $$
  E \colon \
  \Omega^1(\Tanil\MA) \injto
  \Unse{(\Tanil\MA)} \hot \Unse{(\Tanil\MA)} \prto
  \Unse{(\Tanil\MA)}
  $$
  (see~\eqref{eq:res1} below) and replace $\Omega^1(\Tanil\MA)$ by
  $\Unse{(\Tanil\MA)} \hot \MA \hot \Unse{(\Tanil\MA)}$.  This extension
  splits naturally as an extension of left \Mp{\Tanil\MA}modules.  Hence
  if~$\VS$ is a left \Mp{\Tanil\MA}module, then
  $$
  E \hot_{\Tanil\MA} \VS \colon \
  \Unse{(\Tanil\MA)} \hot \MA \hot \VS \injto
  \Unse{(\Tanil\MA)} \hot \VS \prto
  \VS
  $$
  is an allowable extension of right \Mp{\Tanil\MA}modules.  We apply this
  to $\VS \defeq \C$ with the zero module structure.  That is, $x \cdot 1 = 0$
  for all $x \in \Tanil\MA$.  We get an extension $\Unse{(\Tanil\MA)} \hot \MA
  \injto \Unse{(\Tanil\MA)} \prto \C$.  The map $\Unse{(\Tanil\MA)} \to \C$ is
  the usual one with kernel $\Tanil\MA$ and the map $\Unse{(\Tanil\MA)} \hot
  \MA \to \Unse{(\Tanil\MA)}$ is induced by the bilinear map $(x, \ma) \mapsto
  x \odot \sigma(\ma)$.  Thus~\eqref{eq:mu1} is an isomorphism.  The
  isomorphism in~\eqref{eq:mu2} is obtained similarly by applying the functor
  $\C \hot_{\Tanil\MA} \blank$ to~$E$.
\end{proof}

\subsection{The \lanilcur category and Goodwillie's theorem}
\label{sec:lanilcur}

The Extension Lemma~\ref{lem:extension}.(ii) implies that \lanilcurs form the
morphism of a category.  We call this category the \emph{\lanilcur category}.
The category whose objects are the complete bornological algebras and whose
morphisms are the bounded homomorphisms is called the \emph{homomorphism
  category}.  We show that $\Tanil$ is a smooth homotopy functor from the
\lanilcur category to the homomorphism category.  Up to smooth homotopy this
functor is fully faithful, that is, smooth homotopy classes of homomorphisms
$\Tanil\MA \to \Tanil\MA[B]$ correspond bijectively to smooth homotopy classes
of \lanilcurs $\MA \to \MA[B]$.  The proof uses only the axioms of
Section~\ref{sec:Tanil} and hence carries over without change to absolutely
continuous homotopies because $\ABC([0,1])$ is tensoring as well.

\begin{lemma}  \label{lem:Tanil_functor_lanilcur}
  Let $l \colon \MA \to \MA[B]$ be a \lanilcur.  There is a unique bounded
  homomorphism $\Tanil l \colon \Tanil\MA \to \Tanil\MA[B]$ making the diagram
  \begin{equation}  \label{eq:Tanil_functorial_lanilcur}
    \begin{gathered}
      \xymatrix{
        {\Tanil\MA} \ar@{-->}[r]^{\Tanil l} &
          {\Tanil\MA[B]} \\
        {\MA} \ar@{.>}[r]^{l} \ar@{.>}[u]^{\sigma_{\MA}} &
          {\MA[B]} \ar@{.>}[u]_{\sigma_{\MA[B]}}
        }
    \end{gathered}
  \end{equation}
  commute, namely $\Tanil l \defeq \LLH{\sigma_{\MA[B]} \circ l}$.  Moreover,
  $\LLH{l} = \tau_{\MA[B]} \circ \Tanil l$.
\end{lemma}

\begin{proof}
  Since $\sigma_{\MA[B]}$ is a \lanilcur, the map $\sigma_{\MA[B]} \circ l$ is
  a \lanilcur as a composition of two \lanilcurs.  Thus $\Tanil l \defeq
  \LLH{\sigma_{\MA[B]} \circ l}$ is a well-defined bounded homomorphism by the
  universal property~\ref{axiom:Tanil_universal} of $\Tanil\MA$.  In fact, the
  uniqueness part of the universal property asserts that $\Tanil l$ is the
  only bounded homomorphism satisfying $\Tanil l \circ \sigma_{\MA} =
  \sigma_{\MA[B]} \circ l$, that is, making the above diagram commute.  The
  composition $\tau_{\MA[B]} \circ \Tanil l$ is a bounded homomorphism with
  $\tau_{\MA[B]} \circ \Tanil l \circ \sigma_{\MA} = l$.  Hence $\LLH{l} =
  \tau_{\MA[B]} \circ \Tanil l$.
\end{proof}

Comparing \eqref{eq:Tanil_functorial_lanilcur}
and~\eqref{eq:Tanil_functorial_abstract}, we see that if~$l$ is a bounded
homomorphism, we recover the usual functoriality of~$\Tanil$ for bounded
homomorphisms.  The universal property~\ref{axiom:Tanil_universal} of
$\Tanil\MA$ can be interpreted as the assertion that the functor~$\Tanil$ is
left adjoint to the inclusion functor from the homomorphism category to the
\lanilcur category.  Next we establish that $\Tanil$ is a smooth homotopy
functor from the \lanilcur category to the homomorphism category.  A smooth
homotopy in the \lanilcur category is a \lanilcur $\MA \to \CINF([0,1];
\MA[B])$.  Smooth homotopy of \lanilcurs is an equivalence relation.

\begin{lemma}  \label{lem:Tanil_homotopy_lanilcur}
  The functor $\Tanil$ is a smooth homotopy functor from the \lanilcur
  category to the homomorphism category.  In fact, any smooth homotopy in the
  \lanilcur category $H \colon \MA \to \CINF([0,1]; \MA[B])$ can be lifted to
  a unique smooth homotopy in the homomorphism category $\hat{H} \colon
  \Tanil\MA \to \CINF([0,1]; \Tanil\MA[B])$ satisfying $\hat{H} \circ
  \sigma_{\MA} = \CINF([0,1]; \sigma_{\MA[B]}) \circ H$.  The same holds for
  absolutely continuous homotopies.
\end{lemma}

\begin{proof}
  Copy the proof of statement~(iii) of the Tensor Lemma~\ref{lem:tensor},
  using that \lanilcurs form a category instead of
  Lemma~\ref{lem:lanilcur_compose_hom}.  Use that $\CINF([0,1])$ and
  $\ABC([0,1])$ are tensoring.
\end{proof}

\begin{theorem}  \label{the:homotopy_nilpotent}
  Let $(i,p) \colon \MA[N] \injto \MA[E] \prto \MA$ be an a-nilpotent
  extension.  Then~$\MA[N]$ is smoothly contractible and $\MA$ is a smooth
  deformation retract of\/~$\MA[E]$ in the \lanilcur category.  Thus
  $\Tanil\MA[N]$ is smoothly contractible and $\Tanil\MA$ is a smooth
  deformation retract of $\Tanil\MA[E]$ in the homomorphism category.  In
  particular, $\Tanil p \colon \Tanil\MA[E] \to \Tanil\MA$ is a smooth
  homotopy equivalence.
\end{theorem}

\begin{proof}
  Since~$\MA[N]$ is a-nilpotent, so is $\CINF([0,1]; \MA[N])$ by the Homotopy
  Axiom~\ref{axiom:homotopy}.  Hence the bounded linear map $H\colon \MA[N]
  \to \CINF([0,1]; \MA[N])$ defined by $H_t(x) \defeq t\cdot x$ is a
  \lanilcur.  Thus~$\MA[N]$ is smoothly contractible in the \lanilcur
  category.
  
  Let $s \colon \MA \to \MA[E]$ be a bounded linear section.  The curvature
  of~$s$ has values in~$\MA[N]$.  Since~$\MA[N]$ is a-nilpotent, $s$ is a
  \lanilcur.  By definition, $p \circ s = \ID$.  Define $H \colon \MA[E] \to
  \CINF([0,1]; \MA[E])$ by $H = 1 \otimes \ID + t \otimes (s \circ p - \ID)$,
  that is, $H_t(x) = (1-t) x + t\cdot (s\circ p)(x)$.  Clearly, $p \circ H =
  \CINF([0,1]; p)$ (that is, $p\circ H_t(x) = p(x)$ for all $x \in \MA[E]$).
  Thus the curvature of~$H$ factors through the a-nilpotent algebra
  $\CINF([0,1]; \MA[N])$, so that~$H$ is a \lanilcur by
  Axiom~\ref{axiom:lanilcur_nilpotent}.  The \lanilcur~$H$ is a smooth
  homotopy between $H_0 = \ID$ and $H_1 = s \circ p$.  Thus~$\MA$ is a smooth
  deformation retract of~$\MA[E]$ in the \lanilcur category.  Upon applying
  the smooth homotopy functor~$\Tanil$, we get that $\Tanil\MA[N]$ is smoothly
  contractible and that $\Tanil\MA$ is a smooth deformation retract of
  $\Tanil\MA[E]$ in the homomorphism category.
\end{proof}

In particular, this applies to the a-nilpotent extension $\Janil\MA \injto
\Tanil\MA \prto \MA$.  Thus $\Janil\MA$ is contractible and $\MA$
and $\Tanil\MA$ are smoothly homotopy equivalent in the \lanilcur category.

Let~$F$ be a homology theory for a-quasi-free algebras, that is, a smooth
homotopy functor from the homomorphism category of a-quasi-free algebras to
the category of Abelian groups.  We can extend~$F$ to a homology theory on
arbitrary complete bornological algebras by $\bar{F}(\MA) \defeq
F(\Tanil\MA)$.  Since $\Tanil$ is a homotopy functor from the \lanilcur
category to the homomorphism category of a-quasi-free algebras, $\bar{F}$ is a
smooth homotopy functor on the \lanilcur category.  Thus
Theorem~\ref{the:homotopy_nilpotent} implies that $\bar{F}(\MA[N]) = 0$ and
$\bar{F}(\MA[E]) \cong \bar{F}(\MA)$ if $\MA[N] \injto \MA[E] \prto \MA$ is an
a-nilpotent extension.  This is a generalization of a theorem of Goodwillie
about periodic cyclic cohomology~\cite{goodwillie85:cyclic}.

This argument applies in particular to analytic cyclic cohomology and yields
the special case of excision for a-nilpotent extensions.  For nilpotent
extensions of Banach algebras and entire cyclic cohomology, this special case
has been proved by Khalkhali~\cite{khalkhali94:nilpotent}.

\medbreak

The functor $\MA \mapsto \Tanil\MA$ from the \lanilcur category to the
homomorphism category is faithful.  It is obviously not fully faithful, that
is, there are bounded homomorphisms $\Tanil\MA \to \Tanil\MA[B]$ not of the
form $\Tanil l$ for a \lanilcur $l \colon \MA \to \MA[B]$.  Nevertheless, the
functor $\Tanil$ is fully faithful up to smooth homotopy.  To make this
precise, we introduce some notation.  Let $\Mor_l(\MA; \MA[B])$ be the set of
\lanilcurs from $\MA$ to~$\MA[B]$ and let $\Mor_h( \Tanil\MA; \Tanil\MA[B])$
be the set of homomorphisms from $\Tanil\MA$ to $\Tanil\MA[B]$.  We do not
endow these sets with a topology.  However, if~$M$ is a smooth manifold, we
can define the set of smooth maps $M \to \Mor_l(\MA;\MA[B])$ as follows.  A
map $f\colon M \to \Mor_l(\MA;\MA[B])$ is called \emph{smooth} if there is a
necessarily unique $\hat{f} \in \Mor_l \bigl( \MA; \CINF(M; \MA[B]) \bigr)$
with $f(m) = \ev_m \circ \hat{f}$ for all $m\in M$, where $\ev_m \colon
\CINF(M) \to \C$ denotes the evaluation at $m \in M$.  We define a similar
smooth structure on $\Mor_h( \Tanil\MA; \Tanil\MA[B])$.

Define $H^{\MA[B]} \colon \Tanil\MA[B] \to \CINF([0,1]; \Tanil\MA[B])$ by
$$
H^{\MA[B]} =
1 \otimes \ID - t \otimes (\ID - \sigma_{\MA[B]} \circ \tau_{\MA[B]}).
$$
Since $\CINF([0,1]; \tau_{\MA[B]}) \circ H^{\MA[B]} = \CINF([0,1];
\tau_{\MA[B]})$, the curvature of~$H^{\MA[B]}$ has values in $\CINF([0,1];
\Janil\MA[B])$.  This algebra is a-nilpotent by the Homotopy
Axiom~\ref{axiom:homotopy}, so that~$H^{\MA[B]}$ is a \lanilcur.

\begin{proposition}  \label{pro:Tanil_fufa}
  The map
  $$
  h \colon \Mor_h( \Tanil\MA; \Tanil\MA[B]) \times [0,1] \to \Mor_h
  (\Tanil\MA; \Tanil\MA[B]),
  \qquad
  (\phi, t) \mapsto \LLH{H_t^{\MA[B]} \circ \phi \circ \sigma_{\MA}}
  $$
  is a natural smooth deformation retraction from $\Mor_h (\Tanil\MA;
  \Tanil\MA[B])$ onto $\Tanil \bigl(\Mor_l (\MA; \MA[B]) \bigr)$.  It
  satisfies $h(\phi, 1)= \Tanil (\tau_{\MA[B]} \circ \phi \circ
  \sigma_{\MA})$.
\end{proposition}

\begin{proof}
  Since $H^{\MA[B]}$ and~$\sigma_{\MA}$ are \lanilcurs, it is evident that~$h$
  is well-defined.  Furthermore, $h(\phi, 0) = \phi$ because $H_0 = \ID$ and
  $$
  h(\phi, 1) =
  \LLH{\sigma_{\MA[B]} \circ \tau_{\MA[B]} \circ \phi \circ \sigma_{\MA}} =
  \Tanil(\tau_{\MA[B]} \circ \phi \circ \sigma_{\MA})
  $$
  by the definition of the functoriality for \lanilcurs in
  Lemma~\ref{lem:Tanil_functor_lanilcur}.
  
  The proof of the Homotopy Axiom~\ref{axiom:homotopy} carries over to
  $\CINF(M)$ without change.  Thus $\CINF(M)$ is tensoring.  This implies
  that~$h$ is a smooth map.  Details are left to the interested reader.
\end{proof}

\section{The X-complex of $\Tanil\MA$ and analytic cyclic cohomology}
\label{sec:X_Tanil_HA}

The X-complex of a (unital) algebra is introduced by Cuntz and Quillen
in~\cite{cuntz95:cyclic}.  It is obtained by truncating the
\Mp{(b,B)}bicomplex in degree~$1$.  Its even cohomology can be interpreted as
the space of traces on~$\MA$ modulo the equivalence relation generated by
smooth homotopy.  For most algebras~$\MA$, the complex $X(\MA)$ is not very
interesting because it ignores all the information in $\HH^n(\MA)$ with $n \ge
2$.  For \emph{quasi-free} algebras, however, there is no such additional
information above degree~$1$ and the X-complex computes the periodic cyclic
cohomology of the algebra.  The X-complex is a homotopy functor for quasi-free
algebras in the following sense.

\begin{proposition}  \label{pro:X_homotopy}
  Let $\MA$ and~$\MA[B]$ be complete bornological algebras.  Let $\phi_0,
  \phi_1 \colon \MA \to \MA[B]$ be AC-homotopic bounded homomorphisms.  Assume
  that~$\MA$ is quasi-free.

  Then the induced maps $X(\phi_t)$, $t = 0,1$, are chain homotopic.  That is,
  there is a bounded linear map $h \colon X(\MA) \to X(\MA[B])$ of degree~$1$
  such that $X(\phi_0) = X(\phi_1) - [\partial, h]$ with $[\partial, h] \defeq
  \partial_{X(\MA[B])} \circ h + h\circ \partial_{X(\MA)}$.
\end{proposition}

This is proved for fine algebras in Sections~7--8 of~\cite{cuntz95:cyclic} by
writing down an explicit formula for~$h$ and checking that it does the trick.
Such a proof carries over immediately to complete bornological algebras.  It
only remains to verify that the map~$h$ written down in~\cite{cuntz95:cyclic}
is bounded.  In Appendix~\ref{app:X_homotopy}, we will do this and in addition
explain how to derive the formula for~$h$.

The \emph{analytic cyclic (co)homology} of~$\MA$ is the (co)homology of the
complex $X(\Tanil\MA)$.  We define the \emph{bivariant analytic cyclic
  cohomology} $\HA^\ast(\MA; \MA[B])$ as the homology of the complex of
bounded linear maps from $X(\Tanil\MA)$ to $X(\Tanil\MA[B])$.  In these
definitions, we can replace $\Tanil\MA$ and $\Tanil\MA[B]$ by arbitrary
universal a-nilpotent extensions of $\MA$ and~$\MA[B]$, respectively.  If
$\MA[N] \injto \MA[R] \prto \MA$ is a universal a-nilpotent extension
of~$\MA$, then~$\MA[R]$ is naturally smoothly homotopy equivalent to
$\Tanil\MA$.  Thus the complexes $X(\Tanil\MA)$ and $X(\MA[R])$ are homotopy
equivalent by Proposition~\ref{pro:X_homotopy} and compute the same
cohomology.

In addition, Proposition~\ref{pro:X_homotopy} implies that the analytic cyclic
cohomology theories are invariant under absolutely continuous homotopies.
They are also stable with respect to the algebra $\Sch(\Hils)$ of trace class
operators on a Hilbert space and similar objects.  Stability is
straightforward to prove using an idea of Elliott, Natsume, and
Nest~\cite{elliott88:cyclic}.  Using that the fine algebras $\C$ and $\C[u,
u^{-1}]$ are a-quasi-free, we construct the Chern-Connes character from
topological \Mpn{K}theory to analytic cyclic homology.

The identification of $\Omega^1(\Tanil\MA)$ in~\eqref{eq:mu3} gives rise to an
isomorphism $X(\Tanil\MA) \cong \Omega_\an\MA$ as \Mpn{\Ztwo}graded
bornological vector spaces.  We compute the boundary $\partial \colon
\Omega_\an\MA \to \Omega_\an\MA$ coming from the boundary in the X-complex
$X(\Tanil\MA)$.  This boundary map~$\partial$ is quite similar to the more
well-known boundary $B + b$ used to define entire cyclic cohomology.  In fact,
Cuntz and Quillen~\cite{cuntz95:cyclic} have shown that, in the situation of
periodic cyclic cohomology, these two boundaries give rise to chain homotopic
complexes.  Making their argument more explicit, we obtain that if~$\MA$ is a
(unital) locally convex algebra endowed with the bounded bornology, then the
analytic cyclic cohomology $\HA^\ast(\MA)$ is naturally isomorphic to the
entire cyclic cohomology $\HE^\ast(\MA)$ in the sense of Connes
\cite{connes88:entire}, \cite{connes94:ncg}.  Hence entire cyclic cohomology
is a special case of analytic cyclic cohomology.

The homotopy invariance and stability of analytic cyclic cohomology specialize
to analogous statements about entire cyclic cohomology.  These special cases
have been handled previously by Khalkhali.  See~\cite{khalkhali97:survey} for
a nice survey.

\subsection{The X-complex of a quasi-free algebra}
\label{sec:X_complex}

Let~$\MA$ be a complete bornological algebra.  Let $\Omega^1\MA /[,] \defeq
\Omega^1\MA / [\Omega^1\MA,\MA] \defeq \Omega^1\MA / b(\Omega^2\MA)$.  This is
the \emph{commutator quotient} of $\Omega^1\MA$.  Let $X_\beta(\MA)$ be the
\Hochschild complex of~$\MA$ truncated in degree~$1$:
$$
0 \longrightarrow
\Omega^1\MA / b(\Omega^2\MA) \overset{b_\ast}{\longrightarrow}
\Omega^0\MA.
$$
Since $b \circ b=0$, the boundary~$b$ descends to a bounded linear map
$b_\ast \colon \Omega^1\MA / [,] \to \Omega^0\MA$.  The subspace
$b(\Omega^2\MA) \subset \Omega^1\MA$ may fail to be closed, so that
$X_\beta(\MA)$ may be non-separated.  However, if~$\MA$ is quasi-free then
$b(\Omega^2\MA)$ is a direct summand in $\Omega^1\MA$ by
Proposition~\ref{pro:Hochschild_contract}, so that $X_\beta(\MA)$ is
separated.

The complex $X_\beta(\MA)$ will be used in the proof of the Excision Theorem.
It is obviously related to \Hochschild homology and thus can be treated by the
tools of homological algebra.  The X-complex $X(\MA)$ is almost equal to
$X_\beta(\MA)$, we only add another boundary $\partial_0 \colon X_0(\MA) \to
X_1(\MA)$ that is defined by $\partial_0 \defeq d_\ast \colon \MA
\overset{d}{\longrightarrow} \Omega^1\MA \longrightarrow \Omega^1\MA/[,]$.
Thus we send~$\ma_0$ to the class of~$d\ma_0$ modulo commutators.  It is clear
that $b_\ast \circ d_\ast = 0$.  We have $d_\ast \circ b_\ast = 0$ because
$$
d b (\opt{\ma_0} d\ma_1) =
d[\opt{\ma_0}, \ma_1] =
b( d\ma_1 d\opt{\ma_0} - d\opt{\ma_0} d\ma_1) =
-bB (\opt{\ma_0} d\ma_1)
$$
for the bounded linear map $B \colon \Omega^1\MA \to \Omega^2\MA$ defined
by $B(\opt{\ma_0} d\ma_1) \defeq d\opt{\ma_0} d\ma_1- d\ma_1 d\opt{\ma_0}$.
Hence $d b(\Omega^1\MA) \subset -bB(\Omega^1\MA) \subset b(\Omega^2\MA)$.
This means that $d_\ast \circ b_\ast = 0$.  Thus $X(\MA)$ is a
\Mpn{\Ztwo}graded complex of bornological vector spaces.

The X-complex is defined in a natural way and thus functorial for bounded
homomorphisms.  A bounded homomorphism $\phi \colon \MA \to \MA[B]$ induces
the chain map $X(\phi) \colon X(\MA) \to X(\MA[B])$,
\begin{equation}  \label{eq:X_functorial}
  \begin{alignedat}{}{2}
    X(\phi)(\ma) &\defeq \phi(\ma), &\qquad &\forall
    \ma \in \MA;
    \\
    X(\phi)(\opt{\ma_0} d\ma_1 \bmod [,]) &\defeq
    \phi\opt{\ma_0} d\phi(\ma_1) \bmod [,], &\qquad &\forall
    \opt{\ma_0} \in \Unse{\MA},\ \ma_1 \in \MA.
  \end{alignedat}
\end{equation}

The cohomology of $X(\MA)$ is interpreted by Cuntz and
Quillen~\cite{cuntz95:cyclic}.  An even cocycle $l \colon X(\MA) \to \C$ is a
bounded linear map $\MA \to \C$ satisfying $l \circ b = 0$, that is a bounded
trace on~$\MA$.  An even coboundary is a trace of the form $l' \circ d$ for a
bounded linear map $l' \colon \Omega^1\MA / [,] \to \C$, that is, a bounded
trace on $\Omega^1\MA$.  The equivalence relation obtained by addition of such
coboundaries is a simple special case of the relation called \emph{cobordism}
by Connes~\cite{connes94:ncg}.  Alternatively, this relation can be generated
by smooth homotopy.  The space of traces on~$\MA$ is a contravariant functor
that is far from being homotopy invariant.  In order to get a smooth
homotopy functor, we should identify traces that are of the form $T \circ
f_t$, where $f \colon \MA \to \CINF([0,1]; \MA[B])$ is a smooth homotopy and
$T \colon \MA[B] \to \C$ is a bounded trace.  If we have such a smooth
homotopy, then the traces $T \circ f_1$ and $T \circ f_0$ differ only by the
coboundary of the trace
$$
\opt{x} dy \mapsto \int_0^1 f_t\opt{x} \frac{\partial}{\partial t} f_t(y) \,dt
$$
on $\Omega^1\MA$.  Conversely, cobordant traces are smoothly homotopic (see
\cite[Proposition~4.2]{cuntz95:cyclic}).  In particular, the cohomology
$H^0\bigl( X(\MA) \bigr)$ is a smooth homotopy functor for arbitrary
algebras~$\MA$.  The odd cohomology is only smoothly homotopy invariant for
quasi-free algebras.

\subsection{Definition and functoriality of analytic cyclic cohomology}
\label{sec:def_HA}

Let~$\MA$ be a complete bornological algebra.  Define the \emph{analytic
cyclic homology and cohomology} $\HA_\ast(\MA)$ and $\HA^\ast(\MA)$, $\ast =
0,1$, as the homology and cohomology of the \Mpn{\Ztwo}graded complex of
bornological vector spaces $X(\Tanil\MA)$.  (For the cohomology, we only
consider bounded linear maps.)  Define the \emph{bivariant analytic cyclic
cohomology} $\HA^\ast(\MA; \MA[B])$ as the homology of the complex of bounded
linear maps $X(\Tanil\MA) \to X(\Tanil\MA[B])$.  In the notation of
Appendix~\ref{app:complexes}, we have
$$
\HA^\ast(\MA; \MA[B]) \defeq
H^\ast \bigl( X(\Tanil\MA); X(\Tanil\MA[B]) \bigr),
\qquad \ast = 0,1.
$$
We will see that $X(\Tanil\C)$ is naturally chain homotopic to the complex
$\C[0]$ with~$\C$ in degree~$0$ and~$0$ in degree~$1$.  Hence $\HA^\ast (\C;
\MA) \cong \HA_\ast(\MA)$ and $\HA^\ast (\MA; \C) \cong \HA^\ast(\MA)$.

The composition of linear maps gives rise to an associative product on the
homology groups
$$
\HA^i(\MA_1; \MA_2) \otimes \HA^j(\MA_0; \MA_1) \mapsto
\HA^{i+j}(\MA_0; \MA_2),
\qquad i,j = 0,1.
$$
A bounded homomorphism $f \colon \MA \to \MA[B]$ induces a bounded chain
map
$$
X(\Tanil f) \colon X(\Tanil\MA) \to X(\Tanil\MA[B])
$$
and thus a class~$[f]$ in $\HA^0(\MA; \MA[B])$ by
Lemma~\ref{lem:lanilcur_compose_hom} and~\eqref{eq:X_functorial}.  Composition
with~$[f]$ makes $\HA^\ast(\MA; \MA[B])$ a bifunctor for bounded homomorphisms
contravariant in the first and covariant in the second variable.  It makes
$\HA^\ast(\MA)$ and $\HA_\ast(\MA)$ into functors as well.

Roughly speaking, elements of $\HA^0(\MA)$ may be viewed as \emph{homotopy
  classes of traces on a-nilpotent extensions of~$\MA$}.  Firstly, we can view
$\HA^0(\MA) \defeq H^0 \bigl( X(\Tanil\MA) \bigr)$ as the space of smooth
homotopy classes of traces on $\Tanil\MA$.  Let $\MA[N] \injto \MA[E] \prto
\MA$ be an a-nilpotent extension of~$\MA$ and let $T \colon \MA[E] \to \C$ be
a trace on~$\MA[E]$.  By the Universal Extension
Theorem~\ref{the:universal_ext}, there is a morphism of extensions $(\xi,
\psi, \ID)$ from the universal a-nilpotent extension $\Janil\MA \injto
\Tanil\MA \to \MA$ to $\MA[N] \injto \MA[E] \prto \MA$.  Furthermore, this
morphism is unique up to smooth homotopy.  Thus the trace $T \circ \psi$ on
$\Tanil\MA$ defines a unique homotopy class of traces on $\Tanil\MA$ and thus
a unique element of $\HA^0(\MA)$.

\subsection{Homotopy invariance and stability}
\label{sec:homotopy_stability}

\begin{theorem}  \label{the:HA_homotopy}
  Let $f_0, f_1 \colon \MA \to \MA[B]$ be AC-homotopic bounded homomorphisms.
  Then $[f_0] = [f_1]$ in $\HA^0(\MA; \MA[B])$.  Thus analytic cyclic homology
  and (bivariant) analytic cyclic cohomology are invariant under absolutely
  continuous homotopies.
\end{theorem}

\begin{proof}
  Let $f_t \colon \MA \to \MA[B]$, $t = 0,1$, be AC-homotopic homomorphisms.
  Lemma~\ref{lem:Tanil_homotopy} implies that the homomorphisms $\Tanil f_t
  \colon \Tanil\MA \to \Tanil\MA[B]$ are AC-homotopic as well.  By the
  Extension Lemma~\ref{lem:extension}.(iii), the algebra $\Tanil\MA$ is
  a-quasi-free and a fortiori quasi-free by Corollary~\ref{cor:aqf_qf}.
  Proposition~\ref{pro:X_homotopy} therefore implies that the induced maps
  $X(\Tanil f_t) \colon X(\Tanil\MA) \to X(\Tanil\MA[B])$ are chain homotopic.
  Thus $[f_0] = [f_1]$ in $\HA^0(\MA; \MA[B])$ as asserted.  Since the
  functoriality of the analytic cyclic theories is given by composition with
  the class~$[f_t]$, it follows that AC-homotopic homomorphisms induce the
  same map on all analytic cyclic (co)homology groups.
\end{proof}

\begin{theorem}  \label{the:uniqueness_X}
  Let $(i,p) \colon \MA[N] \injto \MA[R] \prto \MA$ be a universal a-nilpotent
  extension of~$\MA$.  Then the complex $X(\MA[R])$ is chain homotopic to
  $X(\Tanil\MA)$ in a canonical way.  If $(\xi, \psi, \ID)$ is a morphism of
  extensions from $\Janil\MA \injto \Tanil\MA \prto \MA$ to $\MA[N] \injto
  \MA[R] \prto \MA$, then $X(\psi) \colon X(\Tanil\MA) \to X(\MA[R])$ is a
  homotopy equivalence of complexes.  The class $[X(\psi)]$ of $X(\psi)$ in
  $H^\ast \bigl( X(\Tanil\MA); X(\MA[R]) \bigr)$ does not depend on the choice
  of~$\psi$.
  
  In particular, if~$\MA$ is a-quasi-free itself, then $X(\tau_{\MA}) \colon
  X(\Tanil\MA) \to X(\MA)$ is a homotopy equivalence of complexes.  If
  $\upsilon \colon \MA \to \Tanil\MA$ is a bounded splitting homomorphism,
  then $X(\upsilon)$ is an inverse to $X(\tau_{\MA})$ up to homotopy.
\end{theorem}

\begin{proof}
  By the Uniqueness Theorem~\ref{the:uniqueness}, the extension $\MA[N] \injto
  \MA[R] \prto \MA$ is smoothly homotopy equivalent to $\Janil\MA \injto
  \Tanil\MA \prto \MA$.  In particular, $\MA[R]$ and $\Tanil\MA$ are smoothly
  homotopy equivalent.  That is, there are homomorphisms $\psi \colon
  \Tanil\MA \to \MA[R]$ and $\psi' \colon \MA[R] \to \Tanil\MA$ such that
  $\psi \circ \psi'$ and $\psi' \circ \psi$ are smoothly homotopic to the
  identity.  Furthermore, $\psi$ and~$\psi'$ are unique up to smooth homotopy.
  Since $\MA[R]$ and $\Tanil\MA$ are a-quasi-free, they are quasi-free by
  Corollary~\ref{cor:aqf_qf}.  The homotopy invariance of the
  X-complex~\ref{pro:X_homotopy} therefore applies.  It shows that $X(\psi
  \circ \psi') = X(\psi) \circ X(\psi') \sim \ID$, $X(\psi' \circ \psi) \sim
  \ID$, and that $[X(\psi)]$ is independent of the choice of~$\psi$.
  
  If~$\MA$ is a-quasi-free, then we can apply this to the universal
  a-nilpotent extension $0 \injto \MA \prto \MA$.  This yields that
  $X(\tau_{\MA})$ and $X(\upsilon)$ are inverses of each other up to chain
  homotopy.
\end{proof}

Next we prove stability with respect to the algebra $\Sch(\Hils)$ of trace
class operators on a Hilbert space.  We endow this Banach algebra with the
bounded bornology.  Actually, it is easier to prove a more general result.
The problem with $\Sch(\Hils)$ is that during the proof we have to consider
the bornological tensor product $\Sch(\Hils) \hot \Sch(\Hils) \cong \Sch(\Hils
\prot \Hils)$, but $\Hils \prot \Hils$ is no longer a Hilbert space.
Therefore, we consider the following more general situation.

Let $\VS$ and~$\VS[W]$ be complete bornological vector spaces and let $b
\colon \VS[W] \times \VS \to \C$ be a bounded bilinear map that is not the
zero map.  That is, there are $\vs_0 \in \VS$, $\vs[w]_0 \in \VS[W]$ with
$b(\vs[w]_0, \vs_0) = 1$.  We endow the complete bornological vector space
$\VS \hot \VS[W]$ with the multiplication defined by
$$
(\vs_1 \otimes \vs[w]_1) \cdot (\vs[v]_2 \otimes \vs[w]_2) \defeq
\vs_1 \otimes b(\vs[w]_1, \vs_2) \cdot \vs[w]_2.
$$
It is straightforward to verify that $(S \otimes T) \cdot (S\otimes T) \subset
c\cdot \bipol{(S \otimes T)}$ for the constant $c \defeq \sup |b(T,S)|$.  Thus
the multiplication above is bounded on $\VS \otimes \VS[W]$ and therefore
extends to $\VS \hot \VS[W]$.  It is evidently associative.  Given $S \in
\CBS(\VS)$, $T \in \CBS(\VS[W])$, we rescale $T \in \CBS(\VS[W])$ such that
$|b(S,T)| \le 1$.  Then $\bipol{\bigl((S \otimes T)\cdot (S\otimes T)\bigr)}
\subset \bipol{(S \otimes T)}$.  It follows that $\bigl( \bipol{(S \otimes T)}
\bigr)^\infty \subset \VS \otimes \VS[W]$ is small.  By
Lemma~\ref{lem:tensor_a_nilpotent}, the algebra $\VS \hot \VS[W]$ is
tensoring.

If we let $\VS \defeq \Hils$, $\VS[W] \defeq \Hils'$, and let $b \colon \Hils
\times \Hils' \to \C$ be the natural pairing with the dual, then $\VS \hot
\VS[W] \cong \Sch(\Hils)$ is isomorphic to the algebra of trace class
operators.  The isomorphism sends $\vs \otimes \vs[w]$ to the rank one
operator $\xi \mapsto \vs \cdot \vs[w](\xi)$.  Thus we write $\Sch(b)$ for the
algebra $\VS \hot \VS[W]$ defined above.

On $\Sch(b)$, we have a natural trace $\tr$ defined by $\tr (\vs \otimes
\vs[w]) \defeq b(\vs[w], \vs)$.  This generalizes the usual trace on
$\Sch(\Hils)$.  In addition, we have an analogue of rank one idempotents.  Let
$\vs_0 \in \VS$ and $\vs[w]_0 \in \VS[W]$ be such that $b(\vs_0, \vs[w]_0) =
1$.  Define $p_0 \defeq \vs_0 \otimes \vs[w]_0 \in \Sch(b)$.  Then $p_0^2 =
p_0$ and $\tr p_0 = 1$.

Thus we can define a bounded homomorphism $\iota \colon \C \to \Sch(b)$ by
$\iota(1) \defeq p_0$.  Tensoring with an arbitrary complete bornological
algebra~$\MA$, we obtain a natural bounded homomorphism $\iota_{\MA} \defeq
\ID[{\MA}] \hot \iota \colon \MA \to \MA \hot \Sch(b)$.

\begin{theorem}[Stability Theorem]  \label{the:HA_stable}
  The class of~$\iota_{\MA}$ is an invertible element in $\HA^0\bigl (\MA;
  \MA\hot \Sch(b) \bigr)$.  Thus $\HA_\ast(\MA) \cong \HA_\ast\bigl (\MA \hot
  \Sch(b) \bigr)$, $\HA^\ast(\MA) \cong \HA^\ast\bigl( \MA \hot \Sch(b)
  \bigr)$, and similarly for the bivariant theory.
\end{theorem}

\begin{proof}
  We abbreviate $\Sch \defeq \Sch(b)$.  We construct a chain map $t_{\MA}
  \colon X\bigl( \Tanil(\MA \hot \Sch) \bigr) \to X(\Tanil\MA)$ that is a
  candidate for the inverse of~$[\iota_{\MA}]$.  Since $\Sch$ is tensoring,
  the Tensor Lemma~\ref{lem:tensor}.(iii) yields a natural bounded
  homomorphism $f_{\MA} \colon \Tanil(\MA \hot \Sch) \to (\Tanil\MA) \hot
  \Sch$ that is determined by $f_{\MA} \circ \sigma_{\MA \hot \Sch} =
  \sigma_{\MA} \hot \ID$.  This implies, by the way, that
  $$
  f_{\MA} \bigl( \opt{\ma_0 \otimes T_0} d(\ma_1 \otimes T_1) \dots
  d(\ma_{2n} \otimes T_{2n}) \bigr) =
  \opt{\ma_0} d\ma_1 \dots d\ma_{2n} \otimes
  (\opt{T_0} \cdot T_1 \cdots T_{2n})
  $$
  for $\ma_j \in \MA$, $T_j \in \Sch$.  Thus we get a natural chain map
  $X(f_{\MA}) \colon X\bigl( \Tanil(\MA \hot \Sch) \bigr) \to X\bigl(
  (\Tanil\MA) \hot \Sch \bigr)$.
  
  If~$\MA[B]$ is any complete bornological algebra, we define $\tr_{\MA[B]}
  \colon X(\MA[B] \hot \Sch) \to X(\MA[B])$ by
  $$
  \tr_{\MA[B]}(\ma[b] \otimes T) \defeq  \tr(T) \cdot \ma[b], \qquad
  \tr_{\MA[B]} \bigl( \ma[b]_0 \otimes T_0 d(\ma[b]_1 \otimes T_1) \bigr)
  \defeq \tr (T_0T_1) \ma[b]_0 d\ma[b]_1.
  $$
  The tracial property of $\tr$ easily implies that $\tr_{\MA[B]} \colon
  X(\MA[B] \hot \Sch) \to X(\MA[B])$ is a well-defined chain map.
  
  For $\MA[B] = \Tanil\MA$, this yields a natural bounded chain map $t_{\MA}
  \defeq \tr_{\Tanil\MA} \circ X(f_{\MA}) \colon X\bigl( \Tanil(\MA
  \hot \Sch) \bigr) \to X(\Tanil\MA)$.  We claim that $[t_{\MA}] \in \HA^0(\MA
  \hot \Sch; \MA)$ is inverse to the homology class~$[\iota_{\MA}]$ of the
  chain map $X( \Tanil\iota_{\MA}) \colon X(\Tanil\MA) \to X\bigl( \Tanil (\MA
  \hot \Sch) \bigr)$.
  
  Since~$\iota_{\MA}$ is a homomorphism, one easily checks that $f_{\MA} \circ
  \Tanil\iota_{\MA} = \iota_{\Tanil\MA}$.  Since $\tr(p_0) = 1$ for any rank
  one idempotent, it follows that $t_{\MA} \circ X(\Tanil \iota_{\MA}) = \ID$.
  That is, $[t_{\MA}] \circ [\iota_{\MA}] = 1$ in $\HA^0(\MA; \MA)$.  It
  remains to show that $[\iota_{\MA}] \circ [t_{\MA}] = 1$ in $\HA^0(\MA \hot
  \Sch; \MA \hot \Sch)$.  There are two standard injections $\MA \hot \Sch \to
  (\MA \hot \Sch) \hot \Sch$, namely
  $$
  i_0 \defeq \ID[{\MA}] \hot \ID[{\Sch}] \hot \iota = \iota_{\MA \hot \Sch}
  \colon (\ma \otimes T) \mapsto \ma \otimes T \otimes p_0,
  \quad
  i_1 \defeq \ID[{\MA}] \hot \iota \hot \ID[{\Sch}] =
  \iota_{\MA} \hot \ID[{\Sch}] \colon
  (\ma \otimes T) \mapsto \ma \otimes p_0 \otimes T.
  $$
  We have $[t_{\MA \hot \Sch}] \circ [i_0] = 1$ and $[t_{\MA \hot \Sch}]
  \circ [i_1] = [\iota_{\MA}] \circ [t_{\MA}]$.  The latter equation follows
  formally from the naturality of~$t_{\MA}$.  We will prove that $\iota \hot
  \ID[{\Sch}]$ and $\ID[{\Sch}] \hot \iota$ are smoothly homotopic maps $\Sch
  \to \Sch \hot \Sch$.  It follows that $i_0$ and~$i_1$ are smoothly homotopic
  for any~$\MA$.  Thus $[i_0] = [i_1]$ by Theorem~\ref{the:HA_homotopy}.
  Hence $[\iota_{\MA}] \circ [t_{\MA}] = 1$ as desired.

  It remains to prove that $i_0 = \iota \hot \ID[{\Sch}]$ and $i_1 =
  \ID[{\Sch}] \hot \iota$ are smoothly homotopic.  We define a
  \emph{multiplier}~$\vartheta$ of $\Sch \hot \Sch$ such that conjugation
  $\Ad(\vartheta)$ with~$\vartheta$ implements the coordinate flip $\Sch \hot
  \Sch \to \Sch \hot \Sch$ sending $T_1 \otimes T_2 \mapsto T_2 \otimes T_1$.
  Recall that $\Sch \hot \Sch \cong (\VS \hot \VS[W]) \hot (\VS \hot \VS[W])$
  as a complete bornological vector space.  We define left and right
  multiplication by~$\vartheta$ by
  $$
  \vartheta \cdot (\vs_1 \otimes \vs[w]_1 \otimes \vs_2 \otimes \vs[w]_2)
  \defeq
  \vs_2 \otimes \vs[w]_1 \otimes \vs_1 \otimes \vs[w]_2,
  \qquad
  (\vs_1 \otimes \vs[w]_1 \otimes \vs_2 \otimes \vs[w]_2) \cdot \vartheta
  \defeq
  \vs_1 \otimes \vs[w]_2 \otimes \vs_2 \otimes \vs[w]_1.
  $$
  Left and right multiplication by~$\vartheta$ extend to bounded linear maps
  $\Sch \hot \Sch \to \Sch \hot \Sch$.  Moreover, $\vartheta \cdot (xy) =
  (\vartheta \cdot x) \cdot y$, $(xy) \cdot \vartheta = x \cdot (y \cdot
  \vartheta)$, and $(x \cdot \vartheta) \cdot y = x \cdot (\vartheta \cdot y)$
  for all $x,y \in \Sch \hot \Sch$.  That is, $\vartheta$ is a multiplier of
  $\Sch \hot \Sch$.  In the case of $\Sch(\Hils)$, we have $\Sch(\Hils) \hot
  \Sch(\Hils) \cong \Sch(\Hils \prot \Hils)$ and~$\vartheta$ is the coordinate
  flip $\Hils \prot \Hils \to \Hils \prot \Hils$ viewed as a multiplier of
  $\Sch(\Hils \prot \Hils)$.
  
  We have $\vartheta^2 = 1$, that is, $\vartheta \cdot \vartheta \cdot x = x =
  x\cdot \vartheta \cdot \vartheta$ and $\Ad(\vartheta) (T_1 \otimes T_2)
  \defeq \vartheta \cdot (T_1 \otimes T_2) \cdot \vartheta^{-1} = T_2 \otimes
  T_1$ for all $T_1,T_2 \in \Sch$.  Therefore $i_1 = \Ad(\vartheta) \circ
  i_0$.  There is a natural path of invertible multipliers joining~$\vartheta$
  to the identity, defined by $\vartheta_t \defeq \cos (\pi t/2) \cdot 1 +
  \sin (\pi t/2) \cdot \vartheta$ and $\vartheta_t^{-1} \defeq \cos (\pi t/2)
  \cdot 1 - \sin (\pi t/2) \cdot \vartheta$ for $t \in [0,1]$.  A computation
  using $\vartheta^2 = 1$ shows that $\vartheta_t^{-1} \cdot \vartheta_t =
  \vartheta_t \cdot \vartheta_t^{-1} = 1$ for all~$t$.  Moreover, $\vartheta_0
  = 1$ and $\vartheta_1 = \vartheta$.  The bounded homomorphisms $i_t \defeq
  \Ad(\vartheta_t) \circ i_0$ provide a smooth homotopy between $i_0$
  and~$i_1$ as desired.
\end{proof}

As special cases we get the well-known results that entire cyclic cohomology
is stable with respect to $\Mat[n] \cong \Sch(\C^{\,n})$, the trace class
operators $\Sch(\Hils)$, and the algebra~$\mathfrak{k}$ of matrices with
rapidly decreasing entries.  The latter is obtained by letting $\VS = \VS[W]$
be the space of rapidly decreasing sequences with the pairing $b\bigl( (x_n),
(y_n) \bigr) \defeq \sum x_n \cdot y_n$.

\subsection{Adjoining units}
\label{sec:unitarization}

We show that analytic quasi-freeness is inherited by unitarizations.  As a
consequence, $X\bigl( \Tanil(\Unse{\MA}) \bigr) \cong X(\Tanil\MA) \oplus
\C[0]$ for all complete bornological algebras~$\MA$.  Thus analytic cyclic
cohomology satisfies excision for the particular extension $\MA \injto
\Unse{\MA} \prto \C$.  We need that the algebra~$\C$ is a-quasi-free.  This is
also important for the Chern-Connes character in the even \Mpn{K}theory.

\begin{proposition}  \label{pro:C_aqf}
  The fine algebra~$\C$ of complex numbers is analytically quasi-free.
\end{proposition}

\begin{proof}
  Let $e \in \C$ be the unit element.  A splitting homomorphism $\C \to
  \Tanil\C$ is equivalent to a lifting of~$e$ to an idempotent~$\hat{e}$
  in~$\Tanil\C$.  That is, $\hat{e} \odot \hat{e} = \hat{e}$ and
  $\tau(\hat{e}) = e$.  Our construction of~$\hat{e}$
  follows~\cite{cuntz95:cyclic}.

  Adjoin units to~$\C$ and~$\Tanil\C$ and consider the involution $\gamma =
  1-2e \in \Unse{\C}$.  It suffices to lift~$\gamma$ to an
  involution~$\hat{\gamma}$ in~$\Unse{(\Tanil\C)}$.  That is, $\Unse{\tau}
  (\hat{\gamma}) = \gamma$ and $\hat{\gamma} \odot \hat{\gamma} = 1$.  We make
  the Ansatz $\hat{\gamma} \defeq (1+n) \odot \bigl(1 - 2\sigma(e)\bigr)$ with
  $n \in \Janil\C$.  This is always a lifting for~$\gamma$ and satisfies
  $$
  \hat{\gamma} \odot \hat{\gamma} =
  (1+n)^2 \odot \bigl(1 - 2\sigma(e) \bigr) \odot
  \bigl(1 - 2\sigma(e) \bigr) = 
  (1+n)^2 (1 - 4dede)
  $$
  because $\Unse{(\Tanil\C)}$ is commutative.  Hence~$\hat{\gamma}$ is an
  involution iff $n = (1-4dede)^{-1/2} - 1$.  The function $f \colon z \mapsto
  (1-4z)^{-1/2} - 1$ is holomorphic in a neighbourhood of the origin and
  satisfies $f(0) = 0$.  Since $\Janil\C$ is a-nilpotent, $n = f(dede)$ is
  well-defined by Lemma~\ref{lem:functional_calculus} and solves our problem.
  Thus~$\C$ is a-quasi-free.  Explicitly, the power series of~$f$ is
  $$
  f(z) =
  \sum_{j=1}^\infty \binom{-1/2}{j} (-4z)^j =
  \sum_{j=1}^\infty \frac{(2j)!}{(j!)^2} z^j =
  \sum_{j=1}^\infty \binom{2j}{j} z^j.
  $$
  Thus
  \begin{equation}  \label{eq:lift_idempotent}
    \hat{e} =
    \tfrac{1}{2}(1 - \hat{\gamma}) =
    e + \bigl(e - \tfrac{1}{2}\bigr) f(dede) =
    e + \sum_{j=1}^\infty
    \binom{2j}{j} e(de)^{2j} - \frac{1}{2} \binom{2j}{j} (de)^{2j}.
  \end{equation}
  
  By the way, the binomial coefficients $\binom{2j}{j}$ are always even.  Thus
  the coefficients in \eqref{eq:lift_idempotent} are integers.
\end{proof}

\begin{proposition}  \label{pro:unse_aqf}
  A complete bornological algebra~$\MA$ is a-quasi-free iff~$\Unse{\MA}$ is
  a-quasi-free.
\end{proposition}

\begin{proof}
  Let~$\MA$ be a-quasi-free.  We have to construct a bounded splitting
  homomorphism $\Unse{\MA} \to \Tanil(\Unse{\MA})$.  Since the algebra~$\C$ is
  a-quasi-free by Proposition~\ref{pro:C_aqf}, we can lift the homomorphism
  $\C \to \Unse{\MA}$ sending $1 \mapsto 1_{\Unse{\MA}}$ to a homomorphism $\C
  \to \Tanil(\Unse{\MA})$ by Theorem~\ref{the:aqf_definitions}.(iii).  That
  is, we can lift~$1_{\Unse{\MA}}$ to an idempotent $e \in
  \Tanil(\Unse{\MA})$.  Thus $e \odot e = e$.  We cut down
  $\Tanil(\Unse{\MA})$ by this idempotent.  Since $\tau_{\MA}(e) =
  1_{\Unse{\MA}}$, we get an allowable extension of complete bornological
  algebras
  $$
  e \odot \Janil(\Unse{\MA}) \odot e \injto
  e \odot \Tanil(\Unse{\MA}) \odot e \prto
  \Unse{\MA}.
  $$
  A bounded linear section for this extension is $x \mapsto e \odot
  \sigma_{\Unse{\MA}}(x) \odot e$.  Since~$\MA$ is a-quasi-free, another
  application of Theorem~\ref{the:aqf_definitions}.(iii) allows us to lift the
  natural homomorphism $\MA \to \Unse{\MA}$ to a bounded homomorphism
  $\upsilon' \colon \MA \to e \odot \Tanil(\Unse{\MA}) \odot e$.  We
  extend~$\upsilon'$ to a bounded linear section $\upsilon \colon \Unse{\MA}
  \to \Tanil(\Unse{\MA})$ by putting $\upsilon(1) \defeq e$ and $\upsilon(x)
  \defeq \upsilon'(x)$ for $x \in \MA$.  This map~$\upsilon$ is a homomorphism
  because $e \odot \upsilon'(x) = \upsilon'(x) \odot e = \upsilon'(x)$ for all
  $x \in \MA$.  Thus~$\Unse{\MA}$ is a-quasi-free.

  Conversely, assume that~$\Unse{\MA}$ is a-quasi-free.  We have to construct
  a bounded splitting homomorphism $\MA \to \Tanil\MA$.
  Theorem~\ref{the:aqf_definitions}.(iii) provides a bounded splitting
  homomorphism $\upsilon \colon \Unse{\MA} \to \Unse{(\Tanil\MA)}$ for the
  a-nilpotent extension
  $$
  (i, \Unse{\tau}) \colon \Janil\MA \injto \Unse{(\Tanil\MA)} \prto \Unse{\MA}.
  $$
  Since $\Unse{\tau} \circ \upsilon(x) = x \in \MA$ for all $x \in \MA$, the
  restriction of~$\upsilon$ to~$\MA$ has values in $\Tanil\MA \subset
  \Unse{(\Tanil\MA)}$.  Therefore, $\upsilon|_{\MA}$ is a bounded splitting
  homomorphism $\MA \to \Tanil\MA$ and~$\MA$ is a-quasi-free.
\end{proof}

\begin{theorem}  \label{the:HA_unitarization}
  Let~$\MA$ be a complete bornological algebra.  Then $\Unse{(\Tanil\MA)}$ and
  $\Tanil(\Unse{\MA})$ are both analytically quasi-free and smoothly homotopy
  equivalent.  The complex $X\bigl( \Tanil(\Unse{\MA}) \bigr)$ is naturally
  chain homotopic to the complex $X(\Tanil\MA) \oplus \C[0]$.  Thus there are
  natural isomorphisms $\HA^\ast(\Unse{\MA}) \cong \HA^\ast(\MA) \oplus
  \C[0]$, $\HA_\ast(\Unse{\MA}) \cong \HA_\ast(\MA) \oplus \C[0]$,
  $\HA^\ast(\MA[B]; \Unse{\MA}) \cong \HA^\ast(\MA[B]) \oplus \HA^\ast(\MA[B];
  \MA)$, and $\HA^\ast(\Unse{\MA}; \MA[B]) \cong \HA_\ast(\MA[B]) \oplus
  \HA^\ast(\MA; \MA[B])$.
\end{theorem}

\begin{proof}
  By Proposition~\ref{pro:unse_aqf}, the extension $\Janil\MA \injto
  \Unse{(\Tanil\MA)} \prto \Unse{\MA}$ is a universal a-nilpotent extension
  of~$\Unse{\MA}$.  The Uniqueness Theorem~\ref{the:uniqueness} implies that
  $\Unse{(\Tanil\MA)}$ and $\Tanil(\Unse{\MA})$ are naturally smoothly
  homotopy equivalent.  Theorem~\ref{the:uniqueness_X} implies that $X( \Tanil
  \bigl( \Unse{\MA}) \bigr)$ is naturally chain homotopic to $X(\Unse
  {(\Tanil\MA)} )$.
  
  We claim that $X(\Unse{\MA[B]}) \cong X(\MA[B]) \oplus \C[0]$ for any
  complete bornological algebra~$\MA[B]$.  Clearly, $X_0(\Unse{\MA[B]}) =
  \Unse{\MA[B]} \cong \MA[B] \oplus \C$.  Let $1 \in \Unse{\MA[B]}$ be the
  unit and let~$V$ be the linear span of terms $xd1$ and $dx - 1dx$ with $x
  \in \Unse{\MA[B]}$.  We will prove that $b \bigl( \Omega^2(\Unse{\MA})
  \bigr) = b(\Omega^2\MA) + V$.  Since $\Omega^1(\Unse{\MA[B]}) \cong
  \Omega^1\MA[B] \oplus V$, this implies that $X_1(\MA[B]) \cong
  X_1(\Unse{\MA[B]})$.  The claim $X(\Unse{\MA[B]}) \cong X(\MA[B]) \oplus
  \C[0]$ follows.  Since
  $$
  b(1 d1dx - d1dx) = dx - 1 dx, \qquad
  b( xd1d1) = xd1
  $$
  for all $x \in \Unse{\MA[B]}$, we have $V \subset b \bigl(
  \Omega^2(\Unse{\MA}) \bigr)$.  Computing the commutators $b(dxd1) = b(d1dx)
  = xd1 - (dx - 1dx)$, $b(1dxdy - dxdy) = d(xy) - 1 d(xy)$, $b(xd1dy) = yx
  d1$, $b(xdyd1) = xy d1$, we find that $b\bigl( \Omega^2(\Unse{\MA[B]})
  \bigr) \subset V + b(\Omega^2\MA)$.
  
  Thus $X\bigl( \Tanil (\Unse{\MA}) \bigr) \sim X(\Unse {(\Tanil\MA)} ) \cong
  X(\Tanil\MA) \oplus \C[0]$.  The isomorphisms of the associated
  (co)\penalty10000\hspace*{0pt}homology groups listed above follow
  immediately.
\end{proof}

\subsection{The Chern-Connes character in \Mpn{K}theory}
\label{sec:Chern_Ktheory}

The definition of the Chern-Connes character in the even \Mpn{K}theory is an
application of the analytic quasi-freeness of the complex numbers~$\C$.
Similarly, the quasi-freeness of the algebra of Laurent polynomials yields the
Chern-Connes character in odd \Mpn{K}theory.  Our construction of the
Chern-Connes character follows the construction of Cuntz and
Quillen~\cite{cuntz95:cyclic} for periodic cyclic homology.

Idempotents in an algebra~$\MA$ are the same thing as (bounded) homomorphisms
$\C \to \MA$.  An idempotent~$e$ corresponds to the bounded homomorphism
$\vec{e} \colon \lambda \mapsto \lambda \cdot e$.  Let $\MA[N] \injto \MA[E]
\prto \MA$ be an a-nilpotent extension and let $e \in \MA$ be an idempotent.
Since~$\C$ is a-quasi-free by Proposition~\ref{pro:C_aqf}, we can apply
Theorem~\ref{the:aqf_definitions} to lift the homomorphism $\vec{e}$ to a
bounded homomorphism $\C \to \MA[E]$.  That is, we can lift~$e$ to an
idempotent $\hat{e} \in \MA[E]$.  Furthermore, the uniqueness assertion of
Theorem~\ref{the:aqf_definitions} implies that~$\hat{e}$ is unique up to
smooth homotopy.  Two idempotents in~$\MA[E]$ are smoothly homotopic iff they
can be joined by a smooth path of idempotents iff the corresponding
homomorphisms $\C \to \MA[E]$ are smoothly homotopic.  Furthermore, the
Universal Extension Theorem~\ref{the:universal_ext} implies that liftings of
smoothly homotopic idempotents remain smoothly homotopic.  Thus the projection
$\MA[E] \to \MA$ induces a bijection between smooth homotopy classes of
idempotents in $\MA[E]$ and~$\MA$.

In particular, we can apply this to the extension $\Janil\MA \injto \Tanil\MA
\prto \MA$ to lift idempotents in~$\MA$ to idempotents in $\Tanil\MA$.  In
this case, there is a canonical way to lift idempotents.  Let $f \in \MA$ be
an idempotent.  The corresponding homomorphism $\vec{f} \colon \C \to \MA$
induces a natural homomorphism $\Tanil \vec{f} \colon \Tanil\C \to \Tanil\MA$
by Lemma~\ref{lem:lanilcur_compose_hom}.  We compose this with the splitting
homomorphism $\upsilon \colon \C \to \Tanil\C$ constructed in the proof of
Proposition~\ref{pro:C_aqf}.  This yields a homomorphism $\C \to \Tanil\MA$
lifting~$\vec{f}$.  The corresponding idempotent~$\hat{f}$ in $\Tanil\MA$ is
obtained by evaluating $\Tanil \vec{f}$ on the idempotent $\hat{e} \in
\Tanil\C$ described in~\eqref{eq:lift_idempotent}.  Thus
$$
\hat{f} =
f + \sum_{j=1}^\infty
\binom{2j}{j} \bigl(f - \tfrac{1}{2}\bigr) \cdot (df)^{2j}.
$$
Since $\hat{f}^2 = \hat{f}$, we have $\partial_0(\hat{f}) = 0$.
Thus~$\hat{f}$ defines a class $\chern(f) \in H_0\bigl( X(\Tanil\MA) \bigr) =
\HA_0(\MA)$.

We have an alternative description of $\chern(f)$.  The homomorphism $\upsilon
\colon \C \to \Tanil\C$ constructed in the proof of
Proposition~\ref{pro:C_aqf} induces a homotopy equivalence of complexes $X(\C)
\congto X(\Tanil\C)$.  Composition with $X(\upsilon)$ therefore gives rise to
an isomorphism $\HA^\ast(\C; \MA) \congto H^\ast \bigl( X(\C); X(\Tanil\MA)
\bigr)$.  Identifying $X(\C) \cong \C[0]$, we obtain an isomorphism $H_0\bigl(
X(\C); X(\Tanil\MA) \bigr) \cong H_0\bigl( X(\Tanil\MA) \bigr) \cong
\HA_0(\MA)$.  The homology class $\chern(f) \in \HA_0(\MA)$ is the image of
the class $[\vec{f}] \in \HA_0(\C; \Tanil\MA)$ under this isomorphism.  This
definition shows immediately that $\chern(f_0) = \chern(f_1)$ if $f_0$
and~$f_1$ are AC-homotopic.

Before we can write down the Chern-Connes character properly, we should define
what we mean by \Mpn{K}theory.  We can use algebraic \Mpn{K}theory, but this
is not really the right domain of definition for the Chern-Connes character.
Rather, we should use a \emph{topological \Mpn{K}theory} for complete
bornological algebras.  The definition of Phillips~\cite{phillips91:K} of a
topological \Mpn{K}theory for \Frechet algebras carries over immediately to
complete bornological algebras.  However, it is not clear whether the
resulting theory has good homological properties (like excision).

Let~$\MA$ be a complete bornological algebra.  Let~$\mathfrak{k}$ be the
algebra of matrices $(x_{i,j})_{i,j \in \N}$ with rapidly decreasing entries.
This algebra is used by Phillips~\cite{phillips91:K} and
Cuntz~\cite{cuntz97:bivariant}.  We construct the algebra $\tilde{\MA} \defeq
\Mat[2] \hot \Unse{(\MA \hot \mathfrak{k})}$.  That is, we first
stabilize~$\MA$ by~$\mathfrak{k}$, then adjoin a unit to the stabilization, and
finally take the \Mp{2 \times 2}matrices over this unitarization.  Let
$I(\MA)$ be the set of all idempotents $e \in \tilde{\MA}$ for which $e -
\bigl( \begin{smallmatrix} 1 & 0 \\ 0 & 0 \end{smallmatrix} \bigr) \in \MA
\hot \mathfrak{k}$.  Two elements of $I(\MA)$ are smoothly homotopic iff they
can be joined by a smooth path inside $I(\MA)$.  The \emph{topological
  \Mpn{K}theory} $K_0(\MA)$ of~$\MA$ is defined as the set of smooth homotopy
classes of elements in $I(\MA)$.  There is a natural addition on $K_0(\MA)$
that makes it an Abelian group.  This is the topological \Mpn{K}theory studied
by Phillips for locally multiplicatively convex \Frechet algebras.

We define the Chern-Connes character $\chern \colon K_0(\MA) \to \HA_0(\MA)$
as follows.  Let $e \in I(\MA)$.  This defines $\chern(e) \in
\HA_0(\tilde{\MA})$ as above.  The stability of $\HA_0$ with respect to
$\Mat[2]$ and~$\mathfrak{k}$ (Theorem~\ref{the:HA_stable}) and
Theorem~\ref{the:HA_unitarization} yield a natural isomorphism
$$
\HA_0(\tilde{\MA}) \cong
\HA_0( \Unse{(\MA \hot \mathfrak{k})} ) \cong
\HA_0(\MA \hot \mathfrak{k}) \oplus \C[0] \cong
\HA_0(\MA) \oplus \C[0].
$$
We define $\chern(e) \in \HA_0(\MA)$ as the component of $\chern(e) \in
\HA_0(\tilde{\MA})$ in $\HA_0(\MA)$.  By the way, $e - \bigl(
\begin{smallmatrix} 1 & 0 \\ 0 & 0 \end{smallmatrix} \bigr) \in \MA \hot
\mathfrak{k}$ implies that the component in $\C[0]$ is always~$1$.  This
procedure defines a natural group homomorphism $K_0(\MA) \to \HA_0(\MA)$.

Phillips's definition of $K_1(\MA)$ is as follows.  Let $U(\MA)$ be the set of
elements $x \in \MA \hot \mathfrak{k}$ such that $1 + x \in \Unse{(\MA \hot
  \mathfrak{k})}$ is invertible.  The smooth homotopy classes of elements in
$U(\MA)$ form an Abelian group $K_1(\MA)$.  We call an element $x \in \MA$
\emph{invertible$-1$} iff $1 + x \in \Unse{\MA}$ is invertible.  Let~$\MA[U]$
be the universal algebra generated by an invertible$-1$ element.  That is,
$\MA[U]$ is the universal algebra on two generators $a,b$ satisfying $(1+b)
\cdot (1+a) = (1+a) \cdot (1+b) = 1$, that is, $a+b+ab = a+b+ba = 0$.  We can
realize~$\MA[U]$ concretely as a maximal ideal
$$
\MA[U] \cong \{ f\in \C[u, u^{-1}] \mid f(1) = 0 \}
$$
in the fine algebra $\C[u,u^{-1}]$ of Laurent polynomials.  The generators
$a,b$ correspond to $u-1$ and $u^{-1}-1$, respectively.

\begin{proposition}  \label{pro:Laurent_aqf}
  The fine algebra~$\MA[U]$ is analytically quasi-free.
  
  The X-complex of~$\MA[U]$ is chain homotopic to the complex $\C[1]$.  The
  \Mpn{1}form $(1+b) da \bmod [,]$ is a cycle and generates the homology of
  $X(\MA[U])$.
\end{proposition}

\begin{proof}
  By the universal property of~$\MA[U]$, a lifting homomorphism $\MA[U] \injto
  \Tanil\MA[U]$ is determined by liftings $\hat{a}, \hat{b} \in \Tanil\MA[U]$
  of $a$ and~$b$ (that is, $\tau(\hat{a}) = a$, $\tau(\hat{b}) = b$)
  satisfying $(1+\hat{a}) (1+\hat{b}) = (1+\hat{b}) (1+\hat{a}) = 1$.

  Lift~$a$ to $\hat{a} \defeq \sigma(a)$.  Then $\bigl( 1 + \sigma(b) \bigr)
  \odot (1 + \hat{a}) = 1 + \sigma(b) + \sigma(a) + \sigma(ba) - dbda = 1 -
  dbda$.  Similarly, $(1 + \hat{a}) \odot \bigl( 1 + \sigma(b) \bigr) =
  1-dadb$.  Since $\Janil\MA[U]$ is a-nilpotent, both $1 - dbda$ and $1 -
  dadb$ are invertible.  Thus $(1 + \hat{a})$ is invertible with inverse
  $$
  1+ \hat{b} \defeq (1 + \hat{a})^{-1} =
  \bigl( 1 + \sigma(b) \bigr) \odot (1-dadb)^{-1} =
  (1 + b) \sum_{j=0}^\infty (dadb)^{j}.
  $$
  Thus
  \begin{equation}  \label{eq:lift_invertible}
    a \mapsto \hat{a} = a,  \qquad
    b \mapsto \hat{b} = b + \sum_{j=1}^\infty (dadb)^j + b(dadb)^j
  \end{equation}
  defines a lifting homomorphism $\MA[U] \to \Tanil\MA[U]$.
  
  It is clear that $X_0(\MA[U]) = \MA[U] \subset \C[u, u^{-1}]$.  We claim
  that the odd part $X_1(\MA[U]) = \Omega^1\MA[U] / [,]$ is isomorphic to
  $\C[u, u^{-1}]$ via the map $\natural \colon f\,dg \bmod [,] \mapsto f\cdot
  g'$.  Since~$\MA[U]$ is commutative, the derivation rule for~$d$ implies
  that $d(f^n) = n f^{n-1} df + [,]$ for all $f \in \MA[U]$, $n \in \Z$.  It
  follows that $f\,dg = f \cdot g' \,da \bmod [,]$ for all $f \in \C[u,
  u^{-1}]$ and $g \in \MA[U]$.  Moreover,
  $$
  \natural \circ b (\opt{f_0}\, df_1 df_2) =
  \opt{f_0} \cdot f_1 \cdot f_2' - \opt{f_0} \cdot (f_1\cdot f_2)' +
  \opt{f_0} \cdot f_1' \cdot f_2 =
  0
  $$
  for all $f_j \in \MA[U]$.  Thus~$\natural$ is a well-defined map
  $X_1(\MA) \to \C[u, u^{-1}]$.  It is bijective.
  
  The boundary $\partial_1 \colon X_1(\MA) \to X_0(\MA)$ is the zero map,
  whereas the boundary $\partial_0 \colon X_0(\MA) \to X_1(\MA)$ is the map $f
  \mapsto f'$.  The map~$\partial_0$ is injective because $f' = 0$ implies
  that~$f$ is constant and thus $f = f(1) = 0$.  The range of~$\partial_0$ is
  the linear span of~$u^n$, $n \in \Z \setminus \{-1\}$.  Hence $X(\MA[U])$ is
  the sum of a contractible complex and the complex $\C[1]$.  The copy
  of~$\C[1]$ is generated by $u^{-1} da = (1+b)da = da + bda$.
\end{proof}

Let $x \in \MA$ be invertible$-1$.  Then $a \mapsto x$, $b \mapsto (1+x)^{-1} -
1$ defines a homomorphism $\vec{x} \colon \MA[U] \to \MA$.  Let $[\vec{x}] \in
\HA^0(\MA[U]; \MA)$ be its homology class.  By
Proposition~\ref{pro:Laurent_aqf}, the algebra~$\MA[U]$ is a-quasi-free.  A
natural bounded splitting homomorphism $\upsilon \colon \MA[U] \to
\Tanil\MA[U]$ is described by~\eqref{eq:lift_invertible}.
Theorem~\ref{the:uniqueness_X} implies that $X(\upsilon) \colon X(\MA[U]) \to
X(\Tanil\MA[U])$ is a homotopy equivalence.  Let $c_1 \in \C \setminus \{0\}$
be a constant.  Then the map $f \colon \C[1] \to X(\MA[U])$ sending $1 \mapsto
c_1\cdot (1+b)da$ is a homotopy equivalence by the second part of
Proposition~\ref{pro:Laurent_aqf}.  Hence composition with $X(\upsilon) \circ
f$ gives rise to an isomorphism
$$
\HA^0(\MA[U]; \MA) =
H^0\bigl( X(\Tanil\MA[U]); X(\Tanil\MA) \bigr) \cong
H^0\bigl( \C[1]; X(\Tanil\MA) \bigr) \cong
H_1\bigl( X(\Tanil\MA) \bigr) \cong
\HA_1(\MA).
$$
We can now apply this to the class $[\vec{x}] \in \HA^0(\MA[U]; \MA)$ of an
invertible$-1$ element~$x$.  The resulting homology class $\chern(x) \in
\HA_1(\MA)$ is represented by the cycle
$$
X(\Tanil\vec{x}) \circ X(\upsilon) (c_1 (1+b)da \bmod [,]) =
c_1 \sum_{j=0}^\infty (1+x)^{-1} (dxdy)^j dx =
c_1 \sum_{j=0}^\infty (dx+ ydx) (dydx)^j
$$
with $y = (1+x)^{-1} - 1 = \vec{x}(b)$.

In the definition of this map, we left a free constant~$c_1$ because, unlike
in the even case, the obvious choice $c_1 = 1$ has some drawbacks.  Its
advantage is that it gives rise to integer coefficients because all
coefficients in~\eqref{eq:lift_invertible} are integers.  Its disadvantage is
that it does not combine well with Bott periodicity.  The integral of the
\Mpn{1}form $u^{-1} du$ over the unit circle is $2\pi i$ by the Cauchy
integral formula.  Hence the constant $2\pi i$ likes to come up even though we
have integer coefficients.  This cannot be helped altogether.  To combine well
with Bott periodicity and the bivariant Chern-Connes character, the best
choice is the compromise $c_1 = (2\pi i)^{-1/2}$ between $1$ and $(2\pi
i)^{-1}$ (see~\cite{cuntz97:bivariant}).  This is also the choice of
Connes~\cite{connes94:ncg}, when translated from his \Mpn{(d_1, d_2)}bicomplex
to $X(\Tanil\MA)$.  Since we do not consider the bivariant Chern-Connes
character here, any choice for~$c_1$ is good for us.

Alternatively, the Chern-Connes character can be described as follows.  Let
$\hat{x} \in \Tanil\MA$ be any lifting of~$x$.  Then~$\hat{x}$ is
invertible$-1$ and $a \mapsto \hat{x}$ defines a lifting homomorphism $\psi
\colon \MA[U] \to \Tanil\MA$ lifting $\vec{x}$.  We let $\chern(x)$ be the
homology class of the cycle $c_1(1+\hat{x})^{-1} \,D\hat{x} \bmod [,] \in
X_1(\Tanil\MA)$.  This is the image of the cycle $c_1(1+b) da \in X_1(\MA[U])$
under $X(\psi)$ and therefore a cycle.  Its homology class does not depend
on~$\psi$ because~$\psi$ is unique up to smooth homotopy by the Universal
Extension Theorem~\ref{the:universal_ext}.  Thus $X(\psi)$ is unique up to
chain homotopy by Proposition~\ref{pro:X_homotopy} because~$\MA[U]$ is
quasi-free.

We can now extend this construction to a Chern-Connes character $K_1(\MA) \to
\HA_1(\MA)$ as in the even case.  Represent an element of $K_1(\MA)$ by an
invertible$-1$ element in $\MA \hot \mathfrak{k}$.  Take its Chern-Connes
character in $\HA_1(\MA \hot \mathfrak{k})$ and apply stability to obtain an
element of $\HA_1(\MA)$.

For any complex, there is a natural pairing between homology and cohomology by
evaluating linear functionals on elements.  For the complex $X(\Tanil\MA)$,
this yields a natural pairing $\HA_\ast(\MA) \times \HA^\ast(\MA) \to \C$ for
$\ast = 0,1$.  If we combine this with the Chern-Connes character, we get a
natural pairing $K_\ast(\MA) \times \HA^\ast(\MA) \to \C$ for $\ast = 0,1$.

\subsection{The X-complex of $\Tanil\MA$ and entire cyclic cohomology}
\label{sec:X_Tanil}

By Proposition~\ref{pro:Omega_I_Tanil}, we have a natural bimodule
homomorphism $\Omega^1(\Tanil\MA) \cong \Unse{(\Tanil\MA)} \hot \MA \hot
\Unse{(\Tanil\MA)}$.  Thus $\Omega^1(\Tanil\MA) / [,]$ is bornologically
isomorphic to $\Unse{(\Tanil\MA)} \hot \MA \cong \Omega_\an^\odd\MA$ via the
map $\opt{x} (D\ma) \opt{y} \bmod [,] \mapsto \opt{y} \odot \opt{x} d\ma$ for
all $x,y \in \Tanil\MA$, $\ma \in \MA$.  Combining this with the definition
$X_0(\Tanil\MA) = \Tanil\MA = \Omega_\an^\even\MA$, we get a grading
preserving bornological isomorphism between the underlying vector spaces of
$X(\Tanil\MA)$ and $\Omega_\an\MA$.  In the following, we use this isomorphism
to identify $X(\Tanil\MA)$ and $\Omega_\an\MA$.  We write~$\partial$ for the
boundary of $X(\Tanil\MA)$ transported to $\Omega_\an\MA$.  Cuntz and
Quillen~\cite{cuntz95:cyclic} obtain explicit formulas for the boundary
$\Omega\MA \to \Omega\MA$ coming from the boundary $X(\Tens\MA) \to
X(\Tens\MA)$.  Since~$\partial$ is bounded on $\Omega_\an\MA$, it must be the
extension to the completion of the boundary on $\Omega\MA$.  Hence we can copy
the formulas for~$\partial$ in~\cite{cuntz95:cyclic}.  They involve the
\emph{Karoubi operator}~$\kappa$ defined by $\kappa (\omega \,d\ma) \defeq
(-1)^{\deg \omega} d\ma \cdot \omega$.  In Appendix~\ref{app:kappa}, we recall
some identities satisfied by~$\kappa$ proved already by Cuntz and
Quillen~\cite{cuntz95:operators}.
\begin{align}
  \label{eq:partial_odd}
  \partial_1 & = b - (1 + \kappa) d
  \qquad \text{on $\Omega_\an^\odd\MA$};
  \\
  \label{eq:partial_even}
  \partial_0 & = \sum_{j=0}^{2n} \kappa^j d - \sum_{j=0}^{n-1} \kappa^{2j} b
  \qquad \text{on $\Omega^{2n}\MA$}.
\end{align}
To assist those readers not familiar with~\cite{cuntz95:cyclic}, we prove
these identities in Appendix~\ref{app:X_boundary} together with a concrete
proof of the isomorphism $\Omega^1(\Tanil\MA) \cong \Unse{(\Tanil\MA)} \hot
\MA \hot \Unse{(\Tanil\MA)}$.  The operator~$B$ on $\Omega\MA$ is defined by
$\sum_{j=0}^n \kappa^j d$ on $\Omega^n\MA$.  It extends to a bounded operator
$B \colon \Omega_\an\MA \to \Omega_\an\MA$.  This operator corresponds to the
operator~$B$ introduced by Connes.  It satisfies $B b = -b B = 0$ and $B B =
0$ (see Appendix~\ref{app:kappa}).  Thus the two boundaries $(b,B)$ make
$\Omega_\an\MA$ a bicomplex.

In~\cite{cuntz95:cyclic}, an explicit chain homotopy between the complexes
$(\Omega\MA, B+b)$ and $(\Omega\MA, \partial)$ is constructed.  This chain
homotopy involves multiplication by~$n!$ on $\Omega^{2n}\MA$ and
$\Omega^{2n+1}\MA$.  Since~$n!$ grows faster than any exponential, this chain
homotopy does not extend to $\Omega_\an\MA$.  However, if we modify the
definition of the bornology $\CBS_\an$, we can carry over the argument
in~\cite{cuntz95:cyclic}.  Let $n! \CBS_\an$ be the convex bornology on
$\Omega\MA$ generated by the sets $S \cup \bigcup_{n=1}^\infty n! \opt{S}
\opt{dS} (dS)^{2n}$ with $S \in \CBS(\MA)$.  Let $C(\MA)$ be the completion of
$(\Omega\MA, n!  \CBS_\an)$.  Define $c \colon \Z_+ \to \C$ by $c_{2n}\defeq
c_{2n+1}\defeq (-1)^n n!$.  By definition, $c(N) \colon \Omega\MA \to
\Omega\MA$ is a bornological isomorphism $(\Omega\MA, \CBS_\an) \cong
(\Omega\MA, n! \CBS_\an)$ and therefore extends to a bornological isomorphism
$\Omega_\an\MA \cong C(\MA)$.

The operators $B$ and~$b$ extend to bounded operators $C(\MA) \to C(\MA)$.  In
fact, a \emph{homogeneous} operator is bounded with respect to the bornology
$\CBS_\an$ iff it is bounded with respect to the bornology $n!  \CBS_\an$.  A
linear map $l \colon \Omega\MA \to \Omega\MA$ is \emph{homogeneous of
  degree~$k$} iff $l(\Omega^n\MA) \subset \Omega^{n+k}\MA$ for all $n \in \N$.
Since~$c(N)$ is a bornological isomorphism, $l$ is bounded with respect to $n!
\CBS_\an$ iff $c(N) \circ l \circ c(N)^{-1}$ is bounded with respect
to~$\CBS_\an$.  Since~$l$ has degree~$k$, we have
$$
c(N) \circ l \circ c(N)^{-1} = \frac{c(n+k)}{c(n)} \cdot l
\qquad \text{on $\Omega^n\MA$.}
$$
Both $c(n+k)/c(n)$ and $c(n)/c(n+k)$ are of polynomial growth in~$n$.  Thus
rescaling with these numbers is a bornological isomorphism with respect to
$\CBS_\an$.  Hence $c(N) \circ l \circ c(N)^{-1}$ is bounded with respect
to~$\CBS_\an$ iff~$l$ is bounded with respect to~$\CBS_\an$.

Consequently, $B+b$ extends to a bounded linear map $C(\MA) \to C(\MA)$ and
makes $(C(\MA), B+b)$ a \Mpn{\Ztwo}graded complex of complete bornological
vector spaces.

\begin{proposition}  \label{pro:X_equals_C}
  The complexes $X(\Tanil\MA)$ and $(C(\MA), B+b)$ are chain homotopic.
\end{proposition}

\begin{proof}
  The proof follows Cuntz and Quillen~\cite{cuntz95:cyclic}.  We only have to
  verify that certain maps are bounded.  This is straightforward if we write
  down explicit formulas for these maps.  Some details are relegated to
  Appendix~\ref{app:X_boundedness}.

  Let $\delta \colon \Omega_\an\MA \to \Omega_\an\MA$ be the boundary that
  corresponds to $B+b$ under the isomorphism $c(N) \colon \Omega_\an\MA \to
  C(\MA)$.  That is, $\delta = c(N)^{-1} \circ (B+b) \circ c(N)$.  Using the
  definition $c(2n) = c(2n+1) = (-1)^n n!$, we easily compute that
  \begin{equation}  \label{eq:def_delta}
    \delta = B - nb \quad \text{on $\Omega^{2n}\MA$,}  \qquad
    \delta = -B/(n+1) + b \quad \text{on $\Omega^{2n+1}\MA$.}
  \end{equation}
  We have to show that the complexes $(\Omega_\an\MA, \delta)$ and
  $(\Omega_\an\MA, \partial)$ are chain homotopic.
  
  The operator~$\kappa$ satisfies a polynomial identity $(\kappa^n-1)
  (\kappa^{n+1}-1) = 0$ on $\Omega^n\MA$ for all $n \in \N$.  Hence~$\kappa^2$
  satisfies such an identity, too.  It follows that the restriction
  of~$\kappa^2$ to $\Omega^n\MA$ has discrete spectrum.  One of the
  eigenvalues is~$1$.  This eigenvalue has multiplicity~$2$, that is, the
  corresponding eigenspace is equal to the kernel of $(\kappa^2 - 1)^2$.
  Let~$P$ be the projection onto the \Mpn{1}eigenspace annihilating all other
  eigenspaces.  Let~$H$ be the operator that is zero on the \Mpn{1}eigenspace
  and the inverse of $1 - \kappa^2$ on all other eigenspaces.  Explicit
  formulas for $P$ and~$H$ show that they are bounded with respect
  to~$\CBS_\an$ (see Appendix~\ref{app:X_boundedness}).  Any operator
  commuting with~$\kappa$ commutes with~$P$ and~$H$.  Thus $P$ and~$H$ are
  chain maps with respect to $\partial$ and~$\delta$.
  
  The restriction of~$\partial$ to the \Mpn{1}eigenspace of~$\kappa^2$ is
  equal to~$\delta$.  If the eigenvalue~$1$ had no multiplicity, we could
  derive this immediately by replacing~$\kappa^2$ with~$1$ in
  \eqref{eq:partial_odd} and~\eqref{eq:partial_even}.  Since $\kappa^2 - 1
  \neq 0$ on $P(\Omega\MA)$, we have to be a bit more careful.  See
  Appendix~\ref{app:X_boundedness} or the proof of
  \cite[Lemma~6.1]{cuntz95:cyclic}.

  The proof will therefore be finished if we show that~$P$ is chain homotopic
  to the identity with respect to the boundaries $\partial$ and~$\delta$.
  Since $\ID - P = H \circ (\ID - \kappa^2)$ and~$H$ is a chain map, this
  follows if $\ID - \kappa^2$ is chain homotopic to~$0$.  We write~$\sim$ for
  chain homotopy.  Equation~\eqref{eq:kappa_i} implies that
  $$
  [\delta, c(N)^{-1} d c(N)] =
  c(N)^{-1} [b+B, d] c(N) =
  c(N)^{-1} (db +bd) c(N) =
  c(N)^{-1} (1 - \kappa) c(N) =
  1 - \kappa
  $$
  because $1 - \kappa$ is homogeneous of degree~$0$ and therefore commutes
  with~$c(N)$.  Thus $\kappa^2 \sim \ID \circ \ID = \ID$ with respect to the
  boundary~$\delta$.  We have $(b- (1+\kappa) d)^2 = - (1+\kappa) (bd+db) = -
  (1+\kappa) (1-\kappa) = \kappa^2 - 1$ by~\eqref{eq:kappa_i}.  However,
  $\partial \circ \partial = 0$ implies easily that $[\partial, b-(1+\kappa)d]
  = (b - (1+\kappa)d)^2$.  For example, if we evaluate $[\partial, b-
  (1+\kappa)d]$ on $\Omega^\even\MA$, then we get $\partial_1 \circ (b -
  (1+\kappa)d) + (b - (1+\kappa)d) \circ \partial_0 = (b - (1+\kappa)d)^2 +
  \partial_1 \circ \partial_0$.  Thus $\kappa^2 \sim \ID$ with respect to the
  boundary~$\partial$.
\end{proof}

By the way, since $\kappa \sim \ID$ with respect to~$\delta$, this holds with
respect to the boundary~$\partial$ as well.

\begin{theorem}  \label{the:entire_equals_analytic}
  Let~$\MA$ be a locally convex algebra.  Endow~$\MA$ with the bounded
  bornology.  The analytic cyclic cohomology $\HA^\ast(\MA)$ is naturally
  isomorphic to the entire cyclic cohomology $\HE^\ast(\MA)$.
\end{theorem}

\begin{proof}
  This follows from Proposition~\ref{pro:X_equals_C}.  We only have to examine
  the dual complex $C(\MA)'$ of bounded linear maps $C(\MA) \to \C$.  By the
  universal property of completions, bounded linear maps $C(\MA) \to \C$
  correspond to bounded linear maps $(\Omega\MA, n! \CBS_\an) \to \C$.  Since
  the sets $\bigcup n! \opt{S} \opt{dS} (dS)^{2n}$ generate the bornology $n!
  \cdot \CBS_\an$, the bounded linear functionals on $(\Omega\MA, n!
  \CBS_\an)$ are those linear maps $\Omega\MA \to \C$ that remain bounded on
  all sets of the form $\bigcup n! \opt{S} \opt{dS} (dS)^{2n}$.  Identifying
  $\Omega\MA \cong \sum_{n=0}^\infty \Omega^n\MA$ and $\Omega^n\MA \cong
  \Unse{\MA} \hot \MA^{\hot n}$, we find that $C(\MA)'$ can be identified with
  the space of families $(\phi_n)_{n \in \Z_+}$ of \Mpn{n+1}linear maps
  $\phi_n \colon \Unse{A} \times A^n \to \C$ satisfying the ``entire growth
  condition''
  \begin{equation}
    \label{eq:entire_growth}
    |\phi_n (\opt{\ma_0}, \ma_1, \dots, \ma_n)| \le \const(S) / [n/2]!
    \quad\text{for all $\opt{\ma_0} \in \opt{S}$, $\ma_1, \dots, \ma_n \in S$}
  \end{equation}
  for all $S \in \CBS(\MA)$.  Here $[n/2] \defeq k$ if $n = 2k$ or $n = 2k+1$
  and $\const(S)$ is a constant depending on~$S$ but not on~$n$.  The boundary
  on $C(\MA)'$ is simply composition with $B + b$.
  
  Thus $C(\MA)'$ is the complex used by Khalkhali~\cite{khalkhali94:entire} to
  define entire cyclic cohomology for non-unital (Banach) algebras.  Connes's
  original definition~\cite{connes88:entire} of entire cyclic cohomology is
  restricted to unital algebras and uses a slightly different complex.
  However, for unital algebras the complex used by Connes is homotopy
  equivalent to the complex used by Khalkhali.  The proof is written down
  in~\cite{khalkhali94:entire} for Banach algebras, but carries over to
  locally convex algebras without difficulty.
\end{proof}

\section{Excision in analytic cyclic cohomology}
\label{sec:excision}

Throughout this section, we consider an allowable extension
\begin{equation}  \label{eq:extension}
  \xymatrix{
    {0} \ar[r] &
      {\MA[K]} \ar[r]^{i} &
        {\MA[E]} \ar[r]^{p} &
          {\MA[Q]} \ar[r] \ar@/^/@{.>}[l]^{s} &
            {0}
    }
\end{equation}
of complete bornological algebras with a bounded linear section~$s$.

\begin{theorem}[Excision Theorem]  \label{the:excision_analytic}
  Let~$\MA$ be a further complete bornological algebra.  There are natural six
  term exact sequences in both variables
  \begin{equation}
    \label{eq:excision_analytic}
    \begin{gathered}
      \xymatrix{
        {\HA^0(\MA; \MA[K])} \ar[r]^{i_\ast} &
          {\HA^0(\MA; \MA[E])} \ar[r]^{p_\ast} &
            {\HA^0(\MA; \MA[Q])} \ar[d]^{\partial_\ast} \\
        {\HA^1(\MA; \MA[Q])} \ar[u]^{\partial_\ast} &
          {\HA^1(\MA; \MA[E])} \ar[l]^{p_\ast} &
            {\HA^1(\MA; \MA[K])} \ar[l]^{i_\ast}
        }
    \displaybreak[0] \\
      \xymatrix{
        {\HA^0(\MA[Q]; \MA)} \ar[r]^{p^\ast} &
          {\HA^0(\MA[E]; \MA)} \ar[r]^{i^\ast} &
            {\HA^0(\MA[K]; \MA)} \ar[d]^{\partial^\ast} \\
        {\HA^1(\MA[K]; \MA)} \ar[u]^{\partial^\ast} &
          {\HA^1(\MA[E]; \MA)} \ar[l]^{i^\ast} &
            {\HA^1(\MA[Q]; \MA)} \ar[l]^{p^\ast}
        }
    \end{gathered}
  \end{equation}
  
  The maps $i_\ast,p_\ast,\partial_\ast,i^\ast,p^\ast,\partial^\ast$ are the
  (signed) composition products with $[i] \in \HA^0(\MA[K]; \MA[E])$, $[p]\in
  \HA^0(\MA[E]; \MA[Q])$, and a certain $\partial \in \HA^1(\MA[Q]; \MA[K])$
  naturally associated to the extension.
\end{theorem}

Naturality means that homomorphisms in the \Mpn{\MA}variable and, more
importantly, morphisms of extensions give rise to commutative diagrams of
six term exact sequences.

In the special case $\MA = \C$, we get excision in analytic cyclic homology
and cohomology.  Specializing further to locally convex algebras with the
bounded bornology, we get excision for the entire cyclic cohomology of locally
convex algebras by Theorem~\ref{the:entire_equals_analytic}.

An analogous excision theorem for the bivariant periodic cyclic cohomology of
algebras without additional structure is due to Cuntz and
Quillen~\cite{cuntz97:excision}.  It was later extended by Cuntz to locally
multiplicatively convex topological algebras~\cite{cuntz97:excision_top}.
These proofs are based on the proof of excision in \Hochschild and cyclic
(co)homology by Wodzicki~\cite{wodzicki89:excision} for a special class of
extensions where the kernel is ``H-unital''.  Recently, Puschnigg has been
able to prove excision for entire cyclic cohomology along these
lines~\cite{puschnigg98:excision}.  However, his proof works only for
tensoring algebras.  Thus it works for Banach algebras but not for fine
algebras.

\subsection{Outline of the proof}
\label{sec:excision_outline}

Define
\begin{equation}
  \label{eq:defRel}
  X(\Tanil\MA[E]: \Tanil\MA[Q]) \defeq
  \Ker\bigl( X(\Tanil p) \colon X(\Tanil\MA[E]) \to X(\Tanil\MA[Q]) \bigr).
\end{equation}
Since $X(\Tanil p) \circ X(\Tanil i) = 0$, we can view $X(\Tanil i)$ as a
chain map $\varrho \colon X(\Tanil\MA[K]) \to X(\Tanil\MA[E]: \Tanil\MA[Q])$.
If~$\varrho$ is a homotopy equivalence of complexes, then we easily get the
Excision Theorem by applying the long exact homology sequence to the extension
of complexes $X(\Tanil\MA[E]: \Tanil\MA[Q]) \injto X(\Tanil\MA[E]) \prto
X(\Tanil\MA[Q])$.

The main new idea of the proof is to use the \emph{left} ideal~$\LL$ generated
by $\MA[K] \subset \Tanil\MA[E]$.  The proof of excision in periodic cyclic
cohomology by Cuntz and Quillen~\cite{cuntz97:excision} works instead with the
\emph{two-sided} ideal~$\II$ generated by~$\MA[K]$ and powers of~$\II$.  The
ideal~$\II$ is neither quasi-free nor H-unital, it only has these properties
in a certain approximate sense.  It turns out, however, that~$\LL$ is much
better behaved than~$\II$.  The two main steps in the proof of excision are
the following theorems:

\begin{theorem}  \label{the:excision_step_I}
  The chain map $\psi \colon X(\LL) \to X(\Tanil\MA[E]: \Tanil\MA[Q])$ induced
  by the inclusion $\LL \subset \Tanil\MA[E]$ is a homotopy equivalence
  $X(\LL) \sim X(\Tanil\MA[E]: \Tanil\MA[Q])$.
\end{theorem}

\begin{theorem}  \label{the:excision_step_II}
  The algebra~$\LL$ is analytically quasi-free.
\end{theorem}

Using these two theorems, the proof of excision is finished as follows.  The
natural projection $\tau_{\MA[E]} \colon \Tanil\MA[E] \to \MA[E]$ restricts to
a projection $\tilde{\tau} \colon \LL \to \MA[K]$ with bounded linear section
$\sigma_{\MA[E]} |_{\MA[K]}$.  The kernel $\MA[N] \defeq \Ker \tilde{\tau}$ is
contained in $\Ker \tau_{\MA[E]} = \Janil\MA[E]$ and therefore a-nilpotent by
Lemma~\ref{lem:nilpotent_inherit_triv}.  By
Theorem~\ref{the:excision_step_II}, $\MA[N] \injto \LL \prto \MA[K]$ is a
universal a-nilpotent extension.  Let $j \colon \Tanil\MA[K] \to \LL$ be the
inclusion map, then $(j|_{\Janil\MA[K]}, j, \ID)$ is a morphism of extensions
from $\Janil\MA[K] \injto \Tanil\MA[K] \prto \MA[K]$ to $\MA[N] \injto \LL
\prto \MA[K]$.  Theorem~\ref{the:uniqueness_X} implies that the induced chain
map $X(j) \colon X(\Tanil\MA[K]) \to X(\LL)$ is a homotopy equivalence.
Combining this with the homotopy equivalence $\psi \colon X(\LL) \sim
X(\Tanil\MA[E] : \Tanil\MA[Q])$ of Theorem~\ref{the:excision_step_I}, we get
that $\varrho = \psi \circ X(j) \colon X(\Tanil\MA[K]) \to X(\Tanil\MA[E] :
\Tanil\MA[Q])$ is a homotopy equivalence.  Theorem~\ref{the:excision_analytic}
follows.

Theorem~\ref{the:excision_step_I} is proved by homological algebra.  We write
down an allowable resolution
$$
\VS[P]_\bullet\colon \
0 \longrightarrow
\VS[P]_1 \longrightarrow
\VS[P]_0 \longrightarrow
\Unse{\LL}
$$
with
$$
\VS[P]_0 \defeq \Unse{(\Tanil\MA[E])} \hot\LL + \Unse{\LL}\hot \Unse{\LL}
\subset \Unse{(\Tanil\MA[E])} \hot \Unse{(\Tanil\MA[E])},
\qquad
\VS[P]_1 \defeq \Unse{(\Tanil\MA[E])}\,D\LL.
$$
We prove that $\VS[P]_0$ and~$\VS[P]_1$ are free \Mpn{\LL}bimodules, using
some isomorphisms collected in Section~\ref{sec:excision_isomorphisms}.  These
isomorphisms assert in particular that $\Unse{(\Tanil\MA[E])}$ is a free left
\Mp{\LL}module.  We show that the commutator quotient complex
$\VS[P]_\bullet/[,]$ is isomorphic to the direct sum $X_\beta(\Tanil\MA[E]:
\Tanil\MA[Q]) \oplus \C[0]$ with
$$
X_\beta(\Tanil\MA[E]: \Tanil\MA[Q]) \defeq
\Ker \bigl(
  p_\ast \colon X_\beta(\Tanil\MA[E]) \prto X_\beta(\Tanil\MA[Q])
\bigr).
$$
In addition, the comparison theorem for projective resolutions implies that
$\VS[P]_\bullet/[,]$ is chain homotopic to $X_\beta(\LL) \oplus \C[0]$.
Indeed, it is elementary to prove $\VS[P]_\bullet/[,] \cong X_\beta(\LL)
\oplus \C[0] \oplus \VS[C]'_\bullet$ with a certain contractible
complex~$\VS[C]'_\bullet$.  Thus $X_\beta(\Tanil\MA[E]: \Tanil\MA[Q]) \cong
X_\beta(\LL) \oplus \VS[C]'_\bullet$.  Theorem~\ref{the:excision_step_I}
follows easily.

The existence of the resolution~$\VS[P]_\bullet$ implies already that~$\LL$ is
quasi-free.  To get the much stronger assertion of
Theorem~\ref{the:excision_step_II}, we have to construct a splitting
homomorphism $\upsilon \colon \LL \to \Tanil\LL$.  We write down explicitly a
bounded bilinear map $\triangleright \colon \MA[E] \times \Tanil\LL \to
\Tanil\LL$ and view it as a bounded linear map $\triangleright \colon \MA[E]
\to \Endo(\Tanil\LL)$ by adjoint associativity.  A computation shows
that~$\triangleright$ has a-nilpotent curvature.  This is the hard part of the
construction of~$\upsilon$.  Using the universal property of $\Tanil\MA[E]$ we
can extend~$\triangleright$ to a bounded unital homomorphism $\triangleright
\colon \Unse{(\Tanil\MA[E])} \to \Endo(\Tanil\LL)$.  The splitting
homomorphism~$\upsilon$ is then defined by $\upsilon(x \odot \ma[k]) \defeq
\triangleright(x) (\ma[k])$ for $x \in \Unse{(\Tanil\MA[E])}$ and $\ma[k] \in
\MA[K]$.  That is, the endomorphism~$\triangleright(x)$ is evaluated at the
point $\ma[k] \in \MA[K] \subset \Tanil\LL$.  The definition
of~$\triangleright$ implies easily that~$\upsilon$ is a bounded splitting
homomorphism for $\tau \colon \Tanil\LL \prto \LL$.

\subsection{Linear functoriality of $\Omega_\an$}
\label{sec:Omega_functorial}
The linear section $s \colon \MA[Q] \to \MA[E]$ induces a bounded linear map
$s_L \colon \Omega_\an\MA[Q] \to \Omega_\an\MA[E]$ defined by
\begin{equation}
  \label{eq:Omega_functorial_l}
  s_L(\opt{\ma[q]_0} d\ma[q]_1 \dots d\ma[q]_n) \defeq
  s\opt{\ma[q]_0} ds(\ma[q]_1) \dots ds(\ma[q]_n).
\end{equation}
Evidently, this well-defines a bounded linear map $(\Omega\MA[Q], \CBS_\an)
\to (\Omega\MA[E], \CBS_\an)$.  The map~$s_L$ is its unique extension to the
completions.  If we identify $\Omega_\an\MA[E] \cong X(\Tanil\MA[E])$ and
$\Omega_\an\MA[Q] \cong X(\Tanil\MA[Q])$ as in Section~\ref{sec:X_Tanil}, then
we get a bounded linear map $s_L \colon X(\Tanil\MA[Q]) \to X(\Tanil\MA[E])$
that is a section for $X(\Tanil p)$.  That is, $X(\Tanil p) \circ s_L$ is the
identity on $X(\Tanil\MA[Q])$.  Therefore,
\begin{equation}  \label{eq:ext_relative}
  \xymatrix@+.5cm{
    {X(\Tanil\MA[E]: \Tanil\MA[Q])\;} \ar@{>->}[r] &
      {X(\Tanil\MA[E])} \ar@{->>}[r]^-{X(\Tanil p)} &
        {X(\Tanil\MA[Q])}
    }
\end{equation}
is an allowable extension of complexes.  Let $L_\bullet \defeq X(\Tanil\MA)$.
If~$\varrho$ is a homotopy equivalence, composition with~$[\varrho]$ gives
rise to an isomorphism $H^\ast \bigl( X(\Tanil\MA[E]: \Tanil\MA[Q]); L_\bullet
\bigr) \cong H^\ast \bigl( X(\Tanil\MA[K]); L_\bullet \bigr)$ and similarly
for the second variable.  As in Appendix~\ref{app:complexes}, we
write $H^\ast (K_\bullet; L_\bullet)$ for the homology of the complex
$\Lin(K_\bullet; L_\bullet)$ of bounded linear maps $K_\bullet \to L_\bullet$.
Hence Theorem~\ref{the:excision_analytic} follows from the Long Exact Homology
Sequence~\ref{the:long_exact_homology} applied to the allowable
extension~\eqref{eq:ext_relative}.  The class $\partial \in \HA^1(\MA[Q];
\MA[K])$ is represented by $[\varrho]^{-1} \circ \partial_0$, where
$\partial_0 \in H^1\bigl( X(\Tanil\MA[Q]); X(\Tanil\MA[E]: \Tanil\MA[Q])
\bigr)$ is the connecting map of the allowable
extension~\eqref{eq:ext_relative} and $[\varrho]^{-1} \in H^0\bigl(
X(\Tanil\MA[E]: \Tanil\MA[Q]); X(\Tanil\MA[K]) \bigr)$ is the inverse
of~$[\varrho]$.  The naturality of~$\partial$ is evident.  Thus the Excision
Theorem~\ref{the:excision_analytic} follows if~$\varrho$ is a homotopy
equivalence.

We will need a ``right handed'' version $s_R \colon \Omega_\an\MA[Q] \to
\Omega_\an\MA[E]$ of~\eqref{eq:Omega_functorial_l} defined by
\begin{equation}
  \label{eq:Omega_functorial_r}
  s_R(d\ma[q]_1 \dots d\ma[q]_n \opt{\ma[q]_{n+1}}) \defeq
  ds(\ma[q]_1) \dots ds(\ma[q]_n) s\opt{\ma[q]_{n+1}}.
\end{equation}
To prove that~$s_R$ is well-defined and bounded we use a symmetry of
$\Omega_\an\MA$ for all complete bornological algebras~$\MA$.  The
\emph{opposite algebra}~$\MA^\opp$ of~$\MA$ is defined as follows.  As a
vector space, $\MA^\opp = \MA$, but the multiplication~$\cdot_\opp$ is given
by $\ma_1 \cdot_\opp \ma_2 \defeq \ma_2 \cdot \ma_1$.

\begin{proposition}  \label{pro:left_right}
  Let~$\MA$ be a complete bornological algebra.  The map
  $$
  \nu( \opt{\ma_0} d\ma_1 \dots d\ma_n) \defeq
  d\ma_n \dots d\ma_1 \cdot \opt{\ma_0}
  $$
  defines a natural isomorphism of bornological DG-algebras $\nu \colon
  \Omega_\an(\MA^\opp) \to (\Omega_\an\MA)^\opp$.  The restriction to the
  even part yields a natural isomorphism $\nu \colon \Tanil(\MA^\opp) \to
  (\Tanil\MA)^\opp$.
\end{proposition}

\begin{proof}
  Since the operator~$b'$ is bounded, we have $\nu(\opt{S} (dS)^\infty)
  = (dS)^\infty \cdot \opt{S} = b'\bigl( (dS)^\infty \bigr) \in
  \CBS_\an$ for all $S \in \CBS(\MA)$.  Thus~$\nu$ is bounded.  Since $\nu
  \circ \nu = \ID$, it is a bornological isomorphism.  It is trivial
  to verify that~$\nu$ is multiplicative, when considered as a map
  $\Omega_\an(\MA^\opp) \to (\Omega_\an\MA)^\opp$.
\end{proof}

Since $s_R = \nu \circ s_L \circ \nu^{-1}$, boundedness of~$s_L$ implies
boundedness of~$s_R$.

\begin{digression}
  The isomorphism $\Tanil (\MA^\opp) \cong (\Tanil \MA)^\opp$ fits well into
  the formalism of Section~\ref{sec:Tanil}.  It is equivalent to the following
  statements: If~$\MA[N]$ is a-nilpotent, so is $\MA[N]^\opp$; if $l \colon
  \MA \to \MA[B]$ is a \lanilcur, it remains a \lanilcur when considered as a
  map $\MA^\opp \to \MA[B]^\opp$.  However, the existence of the bounded
  linear map $s_L \colon X(\Tanil\MA[Q]) \to X(\Tanil\MA[E])$ apparently
  cannot be deduced from properties of a-nilpotent algebras or \lanilcurs.
\end{digression}

\subsection{Some isomorphisms}
\label{sec:excision_isomorphisms}

Let $\LL \subset \Tanil\MA[E]$ be the closed left ideal generated by
$\sigma(\MA[K]) \subset \Tanil\MA[E]$, where $\sigma \colon \MA[E] \to
\Tanil\MA[E]$ is the natural linear map.  Thus~$\LL$ is the image of
$\Unse{(\Tanil\MA[E])} \hot \MA[K]$ under the isomorphism $\mu_{\ref{eq:mu1}}
\colon \Unse{(\Tanil\MA[E])} \hot \MA[E] \congto \Tanil\MA[E]$ in
Proposition~\ref{pro:Omega_I_Tanil}.  Consequently, the linear map
\begin{equation}
  \label{eq:mu5}
  \mu_{\ref{eq:mu5}} \colon \Unse{(\Tanil\MA[E])} \hot \MA[K] \congto \LL,
  \qquad
  x \otimes \ma[k] \mapsto x \odot \sigma(\ma[k]),
\end{equation}
is a bornological isomorphism.  An easy calculation shows that $\opt{\ma[e]_0}
d\ma[e]_1 \dots d\ma[e]_{2n} \in \LL$ iff $\ma[e]_{2n} \in \MA[K]$.

Let $d\ma[e]_1 \dots d\ma[e]_{2n} \opt{\ma[e]_{2n+1}}$ be a monomial in right
handed standard form.  Splitting each entry as $\ma[e]_j = \ma[k]_j +
s(\ma[q]_j)$ with $\ma[k]_j \in \MA[K]$ and $\ma[q]_j \in \MA[Q]$, we can
write $d\ma[e]_1 \dots d\ma[e]_{2n} \opt{\ma[e]_{2n+1}}$ as a sum of
$2^{2n+1}$ \emph{right polarized monomials}, that is, monomials $d\ma[e]_1'
\dots d\ma[e]_{2n}' \opt{\ma[e]_{2n+1}'}$ with the additional property that
either $\ma[e]_j' \in \MA[K]$ or $\ma[e]_j' \in s(\MA[Q])$.  If $S \in
\CBS(\MA[E])$, then there are small sets $S_{\MA[K]} \in \CBS(\MA[K])$ and
$S_{\MA[Q]} \in \CBS(\MA[Q])$ such that $S \subset \frac{1}{2}
\bigl(S_{\MA[K]} + s(S_{\MA[Q]}) \bigr)$.  If $\ma[e]_j \in S$ for all~$j$,
then $d\ma[e]_1 \dots d\ma[e]_{2n} \opt{\ma[e]_{2n+1}}$ is a convex
combination of terms in $(dS_{\MA[E]})^\even \opt{S_{\MA[E]}}$ with
$S_{\MA[E]} = S_{\MA[K]} \cup s(S_{\MA[E]}) \subset \MA[K] \cup s(\MA[Q])$.
We call small subsets of the form $(dS_{\MA[E]})^\even \opt{S_{\MA[E]}}$
\emph{right polarized}.  Hence the bornology $\CBS_\an^\even$ defining
$\Tanil\MA[E]$ is generated by right polarized small sets.

\begin{lemma}  \label{lem:mu68}
  The following two linear maps are bornological isomorphisms:
  \begin{multline}
    \label{eq:mu6}
    \mu_{\ref{eq:mu6}} \colon \Unse{(\Tanil\MA[Q])} \oplus
    \Unse{(\Tanil\MA[E])} \hot \MA[K] \hot \Unse{(\Tanil\MA[Q])} \congto
    \Unse{(\Tanil\MA[E])},
    \\
    \opt{q_1} \oplus \opt{x} \otimes \ma[k] \otimes \opt{q_2} \mapsto
    s_L\opt{q_1} + \opt{x} \odot \sigma(\ma[k]) \odot s_L\opt{q_2},
  \end{multline}
  \begin{multline}
    \label{eq:mu8}
    \mu_{\ref{eq:mu8}} \colon \Unse{(\Tanil\MA[Q])} \oplus
    \Unse{(\Tanil\MA[Q])} \hot \MA[K] \hot \Unse{(\Tanil\MA[E])} \congto
    \Unse{(\Tanil\MA[E])},
    \\
    \opt{q_1} \oplus \opt{q_2} \otimes \ma[k] \otimes \opt{x} \mapsto
    s_R\opt{q_1} + s_R\opt{q_2} \odot \sigma(\ma[k]) \odot \opt{x}.
  \end{multline}
\end{lemma}

\begin{proof}
  It suffices to prove that $\mu_{\ref{eq:mu8}}$ is a bornological
  isomorphism.  We obtain the other map $\mu_{\ref{eq:mu6}}$ by exchanging
  left and right as in Proposition~\ref{pro:left_right}: $\mu_{\ref{eq:mu6}}
  = \nu \circ \mu_{\ref{eq:mu8}} \circ \nu^{-1}$.  We write down explicitly
  the inverse of $\mu_{\ref{eq:mu8}}$ and verify that it is bounded.  We
  will need some detailed estimates about this inverse for the proof
  that~$\LL$ is a-quasi-free.
  
  Let $\omega = d\ma[e]_1 \dots d\ma[e]_{2n} \opt{\ma[e]_{2n+1}} $ be a right
  polarized monomial.  We write $\omega = \mu_{\ref{eq:mu8}}(\cdots)$ and thus
  derive a formula for the inverse~$\mu_{\ref{eq:mu8}}^{-1}$.  If $\ma[e]_j =
  s\ma[q]_j \in s(\MA[Q])$ for all~$j$, then $\omega = \mu_{\ref{eq:mu8}} (
  d\ma[q]_1 \dots d\ma[q]_{2n} \opt{\ma[q]_{2n+1}} )$.  In particular, this
  applies to $1 \in \Unse{(\Tanil\MA[E])}$.  Otherwise, some~$\ma[e]_j$ is
  in~$\MA[K]$.  Pick the first~$j$ with $\ma[e]_j = \ma[k]_j \in \MA[K]$.
  Hence $\ma[e]_k = s(\ma[q]_k) \in s(\MA[Q])$ for $k = 1, \dots, j-1$.
  If~$j$ is even, then
  \begin{multline}  \label{eq:mu8_inverse_even}
    d\ma[e]_1 \dots d\ma[e]_{2n} \opt{\ma[e]_{2n+1}} =
    ds\ma[q]_1 \dots ds\ma[q]_{j - 2} \odot
    (s\ma[q]_{j-1} \cdot \ma[k]_j - s\ma[q]_{j-1} \odot \ma[k]_j) \odot
    d\ma[e]_{j+1} \dots d\ma[e]_{2n} \opt{\ma[e]_{2n+1}}
    \\ =
    \mu_{\ref{eq:mu8}} \bigl(
      d\ma[q]_1 \dots d\ma[q]_{j-2} \otimes (s\ma[q]_{j-1}\cdot \ma[k]_j)
      \otimes d\ma[e]_{j+1} \dots d\ma[e]_{2n} \opt{\ma[e]_{2n+1}}
    \\
    - d\ma[q]_1 \dots d\ma[q]_{j-2} \: \ma[q]_{j-1} \otimes \ma[k]_j \otimes
      d\ma[e]_{j+1} \dots d\ma[e]_{2n} \opt{\ma[e]_{2n+1}}
    \bigr).
  \end{multline}
  If~$j$ is odd and $j \neq 2n+1$, then
  \begin{multline}  \label{eq:mu8_inverse_odd}
    d\ma[e]_1 \dots d\ma[e]_{2n} \opt{\ma[e]_{2n+1}} =
    ds\ma[q]_1 \dots ds\ma[q]_{j - 1} \odot
    (\ma[k]_j \cdot \ma[e]_{j+1} - \ma[k]_j \odot \ma[e]_{j+1}) \odot
    d\ma[e]_{j+2} \dots d\ma[e]_{2n} \opt{\ma[e]_{2n+1}}
    \\ =
    \mu_{\ref{eq:mu8}} \bigl(
      d\ma[q]_1 \dots d\ma[q]_{j-1} \otimes (\ma[k]_j \cdot \ma[e]_{j+1})
      \otimes d\ma[e]_{j+2} \dots d\ma[e]_{2n} \opt{\ma[e]_{2n+1}}
    \\
    - d\ma[q]_1 \dots d\ma[q]_{j-1} \otimes \ma[k]_j \otimes
      \ma[e]_{j+1} \odot d\ma[e]_{j+2} \dots d\ma[e]_{2n} \opt{\ma[e]_{2n+1}}
    \bigr).
  \end{multline}
  Finally, if $j = 2n+1$, then
  \begin{equation}  \label{eq:mu8_inverse_final}
    d\ma[e]_1 \dots d\ma[e]_{2n} \: \ma[k]_{2n+1} =
    \mu_{\ref{eq:mu8}}(
      d\ma[q]_1 \dots d\ma[q]_{2n} \otimes \ma[k]_{2n+1} \otimes 1
    ).
  \end{equation}
  
  Of course, we define $\mu_{\ref{eq:mu8}}^{-1}(\omega)$ to be the pre-image
  of~$\omega$ described above.  Writing elements of $\Tens\MA[E]$ as linear
  combinations of right polarized monomials, we obtain a linear map
  $\mu_{\ref{eq:mu8}}^{-1} \colon \Tens\MA[E] \to \Unse{(\Tens\MA[Q])} \oplus
  \Unse{(\Tens\MA[Q])} \otimes \MA[K] \otimes \Unse{(\Tens\MA[E])}$.  By
  construction, this map satisfies $\mu_{\ref{eq:mu8}} \circ
  \mu_{\ref{eq:mu8}}^{-1} = \ID$.  It is easy to verify that
  $\mu_{\ref{eq:mu8}}^{-1} \circ \mu_{\ref{eq:mu8}} = \ID$.  Equations
  \eqref{eq:mu8_inverse_even}--\eqref{eq:mu8_inverse_final} imply that
  $\mu_{\ref{eq:mu8}}^{-1}$ maps right polarized small sets to small sets.
  Hence $\mu_{\ref{eq:mu8}}^{-1}$ extends to a bounded linear map
  $\mu_{\ref{eq:mu8}}^{-1} \colon \Tanil\MA[E] \to \Unse{(\Tanil\MA[Q])}
  \oplus \Unse{(\Tanil\MA[Q])} \otimes \MA[K] \otimes \Unse{(\Tanil\MA[E])}$.
  The equalities $\mu_{\ref{eq:mu8}} \circ \mu_{\ref{eq:mu8}}^{-1} = \ID$ and
  $\mu_{\ref{eq:mu8}}^{-1} \circ \mu_{\ref{eq:mu8}} = \ID$ carry over to the
  completions.  Thus $\mu_{\ref{eq:mu8}}^{-1}$ is the inverse
  of~$\mu_{\ref{eq:mu8}}$ as asserted.

  For later application, we describe the image of a right polarized
  subset of $\Unse{(\Tanil\MA[E])}$ under $\mu_{\ref{eq:mu8}}^{-1}$ more
  explicitly.  Let $S_{\MA[K]} \in \CBS(\MA[K])$, $S_{\MA[Q]} \in
  \CBS(\MA[Q])$, $S_{\MA[E]} \defeq S_{\MA[K]} \cup s(S_{\MA[Q]})$.  We write
  $A \bipol{\subset} B$ to denote that $A \subset \bipol{B}$ and get
  \begin{multline}
    \label{eq:mu8_inverse_small}
    \mu_{\ref{eq:mu8}}^{-1} \bigl( (dS_{\MA[E]})^\even \opt{S_{\MA[E]}}
    \bigr) \bipol{\subset}
    (dS_{\MA[Q]})^\even \opt{S_{\MA[Q]}}
    \\ \cup
    (dS_{\MA[Q]})^\even \opt{S_{\MA[Q]}} \otimes
    \opt{s(S_{\MA[Q]})} \cdot S_{\MA[K]} \cdot \opt{S_{\MA[E]}} \otimes
    2 \opt{S_{\MA[E]}} \odot (dS_{\MA[E]})^\even \opt{S_{\MA[E]}}.
  \end{multline}
  The constant factor~$2$ converts sums of two terms into convex combinations.
  It is important that in the very first tensor factor, we get nothing bigger
  than $(dS_{\MA[Q]})^\even \opt{S_{\MA[Q]}}$.
\end{proof}

Combining \eqref{eq:mu6} and~\eqref{eq:mu5}, we obtain the bornological
isomorphism
\begin{equation}  \label{eq:mu4}
  \mu_{\ref{eq:mu4}} \colon
  \Unse{\LL} \hot \Unse{(\Tanil\MA[Q])} \congto
  \Unse{(\Tanil\MA[E])}, \qquad
  \opt{l} \otimes \opt{q} \mapsto \opt{l} \odot s_L\opt{q}.
\end{equation}
Evidently, $\mu_{\ref{eq:mu4}}( l_1 \odot \opt{l_2} \otimes \opt{q_3}) = l_1
\odot \mu_{\ref{eq:mu4}}(\opt{l_2} \otimes \opt{q_3})$, that is,
$\mu_{\ref{eq:mu4}}$ is a left \Mpn{\LL}module map with respect to the obvious
left module structure.  Thus $\Unse{(\Tanil\MA[E])}$ is a free left
\Mpn{\LL}module.

Restricting $\mu_{\ref{eq:mu8}}$ to $\Unse{(\Tanil\MA[Q])} \hot \MA[K] \hot
\Unse{\LL}$, we get a bornological isomorphism
\begin{equation}
  \label{eq:mu9}
  \mu_{\ref{eq:mu9}} \colon
  \Unse{(\Tanil\MA[Q])} \hot \MA[K] \hot \Unse{\LL} \congto
  \LL, \qquad
  \opt{q} \otimes \ma[k] \otimes \opt{l} \mapsto
  s_R\opt{q} \odot \sigma(\ma[k]) \odot \opt{l}.
\end{equation}
To prove that~\eqref{eq:mu9} is a bornological isomorphism, tensor both sides
of~\eqref{eq:mu8} on the right with~$\MA[K]$ and simplify using~\eqref{eq:mu5}
and $\Unse{\LL} \cong \C \oplus \LL$.  The map $\mu_{\ref{eq:mu9}}$ is defined
so that $\mu_{\ref{eq:mu9}} \bigl( \opt{q} \otimes \ma[k] \otimes (\opt{l_1}
\odot l_2) \bigr) = \mu_{\ref{eq:mu9}} ( \opt{q} \otimes \ma[k] \otimes
\opt{l_1}) \odot l_2$.  That is, $\mu_{\ref{eq:mu9}}$ is an isomorphism of
right \Mpn{\LL}modules.  Consequently, $\LL$ is a free right \Mpn{\LL}module.

Let $\II \defeq \Ker (\Tanil p \colon \Tanil\MA[E] \to \Tanil\MA[Q])$.  Since
$\II \cong \Unse{(\Tanil\MA[E])} / s_R\Unse{(\Tanil\MA[Q])}$, restriction
of $\mu_{\ref{eq:mu8}}$ gives rise to a bornological isomorphism
\begin{equation}
  \label{eq:mu12}
  \mu_{\ref{eq:mu12}} \colon
  \Unse{(\Tanil\MA[Q])} \hot \MA[K] \hot \Unse{(\Tanil\MA[E])} \congto
  \II, \qquad
  \opt{q} \otimes \ma[k] \otimes \opt{x} \mapsto
  s_R\opt{q} \odot \ma[k] \odot \opt{x}.
\end{equation}
This map is evidently a right \Mpn{\Tanil\MA[E]}module homomorphism, that is,
$\mu_{\ref{eq:mu12}} \bigl( \opt{q} \otimes \ma[k] \otimes (\opt{x_1} \odot
x_2) \bigr) = \mu_{\ref{eq:mu12}} ( \opt{q} \otimes \ma[k] \otimes \opt{x_1})
\odot x_2$.  Symmetrically, \eqref{eq:mu6} yields $\II \cong
\Unse{(\Tanil\MA[E])} \hot \MA[K] \hot \Unse{(\Tanil\MA[Q])}$.  Using
also~\eqref{eq:mu5}, we get a bornological isomorphism
\begin{equation}
  \label{eq:mu11}
  \mu_{\ref{eq:mu11}} \colon
  \LL \hot \Unse{(\Tanil\MA[Q])} \congto \II, \qquad
  l \otimes \opt{q} \mapsto l \odot s_L\opt{q}.
\end{equation}
Combining further with~\eqref{eq:mu4}, we get a bornological isomorphism
\begin{equation}
  \label{eq:mu7}
  \mu_{\ref{eq:mu7}} \colon
  \LL \hot_{\Unse{\LL}} \Unse{(\Tanil\MA[E])} \congto \II, \qquad
  l \otimes \opt{x} \mapsto l \odot \opt{x}.
\end{equation}

\subsection{A free resolution of~$\Unse{\LL}$}
\label{sec:excision_step_I}

This section contains the proof of Theorem~\ref{the:excision_step_I}.  Let
$$
X_\beta(\Tanil\MA[E]: \Tanil\MA[Q]) \defeq
\Ker \bigl( p_\ast\colon X_\beta(\Tanil\MA[E]) \prto
X_\beta(\Tanil\MA[Q]) \bigr).
$$
The isomorphism $X_\beta(\Tanil\MA[E]) \cong \Omega_\an\MA[E]$ maps
$X_\beta(\Tanil\MA[E]: \Tanil\MA[Q])$ to the complex
$$
\Unse{(\Tanil\MA[E])} \,D\MA[K] \oplus \II\,Ds(\MA[Q])
\overset{\partial_1}{\longrightarrow}
\II.
$$

Let~$\MA$ be a complete bornological algebra.  (We will use the cases $\MA =
\Tanil\MA[E]$ and $\MA = \LL$).  The following is an allowable extension of
\Mpn{\MA}bimodules:
\begin{equation}
  \label{eq:res1}
  \begin{gathered}
    \xymatrix{
      {\VS[B]^{\MA}_\bullet \colon \ }
        {\Omega^1(\MA)\;} \ar@<1ex>@{>->}[r]^-{\alpha_1} &
          {\Unse{\MA} \hot \Unse{\MA}} \ar@<1ex>@{->>}[r]^-{\alpha_0}
            \ar@<1ex>@{.>>}[l]^-{h_1} &
            {\;\Unse{\MA}} \ar@<1ex>@{>.>}[l]^-{h_0}
      } \\
    \begin{split}
      \alpha_1(\opt{x} (Dy) \opt{z}) &=
      (\opt{x} \odot y) \otimes \opt{z} - \opt{x} \otimes (y \odot \opt{z});
      \\
      \alpha_0(\opt{x} \otimes \opt{y}) &= \opt{x} \cdot \opt{y};
    \end{split} \qquad
    \begin{split}
      h_0(\opt{x}) &= 1 \otimes \opt{x}; \\
      h_1(\opt{x} \otimes \opt{y}) &= (D\opt{x}) \opt{y}.
    \end{split}
  \end{gathered}
\end{equation}
Evidently, $[\alpha, h] = \ID$, $h \circ h=0$, $\alpha \circ \alpha = 0$.
Thus~$h_\bullet$ is a natural contracting homotopy of the complex
$(\VS[B]^{\MA}_\bullet, \alpha_\bullet)$.  The complex $\VS[B]^{\MA}_\bullet$
is a projective resolution iff~$\MA$ is quasi-free.  The commutator quotient
complex $\VS[B]_\bullet^{\MA}/[,] \defeq \VS[B]_\bullet^{\MA}/
[\VS[B]_\bullet^{\MA},\MA]$ is isomorphic to the complex $X_\beta(\MA) \oplus
\C$.

Define a subcomplex $\VS[P]_\bullet \subset \VS[B]^{\Tanil\MA[E]}_\bullet$ as
follows:
\begin{align}
  \label{eq:PzeroDef}
  \VS[P]_0 &\defeq
  \Unse{(\Tanil\MA[E])} \hot \LL + \Unse{\LL} \hot \Unse{\LL}
  \subset \Unse{(\Tanil\MA[E])} \hot \Unse{(\Tanil\MA[E])};
  \\
  \label{eq:PoneDef}
  \VS[P]_1 &\defeq
  \Unse{(\Tanil\MA[E])} \,D\LL \subset \Omega^1(\Tanil\MA[E]).
\end{align}

Since~$\LL$ is a left ideal in $\Tanil\MA[E]$, we have $\alpha_1 (\VS[P]_1)
\subset \VS[P]_0$ and $\alpha_0 (\VS[P]_0) \subset \Unse{\LL}$.  Thus
$\VS[P]_1 \to \VS[P]_0 \to \Unse{\LL}$ is a subcomplex of
$\VS[B]^{\Tanil\MA[E]}_\bullet \to \Unse{(\Tanil\MA[E])}$.  The subspaces
$\VS[P]_0$ and~$\VS[P]_1$ are \Mpn{\LL}sub-bimodules
of~$\VS[B]^{\Tanil\MA[E]}_0$ and~$\VS[B]^{\Tanil\MA[E]}_1$, respectively.
Moreover, $h_0(\Unse{\LL}) \subset \VS[P]_0$, $h_1(\VS[P]_0) \subset
\VS[P]_1$, so that~$h_\bullet$ is a contracting homotopy of $\VS[P]_1 \to
\VS[P]_0 \to \Unse{\LL}$.  Consequently, $\VS[P]_1 \to \VS[P]_0 \to
\Unse{\LL}$ is an allowable resolution of~$\Unse{\LL}$ by \Mpn{\LL}bimodules.

\begin{theorem}  \label{the:P_free}
  The \Mpn{\LL}bimodules $\VS[P]_0$ and~$\VS[P]_1$ are free.
  Thus~$\VS[P]_\bullet$ is an allowable free \Mpn{\LL}bimodule resolution
  of~$\Unse{\LL}$.  The algebra~$\LL$ is quasi-free.
  
  The inclusion $\VS[P]_\bullet \subset \VS[B]^{\Tanil\MA[E]}_\bullet$ induces
  a natural chain map $\phi_\bullet \colon \VS[P]_\bullet/[,] \to
  \VS[B]^{\Tanil\MA[E]}_\bullet/[,] \cong X_\beta(\Tanil\MA[E]) \oplus \C[0]$
  that is an isomorphism onto $X_\beta(\Tanil\MA[E] : \Tanil\MA[Q]) \oplus
  \C[0]$.
  
  The natural map $\psi \colon X(\LL) \to X(\Tanil\MA[E]: \Tanil\MA[Q])$ is
  split injective and $X (\Tanil\MA[E]: \Tanil\MA[Q]) = \psi\bigl( X(\LL)
  \bigr) \oplus \VS[C]'_\bullet$ with a contractible complex $\VS[C]'_\bullet$
  given by $\VS[C]'_0 \defeq [\LL, s_L(\Tanil\MA[Q])]$ and $\VS[C]'_1 \defeq
  \LL \, Ds_L(\Tanil\MA[Q])$.  Thus~$\psi$ is a homotopy equivalence $X(\LL)
  \sim X(\Tanil\MA[E]: \Tanil\MA[Q])$.
\end{theorem}

\begin{proof}
  Equation~\eqref{eq:mu4} and $\Unse{(\Tanil\MA[E])} \hot \LL \cap
  \Unse{\LL} \hot \Unse{\LL} = \Unse{\LL} \hot \LL$ imply that
  \begin{multline}
    \label{eq:mu21}
    \mu_{\ref{eq:mu21}} \colon
    (\Unse{\LL} \hot \Unse{\LL}) \oplus
    (\Unse{\LL} \hot \Tanil\MA[Q] \hot \LL) \congto \VS[P]_0, 
    \\
    \opt{l_1} \otimes \opt{l_2} \oplus \opt{l_3} \otimes q_4 \otimes l_5
    \mapsto \opt{l_1} \otimes \opt{l_2} +
    \bigl( \opt{l_3} \odot s_L(q_4) \bigr) \otimes l_5,
  \end{multline}
  is a bornological isomorphism.  Since~$\LL$ is a free right
  \Mpn{\LL}bimodule by~\eqref{eq:mu9}, it follows that~$\VS[P]_0$ is a free
  \Mpn{\LL}bimodule.  We claim that
  \begin{equation}  \label{eq:PoneCompute}
    \VS[P]_1 =
    \Omega^1(\Tanil\MA[E]) \odot \MA[K] + \Unse{(\Tanil\MA[E])} \,D\MA[K].
  \end{equation}
  It is evident that $\Unse{(\Tanil\MA[E])} \,D\MA[K] \subset \VS[P]_1$.
  Moreover,
  $$
  \opt{x_0} (Dx_1) \odot \ma[k] =
  \opt{x_0} \,D(x_1\odot \ma[k]) - \opt{x_0} \odot x_1 \,D\ma[k] \in \VS[P]_1
  \quad \forall
  \opt{x_0} \in \Unse{(\Tanil\MA[E])},\ x_1 \in \Tanil\MA[E],\
  \ma[k] \in \MA[K].
  $$
  Thus $\Omega^1(\Tanil\MA[E]) \odot \MA[K] \subset \VS[P]_1$.  Conversely,
  $$
  x_0\, D(x_1 \odot \ma[k]) =
  x_0 (Dx_1) \odot \ma[k] + x_0 \odot x_1 \,D\ma[k] \in
  \Omega^1(\Tanil\MA[E]) \odot \MA[K] + \Unse{(\Tanil\MA[E])} \,D\MA[K].
  $$
  This implies $\VS[P]_1 \subset \Omega^1(\Tanil\MA[E]) \odot \MA[K] +
  \Unse{(\Tanil\MA[E])} \,D\MA[K]$ by~\eqref{eq:mu5} and finishes the proof
  of~\eqref{eq:PoneCompute}.
  
  The inverse~$\mu_{\ref{eq:mu3}}^{-1}$ of the isomorphism
  $\Omega^1(\Tanil\MA[E]) \cong \Unse{(\Tanil\MA[E])} \hot \MA[E] \hot
  \Unse{(\Tanil\MA[E])}$ in~\eqref{eq:mu3} maps $\Omega^1(\Tanil\MA[E]) \odot
  \MA[K]$ onto $\Unse{(\Tanil\MA[E])} \hot \MA[E] \hot (\Unse{(\Tanil\MA[E])}
  \odot \MA[K]) \cong \Unse{(\Tanil\MA[E])} \hot \MA[E] \hot \LL$ and maps
  $\Unse{(\Tanil\MA[E])} \,D\MA[K]$ onto $\Unse{(\Tanil\MA[E])} \hot \MA[K]
  \hot 1$.  Thus
  \begin{multline}
    \label{eq:mu22}
    \mu_{\ref{eq:mu22}} \colon
    \bigl(\Unse{(\Tanil\MA[E])} \hot \MA[K] \hot \Unse{\LL} \bigr) \oplus
    \bigl(\Unse{(\Tanil\MA[E])} \hot \MA[Q] \hot \LL \bigr) \congto \VS[P]_1,
    \\
    \opt{x_1} \otimes \ma[k] \otimes \opt{l_2} +
    \opt{l_3} \otimes \ma[q] \otimes l_4 \mapsto
    \opt{x_1} (D\ma[k]) \opt{l_2} + \opt{l_3} (Ds\ma[q]) l_4
  \end{multline}
  is a bornological isomorphism.  Since $\Unse{(\Tanil\MA[E])}$ is a free left
  \Mpn{\LL}module by~\eqref{eq:mu4} and~$\LL$ is a free right \Mpn{\LL}module
  by~\eqref{eq:mu9}, $\VS[P]_1$ is a free \Mpn{\LL}bimodule.  Thus
  $\VS[P]_\bullet$ is an allowable free \Mpn{\LL}bimodule resolution
  of~$\Unse{\LL}$.  $\LL$ is quasi-free because it satisfies
  condition~\eqref{qf_ix} of Definition~\ref{deflem:quasi_free}.
  
  The inclusion $\VS[P]_\bullet \injto \VS[B]^{\Tanil\MA[E]}_\bullet$ is an
  \Mpn{\LL}bimodule homomorphism and therefore descends to a chain map
  $\phi_\bullet \colon \VS[P]_\bullet/[,] \to \VS[B]^{\Tanil\MA[E]} /
  [\VS[B]^{\Tanil\MA[E]}, \LL] \to \VS[B]^{\Tanil\MA[E]}_\bullet/[,]$.  Using
  \eqref{eq:mu21} and~\eqref{eq:mu11}, we compute that~$\phi_0$ is an
  isomorphism $\VS[P]_0/[,] \congto \C \oplus \II$:
  $$
  \VS[P]_0/[,] \cong
  \Unse{\LL} \oplus (\LL \hot \Tanil\MA[Q]) \cong
  \C \oplus \bigl( \LL \hot \Unse{(\Tanil\MA[Q])} \bigr) \cong
  \C \oplus \II \subset \Unse{(\Tanil\MA[E])}.
  $$
  Equations \eqref{eq:mu22} and~\eqref{eq:mu7} imply that~$\phi_1$ is an
  isomorphism $\VS[P]_1 /[,] \congto \Unse{(\Tanil\MA[E])} \,D\MA[K] \oplus
  \II \,Ds(\MA[Q])$:
  $$
  \VS[P]_1/[,] \cong
  (\Unse{(\Tanil\MA[E])} \hot \MA[K]) \oplus
  (\LL \hot_{\Unse{\LL}} \Unse{(\Tanil\MA[E])} \hot \MA[Q]) \cong
  \Unse{(\Tanil\MA[E])} \hot \MA[K] \oplus
  \II \hot \MA[Q] \cong
  \Unse{(\Tanil\MA[E])} \,D\MA[K] +
  \II \,Ds(\MA[Q]) \subset \Omega_\an^\odd\MA[E].
  $$
  Thus~$\phi_\bullet$ is an isomorphism of chain complexes
  $\VS[P]_\bullet/[,] \congto X_\beta(\Tanil\MA[E]: \Tanil\MA[Q]) \oplus
  \C[0]$ as asserted.
  
  The standard resolution $\VS[B]^{\LL}_\bullet$ is a subcomplex
  of~$\VS[P]_\bullet$.  Let $f_\bullet \colon \VS[B]^{\LL}_\bullet \to
  \VS[P]_\bullet$ be the inclusion.  The comparison theorem implies
  that~$f_\bullet$ is a homotopy equivalence.  Let us make this more explicit.
  Let $\VS[C]_0 \defeq \Unse{\LL} \hot \Tanil\MA[Q] \hot \LL$ and define $g
  \colon \VS[C]_0 \to \VS[P]_0$ by
  $$
  g( \opt{l_1} \otimes q_2 \otimes l_3) \defeq
  \opt{l_1} \odot s_L(q_2) \otimes l_3 - \opt{l_1} \otimes s_L(q_2) \odot l_3.
  $$
  It is evident that $f_0 \oplus g \colon \Unse{\LL} \hot \Unse{\LL} \oplus
  \VS[C]_0 \to \VS[P]_0$ is a bornological isomorphism and that $\alpha_0
  \circ g = 0$.  Thus $\VS[P]_1 = \Ker \alpha_0 \cong \Omega^1\LL \oplus
  \VS[C]_0$.  Let $\VS[C]_1 \defeq \VS[C]_0$ and let the boundary $\VS[C]_1
  \to \VS[C]_0$ be the identity map.  Then $\VS[C]_\bullet$ is a contractible
  complex of \Mpn{\LL}bimodules and $\VS[P]_\bullet \cong \VS[B]^{\LL}_\bullet
  \oplus \VS[C]_\bullet$.  Taking commutators quotients, we obtain
  $$
  X_\beta(\Tanil\MA[E]: \Tanil\MA[Q]) \oplus \C[0] \cong
  \VS[P]_\bullet/[,] \cong
  \VS[B]^{\LL}_\bullet /[,] \oplus \VS[C]_\bullet /[,] \cong
  X_\beta(\LL) \oplus \C[0] \oplus \VS[C]_\bullet /[,].
  $$
  The copies of~$\C[0]$ match and can be left out.  The map $X_\beta(\LL)
  \to X_\beta(\Tanil\MA[E]: \Tanil\MA[Q])$ in this isomorphism is equal
  to~$\psi$.  Thus~$\psi$ is split injective.  Let $\VS[C]'_\bullet$ be the
  image of $\VS[C]_\bullet /[,]$ in $X_\beta(\Tanil\MA[E]: \Tanil\MA[Q])$.  We
  find that $\VS[C]'_0$ is the range of the map $\LL \hot \Tanil\MA[Q] \to
  \Tanil\MA[E]$, $l \otimes q \mapsto [l, s_L(q)] \in \II$ and
  that~$\VS[C]'_1$ is the range of the map $\LL \hot \Tanil\MA[Q] \to
  \Tanil\MA[E]$, $l \otimes q \mapsto \bigl( Ds_L(q) \bigr) l = l \,Ds_L(q)
  \bmod [,]$.  Thus we get the complex~$\VS[C]'_\bullet$ described in the
  theorem.  Evidently, $\partial_1 \colon \VS[C]'_1 \to \VS[C]'_0$ is a
  bornological isomorphism.  Thus~$\partial_0$ vanishes on $\VS[C]'_0$, so
  that $\VS[C]'_\bullet$ is a subcomplex of $X(\Tanil\MA[E]: \Tanil\MA[Q])$ as
  well.  Since~$\psi$ is compatible with~$\partial_0$, we get a direct sum
  decomposition $X(\Tanil \MA[E]: \Tanil\MA[Q]) \cong X(\LL) \oplus
  \VS[C]'_\bullet$.  Since $\VS[C]'_\bullet$ is contractible, the chain
  map~$\psi$ is a homotopy equivalence.
\end{proof}

\begin{digression}
  Let $\Omega(\Tanil\MA[E]: \Tanil\MA[Q])$ be the \textsl{relative \Hochschild
    complex}
  $$
  \Omega(\Tanil\MA[E]: \Tanil\MA[Q]) \defeq
  \Ker \bigl( p_\ast \colon (\Omega(\Tanil\MA[E]), b) \prto
  (\Omega(\Tanil\MA[Q]), b) \bigr).
  $$
  Since all occurring objects are quasi-free, $\Omega(\LL) \sim
  X_\beta(\LL)$ and $\Omega(\Tanil\MA[E]: \Tanil\MA[Q]) \sim
  X_\beta(\Tanil\MA[E]: \Tanil\MA[Q])$.  Theorem~\ref{the:P_free} implies that
  the chain map $j_\ast \colon (\Omega(\LL), b) \injto (\Omega(\Tanil\MA[E]:
  \Tanil\MA[Q]), b)$ induced by the inclusion $\LL \to \Tanil\MA[E]$ is a
  homotopy equivalence.  Thus we get an exact sequence
  $$
  \HH^0(\Tanil\MA[Q]) \injto \HH^0(\Tanil\MA[E]) \to \HH^0(\LL) \to
  \HH^1(\Tanil\MA[Q]) \to \HH^1(\Tanil\MA[E]) \prto \HH^1(\LL).
  $$
  This is some kind of excision in \Hochschild cohomology for the not quite
  exact sequence $\LL \to \Tanil\MA[E] \to \Tanil\MA[Q]$.  For example,
  exactness at $\HH^0(\Tanil\MA[E])$ means that a trace $T \colon \Tanil\MA[E]
  \to \C$ that vanishes on~$\LL$ necessarily vanishes on~$\II$ and thus comes
  from a trace on $\Tanil\MA[Q] = \Tanil\MA[E] / \II$.  Indeed, we have $\II
  \cong \LL \oplus [\LL, s_L(\Tanil\MA[Q])]$ and a trace on $\Tanil\MA[E]$ has
  to vanish on the second summand.
\end{digression}

\subsection{Analytic quasi-freeness of~$\LL$}
\label{sec:excision_step_II}

We consider~$\MA[K]$ as a linear subspace of the algebras $\MA[E]$,
$\Tanil\MA[E]$, $\LL$, and~$\Tanil\LL$, always using the linear maps
$\sigma_{\MA[E]}$ and~$\sigma_{\Tanil\MA[E]}$.  Therefore, we have to
distinguish carefully between the different products in these algebras.  We
write~$\cdot$ for the product on the level of~$\MA[E]$, that is, in $\MA[K]$,
$\MA[E]$, and~$\MA[Q]$; $\odot$ for the Fedosov product on the level of
$\Tanil\MA[E]$, that is, in $\LL$, $\Tanil\MA[E]$, and $\Tanil\MA[Q]$;
and~$\circledcirc$ for the Fedosov product in $\Tanil\LL$.  Let $\VS[G] \defeq
\Unse{(\Tanil\MA[Q])} \hot \MA[K]$.  In this section we always consider
$\VS[G] \subset \LL$ using the map $\opt{q} \otimes \ma[k] \mapsto s_R\opt{q}
\odot \ma[k]$.  Notice that we take~$s_R$ as in~\eqref{eq:mu9}.

In the following computations and definitions, $DX$ is a shorthand for $Dx_2
\dots Dx_{2n}$ with $n \ge 1$ and~$\opt{DX}$ is a shorthand for $Dx_1 \dots
Dx_{2n}$ with $n \ge 0$.  In case $n = 0$, $\opt{DX}$ is the formal unit of
$\Unse{(\Tanil\LL)}$, that is, $x_0 \opt{DX} = x_0$.  Standing conventions are
$x, x_j \in \LL$, $\ma[k], \ma[k]_j \in \MA[K]$, $\ma[q], \ma[q]_j \in \MA[Q]$,
$\vs[g] \in \VS[G]$, $\ma[e], \ma[e]_j \in \MA[E]$ for all $j = 0, 1, \dots$.
Thus~$\opt{x}$ stands for an element of~$\Unse{\LL}$.  By convention,
$D\opt{x} = 0$ if $\opt{x} = 1$.  Notice the subtle difference between
$\opt{DX} = 1$ and $D\opt{x} = 0$.

Let $\xi \colon \LL \hot \Unse{\LL} \to \Omega^1(\LL)$ be the bounded linear
map sending $x_1 \otimes \opt{x_2} \mapsto x_1 \,D\opt{x_2}$.  Thus $\xi( x_1
\otimes 1) = 0$ for all $x_1 \in \LL$.  Recall that the restriction
of~$\mu_{\ref{eq:mu8}}^{-1}$ to~$\LL$ gives rise to a bornological isomorphism
$\mu_{\ref{eq:mu9}}^{-1} \colon \LL \to \VS[G] \hot \opt{\LL} \subset \LL \hot
\opt{\LL}$.  We define a bounded linear map $\alpha \colon \LL \to
\Omega^1(\LL)$ as the composition $\alpha \defeq \xi \circ
\mu_{\ref{eq:mu9}}^{-1}$.  Thus $\alpha|_{\VS[G]} \equiv 0$.  We can also
define~$\alpha$ by
$$
\alpha(\vs[g] \odot \opt{x}) = \vs[g] \,D\opt{x} \qquad
\forall \vs[g] \in \VS[G],\ \opt{x} \in \Unse{\LL}.
$$
Here $D\opt{x} = 0$ if $\opt{x} = 1 \in \Unse{\LL}$.  Observe that
$\sigma(\MA[K]) \subset \VS[G]$, so that
\begin{equation}
  \label{eq:alphaKx}
  \alpha(\ma[k] \odot \opt{x}) = \ma[k] \,D\opt{x}
  \qquad \forall \ma[k] \in \MA[K],\ \opt{x} \in \Unse{\LL}.
\end{equation}
Moreover, $\alpha(\vs[g] \odot \opt{x_0} \odot x_1) = \vs[g] \,D (\opt{x_0}
\odot x_1) = \alpha(\vs[g] \odot \opt{x_0}) \odot x_1 + \vs[g] \odot \opt{x_0}
\,Dx_1$ implies that
\begin{equation}
  \label{eq:alphaMult}
  \alpha(x_0 \odot x_1) = \alpha(x_0) \odot x_1 + x_0 \,Dx_1 \qquad
  \forall x_0, x_1 \in \LL.
\end{equation}

By adjoint associativity, bounded linear maps $\MA[E] \to \Endo(\Tanil\LL)$
are the same thing as bounded bilinear maps $\MA[E] \times \Tanil\LL \to
\Tanil\LL$.  Thus we can define a bounded linear map $\triangleright \colon
\MA[E] \to \Endo(\Tanil\LL)$ as follows:
\begin{align}
  \label{eq:triOpen}
  \ma[e] \triangleright x_0 \,\opt{DX} &\defeq
  \ma[e] \odot x_0 \,\opt{DX} - D\alpha( \ma[e] \odot x_0) \,\opt{DX};
  \\
  \label{eq:triClosed}
  \ma[e] \triangleright Dx_1\,DX &\defeq
  \alpha(\ma[e] \odot x_1) \,DX.
\end{align}
The linear map $\triangleright \colon \MA[E] \to \Endo(\Tanil\LL)$
sends~$\ma[e]$ to the endomorphism $X \mapsto e \triangleright X$ of
$\Tanil\LL$.

\begin{lemma}  \label{lem:delta_lanilcur}
  $\triangleright \colon \MA[E] \to \Endo(\Tanil\LL)$ is a \lanilcur.
\end{lemma}

Before proving this, let us show how to construct a splitting homomorphism
$\LL \to \Tanil\LL$ out of~$\triangleright$.

\begin{lemma}  \label{lem:alpha_to_splitting}
  Extend~$\triangleright$ to a bounded unital homomorphism $\triangleright
  \colon \Unse{(\Tanil\MA[E])} \to \Endo(\Tanil\LL)$.  Define a bounded linear
  map $\upsilon \colon \LL \to \Tanil\LL$ by $\upsilon (y \odot \ma[k]) \defeq
  y \triangleright \ma[k]$ for $y \in \Unse{(\Tanil\MA[E])}$, $\ma[k] \in
  \MA[K]$.  Then~$\upsilon$ is actually a homomorphism and $\tau_{\LL} \circ
  \upsilon = \ID[\LL]$.
\end{lemma}

\begin{proof}
  Since~$\triangleright$ is a \lanilcur, it can be extended to a bounded
  homomorphism $\triangleright \colon \Tanil\MA[E] \to \Endo(\Tanil\LL)$.
  Extend this further to a unital bounded homomorphism $\triangleright \colon
  \Unse{(\Tanil\MA[E])} \to \Endo(\Tanil\LL)$.  By~\eqref{eq:mu5} there is a
  bornological isomorphism $\LL \cong \Unse{(\Tanil\MA[E])} \hot \MA[K]$.
  Thus~$\upsilon$ is a well-defined bounded linear map.  It remains to prove
  that~$\upsilon$ is multiplicative and that $\tau \circ \upsilon = \ID$.  It
  is clear from the definition \eqref{eq:triOpen} and~\eqref{eq:triClosed}
  that $\Janil\LL \subset \Tanil\LL$ is invariant under $\triangleright
  (\MA[E])$ and thus under $\triangleright \Unse{(\Tanil\MA[E])}$.  The
  induced left \Mp{\Tanil\MA[E]}action on $\Tanil\LL / \Janil\LL \cong \LL$
  coincides with the usual left \Mpn{\odot}multiplication
  by~\eqref{eq:triOpen}.  Thus $\tau(y \triangleright X) = y \odot \tau(X)$
  for all $y \in \Unse{(\Tanil\MA[E])}$, $X \in \Tanil\LL$.  Consequently,
  $\tau \circ \upsilon = \ID[\LL]$.
  
  To show that~$\upsilon$ is a splitting homomorphism, we prove first that
  $\triangleright(\ma[e])$ is a left multiplier for all $\ma[e] \in \MA[E]$.
  The bornological isomorphism~\eqref{eq:mu2} applied to $\Tanil\LL$ implies
  that $\Tanil\LL \cong \LL \hot \Unse{(\Tanil\LL)}$ as a right
  \Mpn{\LL}module.  Thus in order to show that $\triangleright(\ma[e])$ is a
  left multiplier, it suffices to prove that $\ma[e] \triangleright (x_0
  \circledcirc X) = (\ma[e] \triangleright x_0) \circledcirc X$ for all $x_0
  \in \LL$, $X \in \Tanil\LL$.  This follows if $\ma[e] \triangleright Dx_0
  \,Dx_1 = \ma[e] \triangleright (x_0 \odot x_1) - (\ma[e] \triangleright x_0)
  \circledcirc x_1$.  Since~$D$ is a graded derivation, \eqref{eq:alphaMult}
  implies
  $$
  D\alpha(\ma[e] \odot x_0 \odot x_1) =
  D\alpha(\ma[e] \odot x_0) \odot x_1 -
  \alpha(\ma[e] \odot x_0) \,Dx_1 +
  D(\ma[e] \odot x_0)\, Dx_1.
  $$
  It follows that
  \begin{multline*}
  \ma[e] \triangleright (x_0 \odot x_1) -
  (\ma[e] \triangleright x_0) \circledcirc x_1 =
  D(\ma[e] \odot x_0)\, Dx_1 -
  D\alpha(\ma[e] \odot x_0 \odot x_1) +
  D\alpha(\ma[e] \odot x_0) \circledcirc x_1 \\ =
  \alpha(\ma[e] \odot x_0) \,Dx_1 =
  \ma[e] \triangleright Dx_0\,Dx_1
  \end{multline*}
  as desired.  Thus $\triangleright(\ma[e])$ is a left multiplier.  In
  particular, \eqref{eq:alphaKx} implies
  $$
  \ma[k] \triangleright (x_0 \circledcirc X) =
  (\ma[k] \triangleright x_0) \circledcirc X =
  (\ma[k] \odot x_0 - D\ma[k] \,Dx_0) \circledcirc X =
  \ma[k] \circledcirc (x_0 \circledcirc X).
  $$
  Hence $\ma[k] \triangleright X = \ma[k] \circledcirc X$ for all $X \in
  \Tanil\LL$, $\ma[k] \in \MA[K]$.

  Since $\triangleright(\ma[e])$ is a left multiplier for all
  $\ma[e] \in \MA[E]$, it follows that $\triangleright(y)$ is a left multiplier
  for all $y \in \Unse{(\Tanil\MA[E])}$ because the left multipliers form a
  bornologically closed unital subalgebra of $\Endo(\Tanil\LL)$.  Therefore,
  \begin{multline*}
    \upsilon( y_1 \odot \ma[k]_1 \odot y_2 \odot \ma[k]_2) =
    (y_1 \odot \ma[k]_1 \odot y_2) \triangleright \ma[k]_2 =
    y_1 \triangleright \ma[k]_1 \triangleright y_2 \triangleright \ma[k]_2 =
    y_1 \triangleright \bigl( \ma[k]_1 \circledcirc
    (y_2 \triangleright \ma[k]_2) \bigr) \\ {} =
    (y_1 \triangleright \ma[k]_1) \circledcirc (y_2 \triangleright \ma[k]_2) =
    \upsilon(y_1\odot \ma[k]_1) \circledcirc \upsilon (y_2 \odot \ma[k]_2)
  \end{multline*}
  for all $y_1, y_2 \in \Unse{(\Tanil\MA[E])}$, $\ma[k]_1, \ma[k]_2 \in
  \MA[K]$.  Hence~$\upsilon$ is multiplicative.
\end{proof}

It remains to show that~$\triangleright$ is a \lanilcur.  Extend~$\alpha$ to a
bounded linear map again called $\alpha \colon \Omega_\an\LL \to
\Omega_\an\LL$ of degree~$1$ by $\alpha(x_0 \,\opt{DX}) \defeq \alpha(x_0)
\,\opt{DX}$ for all $x_0 \in \LL$ and $\alpha(DX) \defeq 0$.  A routine
calculation yields
\begin{align}
  \label{eq:OmtriOpen}
  \omega_\triangleright(\ma[e]_1, \ma[e]_2) (x_0\,\opt{DX})
  &=
  (\ma[e]_1 \ma[e]_2) \odot x_0\,\opt{DX} -
  D\alpha\bigl( (\ma[e]_1 \ma[e]_2) \odot x_0\bigr) \,\opt{DX}
  \\ \notag
  & \qquad {} -
  \ma[e]_1 \triangleright \bigl( \ma[e]_2 \odot x_0\,\opt{DX} -
  D\alpha(\ma[e]_2 \odot x_0)\,\opt{DX} \bigr)
  \\ \notag
  &=
  d\ma[e]_1 d\ma[e]_2\, x_0 \,\opt{DX} -
  D\alpha(d\ma[e]_1 d\ma[e]_2\,x_0) \,\opt{DX}
  \\ \notag
  & \qquad {} +
  \alpha\bigl( \ma[e]_1 \odot \alpha( \ma[e]_2 \odot x_0) \bigr) \,\opt{DX}.
\displaybreak[0]  \\
  \label{eq:OmtriClosed}
  \omega_\triangleright (\ma[e]_1, \ma[e]_2) (Dx_1\,DX)
  &=
  \alpha\bigl( (\ma[e]_1 \ma[e]_2) \odot x_1 \bigr) \,DX -
  \ma[e]_1 \triangleright \alpha (\ma[e]_2 \odot x_1) \,DX
  \\ \notag
  &=
  \alpha\bigl( (\ma[e]_1 \ma[e]_2) \odot x_1 \bigr) \,DX -
  \ma[e]_1 \odot \alpha (\ma[e]_2 \odot x_1 \bigr) \,DX
  \\ \notag
  & \qquad {} +
  D\alpha\bigl( \ma[e]_1 \odot \alpha (\ma[e]_2 \odot x_1) \bigr) \,DX.
\end{align}

Call a subset $S \subset \MA$ of a bornological algebra \emph{bornopotent}
iff~$S^\infty$ is small.  To show that~$\triangleright$ is a \lanilcur, we
have to check that $\omega_\triangleright (S,S)$ is bornopotent for all $S \in
\CBS(\MA[E])$.  The following criterion allows us to avoid computing powers of
$\omega_\triangleright (\ma[e]_1, \ma[e]_2)$ explicitly.

\begin{lemma}  \label{lem:bornopotent_Endo}
  Let~$\VS$ be a bornological vector space and $S \subset \Endo(\VS)$.
  Then~$S^\infty$ is small iff for each $T \in \CBS(\VS)$, there is a small
  set $T' \in \CBS(\VS)$ that contains~$T$ and is \Mpn{S}invariant.
  
  In other words, $S$ is bornopotent iff the \Mpn{S}invariant small sets are
  cofinal in $\CBS(\VS)$.
\end{lemma}

\begin{proof}
  If~$S^\infty$ is small, that is, equibounded, then $S^\infty (T) \cup T \in
  \CBS(\VS)$ by equiboundedness.  This small set is \Mpn{S}invariant by
  definition.  Conversely, if $T' \supset T$ is \Mpn{S}invariant, then
  $S^\infty(T) \subset T'$.  Hence if~$T'$ is small, so is $S^\infty(T)$.
  This means that~$S^\infty$ is equibounded.
\end{proof}

Let $S \in \CBS(\MA[E])$.  We are going to write down invariant sets for
$\omega_\triangleright (S,S)$ that are of the form
$$
F \defeq \coco{\bigl(F_0 (DF_\infty)^\even \cup DF_0
(DF_\infty)^\odd \bigr)} \in \CBS(\Tanil\LL)
$$
with suitable $F_0, F_\infty \in \CBS(\LL)$.  Recall that $A
\bipol{\subset} B$ stands for $A \subset \bipol{B}$.  If we have
\begin{equation}  \label{eq:FinvCondI}
  \begin{gathered}
    3dSdS\, F_0 \bipol{\subset} F_0,
    \\
    3\alpha( dSdS\,F_0) \bipol{\subset} F_0 \,DF_\infty,
    \\
    3\alpha \bigl( S \odot \alpha(S \odot F_0) \bigr) \bipol{\subset}
    F_0 \,DF_\infty \,DF_\infty,
  \end{gathered} \qquad
  \begin{gathered}
    3\alpha \bigl( (S \cdot S) \odot F_0 \bigr) \bipol{\subset}
    F_0 \,DF_\infty,
    \\
    3S \odot \alpha(S \odot F_0) \bipol{\subset} F_0 \,DF_\infty,
  \end{gathered}
\end{equation}
then $\omega_\triangleright (S,S)$ maps $F_0 (DF_\infty)^\even$ and $DF_0
(DF_\infty)^\odd$ into~$F$.  Use \eqref{eq:OmtriOpen}
and~\eqref{eq:OmtriClosed} and that~$F$ is a disk.  Sums of three terms are
converted into convex combinations in the obvious way.  Therefore,
\eqref{eq:FinvCondI} implies that~$F$ is invariant under
$\omega_\triangleright (S,S)$.

First we simplify~\eqref{eq:FinvCondI}.  If $\alpha (2S \odot F_0)
\bipol{\subset} F_0 \,DF_\infty$, then $3\alpha \bigl(S \odot \alpha (S \odot
F_0) \bigr) \subset F_0 \,DF_\infty \,DF_\infty$.  Furthermore, if both
$d(2S)d(2S) \,F_0 \bipol{\subset} F_0$ and $\alpha(F_0) \bipol{\subset} F_0
\,DF_\infty$, then also $3\alpha (dSdS\, F_0) \bipol{\subset} F_0
\,DF_\infty$.  Choose $S^{(2)} \in \CBS(\MA[E])$ containing $2S \cup 4S \cdot
S$.  Equation~\eqref{eq:FinvCondI} follows if
\begin{equation}  \label{eq:FinvCondII}
    dS^{(2)} dS^{(2)}\, F_0 \subset F_0, \qquad
    \opt{S^{(2)}} \odot \alpha (\opt{S^{(2)}} \odot F_0) \bipol{\subset}
    F_0 \,DF_\infty.
\end{equation}

Let $S^{(4)} = 2S^{(2)} \cup 2S^{(2)} \cdot S^{(2)}$.  Choose small disks
$S_{\MA[K]}' \in \CBS_d(\MA[K])$, $S_{\MA[Q]}' \in \CBS_d(\MA[Q])$ such
that~$S^{(4)}$ is contained in the disked hull of $S_{\MA[E]}' \defeq
S_{\MA[K]}' \cup s(S_{\MA[Q]}')$ and such that $\omega_s (S_{\MA[Q]}',
S_{\MA[Q]}')\subset S_{\MA[K]}'$.  To avoid notational clutter, we sometimes
write $S_{\MA[Q]}'$ for $s(S_{\MA[Q]}') \subset \MA[E]$.  Let
\begin{align*}
  S_{\MA[K]}''
  & \defeq
  \opt{S'_{\MA[E]}} \cdot S'_{\MA[K]} \cdot \opt{S'_{\MA[E]}}
  \in \CBS(\MA[K]);
  \\
  F_0
  & \defeq (dS^{(2)})^\even \odot \opt{S^{(2)}} \odot
  (dS_{\MA[Q]}')^\even \odot \opt{S_{\MA[Q]}'} \odot
  S_{\MA[K]}''
  \in \CBS(\LL).
\end{align*}
We claim that this~$F_0$ satisfies~\eqref{eq:FinvCondII} for
suitable~$F_\infty$.  The condition $dS^{(2)}dS^{(2)}\, F_0 \subset F_0$ is
fulfilled by construction.  The second condition in~\eqref{eq:FinvCondII}
follows if we show, for suitably big~$F_\infty$,
\begin{equation}
  \label{eq:FinvCondIII}
  \alpha(S^{(2)} \odot F_0)
  \bipol{\subset}
  (dS_{\MA[Q]}')^\even \odot \opt{S_{\MA[Q]}'} \odot S_{\MA[K]}'' \,DF_\infty
\end{equation}
Observe that $S^{(2)} \odot (dS^{(2)})^\even \odot \opt{S^{(2)}}$ is contained
in the disked hull of $(dS^{(4)})^\even \odot \opt{S^{(4)}}$.  This follows by
bringing an expression $\ma[e]_0 \odot d\ma[e]_1 \dots d\ma[e]_{2n}
\opt{\ma[e]_{2n+1}}$ into right handed standard form.  Hence
we can replace~\eqref{eq:FinvCondIII} by
\begin{equation}  \label{eq:FinvCondIV}
  \alpha \bigl(
    (dS'_{\MA[E]})^\even \odot \opt{S'_{\MA[E]}} \odot
    (dS'_{\MA[Q]})^\even \odot \opt{S'_{\MA[Q]}} \odot
    S''_{\MA[K]}
  \bigr)
  \bipol{\subset}
  (dS'_{\MA[Q]})^\even \opt{S'_{\MA[Q]}} \odot S''_{\MA[K]} \,DF_\infty.
\end{equation}
Now we have simplified \eqref{eq:FinvCondI} sufficiently and
prove~\eqref{eq:FinvCondIV}.  Choose $y_1 \in (dS'_{\MA[E]})^\even \odot
\opt{S'_{\MA[E]}}$ and $y_2 \in (dS'_{\MA[Q]})^\even \opt{S'_{\MA[Q]}} \odot
S_{\MA[K]}''$.  We assume first that $y_1 \in \II$.

By definition, $\alpha = \xi \circ \mu_{\ref{eq:mu9}}^{-1}$.  The map
$\mu_{\ref{eq:mu12}}^{-1}$ extends $\mu_{\ref{eq:mu9}}^{-1}$ to a right
\Mp{\Tanil\MA[E]}module homomorphism $\II \to \VS[G] \hot
\Unse{(\Tanil\MA[E])}$.  Thus we have
$$
\mu_{\ref{eq:mu9}}^{-1} (y_1 \odot y_2) =
\mu_{\ref{eq:mu12}}^{-1} (y_1) \odot y_2
$$
if $y_1 \in \II$.  Since $\mu_{\ref{eq:mu12}}^{-1}$ is equal to the
restriction of $\mu_{\ref{eq:mu8}}^{-1}$ to~$\II$, we can apply the
estimate~\eqref{eq:mu8_inverse_small} to obtain
$$
\mu_{\ref{eq:mu9}}^{-1} (y_1 \odot y_2) \bipol{\in}
(dS'_{\MA[Q]})^\even \opt{S'_{\MA[Q]}} \odot
(\opt{S'_{\MA[E]}} \cdot S'_{\MA[K]} \cdot \opt{S'_{\MA[E]}}) \otimes
F_\infty' =
(dS'_{\MA[Q]})^\even \opt{S'_{\MA[Q]}} \odot S''_{\MA[K]} \otimes F_\infty'.
$$
with $F_\infty' = 2 \opt{S'_{\MA[E]}} \odot (dS'_{\MA[E]})^\even
\opt{S'_{\MA[E]}} \odot (dS_{\MA[Q]}')^\even \odot \opt{S_{\MA[Q]}'} \odot
S_{\MA[K]}''$ depending on $S'_{\MA[E]}$ and~$S''_{\MA[K]}$.  Here $a
\bipol{\in} B$ means that $a \in \bipol{B}$.  Composing with~$\xi$, we get the
desired conclusion in the special case $y_1 \in \II$ for any small
set~$F_\infty$ satisfying $F_\infty \supset F_\infty'$.

If not $y_1 \in \II$, then $y_1 \in s_R\Unse{(\Tanil\MA[Q])}$, because~$y_1$
is a right polarized monomial.  We claim that
\begin{multline}  \label{eq:TQ_odot_G}
  y_1 \odot y_2 \bipol{\in} \VS[G] +
  (dS'_{\MA[Q]})^\odd dS'_{\MA[K]} \odot
  (2dS'_{\MA[Q]})^\even \opt{S'_{\MA[Q]}} \odot 4S''_{\MA[K]}
  \\ \cup
  (dS'_{\MA[Q]})^\even dS'_{\MA[K]} dS'_{\MA[Q]} \odot
  (2dS'_{\MA[Q]})^\even \opt{S'_{\MA[Q]}} \odot 4S''_{\MA[K]}
  \cup
  (dS'_{\MA[Q]})^\even S'_{\MA[K]} \odot 4S''_{\MA[K]}
\end{multline}
Before we prove this, we show how this completes the proof
of~\eqref{eq:FinvCondIV}.  Summands in~$\VS[G]$ are annihilated by~$\alpha$
and can be ignored.  The remaining terms are of the form $y_1' \odot y_2'$
with $y_1' \in \II \cap (dS'_{\MA[E]})^\even \opt{S'_{\MA[E]}}$ and $y_2' \in
(2dS'_{\MA[Q]})^\even \opt{S'_{\MA[Q]}} \odot 4S''_{\MA[K]}$.  Thus we are
almost in the case $y_1 \in \II$ considered above, up to factors of~$2$ that
occur only in~$y_2'$.  However, these factors go into the second tensor factor
and can be accommodated by making~$F'_\infty$ bigger.  Consequently, we
obtain~\eqref{eq:FinvCondI} for suitable~$F_\infty$.  To
prove~\eqref{eq:TQ_odot_G}, write
$$
y_1 = ds\ma[q]_1 \dots ds\ma[q]_{2n} \opt{s\ma[q]_{2n+1}}, \qquad
y_2 = ds\ma[q]_{2n+2} \dots ds\ma[q]_{2n+2m-1} \opt{s\ma[q]_{2n+2m}} \odot
\ma[k]
$$
with $s\ma[q]_j \in S'_{\MA[Q]}$ for all~$j$ and $\ma[k] \in S''_{\MA[K]}$.
Since $dy_1 dy_2 \in \VS[G]$, we have $y_1 \odot y_2 \in y_1 \cdot y_2 +
\VS[G]$.  Bringing $y_1 \cdot y_2$ into right handed standard form, we get a
sum of monomials
$$
\pm ds\ma[q]_1 \dots ds\ma[q]_{j-1} d(s\ma[q]_j s\ma[q]_{j+1})
ds\ma[q]_{j+2} \dots ds\ma[q]_{2n+2m-1} \opt{s\ma[q]_{2n+2m}}
\odot \ma[k], \qquad j = 2n+1, \dots, 2n+2m-2,
$$
and a final term $ds\ma[q]_1 \dots ds\ma[q]_{2n+2m-2} \cdot (s\ma[q]_{2n+2m-1}
\opt{s\ma[q]_{2n+2m}}) \odot \ma[k]$.  We can replace $s\ma[q]_j s\ma[q]_{j+1}$
by $-\omega_s(\ma[q]_j, \ma[q]_{j+1}) = - s(\ma[q]_j \ma[q]_{j+1}) + s\ma[q]_j
s\ma[q]_{j+1}$ because the difference gives rise to terms in~$\VS[G]$.  We have
assumed that $\omega_s(S'_{\MA[Q]}, S'_{\MA[Q]}) \subset S'_{\MA[K]}$.  Hence
we get a sum of terms in $(dS'_{\MA[Q]})^{j-1} dS'_{\MA[K]}
(dS'_{\MA[Q]})^{2n+2m-1 - j} \opt{S'_{\MA[Q]}} \odot S''_{\MA[K]}$ for all $j =
2n+1, \dots, 2n+2m-2$ and a final term in $(dS'_{\MA[Q]})^{2n+2m-2} \cdot
S'_{\MA[K]} \odot S''_{\MA[K]}$.  Finally, we convert this sum into a convex
combination by decorating each summand with factors of~$2$.  It is easy to
check that~\eqref{eq:TQ_odot_G} allows us sufficiently many factors because
geometric series converge.

Hence if we choose~$F_\infty$ suitably big, then~$F$ is invariant under
$\omega_\triangleright(S,S)$.  Since we can choose the small sets $S^{(2)}$,
$S_{\MA[K]}'$, $S_{\MA[Q]}'$, and~$F_\infty$ arbitrarily big, it follows that
invariant sets of the form~$F$ are cofinal in $\CBS(\Tanil\LL)$.  Thus
$\omega_\triangleright (S,S)$ is bornopotent for all $S \in \CBS(\MA[E])$ by
Lemma~\ref{lem:bornopotent_Endo}.  The map~$\triangleright$ is a \lanilcur.
The splitting homomorphism~$\upsilon$ in Lemma~\ref{lem:alpha_to_splitting}
can be defined.  The algebra~$\LL$ is analytically quasi-free.  The proof of
the Excision Theorem~\ref{the:excision_analytic} is complete.

\section{The Chern-Connes character in \Mpn{K}homology}
\label{sec:Chern_character}

The existence of a Chern character for suitable Fredholm modules was one of
the main motivations for the introduction of cyclic
cohomology~\cite{connes85:ncg}.  In cyclic cohomology, we can only obtain a
character for finitely summable Fredholm modules.  Connes enlarged the
periodic cyclic cohomology to entire cyclic cohomology in order to accommodate
the characters of \Mpn{\theta}summable Fredholm
modules~\cite{connes88:entire}.  However, the assumption of
\Mpn{\theta}summability is still a serious restriction.

Khalkhali~\cite{khalkhali94:connections} proves that $\HE^0(\MA) = \HP^0(\MA)
= \HH^0(\MA)$ and $\HE^1(\MA) = \HP^1(\MA) = 0$ if~$\MA$ is an \emph{amenable}
Banach algebra.  This applies, in particular, to nuclear \Cstar{}algebras
because these are amenable Banach algebras \cite{helemskii89:homology}.  Since
there are no bounded traces on a stable \Cstar{}algebra, we have $\HE^\ast
(\MA) = 0$ for all stable nuclear \Cstar{}algebras.  We get more reasonable
results if we endow \Cstar{}algebras with the precompact bornology $\COMP$,
that is, consider $\HA^\ast\bigl( (\MA, \COMP) \bigr)$.  We will construct a
Chern-Connes character $\chern \colon \KK_\ast(\MA,\C) \to
\HA^\ast(\MA,\COMP)$.  We write $K^\ast(\MA) \defeq \KK_\ast(\MA,\C)$ and call
this theory \emph{\Mpn{K}homology}, following Kasparov.  Thus $\HA^\ast(\MA,
\COMP)$ is large enough to accommodate the characters of arbitrary Fredholm
modules.  Alternatively, we can replace the precompact bornology by the fine
bornology $\FINE$.  The resulting theory $\HA^\ast(\MA, \FINE)$ is considered
already by Connes~\cite{connes94:ncg}.

There is a natural \emph{index pairing} $K_\ast(\MA) \times K^\ast(\MA) \to
\Z$ for $\ast = 0,1$.  The Chern-Connes characters in \Mpn{K}homology and
\Mpn{K}theory are compatible with the index pairing in the usual sense:
\begin{equation}  \label{eq:chern_compatible}
  \5{\chern(x)}{\chern(y)} = \5{x}{y}
  \qquad \forall x\in K_\ast(\MA),\ y \in K^\ast(\MA).
\end{equation}

The character that we construct here agrees with the usual Chern-Connes
character for finitely summable Fredholm modules.  I expect that it agrees
with the character defined for \Mpn{\theta}summable Fredholm
modules~\cite{connes94:ncg}.  However, I have not checked this.

Our construction of the Chern-Connes character is based on Connes's formulas
for the character of a finitely summable Fredholm module.  Consider a
\Mpn{1}summable even Fredholm module over an algebra~$\MA$.  Its character can
be described by cyclic cocycles $\tau_{2n} \colon \Omega^{2n}\MA \to \C$ for
any $n \in \Z_+$.  Since the cycles~$\tau_{2n}$ represent the same element of
periodic cyclic cohomology, there are bounded linear functionals $h_{2n+1}
\colon \Omega^{2n+1}\MA \to \C$ such that $h_{2n+1} \circ \partial =
\tau_{2n+2} - \tau_{2n}$.  Thus we can write~$\tau_0$ as a boundary: $\tau_0 =
\sum_{n=0}^\infty h_{2n+1} \circ \partial$.  Of course, $\sum_{n=0}^\infty
h_{2n+1}$ is unbounded because, in general, $[\tau_0] \neq 0$ in
$\HP^\ast(\MA)$.  We apply a suitable cut-off to each~$h_{2n+1}$, replacing it
by~$\tilde{h}_{2n+1}$.  Then $\tilde{\tau} = \tau_0 - \tilde{h}_{\ast} \circ
\partial$ should be the Chern-Connes character.  We automatically have
$\tilde{\tau} \circ \partial = 0$ and only have to control the growth
of~$\tilde{\tau}$ to check that it is a bounded linear functional.

The appropriate cut-off is achieved as follows.  Let~$(p_n)$ be an increasing
sequence of projections on~$\Hils$ that commute with~$F$ and satisfy $\rank
p_n = n$ and $\lim p_n = 1$ strongly.  Let $p_n^\bot = 1 - p_n$.  Up to
normalization constants, $h_{2n+1}(\opt{\ma_0} d\ma_1\dots d\ma_{2n+1}) = \tr(
\gamma F \opt{\ma_0} [F,\ma_1] \dots [F,\ma_{2n+1}])$.  We put
\begin{multline*}
  \tilde{h}_{2n+1}(\opt{\ma_0} d\ma_1 \dots d\ma_{2n+1}) =
  h_{2n+1}\bigl(\opt{p_n^\bot \ma_0 } d(p_n^\bot \ma_1) \dots
  d(p_n^\bot \ma_{2n+1}) \bigr)
  \\ =
  \tr( \gamma F p_n^\bot \opt{\ma_0} p_n^\bot [F,\ma_1] p_n^\bot [F,\ma_2]
  p_n^\bot \dots [F,\ma_{2n+1}]).
\end{multline*}
This does the trick.  We compute that $\tau_0 - \sum \tilde{h}_{2n+1} \circ
\partial$ is a sum of traces of products where one of the factors is~$p_n$ and
there are, in degree~$2n$, about~$2n$ factors of the form $p_n^\bot
[F,\ma_j]$.  Since the operators~$p_n$ have finite rank, we only have to take
traces of finite rank operators.  This is well-defined without any summability
condition on the commutators $[F,\ma_{2n+1}]$.  In addition, if $S\subset
\Comp(\Hils)$ is precompact, then $p_n^\bot S \to 0$ uniformly in the operator
norm.  If we fix a precompact set $S\subset \MA$, then $[F,S]$ is a precompact
subset of $\Comp(\Hils)$.  Thus the operator norm of a product
$$
\| p_n^\bot [F,\ma_1] p_n^\bot [F,\ma_2] \cdots p_n^\bot [F,\ma_{2n}]
\|_\infty, \qquad
\ma_1,\dots,\ma_{2n}\in S,
$$
can be estimated by $C_\epsilon \epsilon^{2n}$ for all $\epsilon>0$, with
constants~$C_\epsilon$ depending only on~$\epsilon$ and~$S$.  This produces
the exponential decay needed for analytic cocycles.

Actually, it is not necessary that~$\Hils$ is a Hilbert space.  In the actual
construction, we replace~$\Hils$ by a separable Banach space that has
Grothendieck's metric approximation property.  Also, the operators~$p_n$ need
not be projections.  The precise requirements are listed below.

Unfortunately, these elementary techniques only work for \Mpn{K}homology, not
for the bivariant case.  In the case of \Mpn{K}homology, we have to construct
bounded linear functionals with values in~$\C$.  In the bivariant situation,
the range should be the complicated complex $X(\Tanil\MA[B])$.  I expect that
there is no bivariant Chern character from $\KK_\ast(\MA,\MA[B])$ to
$\HA^\ast(\MA; \MA[B])$.  I have some hope that analytic cyclic cohomology
$\HA^\ast(\MA,\COMP)$ is well-behaved on the category of \Cstar{}algebras,
that is, invariant under continuous homotopies and stable with respect to the
spatial tensor product with the compact operators $\Comp(\Hils)$.  However, I
have not been able to prove these assertions.

There is a bivariant Chern character in Puschnigg's \emph{local} analytic
cyclic cohomology~\cite{puschnigg98:cyclic}.  This follows from the following
crucial property of the local theory: If $\MA \subset \MA[B]$ is a dense
``smooth'' (see~\cite{blackadar91:smooth}) subalgebra of a nuclear
\Cstar{}algebra~$\MA[B]$, then the inclusion $\MA \to \MA[B]$ is an invertible
element in the local theory $H^\ast (\MA; \MA[B])$.  Since $\CINF([0,1];
\MA[B]) \subset \NBC([0,1]; \MA[B])$ is a dense smooth subalgebra, the local
theory does not distinguish between $\CINF([0,1]; \MA[B])$ and $\NBC([0,1];
\MA[B])$.  Thus smooth homotopy invariance implies continuous homotopy
invariance.  Since $\Sch(\Hils) \prot \MA[B] \subset \Comp(\Hils)
\otimes_{\mathrm{spatial}} \MA[B]$ is a dense smooth subalgebra, stability for
projective tensor products with $\Sch(\Hils)$ implies \Cstar{}stability.
Excision carries over to the local theory.  Thus the local theory, when
restricted to \Cstar{}algebras, has all the desirable homological properties.

\subsection{The finitely summable case}
\label{sec:Chern_finite}

The construction of the character for Fredholm modules without summability
condition depends on detailed formulas for the Chern-Connes character in the
finitely summable case.  As a preparation, we therefore repeat Connes's
construction of the Chern-Connes character in the finitely summable case
\cite{connes85:ncg}, \cite{connes94:ncg}.  The main difference is that we use
different normalization constants.  Furthermore, we work in a slightly more
general situation, considering representations on (almost) arbitrary complete
bornological vector spaces, not just on Hilbert space.  This creates no
additional difficulty.

Let~$\VS$ be a complete bornological vector space.  Let $\Fin(\VS) \defeq \VS
\otimes \VS'$ be the algebra of bounded finite rank operators on~$\VS$.  An
elementary tensor $\vs \otimes \vs' = \ket{\vs} \bra{\vs'}$ corresponds to the
rank one endomorphism $\VS \ni \vs_2 \mapsto \vs\cdot \vs'(\vs_2) \in \VS$
of~$\VS$.  Clearly, the resulting homomorphism $\Fin(\VS) \to \Endo(\VS)$ is a
bijection onto the ideal of finite rank operators.  Thus $\Fin(\VS)$ is a
bimodule over $\Endo(\VS)$ with $T\cdot \ket{\vs} \bra{\vs'} = \ket{T\vs}
\bra{\vs'}$ and $\ket{\vs}\bra{\vs'} T = \ket{\vs} \bra{T'\vs'}$, where
$T'\colon \VS' \to \VS'$ denotes the transpose of~$T$.  The trace $\tr$
defined by $\tr\, \ket{\vs} \bra{\vs'} \defeq \vs'(\vs)$ is a
\Mp{\Endo(\VS)}bimodule trace in the sense that $\tr (Tx) = \tr (xT)$ for all
$x\in \Fin(\VS)$, $T\in \Endo(\VS)$.  It suffices to verify this relation on
elementary tensors $x = \ket{\vs} \bra{\vs'}$, where it is equivalent to the
definition $\5{\vs'}{T\vs} = \5{T'\vs'}{\vs}$ of the transpose.

We now complete $\Fin(\VS)$ to the complete bornological algebra $\Sch(\VS)
\defeq \VS \hot \VS'$ of ``nuclear endomorphisms'' of~$\VS$.  The bimodule
structure $\Endo(\VS) \times \Fin(\VS) \times \Endo(\VS) \to \Fin(\VS)$
extends uniquely to a bounded trilinear map $\Endo(\VS) \times \Sch(\VS)
\times \Endo(\VS) \to \Sch(\VS)$ making $\Sch(\VS)$ a \Mp{\Endo(\VS)}bimodule.
The trace $\tr\colon \Fin(\VS) \to \C$ extends uniquely to a bounded linear
functional $\tr \colon \Sch(\VS) \to \C$ that is a \Mp{\Endo(\VS)}bimodule
trace on $\Sch(\VS)$ in the sense that
\begin{equation}  \label{eq:tr_trace}
  \tr(Tx) = \tr (xT) \qquad \forall x\in\Sch(\VS),\ T\in \Endo(\VS).
\end{equation}
The inclusion $\Fin(\VS) \subset \Endo(\VS)$ extends to a bounded homomorphism
$\natural \colon \Sch(\VS) \to \Endo(\VS)$.  In general~$\natural$ need not be
injective.  However, $\natural$ is injective if~$\VS$ is a Banach space having
Grothendieck's approximation property.  Hence for our later applications we
can safely assume that~$\natural$ is injective.  We make this assumption and
consider $\Sch(\VS) \subset \Endo(\VS)$ because it simplifies the definitions.

We can now define \Mpn{1}summable Fredholm modules.  In the odd case, we need
an operator $F\in \Endo(\VS)$ satisfying $F^2 = \ID$.  Let
$$
\MA[B]_\odd^0 \defeq \MA[B]_\odd^0(\VS,F) \defeq
\{ T\in \Endo(\VS) \mid [F,T] \in \Sch(\VS) \}.
$$
We write $\MA[B]_\odd^0$ if $(\VS,F)$ is clear from the context.  Since
taking commutators is a derivation, $\MA[B]_\odd^0$ is a subalgebra of
$\Endo(\VS)$.  We endow $\MA[B]_\odd^0$ with the coarsest bornology for which
the inclusion $\MA[B]_\odd^0 \to \Endo(\VS)$ and the derivation $\MA[B]_\odd^0
\to \Sch(\VS)$ sending $T\mapsto [F,T]$ are bounded.

This bornology makes $\MA[B]_\odd^0$ a complete bornological algebra.  The
easiest way to see this is to decompose~$\VS$ into eigenspaces for~$F$.
Let~$\VS_{\pm}$ be the range of the idempotent $\frac{1}{2}(1\pm F)$.  Then
$\VS = \VS_+ \oplus \VS_-$.  Write an endomorphism $T\in \Endo(\VS)$ as a
\Mp{2\times 2}matrix with respect to this decomposition.  Then $T\in
\MA[B]_\odd^0$ iff the off-diagonal terms are nuclear operators, that is,
elements of $\Sch(\VS_+,\VS_-)$ and $\Sch(\VS_-,\VS_+)$, respectively.  It
follows that $\MA[B]_\odd^0 \cong \Endo(\VS_+) \oplus \Endo(\VS_-) \oplus
\Sch(\VS_+,\VS_-) \oplus \Sch(\VS_-,\VS_+)$ as a bornological vector space.
This is evidently complete.

Let~$\MA$ be a complete bornological algebra.  A \emph{\Mpn{1}summable odd
  \Mp{(\VS,F)}Fredholm module} over~$\MA$ is, by definition, a bounded
homomorphism $\phi\colon \MA \to \MA[B]_\odd^0(\VS,F)$.

In the even case, we add an operator $\gamma \in \Endo(\VS)$ satisfying
$F\gamma + \gamma F = 0$ and $\gamma^2 = \ID$.  Let
$$
\MA[B]_\even^0 \defeq \MA[B]_\even^0(\VS,F,\gamma) \defeq
\{ T\in \Endo(\VS) \mid [F,T] \in \Sch(\VS),\ [\gamma,T] = 0 \}.
$$
We write $\MA[B]_\even^0$ if $(\VS,F,\gamma)$ is clear from the context.
It is easy to see that $\MA[B]_\even^0$ is a closed subalgebra of
$\MA[B]_\odd^0$.  We endow it with the subspace bornology.  This is the
coarsest bornology making the inclusion to $\Endo(\VS)$ and the derivation $T
\mapsto [F,T]$ to $\Sch(\VS)$ bounded.

Let~$\MA$ be a complete bornological algebra.  A \emph{\Mpn{1}summable even
  \Mp{(\VS,F,\gamma)}Fredholm module} over~$\MA$ is, by definition, a bounded
homomorphism $\phi\colon \MA \to \MA[B]_\even^0(\VS,F,\gamma)$.

We now construct the Chern-Character for even and odd Fredholm modules.  We
begin with some observations that are relevant for both cases.  Let $\delta
\colon \MA[B]_\odd^0 \to \Sch(\VS)$ be the derivation $\delta(T) \defeq
\frac{i}{2} [F,T]$.  The normalization constant~$i/2$ insures that certain
normalization constants later will be just~$1$ instead of~$(-4)^{n}$.  The
derivation~$\delta$ induces a bounded homomorphism $\psi \colon
\Omega\MA[B]_\odd^0 \to \MA[B]_\odd^0$,
$$
\psi( \opt{x_0} dx_1 \dots dx_n) \defeq
\opt{x_0} \cdot \delta x_1 \cdots \delta x_n.
$$
The restriction of~$\psi$ to $\sum_{j=1}^\infty \Omega^j \MA[B]_\odd^0$ is
a bounded linear map $\psi \colon \Omega^{\ge 1} \MA[B]_\odd^0 \to \Sch(\VS)$.

A computation using $F^2 = 1$ shows that~$F$ anti-commutes with~$\delta(x)$
for all $x \in \MA[B]_\odd^0$.  Hence
\begin{multline*}
  F \psi(\opt{x_0} dx_1 \dots dx_n) -
  (-1)^n \psi(\opt{x_0} dx_1 \dots dx_n) F =
  (F \opt{x_0} - \opt{x_0} F) \delta x_1 \dots \delta x_n
  \\ =
   \tfrac{2}{i} \cdot \delta\opt{x_0} \delta x_1 \dots \delta x_n =
   \tfrac{2}{i} \psi(d\opt{x_0} dx_1 \dots dx_n)
\end{multline*}
for all $x_0,\dots,x_n\in \MA[B]_\odd^0$.  More succinctly,
\begin{equation}
  \label{eq:F_commutator_d}
  F \psi(\omega) - (-1)^n \psi(\omega) F = \tfrac{2}{i} \psi(d\omega) \qquad
  \forall \omega\in \Omega^n\MA[B]_\odd^0.
\end{equation}

We now restrict attention to the even case.  For $n \in \Z_+$, we define
$\tau_{2n}\colon \Omega^{2n}\MA[B]_\even^0 \to \C$ by
\begin{equation}  \label{eq:tau_even_def}
  \tau_{2n} (\omega) \defeq
  c_{2n} \tr\bigl( iF \gamma \psi(d\omega) \bigr)
\end{equation}
for all $\omega \in \Omega^{2n}\MA[B]_\even^0$ with certain normalization
constants
$$
c_{2n} \defeq \frac{2n(2n-2)(2n-4) \cdots 2}{(2n-1)(2n-3)(2n-5)\cdots 1}.
$$
Thus $c_0=1$ and $c_{2n+2} = (2n+2)/(2n+1) c_{2n}$ for all $n\in\Z_+$.  For
$n\ge1$, we can use~\eqref{eq:F_commutator_d} and~\eqref{eq:tr_trace} to
rewrite the definition of~$\tau_{2n}$ as
\begin{equation}  \label{eq:tau_even_rewrite}
  \tau_{2n} (\omega) =
  -\frac{c_{2n}}{2}
  \tr\bigl( F \gamma (F\psi(\omega) - \psi(\omega) F) \bigr) =
  \frac{c_{2n}}{2}
  \tr\bigl(  2F^2 \gamma \psi(\omega) \bigr) =
  c_{2n} \tr\bigl( \gamma \psi(\omega) \bigr).
\end{equation}
Equation~\eqref{eq:tr_trace} applies because $\psi(\omega) \in \Sch(\VS)$.  Up
to normalization constants this is Connes's formula for the character of a
\Mp{2n}summable even Fredholm module \cite[p.~293]{connes94:ncg}.

\begin{lemma}  \label{lem:tau_clgrt_even}
  The functionals~$\tau_{2n}$ are closed graded traces for all $n\in\Z_+$.
  That is, $\tau_{2n} \circ d = 0$ and~$\tau_{2n}$ vanishes on graded
  commutators.  The latter condition implies (in fact, is equivalent to)
  $\tau_{2n} \circ \kappa = \tau_{2n}$ and $\tau_{2n}\circ b = 0$ for
  all $n\in\Z_+$.
\end{lemma}

\begin{proof}
  By definition, $\tau_{2n} \circ d = 0$.  Since~$\gamma$ anti-commutes
  with~$F$ and commutes with $\MA[B]_\even^0$, we have $\gamma \psi(\omega) =
  (-1)^{\deg \omega} \psi(\omega) \gamma$ for all homogeneous $\omega \in
  \Omega \MA[B]_\even^0$.  In addition, \eqref{eq:F_commutator_d} implies that
  $F\psi(d\omega) + (-1)^{\deg \omega} \psi(d\omega) F = \frac{2}{i}
  \psi(dd\omega) = 0$.  Consequently, $F\gamma$ commutes with $\psi(d\omega)$
  for all $\omega \in \Omega\MA[B]_\even^0$.  Let $\omega \in
  \Omega^k\MA[B]_\even^0$, $\omega'\in \Omega^{2n-k}\MA[B]_\even^0$, then
  \eqref{eq:tr_trace} implies
  \begin{multline*}
    \tau_{2n}([\omega,\omega']) =
    c_{2n} \tr\bigl( iF \gamma
      \psi( d(\omega \omega' - (-1)^{k} \omega' \omega)) \bigr)
    \\ =
    ic_{2n} \tr\bigl(
        F\gamma \psi(d\omega) \psi(\omega')
      + (-1)^k F\gamma \psi(\omega) \psi(d\omega')
      - (-1)^k F\gamma \psi(d\omega') \psi(\omega)
      - F\gamma \psi(\omega') \psi(d\omega)
    \bigr)
    \\ =
    ic_{2n} \tr \bigl(
      [F\gamma, \psi(d\omega)] \cdot \psi(\omega')
    - (-1)^k [F\gamma, \psi(d\omega')] \cdot \psi(\omega)
    \bigr)
    = 0.
  \end{multline*}
  Thus~$\tau_{2n}$ is a closed graded trace.  Since $b(\omega)$ and
  $(1-\kappa)(\omega)$ are graded commutators by definition, it follows that
  $\tau_{2n+1} \circ b = 0$ and $\tau_{2n+1} \circ \kappa = \tau_{2n+1}$.
\end{proof}

The cyclic cocycles~$\tau_{2n}$ should define the same element of periodic
cyclic cohomology.  In fact, we can write down explicit linear functionals
$h_{2n+1} \colon \Omega^{2n+1}\MA[B]_\even^0 \to \C$ such that $h_{2n+1} \circ
\partial = \tau_{2n} - \tau_{2n+2}$ for all $n\in\Z_+$.  Here~$\partial$
denotes the boundary in the X-complex.  Since~$h_{2n+1}$ lives on the odd part
only, the relevant part of~$\partial$ is given by~\eqref{eq:partial_even}.
Define
$$
h_{2n+1}(\omega) = \frac{c_{2n}}{2n+1} \tr\bigl(iF \gamma \psi(\omega) \bigr)
$$
for all $\omega\in \Omega^{2n+1}\MA[B]_\even^0$.

\begin{lemma}  \label{lem:tau_h_even}
  We have $h_{2n+1} \circ \kappa = h_{2n+1}$ and $h_{2n+1} \circ \partial =
  \tau_{2n} - \tau_{2n+2}$ for all $n\in\Z_+$.
\end{lemma}

\begin{proof}
  If $x\in \MA[B]^0_\even$ and $\omega \in \Omega\MA[B]_\even^0$, then
  $$
  \tr \bigl(F \gamma \psi( dx\cdot \omega)\bigr) =
  \tr \bigl(F \gamma \delta x \cdot \psi(\omega)\bigr) =
  \tr \bigl(\delta x \cdot F \gamma \psi(\omega) \bigr) =
  \tr \bigl(F \gamma \psi(\omega) \delta x  \bigr) =
  \tr \bigl(F \gamma \psi(\omega dx) \bigr)
  $$
  because $F \gamma$ and~$\delta x$ commute.  This implies $h_{2n+1} \circ
  \kappa = h_{2n+1}$.  The boundary~$\partial$ from even to odd degrees is
  given by~\eqref{eq:partial_even}.  Since $h_{2n+1} \circ \kappa = h_{2n+1}$,
  the desired equation $h_{2n+1} \circ \partial = \tau_{2n} - \tau_{2n+2}$ is
  equivalent to $(2n+1) h_{2n+1} \circ d = \tau_{2n}$ and $(n+1) h_{2n+1}
  \circ b = \tau_{2n+2}$.  The first equality is the definition
  of~$\tau_{2n}$.  For the second equation, let $\omega \in \Omega^{2n+1}
  \MA[B]_\even^0$, $x\in \MA[B]_\even^0$ and compute
  \begin{multline*}
    \tr\bigl(iF \gamma \cdot \psi\circ b(\omega dx) \bigr) =
    \tr\bigl(iF \gamma \psi(-[\omega,x]) \bigr) =
    \tr\bigl(iF \gamma x \psi(\omega) - iF \gamma \psi(\omega) x \bigr)
    \\ =
    \tr\bigl(iF x \gamma \psi(\omega) - \gamma \psi(\omega) x iF \bigr)
    =
    \tr\bigl(\gamma \psi(\omega) [iF,x] \bigr) =
    2 \tr\bigl(\gamma \psi(\omega) \delta x \bigr) =
    2 \tr\bigl(\gamma \psi(\omega dx)\bigr).
  \end{multline*}
  Decorated with appropriate normalization constants, this means $(n+1)
  h_{2n+1} \circ b = \tau_{2n+2}$.
\end{proof}

Consequently, we have $[\tau_{2n}] = [\tau_0]$ in $\HP^0(\MA[B]_\even^0)$ for
all $n\in\Z_+$.  To check the normalization constant~$c_0$, we compute the
functional $K_0(\MA[B]_\even^0) \to \C$ induced by $[\tau_{2n}] = [\tau_{0}]$.
Actually, it suffices to compute $\5{\tau_0}{\chern(p)}$ for a single rank one
idempotent $p \in \MA[B]_\even^0$ with $\gamma p = p \gamma = p$.  We obtain
the desired result
$$
\tau_0\bigl( \chern(p)\bigr) =
\tau_0(p) =
c_0 \tr (iF \gamma [iF/2,p]) =
\tfrac{1}{2} \tr\bigl( \gamma (p -  Fp F) \bigr) =
\tr(p) =
1.
$$

Let us now turn to the odd case.  Recall that we left a normalization
constant~$c_1$ in the Chern-Connes character $K_1 \to \HA_1$.  In order to get
the right index pairing, we have to compensate this constant in $K$-homology.
Let $\tilde{c}_1 = c_1^{-1}$.  For $n\in\Z_+$, define $\tau_{2n+1}\colon
\Omega^{2n+1} \MA[B]_\odd^0 \to \C$ by
\begin{equation}  \label{eq:tau_odd_def}
  \tau_{2n+1} (\omega) \defeq
  \tilde{c}_1 \tr\bigl( F \psi(d\omega) \bigr) =
  \tilde{c}_1 \tfrac{i}{2} \tr\bigl( F (F\psi(\omega) +
  \psi(\omega) F) \bigr) =
  i\tilde{c}_1 \tr\bigl( \psi(\omega) \bigr)
\end{equation}
for all $\omega\in \Omega^{2n+1}\MA[B]_\odd^0$.  The reformulations
in~\eqref{eq:tau_odd_def} use \eqref{eq:tr_trace}
and~\eqref{eq:F_commutator_d}.  Up to normalization constants, this is
Connes's formula for the character of an odd Fredholm module
\cite[p.~293]{connes94:ncg}.

\begin{lemma}  \label{lem:tau_clgrt_odd}
  The functionals~$\tau_{2n+1}$ are closed graded traces for all $n\in\Z_+$.
  That is, $\tau_{2n+1} \circ d = 0$ and~$\tau_{2n+1}$ vanishes on graded
  commutators.  The latter condition implies (in fact, is equivalent to)
  $\tau_{2n+1} \circ \kappa = \tau_{2n+1}$ and $\tau_{2n+1}\circ b = 0$ for
  all $n\in\Z_+$.
\end{lemma}

\begin{proof}
  Equation~\eqref{eq:tau_odd_def} shows that $\tau_{2n+1}\circ d = 0$.
  Equation~\eqref{eq:tr_trace} implies that $\tau_{2n+1} (\omega \cdot
  \omega') = \tau_{2n+1}(\omega' \cdot \omega)$ for all $\omega \in \Omega^k
  \MA[B]_\odd^0$, $\omega'\in \Omega^{2n+1-k} \MA[B]_\odd^0$.  Since either
  $k$ or $2n+1-k$ is even, the graded commutator of $\omega$ and~$\omega'$ is
  $\omega\omega' - \omega'\omega$.  Hence~$\tau_{2n+1}$ vanishes on all graded
  commutators.
\end{proof}

The cyclic cocycles~$\tau_{2n+1}$ should define the same element of periodic
cyclic cohomology.  In fact, we can write down explicit linear functionals
$h_{2n}\colon \Omega^{2n}\MA[B]_\odd^0 \to \C$ such that $h_{2n} \circ
\partial = \tau_{2n-1} - \tau_{2n+1}$ for all $n\in\N$.  Here~$\partial$
denotes the boundary in the X-complex.  Since~$h_{2n}$ lives on the even part,
$\partial$ is given by~\eqref{eq:partial_odd} and can be replaced by $b -
(1+\kappa)d$.  Define
$$
h_{2n}(\omega) = -\tfrac{1}{2} \tilde{c}_1 \tr\bigl(F \psi(\omega) \bigr).
$$
This is well-defined for all $n\ge 1$, but not for $n=0$.

\begin{lemma}  \label{lem:tau_h_odd}
  We have $h_{2n} \circ \kappa = h_{2n}$ and $h_{2n} \circ \partial =
  \tau_{2n-1} - \tau_{2n+1}$ for all $n\in\N$.
\end{lemma}

\begin{proof}
  If $x\in \MA[B]_\odd$ and $\omega \in \Omega\MA[B]_\odd^0$, then
  $$
  \tr \bigl(F \psi( dx\cdot \omega)\bigr) =
  \tr \bigl(F \delta x \psi(\omega) \bigr) =
  -\tr \bigl(\delta x F \psi(\omega) \bigr) =
  -\tr \bigl(F \psi(\omega) \delta x \bigr) =
  -\tr \bigl(F \psi(\omega dx) \bigr)
  $$
  because $F$ and~$\delta x$ anti-commute.  This implies $h_{2n} \circ
  \kappa = h_{2n}$.  By definition, $h_{2n} \circ d = -\frac{1}{2}
  \tau_{2n-1}$.  Therefore, $h_{2n}\circ (1+\kappa)d = 2 h_{2n} d = -
  \tau_{2n-1}$.  Moreover,
  \begin{multline*}
    h_{2n} \circ b(\omega dx) =
    h_{2n} (\omega\cdot x - x\cdot \omega) =
    -\tfrac{1}{2} \tilde{c}_1
    \tr\bigl(F \psi(\omega) x - F x\psi(\omega) \bigr)
    \\ =
    -i\tilde{c}_1 \tr\bigl(-\psi(\omega) x iF/2 + \psi(\omega) (iF/2)x \bigr)
    =
    -i\tilde{c}_1 \tr(\psi(\omega) \delta x)
    =
    -i\tilde{c}_1 \tr\bigl(\psi(\omega dx) \bigr).
  \end{multline*}
  Thus $h_{2n} \circ b = - \tau_{2n+1}$ and consequently $h_{2n} \circ
  \partial = h_{2n} \circ (b - (1+\kappa)d) = \tau_{2n-1} - \tau_{2n+1}$.
\end{proof}

It follows that $[\tau_{2n+1}] = [\tau_1]$ in $\HP^0(\MA[B]_\odd^0)$ for all
$n\in\Z_+$.  To check our normalization constants, we restrict attention to
the case where~$\VS$ is a separable Hilbert space, $F = F^\ast$, and the
projections $\frac{1}{2}(1\pm F)$ have infinite dimensional ranges $\VS_\pm$.
Thus $\VS_+$ and~$\VS_-$ are isomorphic.  We write elements of $\Endo(\VS)$ as
\Mp{2\times 2}matrices over~$\VS_+$.  Let $S\colon \VS_+ \to \VS_+$ be the
unilateral shift, $S^\ast$ its adjoint and $p \defeq \ID - SS^\ast$.  Thus~$p$
is a rank one projection.  Let
$$
u \defeq
\begin{pmatrix} S^\ast & 0 \\ p & S 
\end{pmatrix}.
$$
This is the bilateral shift and hence unitary.  The compression of~$u$
to~$\VS_+$ is~$S^\ast$ and therefore has index~$1$.  An easy computation shows
that
$$
u^{-1} [F,u] = 
\begin{pmatrix} -2p & 0 \\ 0 & 0 
\end{pmatrix}.
$$
Hence we get the desired result
$$
\5{\tau_1}{\chern(u)} =
\tau_1(c_1 u^{-1} du) =
i\tilde{c}_1 c_1 \tr( u^{-1}[iF/2,u] ) =
-\tr(-2p/2) =
1.
$$

\subsection{The metric approximation property and cut-off sequences}
\label{sec:cut_off}

Let~$\VS$ be a Banach space and let $F$ and possibly~$\gamma$ be as above.
The algebra $\Sch(\VS) \defeq \VS \hot \VS'$ is a Banach space with the
projective tensor product norm $\|\blank\|_1$.  Let $\Comp(\VS)$ be the ideal
of compact endomorphisms of~$\VS$.  This is a closed subalgebra of
$\Endo(\VS)$ with respect to the operator norm $\|\blank\|_\infty$ on
$\Endo(\VS)$.  We assume that there is a sequence~$(p_n)$ in $\Sch(\VS)$ with
the following properties:
\begin{enumerate}%
  
\item the trace norms $\| p_n \|_1$ in $\Sch(\VS)$ are of at most exponential
  growth;
  
\item the operator norms $\| p_n \|_\infty$ in $\Comp(\VS)$ are uniformly
  bounded;
  
\item the set of $\vs \in \VS$ with $\lim p_n(\vs) = \vs$ is dense in~$\VS$;
  
\item $p_n$ commutes with $F$ (and~$\gamma$ in the even case) for all
  $n \in \N$.

\end{enumerate}
We call~$(p_n)$ a \emph{cut-off sequence}.  If there is a sequence $(p_n)$ in
$\Sch(\VS)$ that satisfies 2 and~3, then there is a cut-off sequence.  The
first condition can be achieved by repeating each~$p_n$ sufficiently often to
insure that the norms $\|p_n\|_1$ grow slowly.  The fourth condition can be
achieved by averaging.  Replace~$p_n$ by $\frac{1}{2}(p_n + Fp_n F)$ to make
it commute with~$F$.  If necessary replace the result by $\frac{1}{2}(p_n +
\gamma p_n \gamma)$ to make~$p_n$ commute with~$\gamma$.

An equibounded net of operators that converges pointwise on a dense subset
automatically converges pointwise on all of~$\VS$, and the convergence is
uniform on precompact subsets of~$\VS$.  This follows because each precompact
subset of~$\VS$ is contained in the disked hull of a sequence~$(x_n)$
in~$\VS$, where we can achieve that all~$x_n$ are contained in a given dense
subspace of~$\VS$.  Hence condition~3 can be replaced by the equivalent
requirement that $p_n \to \ID$ uniformly on compact subsets of~$\VS$.
Moreover, we can replace the~$(p_n)$ by nearby operators of finite rank.  That
is, if there is a cut-off sequence~$(p_n)$, then there is another cut-off
sequence with $p_n' \in \Fin(\VS)$ for all $n \in \N$.

Consequently, a cut-off sequence exists iff~$\VS$ is separable and has
\emph{Grothendieck's metric approximation
  property}~\cite{grothendieck55:produits}.  Hence $\Fin(\VS)$ is dense in
$\Comp(\VS)$ and the natural map $\Sch(\VS) \to \Endo(\VS)$ is injective.

Let $p_n^\bot = 1 - p_n$.  For a subset~$S$ of a Banach space~$\VS[W]$, let
$\|S\| = \sup \{ \|x\|_{\VS[W]} \mid x\in S\}$.

\begin{lemma}  \label{lem:cutoff_on_compacts}
  Let~$(p_n)$ be a cut-off sequence.  If $S\subset \Sch(\VS)$ is precompact,
  then $\| p_n^\bot S \|_1 \to 0$ for $n \to \infty$.  If $S \subset
  \Comp(\VS)$ is precompact, then $\| p_n^\bot S \|_\infty \to 0$ for $n \to
  \infty$.
\end{lemma}

\begin{proof}
  Let $S\subset \Sch(\VS) = \VS\prot \VS'$ be precompact.  Then there are
  compact disks $S_1 \subset \VS$, $S_2 \subset \VS'$ such that $S \subset
  \coco{(S_1\otimes S_2)}$ by Grothendieck's Theorem~\ref{the:Grothendieck}.
  We have $p_n^\bot \circ S \subset \coco{(p_n^\bot (S_1) \otimes S_2)}$.
  Since $p_n^\bot$ converges to~$0$ uniformly on compact subsets, it follows
  that $\|p_n^\bot (S_1)\|_{\VS} \to 0$ for $n \to \infty$.  Consequently,
  $\|p_n^\bot S\|_1 \to 0$ as asserted.
  
  Let $S \subset \Comp(\VS)$ be precompact.  Let $B_{\VS} \subset \VS$ and
  $B_{\Comp} \subset \Comp(\VS)$ be the unit balls.  We claim that $S(B_{\VS})
  \subset \VS$ is precompact.  Given any $\epsilon>0$, we have $S \subset F +
  \epsilon B_{\Comp}$ for a suitable finite set $F \subset \Comp(\VS)$.  By
  definition, a compact operator maps~$B_{\VS}$ to a precompact set.  Thus
  $F(B_{\VS}) \subset F' + \epsilon B_{\VS}$ for a suitable finite set $F'
  \subset \VS$.  Furthermore, $\epsilon B_{\Comp}$ maps~$B_{\VS}$ into
  $\epsilon B_{\VS}$.  Consequently, $S(B_{\VS}) \subset F' + 2\epsilon
  B_{\VS}$.  Since this holds for all $\epsilon>0$, $S(B_{\VS})$ is
  precompact.
  
  The sequence of operators~$(p_n^\bot)$ converges to~$0$ uniformly on the
  precompact subset $S(B_{\VS})$.  Hence $\| p_n^\bot \circ S(B_{\VS})
  \|_{\VS} \to 0$ for $n \to \infty$.  Thus $\| p_n^\bot \circ S\|_\infty \to
  0$ for $n\to \infty$.
\end{proof}

\subsection{Construction of the Chern-Connes character}
\label{sec:Chern_homology}

After these preparations, we can define the Chern-Connes character for
Fredholm modules without summability restrictions.  Let~$\VS$ be a separable
Banach space with Grothendieck's metric approximation property.  Given $F\in
\Endo(\VS)$ with $F^2 = \ID$, we define
$$
\MA[B]_\odd \defeq
\MA[B]_\odd(\VS,F) \defeq
\{ T \in \Endo(\VS) \mid [F,T] \in \Comp(\VS) \}.
$$
It is easy to see that~$\MA[B]_\odd$ is a closed subalgebra of $\Endo(\VS)$
and thus a Banach algebra.  We endow~$\MA[B]_\odd$ with the \emph{precompact}
bornology.  Given $F,\gamma \in \Endo(\VS)$ with $F^2 = \gamma ^2 = \ID$,
$F\gamma + \gamma F = 0$, we define
$$
\MA[B]_\even \defeq
\MA[B]_\even(\VS,F,\gamma) \defeq
\{ T \in \Endo(\VS) \mid [F,T] \in \Comp(\VS),\ [\gamma,T] = 0 \}.
$$
This is a closed subalgebra of $\Endo(\VS)$.  We endow $\MA[B]_\even$ with
the precompact bornology.

By definition, an \emph{odd \Mp{(\VS,F)}module} over a complete bornological
algebra~$\MA$ is a bounded homomorphism $\phi\colon \MA \to
\MA[B]_\odd(\VS,F)$.  An \emph{even \Mp{(\VS,F,\gamma)}module} over a complete
bornological algebra~$\MA$ is a bounded homomorphism $\phi\colon \MA \to
\MA[B]_\even(\VS,F,\gamma)$.

To construct the Chern-Connes character $\chern(\phi) \in \HA^\ast(\MA)$ of an
even \Mp{(\VS,F,\gamma)}Fredholm module $\phi\colon \MA \to \MA[B]_\even$, we
need an even cocycle $\tilde{\tau}_\even \colon X(\Tanil\MA[B]_\even) \to \C$.
Given~$\tilde{\tau}_\even$, we can define
$$
\chern(\phi) \defeq
\phi^\ast([\tilde{\tau}_\even]) =
[\tilde{\tau}_\even]\circ [\phi] =
[\tilde{\tau}_\even \circ X(\Tanil\phi)].
$$
This is automatically a natural transformation from even
\Mp{(\VS,F,\gamma)}modules to $\HA^0(\MA)$.  Furthermore, any natural
transformation must be of this form with $[\tilde{\tau}_\even] =
\chern(\ID[{\MA[B]_\even}])$.

\begin{lemma}
  The algebra $\MA[B]_\odd^0(\VS,F,\gamma)$ is dense in
  $\MA[B]_\odd(\VS,F,\gamma)$.  The algebra $\MA[B]_\even^0(\VS,F,\gamma)$ is
  dense in $\MA[B]_\even(\VS,F,\gamma)$.
\end{lemma}

\begin{proof}
  We only consider the odd case.  The even case is similar.  Let~$\VS_{\pm 1}$
  be the range of the idempotent $\frac{1}{2}(1\pm F)$.  Then $\VS = \VS_1
  \oplus \VS_{-1}$.  Write $x \in \Endo(\VS)$ as a
  \Mp{2\times2}matrix $x = \bigl(\begin{smallmatrix} a & b \\ c & d
  \end{smallmatrix}\bigr)$ with respect to this decomposition.  Computing
  commutators with~$F$ shows that $x\in \MA[B]_\odd$ iff the off-diagonal
  terms $b$ and~$c$ are compact, and $x\in \MA[B]_\odd^0$ iff $b$ and~$c$ are
  nuclear.  Since~$\VS$ has Grothendieck's approximation property, $\Sch(\VS)$
  is dense in $\Comp(\VS)$.  Thus $\MA[B]_\odd^0$ is dense in $\MA[B]_\odd$.
\end{proof}

Hence $\Omega \MA[B]_\even^0$ is dense in $\Omega_\an \MA[B]_\even$.  A
bounded cocycle $\tilde{\tau}_\even \colon X(\Tanil \MA[B]_\even) \to \C$ is
therefore equivalent to a linear map $\Omega \MA[B]_\even^0 \to \C$ that is
bounded on $\opt{S} (dS)^\infty$ for all $S \in \CBS(\MA[B]_\even)$ satisfying
$S \subset \MA[B]_\even^0$.  Therefore, it suffices to work in the algebra
$\MA[B]_\even^0$.

The odd case is quite similar: A natural transformation from odd
\Mp{(\VS,F)}modules over~$\MA$ to $\HA^1(\MA)$ is equivalent to an odd cocycle
$\tilde{\tau}_\odd \colon X(\Tanil\MA[B]_\odd) \to \C$.  It suffices to
construct the restriction of $\tilde{\tau}_\odd$ to $\Omega \MA[B]_\odd^0$ and
verify that it is bounded in the appropriate sense.

In the even case, define
$$
\tilde{h}_{2n-1}(\opt{x_0} dx_1 \dots dx_{2n-1}) \defeq
h_{2n-1}\bigl(
  \opt{p_n^\bot x_0} d(p_n^\bot x_1) \dots d(p_n^\bot x_{2n-1})
\bigr)
$$
for all $x_0,\dots,x_{2n-1} \in \MA[B]_\even^0$, $n \in \N$.  Here and in
the following, $\opt{p_n^\bot x_0}$ denotes $1 \in \Unse{(\MA[B]_\odd^0)}$ if
$\opt{x_0} = 1$.  Let $\tilde{h}_\odd \defeq \sum_{n=1}^\infty
\tilde{h}_{2n-1} \colon \Omega \MA[B]_\odd^0 \to \C$ and
$$
\tilde{\tau}_\even \defeq \tau_0 - \tilde{h}_\odd \circ \partial.
$$

In the odd case, define
$$
\tilde{h}_{2n}(\opt{x_0} dx_1\dots dx_{2n}) \defeq
h_{2n}\bigl(
  \opt{ p_n^\bot x_0} d(p_n^\bot x_1) \dots d(p_n^\bot x_{2n})
\bigr)
$$
for all $x_0,\dots,x_{2n} \in \MA[B]_\odd^0$, $n\in\N$.
Let $\tilde{h}_\even \defeq \sum_{n=1}^\infty \tilde{h}_{2n}\colon \Omega
\MA[B]_\odd^0 \to \C$ and
$$
\tilde{\tau}_\odd \defeq \tau_1 - \tilde{h}_\even \circ \partial.
$$
By construction, $\tilde{\tau}_\even \circ \partial = \tau_0 \circ \partial
= 0$ and $\tilde{\tau}_\odd \circ \partial = \tau_1 \circ \partial = 0$.

\begin{proposition}  \label{pro:chern_exists}
  The functionals $\tilde{h}_\even \colon \Omega \MA[B]_\even^0\to \C$ and
  $\tilde{h}_\odd \colon \Omega \MA[B]_\odd^0\to \C$ extend uniquely to
  bounded linear functionals on $\Omega_\an(\MA[B]_\even^0)$ and
  $\Omega_\an(\MA[B]_\odd^0)$, respectively.  Thus $[\tilde{\tau}_\even] =
  [\tau_0]$ in $\HA^0(\MA[B]_\even^0)$ and $[\tilde{\tau}_\odd] = [\tau_1]$ in
  $\HA^1(\MA[B]_\odd^0)$.
  
  The functionals $\tilde{\tau}_\even \colon \Omega \MA[B]_\even^0 \to \C$ and
  $\tilde{\tau}_\odd \colon \Omega \MA[B]_\odd^0 \to \C$ extend uniquely to
  bounded linear functionals $\tilde{\tau}_\even \colon
  \Omega_\an(\MA[B]_\even) \to \C$ and $\tilde{\tau}_\odd \colon
  \Omega_\an(\MA[B]_\odd) \to \C$, respectively.
  
  These extensions define elements $[\tilde{\tau}_\even] \in
  \HA^0(\MA[B]_\even)$ and $[\tilde{\tau}_\odd] \in \HA^1(\MA[B]_\odd)$,
  respectively, that do not depend on the choice of the cut-off
  sequence~$(p_n)$.
\end{proposition}

\begin{proof}
  We will give the details of the proof in the odd case.  The even case is
  similar but slightly more complicated because the relevant part of the
  boundary~$\partial$ involves $\sum_{j=0}^{n-1} \kappa^{2j}b$.  Since we no
  longer have $\tilde{h}_{2n+1} \circ \kappa = \tilde{h}_{2n+1}$, this gives
  rise to complicated expressions.
  
  Let $S\in \CBS(\MA[B]_\odd^0)$.  Thus $T\defeq \frac{1}{2i} [F,S] \subset
  \Sch(\VS)$ is precompact.  Since~$S$ and~$\{p_n\}$ are equibounded, there is
  a constant $C \in \left[1, \infty\right[$ such that $C \ge \| S \|_\infty
  \|p_n^\bot\|_\infty$ for all $n \in \N$.  Lemma~\ref{lem:cutoff_on_compacts}
  implies that $\epsilon_{1,n} \defeq \| p_n^\bot T \|_1 \to 0$ for $n \to
  \infty$.  Since~$F$ commutes with~$p_n^\bot$, we have $\delta(p_n^\bot x) =
  p_n^\bot \delta x \in p_n^\bot T$ for all $x\in S$.  Thus
  \begin{multline*}
    |\tilde{h}_{2n}(\opt{x_0} dx_1\dots dx_{2n}) | =
    |\tr(\opt{p_n^\bot x_0} (p_n^\bot \delta x_1) \dots
    (p_n^\bot \delta x_{2n}) |
    \\ \le
    \| \opt{p_n^\bot x_0} (p_n^\bot \delta x_1) \dots
    (p_n^\bot \delta x_{2n}) \|_1 \le 
    C \cdot \epsilon_{1,n}^{2n}
  \end{multline*}
  if $x_0,\dots,x_{2n} \in S$.  Since $\epsilon_{1,n} < 1$ for sufficiently
  large~$n$, the set of numbers $C \cdot \epsilon_{1,n}^{2n}$ is bounded.
  Thus~$\tilde{h}_\even$ is bounded on the set $\opt{S}(dS)^\even$.  By the
  universal property of the completion, $\tilde{h}_\even$ extends uniquely to
  a bounded operator $\Omega_\an(\MA[B]_\odd^0) \to \C$.

  The same argument works for $\tilde{h}_\odd$ instead of $\tilde{h}_\even$.

  To estimate the growth of~$\tilde{\tau}_\odd$, we want to use the equation
  $h_\even \circ \partial = \tau_1$ proved in Lemma~\ref{lem:tau_h_odd}.
  Define a linear map $l\colon \Omega \MA[B]_\odd^0 \to \Omega \MA[B]_\odd^0$
  by
  $$
  l(\opt{x_0} dx_1 \dots dx_k) \defeq
  \opt{p_n^\bot x_0} d(p_n^\bot x_1) \dots d(p_n^\bot x_k),
  $$
  for $k=2n$ and $k=2n+1$ (in the even case, it is better to
  take~$p_n^\bot$ for $k = 2n$ and $k = 2n-1$).  By convention, $p_0 \defeq
  0$.  Since $\tilde{h}_{2n} = h_{2n} \circ l$, we get
  \begin{displaymath}
    \tilde{\tau}_\odd =
    \tau_1 - \tilde{h}_\even\circ \partial =
    \tau_1 - h_\even \circ \partial\circ l + h_\even\circ [\partial,l]
    =
    \tau_1 - \tau_1 \circ l + h_\even \circ [b - (1+\kappa)d,l].
  \end{displaymath}
  Evidently, $\kappa \circ
  l \circ d = l \circ \kappa \circ d$ and hence
  \begin{equation}  \label{eq:tau_reformulate}
    \tilde{\tau}_\odd =
    \tau_1 - \tau_1 \circ l + h_\even \circ [b,l] - 2h_\even \circ [d,l].
  \end{equation}
  
  In the even case, we find similarly that $\tilde{\tau}_\even = h_\even \circ
  [\partial,l] + \tau_0 - \tau_0 \circ l$ with $\partial = B -
  \sum_{j=0}^{n-1} \kappa^{2j} b$ on $\Omega^{2n} \MA[B]_\even^0$.  Since~$l$
  does not commute with~$\kappa$, this gives rise to complicated sums.  These
  can be handled.  It is, however, more convenient to replace~$\partial$ by
  the boundary~$\delta$ given by $\delta = B - nb$ on
  $\Omega^{2n}\MA[B]_\even^0$ and $\delta = -B/(n+1) + b$ on
  $\Omega^{2n+1}\MA[B]_\even^0$.  The boundaries $\partial$ and~$\delta$ give
  rise to chain homotopic complexes.  Hence it is admissible to replace
  $\partial$ by~$\delta$.
  
  The summands $\tau_1 - \tau_1 \circ l$ or $\tau_0 - \tau_0 \circ l$ are
  harmless.  In fact, $\tau_1 - \tau_1 \circ l$ is zero because we declared
  $p_0 = 0$.  In the even case $\tau_0 - \tau_0 \circ l$ is not zero unless
  $p_1 = 0$.  Even if $p_1 \neq 0$, it is easy to prove that $\tau_0 - \tau_0
  \circ l$ is bounded.  Hence we concentrate on the other two summands
  in~\eqref{eq:tau_reformulate}.  Let $S \in \CBS(\MA[B]_\odd)$ and $S \subset
  \MA[B]_\odd^0$.  It suffices to verify that $\|\psi \circ [d,l]\|_1$ and
  $\|\psi \circ [b,l]\|_1$ remain bounded on $\opt{S}(dS)^\infty$ because
  $|\tr(X)| \le \|X\|_1$.

  The set $T \defeq \frac{1}{2i} [F,S] \subset \Comp(\VS)$ is precompact.  Let
  $$
  \epsilon_{\infty,n} \defeq \sup_{m\ge n} \| p_m^\bot T \|_\infty.
  $$
  The sequence $(\epsilon_{\infty,n})$ is decreasing by definition and
  converges towards~$0$ for $n \to \infty$ by
  Lemma~\ref{lem:cutoff_on_compacts}.  Fix $n \in \N$.  We compute $\psi \circ
  [d,l]$ on $\opt{S}(dS)^{2n+1}$.  Since~$l$ maps closed forms again to closed
  forms, $[d,l]$ annihilates closed forms $dx_1 \dots dx_{2n+1}$.  We compute
  \begin{multline*}
    \psi\circ [d,l](x_0 dx_1\dots dx_{2n+1}) =
    \psi\bigl(d(p_n^\bot x_0) \dots d(p_n^\bot x_{2n+1}) -
    d(p_{n+1}^\bot x_0) \dots d(p_{n+1}^\bot x_{2n+1}) \bigr)
    \\ =
    \sum_{j=0}^{2n+1} \delta(p_{n+1}^\bot x_0) \dots
    \delta(p_{n+1}^\bot x_{j-1})
    \delta\bigl((p_n^\bot - p_{n+1}^\bot)x_j\bigr)
    \delta (p_n^\bot x_{j+1}) \dots \delta(p_n^\bot x_{2n+1})
    \\ =
    \sum_{j=0}^{2n+1}
    (p_{n+1}^\bot \delta x_0) \dots (p_{n+1}^\bot \delta x_{j-1})
    (p_{n+1} - p_n) \delta x_j
    (p_n^\bot \delta x_{j+1}) \dots (p_n^\bot \delta x_{2n+1}).
  \end{multline*}
  We used that~$p_n$ commutes with~$F$, so that $\delta(p_n^\bot x) = p_n^\bot
  \delta x$.  If $x_0, \dots, x_{2n+1} \in S$, we can estimate
  \begin{multline*}
    \| \psi\circ [d,l](x_0 dx_1\dots dx_{2n+1}) \|_1
    \\ \le
    \sum_{j=0}^{2n+1}
    \| (p_{n+1}^\bot \delta x_0) \|_\infty \cdots
    \| (p_{n+1}^\bot \delta x_{j-1}) \|_\infty \cdot
    \| p_{n+1} - p_n\|_1 \cdot \|\delta x_j\|_\infty
    \| p_n^\bot \delta x_{j+1}\|_\infty \cdots
    \| p_n^\bot \delta x_{2n+1} \|_\infty
    \\ \le
    2(n+1) \epsilon_{\infty,n}^{2n+1} \cdot
    \bigl( \|p_{n+1}\|_1 + \|p_n\|_1 \bigr)\cdot  \|T\|_\infty.
  \end{multline*}
  Since~$(p_n)$ is a cut-off sequence $\|p_n\|_1$ has at most exponential
  growth.  Since $\epsilon_{\infty,n} \to 0$, the term
  $\epsilon_{\infty,n}^{2n+1}$ is $O(\epsilon^n)$ for any $\epsilon>0$.  Hence
  it decays fast enough to compensate the at most exponential growth of
  $2(n+1) \bigl( \|p_n\|_1 + \|p_{n+1}\|_1 \bigr)$.  Therefore, $\| \psi\circ
  [d,l] (\opt{S}(dS)^{2n+1}) \|_1 = O(\epsilon^n)$ for all $\epsilon>0$.  In
  particular, $\psi \circ [d,l]$ remains bounded on $\opt{S}(dS)^\odd$.

  Next, we compute $\psi\circ [b,l]$ on $S(dS)^{2n+1}$:
  \begin{multline*}
    \psi\circ [b,l](x_0 dx_1 \dots dx_{2n+1})
    \\ =
    \sum_{j=0}^{2n} (-1)^j \psi\Bigl(
    p_n^\bot x_0 d(p_n^\bot x_1) \dots d(p_n^\bot x_{j-1})
    d(p_n^\bot x_j p_n^\bot x_{j+1})
    d(p_n^\bot x_{j+2}) \dots d(p_n^\bot x_{2n+1})
    \\ -
    p_n^\bot x_0 d(p_n^\bot x_1) \dots d(p_n^\bot x_{j-1})
    d(p_n^\bot x_j x_{j+1}) d(p_n^\bot x_{j+2}) \dots d(p_n^\bot x_{2n+1})
    \Bigr)
    \\ +
    (-1)^{2n+1} \psi\bigl( p_n^\bot x_{2n+1} (p_n^\bot x_0 - x_0)
    d(p_n^\bot x_1) \dots d(p_n^\bot x_{2n}) \bigr)
    \\ =
    - p_n^\bot x_0 p_n x_1
    (p_n^\bot \delta x_2) \cdots (p_n^\bot \delta x_{2n+1}) +
    p_n^\bot x_{2n+1} p_n x_0 \cdot
    (p_n^\bot \delta x_1) \cdots (p_n^\bot \delta x_{2n})
    \\
    + \sum_{j=1}^{2n} (-1)^{j+1} p_n^\bot x_0
    (p_n^\bot \delta x_1) \dots
    (p_n^\bot \delta x_{j-1}) \cdot
    \bigl(p_n^\bot \delta(x_j p_n x_{j+1}) \bigr) \cdot
    (p_n^\bot \delta x_{j+2}) \dots
    (p_n^\bot \delta x_{2n+1})
  \end{multline*}
  If $x_0,\dots,x_{2n+1} \in S$, we can estimate
  \begin{multline*}
    \|\psi\circ [b,l](x_0 dx_1 \dots dx_{2n+1})\| \le
    2n \|p_n^\bot\|_\infty^2 \|S\|_\infty \|p_n^\bot T\|_\infty^{2n-1}
    \| \delta(Sp_n S) \|_1
    \\ +
    2 \|p_n^\bot\|_\infty \|S\|_\infty^2 \|p_n\|_1 \|p_n^\bot T\|_\infty^{2n}
    =
    O(n \|p_n\|_1 \cdot \epsilon_{\infty,n}^{2n-1} ).
  \end{multline*}
  This is again $O(\epsilon^n)$ for any $\epsilon > 0$ for $n \to \infty$
  because the super-exponential decay of $\epsilon_{\infty,n}^{2n-1}$
  overcompensates the at most exponential growth of $n \|p_n\|_1$.  The
  estimation on $(dS)^\odd$ is similar.  The only difference is that the terms
  in the sum with $j=0$ and $j=2n+1$ drop out.

  Consequently, $\tilde{\tau}_\odd$ remains bounded on sets of the form
  $\opt{S}(dS)^\odd$ with $S\in \CBS(\MA[B]_\odd)$ and $S\subset
  \MA[B]_\odd^0$.  It remains to verify that the choice of the cut-off
  sequence does not matter.

  Assume that~$\tilde{h}'_\even$ is produced by another cut-off
  sequence~$(p_n')$.  Then $\tilde{h}'_\even - \tilde{h}_\even$ extends to a
  bounded linear functional on $\Omega_\an \MA[B]_\odd$.  This is verified in
  the same way as the boundedness of $h_\even \circ [d,l]$.  Therefore, $\tau_1
  - \tilde{h}_\even \circ \partial$ and $\tau_1 - \tilde{h}'_\even \circ
  \partial$ are cohomologous, the difference being the coboundary
  $(\tilde{h}_\even - \tilde{h}'_\even) \circ \partial$.
\end{proof}

We now consider the homotopy invariance and additivity properties of the
Chern-Connes character.  The class $\chern(\phi) = [\phi] \circ
[\tilde{\tau}]$ depends only on the class of the homomorphism~$\phi$ in
$\HA^0(\MA; \MA[B]_\even)$ or $\HA^0(\MA; \MA[B]_\odd)$.
Theorem~\ref{the:HA_homotopy} implies that $\chern(\phi_0) = \chern(\phi_1)$
if $\phi_0, \phi_1 \colon \MA \to \MA[B]_{\dots}$ are AC-homotopic.  These
homotopies leave the operators~$F$ (and~$\gamma$) fixed.  Alternatively, we
can leave the homomorphism $\phi \colon \MA \to \Endo(\VS)$ (and~$\gamma$)
fixed and vary the operator~$F$.  We claim that the Chern-Connes character is
also invariant under this kind of ``operator homotopy''.

\begin{lemma}
  Let $t\mapsto F_t \in \Endo(\VS)$ be a continuous path of operators
  satisfying $F_t^2 = 1$ and $[F_t, \phi(\MA)] \subset \Comp(\VS)$ for all
  $j\in [0,1]$.  Then there is a \emph{smooth} path of invertibles $t\mapsto
  U_t$ such that $U_0 = \ID$, $F_1 = U_1 F_0 U_1^{-1}$, and
  $[U_t,\phi(\MA)]\subset \Comp(\VS)$ for all $t\in [0,1]$.  If $\gamma \in
  \Endo(\VS)$ is a grading operator that anticommutes with~$F_t$ for all $t\in
  [0,1]$, then we can achieve that $U_t$ commutes with~$\gamma$ for all~$t$.
\end{lemma}

\begin{proof}
  Let $\MA[C] \defeq \{ x\in \Endo(\VS) \mid [x,\phi(\MA)]\subset \Comp(\VS)
  \}$.  This is a closed unital subalgebra of $\Endo(\VS)$ and thus a unital
  Banach algebra.  We can replace the path~$F_t$ by the path of idempotents
  $p_t = \frac{1}{2}(1+F_t)$ in~$\MA[C]$.  Standard techniques of
  \Mpn{K}theory~\cite{blackadar86:ktheory} imply that there is a continuous
  path of invertibles $t\mapsto U_t'$ in~$\MA[C]$ such that $U_0' = \ID$ and
  $U_1'F_0 (U_1')^{-1} = F_1$.  It is well-known that $\CINF([0,1]; \MA[C])$
  is dense and closed under holomorphic functional calculus in $\NBC([0,1];
  \MA[C])$.  Thus we can replace the path~$(U_t')$ by a nearby smooth path of
  invertible elements~$(U_t)$.  In addition, we can arrange that $U_0 = \ID$
  and $U_1 = U_1'$.
\end{proof}

Choose a cut-off sequence~$(p_n)$ for the operator~$F_0$ (recall that one of
the conditions for the~$(p_n)$ was to commute with~$F$).  Then the sequence
$(U_1 p_n U_1^{-1})$ is a cut-off sequence for the operator~$F_1$.  The
associated functionals $\tilde{\tau}_0$ and~$\tilde{\tau}_1$ are related by
$$
\tilde{\tau}_1(\opt{x_0} dx_1\dots dx_n) =
\tilde{\tau}_0(\opt{U_1^{-1}x_0U_1} d(U_1^{-1}x_1U_1) \dots d(U_1^{-1}x_nU_1)).
$$
Therefore, the Chern-Connes character of $(\phi,\VS,F_1)$ is equal to the
character of $(\phi_1,\VS,F_0)$ with $\phi_1(\ma) = U_1^{-1} \phi(\ma) U_1$.
Since~$\phi_1$ is smoothly homotopic to~$\phi$ via $\phi_t(\ma) = U_t^{-1}
\phi(\ma) U_t$, the Chern-Connes characters of $(\phi,\VS,F_0)$ and
$(\phi,\VS,F_1)$ are equal.

Consequently, operator homotopic Fredholm modules have the same character.

For $j=1,2$, let $\phi_j\colon \MA \to \Endo(\VS_j)$, $F_j$, and
possibly~$\gamma_j$ be the data describing Fredholm modules.  The \emph{direct
  sum} of these two Fredholm modules is described by $\phi \defeq \phi_1
\oplus \phi_2 \colon \MA \to \Endo(\VS_1 \oplus \VS_2)$, $F \defeq F_1 \oplus
F_2$, and possibly $\gamma \defeq \gamma_1 \oplus \gamma_2$.  It is easy to
verify that $\chern(x_1 \oplus x_2) = \chern(x_1) + \chern(x_2)$.  That is,
the Chern-Connes character is additive.  A Fredholm module is
\emph{degenerate} if $[\phi(\MA), F] = 0$.  Evidently, the Chern-Connes
character of a degenerate Fredholm module is~$0$.

I am not certain what equivalence relation to put on the Fredholm modules over
Banach spaces considered above.  Therefore, we now restrict attention to
\Cstar{}algebras and Kasparov's \Mpn{K}homology based on Fredholm modules over
Hilbert space.

\subsection{\Mpn{K}homology and the index pairing}
\label{sec:Chern_add_homotopy}

The \Mpn{K}homology of \Cstar{}algebras is defined as follows.  Let~$\MA$ be a
separable \Cstar{}algebra and let~$\Hils$ be the separable infinite
dimensional Hilbert space.  An \emph{odd Fredholm module} for~$\MA$ consists
of an operator $F\in \Endo(\Hils)$ and a \Mpn{\ast}homomorphism $\phi\colon
\MA \to \Endo(\Hils)$ such that $F = F^\ast$, $F^2 = 1$, and $[F,\phi(\MA)]
\subset \Comp(\Hils)$.  The \Mpn{K}homology $K^1(\MA)$ of~$\MA$ is defined as
the set of equivalence classes of odd Fredholm modules for a certain
equivalence relation.

Let $(\phi_0,F_0)$ and $(\phi_1,F_1)$ describe odd Fredholm modules.  These
define the same element of \Mpn{K}homology iff there are degenerate odd
Fredholm modules $(\phi_j',F_j')$, $j=0,1$, a \Mpn{\ast}homomorphism
$\psi\colon \MA \to \Endo(\Hils)$, and a norm continuous path $t\mapsto F_t$,
$t\in [0,1]$, such that $(\psi,F_t)$ is an odd Fredholm module for all $t\in
[0,1]$ and $(\psi,F_j)$ is unitarily equivalent to $(\phi_j \oplus \phi_j',
F_j \oplus F_j')$ for $j=0,1$.  This is the right equivalence relation because
the equivalence relation generated by addition of degenerate modules and
operator homotopy is the same as the equivalence relation of
homotopy~\cite{kasparov80:K}.

The definition of the even \Mpn{K}homology $K^0(\MA)$ is essentially the same,
taking even Fredholm modules instead of odd ones.

Let $K_\ast(\MA)$ be the topological \Mpn{K}theory of~$\MA$.  The index pairing
$\mathrm{Ind}\colon K^\ast(\MA) \times K_\ast(\MA) \to \Z$, $\ast=0,1$, is
defined as follows.

In the even case, represent the elements of $K^0(\MA)$ and $K_0(\MA)$ by an
even Fredholm module $(\phi,F,\gamma)$ and an element $a\in \Mat[n](\MA)$ such
that $1+a \in \Unse{(\Mat[n](\MA))}$ is an idempotent.  Let $\phi_n \defeq
\phi\otimes\ID\colon \MA\otimes \Mat[n] \to \Endo(\Hils)\otimes \Mat[n] \cong
\Endo(\Hils \otimes \C^{\,n})$, let $F_n \defeq F\otimes \ID[\C^n] \in
\Endo(\Hils \otimes \C^{\,n})$, and $\gamma_n = \gamma\otimes \ID[\C^n]$.
Then $(\phi_n,F_n,\gamma_n)$ is again an even Fredholm module.  Let
$\Hils[K]^{\pm}$ be the even and odd part of the range of the idempotent
$1+\phi_n(a)$.  The compression of~$F$ to an operator from~$\Hils[K]^{+}$
to~$\Hils[K]^{-}$ is a Fredholm operator.  We let $\operatorname{Ind}\bigl(
[(\phi,F,\gamma)], [a] \bigr)$ be the index of this Fredholm operator.  The
index of a Fredholm operator~$X$ is the difference
$$
\operatorname{Ind}(X) \defeq \dim \Ker X - \dim \Coker X.
$$

In the odd case, represent the elements of $K^1(\MA)$ and $K_1(\MA)$ by an odd
Fredholm module $(\phi,F)$ and an element $a\in \Mat[n](\MA)$ such that $1+a
\in \Unse{(\Mat[n](\MA))}$ is invertible.  Define $\phi_n$ and~$F_n$ as in the
even case to get an odd Fredholm module $(\phi_n,F_n)$.  Let $\Hils^n_+
\subset \Hils \otimes \C^{\,n}$ be the range of the projection $p \defeq
\frac{1}{2}(1+F_n)$.  The compression $p \bigl(1+ \phi_n(a)\bigr) p \colon
\Hils^n_+ \to \Hils^n_+$ is a Fredholm operator.  We let $\operatorname{Ind}
\bigl( [(\phi,F)], [a] \bigr)$ be the index of this Fredholm operator.

Let us give an alternative description of these index pairings.  Let
$(\phi,F,\gamma)$ be an even Fredholm module.  Let~$\Hils_\pm$ be the range of
the projections $\frac{1}{2}(1\pm \gamma)$.  Since~$F$ anticommutes
with~$\gamma$ and $F^2 = 1$, the restriction of~$F$ to~$\Hils_+$ is an
isomorphism onto~$\Hils_-$ and the restriction of~$F$ to~$\Hils_-$ is the
inverse of this isomorphism.  We use these isomorphisms to identify $\Hils_+$
with~$\Hils_-$.  We assume that~$\Hils_+$ is infinite dimensional.  This can
be achieved by adding a degenerate module.

The isomorphism $\Endo(\Hils) \cong \Mat[2]\bigl(\Endo(\Hils_+)\bigr)$ maps
the algebra $\MA[B]_\even = \MA[B]_\even(\Hils,F,\gamma) \subset
\Endo(\Hils)$ to the algebra of diagonal matrices $\bigl(\begin{smallmatrix} x
  & 0 \\ 0 & y
\end{smallmatrix}\bigr)$ with $x - y \in \Comp(\Hils_+)$.  Therefore, we have
a split extension
$$
\Comp(\Hils_+) \injto \MA[B]_\even \prto \Endo(\Hils_+).
$$
Using the long exact sequence in \Mpn{K}theory, we compute that
$K_0(\MA[B]_\even) \cong \Z$, $K_1(\MA[B]_\even) = 0$.  A generator of
$K_0(\MA[B]_\even)$ is the idempotent $E_\even \defeq
\bigl(\begin{smallmatrix} p & 0 \\ 0 & 0 \end{smallmatrix}\bigr)$ with a rank
one projection $p\in \Comp(\Hils_+)$.  We choose an isomorphism $j_\even\colon
K_0(\MA[B]_\even) \to \Z$ that sends $[E_\even] \mapsto 1$.

If $(\phi,F,\gamma)$ is an even Fredholm module, then~$\phi$ is a
\Mpn{\ast}homomorphism $\MA \to \MA[B]_\even$.  It is easy to check that the
index pairing with $(\phi,F,\gamma)$ is equal to the composition
$$
K_0(\MA) \overset{\phi_\ast}{\longrightarrow}
K_0(\MA[B]_\even) \overset{j_\even}{\longrightarrow}
\Z.
$$

In the odd case, we can proceed similarly.  Let $(\phi,F)$ be an odd Fredholm
module.  Let~$\Hils_\pm$ be the range of the projection $\frac{1}{2}(1\pm F)$.
We assume for simplicity that both $\Hils_+$ and~$\Hils_-$ are infinite
dimensional.  The algebra $\MA[B]_\odd = \MA[B]_\odd(\Hils,F)$ can be
described in this decomposition as the set of all \Mp{2\times 2}matrices
$\bigl(\begin{smallmatrix} x & y \\ z & w \end{smallmatrix}\bigr)$ with $y$
and~$z$ compact.  Writing $\Endo/\Comp(\Hils)$ for the Calkin algebra
on~$\Hils$, we therefore get a short exact sequence
$$
\Comp(\Hils) \injto \MA[B]_\odd \prto
\Endo/\Comp (\Hils_+) \oplus \Endo/\Comp (\Hils_-).
$$
The inclusion $\Comp(\Hils) \to \MA[B]_\odd$ induces the zero map on
\Mpn{K}theory because there is an isometry in $\MA[B]_\odd$ with
\Mpn{1}dimensional cokernel.  The long exact sequence for the above extension
shows therefore that $K_0(\MA[B]_\odd) = 0$, $K_1(\MA[B]_\odd) \cong \Z$.  To
write down a generator of $K_1(\MA[B]_\odd)$, let $S_\pm \in \Endo(\Hils_\pm)$
be unilateral shifts (of multiplicity~$1$) and let $p\colon \Hils_- \to
\Hils_+$ be a partial isometry of rank one that maps the cokernel of~$S_-$ to
the cokernel of~$S_+$.  Then
$$
E_\odd \defeq 
\begin{pmatrix}
  S_+^\ast & 0 \\
  p & S_-
\end{pmatrix}
$$
is a unitary element of $\MA[B]_\odd$ that generates $K_1(\MA[B]_\odd)$.
Let $j_\odd\colon K_1(\MA[B]_\odd) \to \Z$ be the isomorphism sending~$E_\odd$
to~$1$.  The compression of~$E_\odd$ to~$\Hils_+$ is the
co-isometry~$S_+^\ast$ and thus has index~$1$.  Hence the index map
$K_1(\MA) \to \Z$ associated to the Fredholm module $(\phi,F)$ is the
composition
$$
K_1(\MA) \overset{\phi_\ast}{\longrightarrow}
K_1(\MA[B]_\odd) \overset{j_\odd}{\longrightarrow}
\Z.
$$

\begin{theorem}  \label{the:character_Cstar}
  Let~$\MA$ be a separable \Cstar{}algebra.  The Chern-Connes character
  descends to a natural additive map $\chern\colon K^\ast(\MA) \to
  \HA^\ast(\MA,\COMP)$ such that $\5{\chern(x)}{\chern(y)} = \5{x}{y} \in \Z$
  for all $x\in K_\ast(\MA)$, $y\in K^\ast(\MA)$, $\ast=0,1$.
\end{theorem}

\begin{proof}
  An element of $K^0(\MA)$ is represented by an even
  \Mp{(\Hils,F,\gamma)}Fredholm module.  The Chern-Connes character of such a
  module was defined in the previous section.  We have already observed that
  it is invariant under operator homotopies, additive, and vanishes on
  degenerate modules.  Hence it descends to a map $\chern\colon K^0(\MA) \to
  \HA^0(\MA)$.  The same applies in the odd case.
  
  It remains to compare the pairings $\5{\chern(x)}{\chern(y)}$ and $\5{x}{y}$
  for $x\in K_\ast(\MA)$, $y\in K^\ast(\MA)$, $\ast=0,1$.  We consider the
  even case.  Represent~$y$ by an even Fredholm module $\phi\colon \MA \to
  \MA[B]_\even(\Hils,F,\gamma)$.  The Chern-Connes character is defined by
  $\chern(y) = \phi^\ast([\tilde{\tau}_\even])$.  Since the Chern-Connes
  character in \Mp{K}theory is natural, we have
  $$
  \5{\chern(x)}{\chern(y)} =
  \5{\phi_\ast\chern(x)}{\tilde{\tau}_\even} =
  \5{\chern(\phi_\ast x)}{\tilde{\tau}_\even}.
  $$
  Comparing with the index pairing, we see that it suffices to show that
  the functional $K_0(\MA[B]_\even) \to \C$ sending $x\mapsto
  \5{\chern(x)}{\tilde{\tau}_\even}$ is equal to~$j_\even$.  Since
  $K_0(\MA[B]_\even)$ is generated by~$[E_\even]$, it suffices to compute
  $\5{\chern(E_\even)}{\tilde{\tau}_\even} = 1$.  Since $E_\even \in
  \MA[B]_\even^0$, we can replace $\tilde{\tau}_\even$ by~$\tau_0$.  We have
  already checked $\5{\chern(E_\even)}{\tau_0} = 1$ in
  Section~\ref{sec:Chern_finite}.

  Similarly, in the odd case the compatibility with the index pairing can be
  reduced to the computation $\5{\chern(E_\odd)}{\tau_1} = 1$ that was done in
  Section~\ref{sec:Chern_finite}.
\end{proof}

\section{Tensoring algebras and the exterior product}
\label{sec:tensoring}

Recall that a complete bornological algebra~$\MA[C]$ is tensoring iff $\MA[C]
\hot \MA[N]$ is a-nilpotent for all a-nilpotent algebras~$\MA[N]$.  In
Lemma~\ref{lem:tensor_a_nilpotent}, we found a sufficient condition
for~$\MA[C]$ to be tensoring involving the bornology of~$\MA[C]$.  This
condition is necessary:

\begin{proposition}  \label{pro:tensoringConditions}
  Let~$\MA[C]$ be a complete bornological algebra.  The following
  conditions are all equivalent to~$\MA[C]$ being tensoring:
  \begin{enumerate}[(i)]
  \item $\MA[C] \hot \Janil\C$ is a-nilpotent;
  \item there is $\MA[B] \neq \{0\}$ such that $\MA[C] \hot \Janil\MA[B]$ is
    a-nilpotent and the natural maps $\MA[C] \otimes \MA[B]^{\otimes j} \to
    \MA[C] \hot \MA[B]^{\hot j}$ are injective for all $j \in \N$;
  \item for all $S \in \CBS(\MA[C])$ there is $\lambda \in \R$ such that
    $(\lambda S)^\infty \in \CBS(\MA[C])$;
  \item $\MA[C]$ is an inductive limit of Banach algebras.
  \end{enumerate}
\end{proposition}

That the natural maps $\MA[C] \otimes \MA[B]^{\otimes j} \to \MA[C] \hot
\MA[B]^{\hot j}$ be injective is a technical assumption to exclude
pathological analysis.  It excludes, for instance, the case $\MA[B]^{\hot n} =
0$ for some $n \in \N$.

\begin{proof}
  Of course, (i) and~(ii) hold if~$\MA[C]$ is tensoring and (i)
  implies~(ii).  Furthermore, (iii) implies that~$\MA[C]$ is tensoring by
  Lemma~\ref{lem:tensor_a_nilpotent}.  Therefore, we will be finished if we
  show that~(ii) implies~(iii) and that (iii) and~(iv) are equivalent.
  
  Since $\MA[B] \neq \{0\}$, there is $x \in \MA[B] \setminus \{0\}$.  Let $S
  \in \CBS(\MA[C])$, then $S \otimes \{ dx dx \} \in \CBS(\MA[C] \hot
  \Janil\MA[B])$.  Since $\MA[C] \hot \Janil\MA[B]$ is a-nilpotent, the set
  $(S \otimes \{ dxdx \})^\infty = \bigcup S^n \otimes (dx)^{2n} \subset
  \MA[C] \hot \Janil\MA[B]$ is small.  That is, it is contained in
  $\coco{(S_{\MA[C]} \otimes \opt{S_{\MA[B]}} (dS_{\MA[B]})^\even)}$ for
  suitable small disks $S_{\MA[C]} \in \CBS_c(\MA[C])$ and $S_{\MA[B]} \in
  \CBS_c(\MA[B])$.  Since the map $\MA[C] \otimes \MA[B]^{\otimes j} \to
  \MA[C] \otimes \MA[B]^{\hot j}$ is injective, this implies that
  $\|y\|_{S_{\MA[C]}} \cdot \| x \|_{S_{\MA[B]}}^{2j} \le 1$ for all $y \in
  S^j$ (Use that the projective tensor product norm on a Banach space
  satisfies $\| a \otimes b \| = \|a \| \cdot \|b \|$.)  Let $\lambda \defeq
  \|x\|^2_{S_{\MA[B]}}$, then $\lambda^j S^j \subset S_{\MA[C]}$.  That is,
  $(\lambda S)^\infty \subset S_{\MA[C]}$.  Thus~(ii) implies~(iii).
  
  Let $(\MA[B]_i)_{i \in I}$ be an inductive system of Banach algebras, let
  $\MA[C] \defeq \indlim {(\MA[B]_i)_{i \in I}}$.  If $S \in \CBS(\MA[C])$,
  then~$S$ is contained in the image of a bounded subset of~$\MA[B]_i$ for
  some $i \in I$.  Hence $\lambda S$ is contained in the image of the unit
  ball of~$\MA[B]_i$ for suitable $\lambda > 0$.  It follows that $(\lambda
  S)^\infty$ is small.  Thus~$\MA[C]$ satisfies~(iii).  Conversely,
  if~$\MA[C]$ satisfies~(iii), then the sets $S \in \CBS_c(\MA[C])$ with $S^2
  \subset \lambda S$ for some $\lambda \in \R$ are cofinal.  If $S^2 \subset
  \lambda S$, then~$\MA[C]_S$ is a Banach algebra with respect to some norm.
  We have $\MA[C] \cong \indlim_{S \in \CBS_c(\MA[C])} \MA[C]_S$ by
  Theorem~\ref{the:bornologies_inductive}.  We can replace $\CBS_c(\MA[C])$ by
  any cofinal subset without changing the limit.  Thus we can write~$\MA[C]$
  as an inductive limit of Banach algebras.
\end{proof}

\begin{theorem}  \label{the:tensoring_Frechet}
  Let~$\MA[C]$ be a \Frechet algebra.  Then the following are equivalent:
  \begin{enumerate}[(i)]%

  \item for all precompact subsets $S \subset \MA[C]$, there is $\lambda > 0$
    such that $(\lambda S)^\infty$ is precompact;
    
  \item for each null-sequence $(x_n)_{n \in \N}$ in~$\MA[C]$, there is $N \in
    \N$ such that $\{x_n \mid n \ge N \}^\infty$ is precompact;
    
  \item there is a neighborhood~$U$ of the origin in~$\MA[C]$ such
    that~$S^\infty$ is precompact for all precompact sets $S \subset U$.

  \end{enumerate}
\end{theorem}

This theorem is due to Puschnigg~\cite{puschnigg96:asymptotic}.  We prove it
in appendix~\ref{app:proof_admissible}.

\begin{lemma}
  Completions, subalgebras, quotients, direct sums, inductive limits, tensor
  products, and allowable extensions of tensoring algebras are again
  tensoring.
\end{lemma}

\begin{proof}
  We only prove the assertion about extensions, the other assertions are
  trivial.
  
  Let $\MA[K] \injto \MA[E] \prto \MA[Q]$ be an extension, possibly without a
  bounded linear section.  Since $\Janil\C$ is nuclear, tensoring with
  $\Janil\C$ preserves extensions.  Thus $\MA[K] \hot \Janil\C \injto \MA[E]
  \hot \Janil\C \prto \MA[Q] \hot \Janil\C$ is again an extension.  (If we had
  bounded linear sections, we would not need the nuclearity of
  $\Janil\MA[C]$.)  If $\MA[K]$ and~$\MA[E]$ are tensoring, then $\MA[K] \hot
  \Janil\C$ and $\MA[Q] \hot \Janil\C$ are a-nilpotent.  The Extension
  Axiom~\ref{axiom:extension} implies that $\MA[E] \hot \Janil\C$ is
  a-nilpotent.  Thus~$\MA[E]$ is tensoring by
  Proposition~\ref{pro:tensoringConditions}.
\end{proof}

\begin{corollary}  \label{cor:tensoring_Tanil}
  $\MA$ is tensoring iff $\Tanil\MA$ is tensoring.
\end{corollary}

\begin{proof}
  If~$\MA$ is tensoring, then so is $\Tanil\MA$ as an extension of~$\MA$ by
  the a-nilpotent and thus tensoring algebra $\Janil\MA$.  Conversely, if
  $\Tanil\MA$ is tensoring, so is its quotient~$\MA$.
\end{proof}

\begin{example}
  The algebra $\prod_{\N} \C$, the product of countably many copies of~$\C$, is
  not tensoring.  Hence a product of countably many tensoring algebras
  may fail to be tensoring.
  
  Consider $\vs \defeq (1, 2, 3, 4, \dots) \in \prod_{\N} \C$.  The
  corresponding one-point set $\{\vs\}$ is small, but $\{c \cdot \vs\}^\infty$
  is not small for any $c>0$.  To see this, choose $m \in \N$ such that $m
  \cdot c>1$.  The $m$th component of $c^n \vs^n$ is $(c \cdot m)^n$.  Since
  the set $(c \cdot m)^n$, $n \in \N$, is unbounded in~$\C$, the set $\{c
  \cdot \vs\}^\infty$ is not small.
\end{example}

Let $\MA[C]_1$ and~$\MA[C]_2$ be tensoring.  Thus $\ID[{\MA[C]_1}] \hot
\sigma_{\MA[C]_2}$ and $\sigma_{\MA[C]_1} \hot \ID[{\Tanil\MA[C]_2}]$ are
\lanilcurs by the Tensor Lemma~\ref{lem:tensor}.  Hence their composition
$\sigma_{\MA[C]_1} \hot \sigma_{\MA[C]_2}$ is a \lanilcur by the Extension
Lemma~\ref{lem:extension}.  However, if $\MA[C]_1$ or~$\MA[C]_2$ is not
tensoring, then $\sigma_{\MA[C]_1} \hot \sigma_{\MA[C]_2}$ is (usually) not a
\lanilcur.  There is no reason to expect $\sigma_1 \hot \sigma_2$ to define an
element of $\HA^\ast(\MA[C]_1 \hot \MA[C]_2, \Tanil\MA[C]_1 \hot
\Tanil\MA[C]_2)$.  Consequently, there is no exterior product in analytic
cyclic cohomology for arbitrary algebras.  However, counterexamples may be
hard to find because the most obvious non-tensoring algebras like $\C[t]$ are
still homotopy equivalent to tensoring algebras.  For tensoring algebras, the
exterior product exists:

\begin{theorem}  \label{the:exterior_product_tensoring}
  Let $\MA[C]_1$ and~$\MA[C]_2$ be tensoring algebras.  Then $X\bigl(
  \Tanil(\MA[C]_1 \hot \MA[C]_2) \bigr)$ is naturally chain homotopic to
  $X\bigl( \Tanil\MA[C]_1 \bigr) \hot X\bigl( \Tanil\MA[C]_2)$.
  
  Thus we obtain an exterior product $\HA^\ast(\MA[C]_1; \MA[C]_2) \times
  \HA^\ast(\MA[C]_3; \MA[C]_4) \to \HA^\ast(\MA[C]_1 \hot \MA[C]_3; \MA[C]_2
  \hot \MA[C]_4)$ if the algebras $\MA[C]_j$, $j= 1, \dots,4$, are
  tensoring.
\end{theorem}

\begin{proof}
  In~\cite{puschnigg98:product}, Puschnigg proves the assertion if $\MA[B]_1$
  and~$\MA[B]_2$ are Banach algebras with the bounded bornology.  This is the
  special case of entire cyclic cohomology for Banach algebras.  Furthermore,
  the chain maps between $X(\Tanil\MA[B]_1) \hot X(\Tanil\MA[B]_2)$ and
  $X\bigl( \Tanil(\MA[B]_1 \hot \MA[B]_2) \bigr)$ and the homotopy operators
  joining the compositions to the identity are \emph{natural}.  We can reduce
  the case of tensoring algebras to this special case by writing a tensoring
  algebra as an inductive limit of Banach algebras.
  
  The functors $\MA[C]_1 \hot \MA[C]_2 \mapsto X(\Tanil\MA[C]_1) \hot
  X(\Tanil\MA[C]_2)$ and $\MA[C]_1 \hot \MA[C]_2 \mapsto X\bigl(
  \Tanil(\MA[C]_1 \hot \MA[C]_2) \bigr)$ commute with inductive limits.  That
  is, if $\MA[C]_j = \indlim {(\MA[C]_{j,i})_{i \in I}}$ for $j = 1,2$, then
  $$
  X(\Tanil\MA[C]_1) \hot X(\Tanil\MA[C]_2) \cong
  \indlim X(\Tanil\MA[C]_{1,i}) \hot X(\Tanil\MA[C]_{2,i}),  \qquad
  X\bigl( \Tanil(\MA[C]_1 \hot \MA[C]_2) \bigr) \cong
  \indlim X\bigl( \Tanil(\MA[C]_{1,i} \hot \MA[C]_{2,i}) \bigr).
  $$
  Thus if we have natural maps between these complexes for Banach algebras, we
  get induced natural maps for arbitrary tensoring algebras.
\end{proof}

\chapter{Periodic cyclic cohomology}
\label{sec:periodic}

The methods used above for analytic cyclic cohomology can also be used to
study periodic cyclic cohomology.  In Section~\ref{sec:pro_algebras}, we carry
over the machinery of universal nilpotent extensions.  However, the
appropriate analogue of the analytic tensor algebra of an algebra~$A$ is no
longer an algebra any more, but the projective system of algebras $\Tpnil A
\defeq \bigl( \Tens A / (\Jens A)^n \bigr)_{n \in \N}$.  Thus to handle
periodic cyclic cohomology properly, we have to define it on the category of
pro-algebras, that is, algebras in the category of projective systems of
``vector spaces''.

The ``vector spaces'' are put in quotation marks because we have great freedom
to choose the entries of our projective systems.  We define the category
$\PRO(\Cat)$ of projective systems over any additive category $\Cat$ in
Section~\ref{sec:pro_systems}.  We can take for $\Cat$ the vector spaces over
any ground field of characteristic zero; or complete locally convex
topological vector spaces; or complete bornological vector spaces.  In these
cases the theory goes through without difficulty.  To formulate results in a
more concrete way, we will stick to complete bornological vector spaces and
leave it to the interested reader to adapt the theory to other cases.

\begin{digression}
  It would be very desirable to take for $\Cat$ the category of abelian
  groups, so that we get the periodic cyclic cohomology of rings without
  additional structure.  Unfortunately, there appears to be no reasonable
  definition of homotopy in this category.  The proof of
  Proposition~\ref{pro:X_homotopy} involves the integration of functions and
  thus division by integers.  This only works if all objects are vector spaces
  over the rational numbers~$\Q$.
\end{digression}

The analogue of the analytic tensor algebra for pro-algebras is the
\emph{periodic tensor algebra} $\Tpnil A$.  This is defined, roughly
speaking, as the projective limit of the pro-algebras $\bigr( \Tens A / (\Jens
A)^k \bigl)_{k \in \N}$.  The algebra $\Tens A / (\Jens A)^k$ is most easily
\emph{defined} by $\Tens A / (\Jens A)^k \defeq A \oplus \Omega^2 A \oplus
\dots \oplus \Omega^{2k-2} A$ with the ``cut off'' Fedosov product as
multiplication.  ``Cut off'' means that all terms of degree at least~$2k$ are
ignored.  We do not define the tensor algebra $\Tens A$ of a pro-algebra here.
Suffice it to say that it always exists, is very complicated, and that the
quotients $\Tens A / (\Jens A)^k$ are isomorphic to the simple objects
described above.  We define $\Omega_p A$ as the projective limit of the
projective system of differential pro-algebras $A \oplus \Omega^1 A \oplus
\cdots \oplus \Omega^k A$ with the ``cut off'' multiplication and
differential.

The appropriate analogue of a-nilpotent algebras is the class of \emph{locally
  nilpotent pro-algebras}.  Let~$N$ be a pro-algebra with multiplication $m_N
\colon N \hot N \to N$.  Let $m_N^k \colon N^{\hot k} \to N$ be the iterated
multiplication map.  We call~$N$ \emph{locally nilpotent} iff for all
pro-linear maps $l \colon N \to C$ whose range is a \emph{constant} pro-vector
space, there is $k \in \N$ such that $l \circ m_N^k = 0$.  A \emph{constant
  pro-vector space} is a complete bornological vector space viewed as a
pro-vector space indexed by a one point set.

We define a class of linear maps with locally nilpotent curvature so that all
the axioms in Section~\ref{sec:Tanil} are satisfied with periodic tensor
algebras, locally nilpotent algebras, and ``\lonilcurs'' instead of analytic
tensor algebras, analytically nilpotent algebras, and \lanilcurs.  There are,
however, some simplifications.  Tensor products $A \hot N$ with locally
nilpotent~$N$ are locally nilpotent without any assumption on~$A$.  Hence the
analogue of the Homotopy Axiom is trivial.  A pro-algebra~$R$ is \emph{locally
  quasi-free} in the sense that there is a bounded splitting homomorphism $R
\to \Tpnil R$ if and only if~$R$ is quasi-free in the \Hochschild
homological sense.

Thus we can carry over all consequences of the axioms in Section
\ref{sec:Tanil}.  In particular, we have the Universal Extension Theorem and
the Uniqueness Theorem.  We define the X-complex of a pro-algebra in the
obvious way and define \emph{periodic cyclic cohomology} using the X-complex
$X(\Tpnil A)$.  The resulting theory is automatically homotopy invariant for
absolutely continuous homotopies and stable with respect to generalized
algebras of trace class operators.  Moreover, we get the homotopy equivalence
$X(\Tpnil A) \sim (\Omega_p A, B+b)$ as in the analytic case.  Actually, the
proof is simpler because rescaling with $n!$ does not matter in $X(\Tpnil
A)$.

The bivariant periodic cyclic cohomology for pro-algebras satisfies excision
for allowable extensions of pro-algebras.  The proof in
Section~\ref{sec:excision} simplifies considerably because the \Hochschild
homological notion of quasi-freeness is already the right one for the
Uniqueness Theorem.  Hence the construction in
Section~\ref{sec:excision_step_II} of the splitting homomorphism $\LL \to
\Tanil\LL$ can be omitted.  Pro-algebras and excision for pro-algebras have
been studied by Valqui~\cite{valqui98:pro} and
Gr\o{}nb\ae{}k~\cite{groenbaek99:bivariant}.  Our emphasis is quite different
from Valqui's because we treat only extensions with a pro-linear section.

In Section~\ref{sec:dimension_estimate} we return to ordinary algebras and
consider excision results in \emph{cyclic (co)\hspace*{0pt}homology}.  The
cyclic homology $\HC_\ast(\MA)$ and cohomology $\HC^\ast(\MA)$ are defined in
Appendix~\ref{app:cyclic} using the \emph{Hodge filtration} $F_n(\MA)$ on the
complex $X(\Tanil\MA)$.  Let $\MA[K] \injto \MA[E] \prto \MA[Q]$ be an
allowable extension of complete bornological algebras.  Recall that $\varrho
\colon X(\Tanil\MA[K]) \to X(\Tanil\MA[E] : \Tanil\MA[Q])$ is the natural map
induced by the inclusion $\MA[K] \to \MA[E]$.  The Excision
Theorem~\ref{the:excision_analytic} asserts that~$\varrho$ is a homotopy
equivalence.  We can estimate how the maps implementing this homotopy
equivalence shift the Hodge filtration:

\begin{theorem}  \label{the:excision_cyclic}
  There are a chain map $f \colon X(\Tanil\MA[E] : \Tanil\MA[Q]) \to
  X(\Tanil\MA[K])$ and a map $h \colon X(\Tanil\MA[E] : \Tanil\MA[Q]) \to
  X(\Tanil\MA[E] : \Tanil\MA[Q])$ of degree~$1$ such that $f \circ \varrho =
  \ID$ on $X(\Tanil\MA[K])$, $h \circ \varrho = 0$, and $\ID - \varrho \circ f
  = [h, \partial]$.  In addition,
  \begin{displaymath}
    f\bigl( F_{3n+2}(\MA[E] : \MA[Q]) \bigr) \subset F_n(\MA[K]), \qquad
    h\bigl( F_{3n+2}(\MA[E] : \MA[Q]) \bigr) \subset F_n(\MA[E] : \MA[Q]).
  \end{displaymath}

  If there is a bounded splitting homomorphism $s \colon \MA[Q] \to \MA[E]$,
  we can achieve the better estimates
  \begin{gather*}
    f\bigl( F_{4n-1} (\MA[E] : \MA[Q]) \bigr) \subset
    F_{2n-1} (\MA[K]), \qquad
    h\bigl( F_{4n-1} (\MA[E] : \MA[Q]) \bigr) \subset
    F_{2n-1} (\MA[E] : \MA[Q]).
    \\
    f\bigl( F_{4n+2} (\MA[E] : \MA[Q]) \bigr) \subset
    F_{2n} (\MA[K]), \qquad
    h\bigl( F_{4n+2} (\MA[E] : \MA[Q]) \bigr) \subset
    F_{2n} (\MA[E] : \MA[Q]).
  \end{gather*}
\end{theorem}

Using the definition of cyclic cohomology in terms of the Hodge filtration in
Appendix~\ref{app:cyclic}, it follows that the chain map~$f$ induces maps
$\HC_n(f) \colon \HC_{3n+2}(\MA[E] : \MA[Q]) \to \HC_n (\MA[K])$.  The
estimates about the chain homotopy~$h$ imply that $\HC_n(f) \circ
\HC_{3n+2}(i) = S^{n+1} \colon \HC_{3n+2} (\MA[K]) \to \HC_{n}(\MA[K])$ and
$\HC_n(i) \circ \HC_n(f) = S^{n+1} \colon \HC_{3n+2} (\MA[E]: \MA[Q]) \to
\HC_n(\MA[E]: \MA[Q])$.  That is, we get the \Mpn{S}operator iterated $n+1$
times.  Dual statements hold in periodic cyclic cohomology.

The connecting map $X(\Tanil\MA[Q]) \to X(\Tanil\MA[K])$ is the composition $f
\circ [\partial, s_L]$, with $[\partial, s_L] \colon X(\Tanil\MA[Q]) \to
X(\Tanil\MA[E]: \Tanil\MA[Q])$ the connecting map of
Theorem~\ref{the:long_exact_homology}.  It maps $F_{3n+3} (\MA[E] : \MA[Q])
\mapsto F_n (\MA[K])$ because $\partial (F_n) \subset F_{n-1}$ for all $n \in
\N$.  Hence the connecting map in periodic cyclic (co)homology restricts to
maps $\HC^n (\MA[K]) \to \HC^{3n+3} (\MA[Q])$ and $\HC_{3n+3} (\MA[Q]) \to
\HC_n(\MA[K])$.  Puschnigg obtains the same estimate for the connecting map
and proves in addition that it is optimal for a class of extensions including
all extensions of the form $\Jens\MA \injto \Tens\MA \prto \MA$.  It follows
that the estimates in Theorem~\ref{the:excision_cyclic} are optimal in the
worst case.  I do not know whether the above estimates in the split case are
optimal.

The proof of Theorem~\ref{the:excision_cyclic} involves lengthy and
complicated bookkeeping.  It is much easier to prove the weaker statement that
$f \bigl( F_{3n + \const}(\MA[E]: \MA[Q]) \bigr) \subset F_n(\MA[K])$ for a
suitable constant $\const$.  Considerable work is required to detect the
cancellation that is responsible for the optimal estimate.

In Section~\ref{sec:permanence} we argue why it is a good idea to study cyclic
cohomology theories for bornological algebras and not for topological
algebras.  Bornological algebras contain both algebras without additional
structure and \Frechet algebras as full subcategories.  The ``bornological''
cyclic cohomology groups are the usual ones on these subcategories.  Thus we
can handle the most important examples of algebras in a common framework.
Furthermore, there are interesting examples of complete bornological algebras
that are not topological algebras, but for which a ``topological'' cyclic
cohomology is needed.

\section{Projective systems}
\label{sec:pro_systems}

Let~$\Cat$ be an additive category, that is, the morphisms in~$\Cat$ are
Abelian groups.  We define the category $\PRO(\Cat)$ of \emph{projective
  systems over~$\Cat$}.  To define algebras, we need an additive tensor
product functor $\Cat \times \Cat \to \Cat$.  We will call the objects of
$\Cat$ ``vector spaces'', the morphisms in~$\Cat$ ``linear maps'', and the
morphisms $A \otimes B \to C$ ``bilinear maps''.

A \emph{projective system over $\Cat$} or briefly \emph{pro-vector space}
consists of a directed set~$I$ of indices, ``linear spaces''~$V_i$ for all $i
\in I$, and ``linear maps'' $\phi_{i, j} \colon V_i \to V_j$ for all $i \ge
j$, $i, j \in I$.  These so-called \emph{structure maps} are assumed to be
compatible in the sense that $\phi_{j,k} \circ \phi_{i,j} = \phi_{i,k}$
whenever $i \ge j \ge k$, $i, j, k \in I$.  Furthermore, $\phi_{i,i} = \ID$
for all $i \in I$.  These are the \emph{objects} of $\PRO(\Cat)$.

The space of \emph{morphisms} between pro-vector spaces $(V_i)_{i \in I}$ and
$(W_j)_{j \in J}$ is the Abelian group defined by
$$
\Mor \bigl( (V_i); (W_j) \bigr) \defeq
\prolim_j \indlim_i \Mor_{\Cat}(V_i; W_j).
$$
We call these morphisms \emph{pro-linear maps}.  That is, a morphism $f
\colon V \to W$ consists of the following data.  For all $j \in J$, an index
$i(j) \in I$ and a ``linear map'' $f_j \colon V_{i(j)} \to W_j$.  These maps
are compatible in the sense that if $j' \ge j$, then there is $i' \ge i(j),
i(j')$ such that the diagram
$$
\xymatrix@+.5cm{
  {V_{i'}} \ar[r]^{\phi_{i', i(j')}^{V}} \ar[rd]^{\phi_{i', i(j)}^{V}} &
    {V_{i(j')}} \ar[r]^{f_{j'}} &
      {W_{j'}} \ar[d]^{\phi_{j', j}^{W}} \\
  & {V_{i(j)}} \ar[r]_{f_j} &
      {W_{j}}
  }
$$
commutes.  On these raw data, an equivalence relation is introduced as
follows.  Collections of ``linear maps'' $f_j \colon V_{i(j)} \to W_j$, $f'_j
\colon V_{i'(j)} \to W_j$, $j \in J$, describe the same morphism iff there is
a function $j \mapsto i''(j)$ from~$J$ to~$I$ such that $i''(j) \ge i'(j),
i(j)$ for all $j \in J$ and the diagram
$$
\xymatrix@+.5cm{
  {V_{i''(j)}} \ar[r]^-{\phi_{i''(j),i'(j)}^{V}}
    \ar[d]^-{\phi_{i''(j),i(j)}^{V}} &
    {V_{i'(j)}} \ar[d]^-{f'_j} \\
  {V_{i(j)}} \ar[r]^-{f_j} &
    {W_j}
  }
$$
commutes.

These definitions are quite complicated and abound in indices.  However, there
is a reasonable way to handle pro-vector spaces by comparing them to
\emph{constant} pro-vector spaces.  This technique is known already to Artin
and Mazur~\cite{artin69:etale}.  A constant pro-vector space is just a
``vector space''~$C$, viewed as pro-vector space that is indexed by a one
point set.  If $V = (V_i)_{i \in I}$ is a pro-vector space, then the map $\ID
\colon V_i \to V_i$ defines a pro-linear map $V \to V_i$.  Conversely, if $l
\colon V \to C$ is a pro-linear map with constant range, then there is $i \in
I$ such that~$l$ factors as $V \to V_i \to C$ for some ``linear map'' $V_i \to
C$.  This is immediate from the definition of a pro-linear map.  The
pro-vector space~$V$ can be characterized uniquely by the class of all
pro-linear maps $V \to C$ with constant range.

The tensor product in~$\Cat$ gives rise to a tensor product of pro-vector
spaces.  We define
\begin{equation}
  \label{eq:def_tensor_pro}
  (V_i)_{i \in I} \otimes (W_j)_{j \in J} \defeq
  (V_i \otimes W_j)_{i \in I, j \in J}.
\end{equation}
The tensor product is indexed by the product $I \times J$ with the product
ordering and structure maps $\phi_{i',i} \otimes \phi_{j',j}$ for $i' \ge i$
and $j' \ge j$.  The tensor product can be characterized as follows.  For each
pair of morphisms $l_V \colon V \to C_V$, $l_W \colon W \to C_W$ with constant
range, there is a unique morphism $l_V \otimes l_W \colon V \otimes W \to C_V
\otimes C_W$ and each morphism $l \colon V \otimes W \to C$ with constant
range factors through a morphism of the form $l_V \otimes l_W$.

A \emph{pro-algebra}~$A$ is a pro-vector space together with a pro-linear map
$m_A \colon A \otimes A \to A$, the multiplication, that is associative in the
usual sense $m_A \circ (m_A \otimes \ID) = m_A \circ (\ID \otimes m_A)$.

A \emph{left module} over a pro-algebra~$A$ is a pro-vector
spaces~$V$ with a morphism of projective systems $m_{AV} \colon A \otimes V
\to V$ that is associative in the sense that $m_{AV} \circ
(m_A \otimes \ID[V]) = m_{AV} \circ (\ID[A] \otimes m_{AV})$.

A homomorphism of pro-algebras $A \to B$ is a morphism of pro-vector spaces $f
\colon A \to B$ such that $m_B \circ (f \otimes f) = f \circ m_A$.  A
homomorphism of \Mpn{A}modules $V \to W$ is a morphism of pro-vector spaces $f
\colon V \to W$ such that $m_{AW} \circ (\ID \otimes f) = f \circ m_{AV}$.

Since there are finite direct sums in~$\Cat$, there are finite direct sums in
$\PRO(\Cat)$ as well, defined by $(V_i)_{i \in I} \oplus (W_j)_{j \in J}
\defeq (V_i \oplus W_j)_{i \in I, j \in J}$.  We define \emph{allowable
  extensions of pro-vector spaces} in the obvious way as diagrams $(i,p)
\colon V' \injto V \prto V''$ such that there are pro-linear maps $s \colon
V'' \to V$, $t \colon V \to V'$ satisfying $t \circ i = \ID[V']$, $p \circ s =
\ID[V'']$ and $i \circ t + s \circ p = \ID[V]$.  Thus $V \cong V' \oplus V''$.

More generally, we can define arbitrary extensions as follows.  Assume that we
know what an extension of ``vector spaces'' is.  If we have no other concept
of extension, we may take allowable extensions.  A diagram $(i,p) \colon V'
\to V \to V''$ of pro-vector spaces with $p \circ i = 0$ is an \emph{extension
  of pro-vector spaces} iff it is a projective limit of extensions of constant
pro-vector spaces.  That is, we have a projective system of extensions $(i_j,
p_j) \colon C'_j \injto C_j \prto C''_j$ of ``vector spaces'' and $V' \cong
\prolim C'_j$, $V \cong \prolim C_j$, $V'' \cong \prolim C''_j$, $i \cong
\prolim i_j$, $p \cong \prolim p_j$.  An allowable extension of pro-vector
spaces is an extension in this sense.

In the category $\PRO(\Cat)$, projective limits and products exist for formal
reasons because we can identify a pro-pro-vector space with a pro-vector space
in a natural way.  If $(V_{i,j})_{j \in J_i}$ is a pro-vector space for all $i
\in I$ and if these pro-vector spaces form a projective system with structure
maps $\phi_{i, i'} \in \Mor \bigl( (V_{i,j})_{j}; (V_{i',j})_j \bigr)$, then
the projective limit is simply $(V_{i,j})_{i \in I, j \in J_i}$ indexed by
$\bigcup_{i \in I} J_i$ with the obvious partial ordering and the obvious
structure maps.

Products and projective limits of pro-algebras are again pro-algebras.  Tensor
products of pro-algebras are pro-algebras.

\section{Periodic cyclic cohomology for pro-algebras}
\label{sec:pro_algebras}

Let~$A$ be a pro-algebra.  We define $\Omega^n A \defeq \Unse{A} \hot A^{\hot
  n}$ as usual, using the tensor product of pro-vector spaces.  The
\emph{unitarization} is defined as the direct sum $\Unse{A} \defeq A \oplus
\C$, with the usual multiplication.  The multiplication of differential forms
gives rise to pro-linear maps $\Omega^n A \hot \Omega^m A \to \Omega^{n+m} A$.
Since the multiplication of differential forms uses only the multiplication
in~$A$, it continues to make sense for pro-algebras.  The differential $d
\colon \Omega^n A \to \Omega^{n+1} A$ is easy to carry over to pro-algebras.

The finite direct sum $\Omega^{\le n} A \defeq A \oplus \Omega^1 A \oplus
\cdots \oplus \Omega^n A$ is a pro-vector space.  We make $\Omega^{\le n} A$ a
differential pro-algebra by cutting off the multiplication and differential,
that is, ignoring all terms in $\Omega^m A$ with $m > n$.  The natural
projection $\Omega^{\le m} A \to \Omega^{\le n} A$ for $m \ge n$ annihilating
$\Omega^k A$ with $k > n$ is a differential homomorphism, that is,
multiplicative and compatible with the differential.  Thus we get a projective
system $(\Omega^{\le n} A)_{n \in \N}$ of differential pro-algebras.  The
\emph{periodic differential envelope} $\Omega_p A$ of~$A$ is the projective
limit of this system.  Thus $\Omega_p A$ is a differential pro-algebra.  The
even part of $\Omega_p A$ endowed with the Fedosov product is the
\emph{periodic tensor algebra} $\Tpnil A$ of~$A$.

The natural projection $\Omega_p A \to \Omega^{\le 0} A = A$ restricts to a
homomorphism $\tau_A \colon \Tpnil A \to A$.  We write $\Jpnil A$ for its
kernel.  Thus $\Jpnil A$ is the projective limit of the pro-algebras $\Janil
_p A/ (\Jpnil A)^n \defeq \Omega^2 A \oplus \dots \oplus \Omega^{2n-2} A$.
The usual formula defines a pro-linear section $\sigma_A \colon A \to \Tpnil
A$ for~$\tau_A$.  It is the map $A \to \Tpnil A$ induced by the natural
inclusions $A \to A \oplus \Omega^2 A \oplus \cdots \oplus \Omega^{2k-2} A$.

To carry over the formalism of Section~\ref{sec:Tanil}, we define
\emph{locally nilpotent} pro-algebras and \emph{\lonilcurs{}}.  The latter is
an abbreviation for \emph{pro-linear maps with locally nilpotent curvature}.

A constant pro-vector space is just an ordinary ``vector space'', viewed as a
projective system in the trivial way.  Let~$A$ be a pro-algebra with
multiplication $m_A \colon A \hot A \to A$.  Let $m_A^n \colon A ^{\hot n} \to
A$ be the iterated multiplication.  We call~$A$ \emph{locally nilpotent} iff
for all pro-linear maps $f \colon A \to C$ with constant range, there is $n
\in \N$ such that $f \circ m_A^n = 0$.  Typical examples are nilpotent
pro-algebras, where $m_A^n = 0$ for some $n \in \N$, and projective limits of
nilpotent pro-algebras.

\begin{axiom}  \label{axiom:Janilp_lonilpotent}
  The pro-algebra $\Jpnil A$ is locally nilpotent for all pro-algebras~$A$.
\end{axiom}

\begin{proof}
  Any pro-linear map $l \colon \Jpnil A \to C$ with \emph{constant} range
  factors through $\Jpnil A / (\Jpnil A)^n$ for suitable~$n$ by the
  concrete construction of projective limits.  The \Mpn{n}fold multiplication
  in the pro-algebra $\Jpnil A / (\Jpnil A)^n$ is zero because a product
  of~$n$ forms of positive even degree has degree at least $2n$ and is
  therefore zero in $\Jpnil A / (\Jpnil A)^n = \Omega^2 A \oplus \dots
  \oplus \Omega^{2n-2} A$.  Thus $l \circ m^n_{\Jpnil A} = 0$.
\end{proof}

A pro-algebra~$R$ is \emph{locally quasi-free} iff there is a splitting
homomorphism $R \to \Tpnil R$.  A \emph{universal locally nilpotent
  extension} is an allowable extension $N \injto R \prto A$ with locally
nilpotent~$N$ and locally quasi-free~$R$.  A morphism of pro-vector spaces $l
\colon A \to B$ between pro-algebras is called a \emph{\lonilcur{}} iff its
curvature $\omega_l \colon A \hot A \to B$ is locally nilpotent in the
following sense: For each pro-linear map $f \colon B \to C$ with constant
range, there is $n \in \N$ such that $f \circ m_B^n \circ \omega_l^{\hot n} =
0$.

\begin{axiom}  \label{axiom:lonilpotent_lonilcur}
  A pro-algebra~$N$ is locally nilpotent if the identity map $\Null(N) \to N$
  is a \lonilcur.
\end{axiom}

\begin{proof}
  The curvature of the identity map $l \colon \Null(N) \to N$ is equal
  to~$-m_N$.  If~$l$ is a \lonilcur, then for any pro-linear map $f \colon N
  \to C$ with constant range there is $k \in \N$ such that $f \circ m_N^k
  \circ \omega_l = 0$, that is, $-f \circ m_N^k \circ m_N = f \circ m_N^{2k} =
  0$.  Thus~$N$ is locally nilpotent.
\end{proof}

\begin{axiom}  \label{axiom:lonilcur_lonilpotent}
  Let $l \colon A \to B$ be a pro-linear map between pro-algebras whose
  curvature factors through a locally nilpotent algebra, that is, $\omega_l =
  f \circ \omega'$ for a homomorphism of pro-algebras $f \colon N \to B$ and a
  pro-linear map $\omega' \colon A \hot A \to N$ with locally nilpotent~$N$.
  Then~$l$ is a \lonilcur.
\end{axiom}

\begin{proof}
  Assume that the pro-linear map $l \colon A \to B$ factors as described
  above.  Let $g \colon B \to C$ be a pro-linear map with constant range.
  Then $g \circ f \circ m_N^k = 0$ for some $k \in \N$ because~$N$ is locally
  nilpotent.  Since~$f$ is a homomorphism, we have $f \circ m_N^k = m_B^k
  \circ f^{\hot n}$.  Hence $0 = g \circ m_N^k \circ f^{\hot n} \circ
  (\omega')^{\hot n} = g \circ m_N^k \circ \omega_l ^{\hot n}$.  This means
  that~$l$ is a \lonilcur.
\end{proof}

\begin{axiom}  \label{axiom:Tanilp_universal}
  Let~$A$ be a pro-algebra.  The pro-algebra $\Tpnil A$ and the bounded
  linear map $\sigma_A \colon A \to \Tpnil A$ have the following universal
  property.  If $l \colon A \to B$ is a \lonilcur into a pro-algebra~$B$, then
  there is a unique homomorphism of pro-algebras $\LLH{l} \colon \Tpnil A
  \to B$ with $\LLH{l} \circ \sigma_A = l$.  Furthermore, any map of the form
  $f \circ \sigma_A$ with a bounded homomorphism $f \colon \Tpnil A \to B$
  is a \lonilcur.
\end{axiom}

\begin{proof}
  Let~$l$ be a \lonilcur.  We have to construct a homomorphism of pro-algebras
  $\LLH{l} \colon \Tpnil A \to B$ satisfying $\LLH{l} \circ \sigma_A = l$.
  Let $\opt{l} \omega_l^k \defeq m_B^{k+1} \circ (\opt{l} \hot \omega_l^{\hot
    k}) \colon \Omega^{2k} A \to B$, where $m_B^{k+1} \colon \Unse{B} \hot
  B^{\hot k} \to B$ is the multiplication in~$B$ and~$\Unse{B}$ and~$\opt{l}$
  is the unital extension of~$l$.  The map $\opt{l} \omega_l^k$ describes the
  map $\opt{x_0} dx_1 \dots dx_{2k} \mapsto l\opt{x_0} \cdot \omega_l(x_1,
  x_2) \cdots \omega_l(x_{2n-1}, x_{2k})$ in the situation of pro-algebras.
  
  Let $f \colon B \to C$ be a pro-linear map with constant range.  We claim
  that there is some $n \in \N$ such that $f \circ \opt{l} \omega_l^k = 0$ for
  all $k \ge n$.  The pro-linear map $f \circ \Unse{m_B} \colon \Unse{B} \hot
  B \to C$ factors as $f = g \circ (f_1 \hot f_2)$ with pro-linear maps $f_1
  \colon \Unse{B} \to C_1$, $f_2 \colon B \to C_2$ with constant range and a
  bounded linear map $g \colon C_1 \hot C_2 \to C$.  Since~$l$ is a \lonilcur,
  there is some $n \in \N$ such that $f_2 \circ \omega_l^n = 0$.  Hence
  $$
  f \circ \opt{l} \omega_l^k =
  f \circ \Unse{m_B} \circ (\opt{l} \omega_l^{k-n} \hot \omega_l^n) =
  g \circ (f_1 \circ \opt{l} \omega_l^{k-n} \hot f_2 \circ \omega_l^n) =
  0.
  $$
  
  To construct~$\LLH{l}$, write $B = (B_i)_{i \in I}$ as a projective system.
  For each $i \in I$, we have the pro-linear map $f_i \colon B \to B_i$ with
  constant range and have found above a number $n_i \in \N$ such that $\opt{l}
  \omega_l^k = 0$ for all $k \ge n_i \in \N$.  We define pro-linear maps
  $\LLH{l}_i \colon \Tpnil A / (\Jpnil A)^{n_i} \to B_i$ by
  $$
  \LLH{l}_i \defeq f_i \circ (
  l \oplus \opt{l} \omega_l^1 \oplus \dots \oplus \opt{l} \omega_l^{n_i-1})
  \colon
  A \oplus \Omega^2 A \oplus \dots \oplus \Omega^{2n_i-2} A \to B \to B_i.
  $$
  The pro-linear maps~$\LLH{l}_i$ form a morphism of projective systems
  from $\bigl( \Tens A / (\Jens A)^n \bigr)_{n \in \N}$ to $(B_i)_{i \in I}$.
  This follows from $f_i \circ \opt{l} \omega_l^k = 0$ for all $k \ge n_i$.
  This morphism of projective systems induces a pro-linear map $\LLH{l} \colon
  \Tpnil A = \prolim {\bigl(\Tens / (\Jens A)^n \bigr)} \to \prolim
  {(B_i)_{i \in I}} = B$.
  
  The pro-linear map~$\LLH{l}$ is defined by the same formula that occurs in
  the universal property of the algebraic tensor algebra $\Tens A$ of an
  algebra in~\eqref{eq:LLH_Tens}.  Since that formula defines a homomorphism,
  the map~$\LLH{l}$ is a homomorphism of pro-algebras.  Details are left to
  the reader.  Moreover, $\LLH{l} \circ \sigma_A = l$ by construction,
  and~$\LLH{l}$ is the only homomorphism $\Tpnil A \to B$ with that
  property.
  
  It remains to show that any linear map of the form $l = f \circ \sigma_A
  \colon A \to B$ with a homomorphism of pro-algebras $f \colon \Tpnil A \to
  B$ is a \lonilcur.  If $g \colon B \to C$ is a pro-linear map with constant
  range, then $g \circ f$ factors through $\Tpnil A / (\Jpnil A)^k$ for
  some $k \in \N$.  Hence $g \circ \omega_l = g \circ f \circ \omega_\sigma$
  factors through $\Jpnil A / (\Jpnil A)^k$.  Since the latter pro-algebra
  is \Mpn{k}nilpotent, it follows that $g \circ \omega_l^k = 0$.
\end{proof}

Hence the periodic tensor algebra, \lonilcurs, and locally nilpotent algebras
are related by the same axioms as in Section~\ref{sec:Tanil}.  It remains to
verify the analogues of the Homotopy Axiom and the Extension Axiom.  The
Homotopy Axiom becomes trivial:

\begin{lemma}  \label{lem:lonilpotent_tensor}
  Let $A$ and~$N$ be pro-algebras.  If~$N$ is locally nilpotent, so is $A \hot
  N$.
\end{lemma}

\begin{proof}
  Let $f \colon A \hot N \to C$ be a pro-linear map with constant range.  It
  factors as $g \circ (f_1 \hot f_2)$ for pro-linear maps with constant range
  $f_1 \colon A \to C_1$, $f_2 \colon B \to C_2$, and a bounded linear map $g
  \colon C_1 \hot C_2 \to C$.  Since~$N$ is locally nilpotent, there is $k \in
  \N$ with $f_2 \circ m_N^k = 0$.  The \Mpn{k}fold multiplication in $A \hot
  N$ is $m_A^k \hot m_N^k$ up to a coordinate flip $\vartheta \colon (A \hot
  N)^{\hot k} \to A^{\hot k} \hot N^{\hot k}$.  Hence $f \circ m_{A \hot N}^k
  = \bigl( (f_1 \circ m_A^k) \hot (f_2 \circ m_N^k) \bigr) \circ \vartheta =
  0$.  Thus $A \hot N$ is locally nilpotent.
\end{proof}

Therefore, we do not have to worry about tensoring algebras and can the
algebra $\C[t]$ of polynomials instead of $\ABC([0,1])$ or $\CINF([0,1])$.
Hence the homotopies constructed in the Universal Extension Theorem are
\emph{polynomial}.

\begin{axiom}[Extension Axiom]  \label{axiom:extension_lonilpotent}
  Let $(\iota, \pi) \colon N' \injto N \prto N''$ be an extension of
  pro-algebras.  Assume that $N'$ and~$N''$ are locally nilpotent.  Then~$N$
  is locally nilpotent.
\end{axiom}

\begin{proof}
  By definition, an extension of pro-algebras is a projective limit of
  extensions.  That is, we have extensions $(\iota_i, \pi_i) \colon N'_i
  \injto N_i \prto N_i''$ of complete bornological vector spaces indexed by a
  directed set~$I$ and morphisms $(\phi_{j,i}', \phi_{j,i}, \phi_{j,i}'')$ for
  $j \ge i$ between these extensions such that $N' = \prolim {(N'_i)_{i \in
      I}}$, $N = \prolim {(N_i)_{i \in I}}$, $N'' = \prolim {(N_i'')_{i \in
      I}}$, $\iota = \prolim {(\iota_i)_{i \in I}}$, $\pi = \prolim
  {(\pi_i)_{i \in I}}$.
  
  We can describe the \Mpn{n}fold multiplication $m^n_N = m^n$ in~$N$ by maps
  $m^n_i \colon N_{j_n(i)}^{\hot n} \to N_i$.  Since the map $\pi \colon N \to
  N''$ is multiplicative, we can achieve (by increasing~$j_n(i)$ if necessary)
  that~$m^n_i$ maps the kernel of $\pi_{j_n(i)}^{\hot n} \colon
  N_{j_n(i)}^{\hot n} \to (N''_{j_n(i)})^{\hot n}$ into~$N'_i$.  The maps on
  the quotients $(N''_{j_n(i)})^{\hot n} \to N''_i$ induced by~$m^n_i$
  describe the \Mpn{n}fold multiplication in~$N''$.  The restrictions to
  $(N'_{j_n(i)})^{\hot n} \to N'_i$ describe the \Mpn{n}fold multiplication
  in~$N'$.  We write $m^n_{j, i} \defeq m^n_i \circ \phi_{j, j(i)} \colon
  N_j^{\hot n} \to N_i$ for all $j \ge j_n(i)$.
  
  We have to show that for all $i \in I$, there are $n \in \N$ and $j \ge
  j_n(i)$ such that $m^n_{j, i} = 0$.  Since~$N'$ is locally nilpotent, there
  are $n' \in \N$ and $i' \in I$, $i' \ge j_{n'}(i)$, such that $m^{n'}_{i',
    i}$ annihilates $(N'_{i'})^{\hot n'}$.  Since~$N''$ is locally nilpotent,
  there are $n'' \in \N$ and $i'' \in I$, $i'' \ge j_{n''}(i')$, such that
  $$
  p_{i'} \circ m^{n''}_{i'', i'} \colon
  N_{i''}^{\hot n''} \to N_{i'} \to N''_{i'}
  $$
  is the zero map.  That is, $N_{i''}^{\hot n''}$ is mapped into~$N'_{i'}$.
  Thus $m^{n'}_{i', i} \circ (m^{n''}_{i'', i'})^{\hot n'} = 0$.  Hence $m^{n'
    \cdot n''}_{j, i} = 0$ for all $j \in I$ with $j \ge i''$ and $j \ge
  j_{n'\cdot n''}(i)$.  Thus~$N$ is locally nilpotent.
\end{proof}

We have verified all axioms of Section~\ref{sec:Tanil}.  Hence the
consequences drawn there carry over.  We only formulate the Uniqueness
Theorem:

\begin{theorem}[Uniqueness Theorem]  \label{the:uniqueness_pro}
  Let~$A$ be a pro-algebra.  Then $\Jpnil A \injto \Tpnil A \prto A$ is a
  universal locally nilpotent extension of~$A$.  Any universal locally
  nilpotent extension $N \injto R \prto A$ is smoothly homotopy equivalent
  relative to~$A$ to $\Jpnil A \injto \Tpnil A \prto A$.  In particular,
  $R$ is smoothly homotopy equivalent to $\Tpnil A$ and~$N$ is smoothly
  homotopy equivalent to $\Jpnil A$.
  
  In addition, any morphism of extensions $(\xi, \psi, \ID)$ between two
  universal locally nilpotent extensions is a smooth homotopy equivalence of
  extensions relative to~$A$.  Any two morphisms of extensions $(\xi, \psi,
  \ID)$ are smoothly homotopic.
\end{theorem}

The X-complex of a pro-algebra is defined as for ordinary algebras by $X_0(A)
\defeq A$ and $X_1(A) \defeq \Omega^1 A / b(\Omega^2 A)$.  We have to declare,
however, what this quotient should mean for pro-algebras.  It can be
characterized by its universal property: $\Omega^1 A / [,] \defeq \Omega^1 A /
b(\Omega^2 A)$ is a projective system of (possibly non-separated) bornological
vector spaces such that the pro-linear maps $\Omega^1 A / [,] \to V$ are in
bijection with pro-linear maps $l \colon \Omega^1 A \to V$ satisfying $l \circ
b = 0$.  If~$A$ is quasi-free in the sense that $\Omega^1 A$ is projective as
a module over~$A$, then $X_1(A)$ is a direct summand in $\Omega^1 A$ and in
particular a projective system of separated spaces.

Homological algebra for modules over a pro-algebra works just as well as for
modules over a complete bornological vector space.  Furthermore, the proof of
Lemma~\ref{deflem:quasi_free} carries over easily to pro-algebras.
Thus if~$R$ is a quasi-free pro-algebra in the sense that~$\Unse{R}$ has an
allowable projective \Mpn{R}bimodule resolution of length~$1$, then there are
homomorphisms $\upsilon_n \colon R \to \Tpnil R / (\Jpnil R)^n$ that
combine to a homomorphism of pro-algebras $R \to \Tpnil R$.  That is, $R$ is
locally quasi-free.  Hence local quasi-freeness and quasi-freeness are the
same concepts.

We define the \emph{bivariant periodic cyclic cohomology} $\HP^\ast( A; B)$ as
the homology of the complex of pro-linear maps $X(\Tpnil A) \to X(\Tpnil
B)$.  Similarly, we define the periodic cyclic (co)homology groups
$\HP^\ast(A)$ and $\HP_\ast(A)$ as the homology groups of the complexes of
pro-linear maps $X(\Tpnil A) \to \C[0]$ and $\C[0] \to X(\Tpnil A)$.  Here
$\C[0]$ is viewed as a constant pro-vector space.

The X-complex continues to be invariant under absolutely continuous homotopies
for quasi-free pro-algebras.  We describe homotopies by homomorphisms of
pro-algebras $A \to B \hot \ABC([0,1])$, where $\ABC([0,1])$ is viewed as a
constant pro-algebra.  The proof of Proposition~\ref{pro:X_homotopy} carries
over without change.  Hence periodic cyclic cohomology is invariant under
absolutely continuous homotopies and stable with respect to generalized
algebras of trace class operators.  The arguments in
Section~\ref{sec:homotopy_stability} carry over.  The algebras $\C$
and~$\MA[U]$ remain quasi-free when viewed as constant pro-algebras.  Hence we
can associate a Chern-Connes character to any homomorphism of pro-algebras $\C
\to A$ or $\MA[U] \to A$ as in Section~\ref{sec:Chern_Ktheory}.  If the
\Mpn{K}theory of pro-algebras is defined appropriately, we can carry over the
Chern-Connes character in \Mpn{K}theory.  Finally, the complex $X(\Tpnil A)$
is chain homotopic to the pro-complex $(\Omega_p A, B+b)$.  The argument in
Section~\ref{sec:X_Tanil} carries over with some simplifications because
rescaling by $[n/2]!$ on $\Omega^n A$ does not matter any more.

We have excision for allowable extensions in both variables for $\HP^\ast( A;
B)$:

\begin{theorem}[Excision Theorem]  \label{the:excision_HP}
  Let~$A$ be a pro-algebra.  Let $(i,p) \colon K \injto E \prto Q$ be an
  extension of pro-algebras with a pro-linear section~$s$.  That is, $s \colon
  Q \to E$ is a pro-linear map satisfying $p \circ s = \ID$.  There are six
  term exact sequences in both variables:
  \begin{gather*}
    \xymatrix{
      {\HP^0(A; K)} \ar[r]^{i_\ast} &
        {\HP^0(A; E)} \ar[r]^{p_\ast} &
          {\HP^0(A; Q)} \ar[d] \\
      {\HP^1(A; Q)} \ar[u] &
        {\HP^1(A; E)} \ar[l]^{p_\ast} &
          {\HP^1(A; K)} \ar[l]^{i_\ast}
      } \displaybreak[0] \\
    \xymatrix{
      {\HP^0(Q; A)} \ar[r]^{p^\ast} &
        {\HP^0(E; A)} \ar[r]^{i^\ast} &
          {\HP^0(K; A)} \ar[d] \\
      {\HP^1(K; A)} \ar[u] &
        {\HP^1(E; A)} \ar[l]^{i^\ast} &
          {\HP^1(Q; A)} \ar[l]^{p^\ast} 
      }
  \end{gather*}
  The maps $i_\ast, i^\ast, p_\ast, p^\ast$ are given by composition with the
  \Mp{\HP}classes $[i]$ and~$[p]$ of the homomorphisms~$i$ and~$p$
  respectively.
\end{theorem}

The proof in Section~\ref{sec:excision} actually simplifies considerably.  Let
$(i,p) \colon K \injto E \prto Q$ be an allowable extension of pro-algebras
with a pro-linear section $s \colon Q \to E$.  As in
Section~\ref{sec:Omega_functorial}, it follows that $X(\Tpnil p) \colon
X(\Tpnil E) \to X(\Tpnil Q)$ is split surjective, so that the long exact
homology sequence applies.  The isomorphisms in
Section~\ref{sec:excision_isomorphisms} carry over easily.  We define the
analogue of~$\LL$ in the obvious way.  We obtain the allowable free
\Mpn{\LL}bimodule resolution $\VS[P]_1 \injto \VS[P]_0 \prto \Unse{\LL}$ as in
Section~\ref{sec:excision_step_I} and compute that the associated commutator
complex is isomorphic to $X_\beta( \Tpnil E : \Tpnil Q) \oplus \C[0]$.  It
follows that~$\LL$ is quasi-free.  Thus the argument in
Section~\ref{sec:excision_step_II} is no longer necessary, we automatically
get a homomorphism of pro-algebras $\upsilon \colon \LL \to \Tpnil \LL$.

\section{Dimension estimates in the Excision Theorem}
\label{sec:dimension_estimate}

We prove Theorem~\ref{the:excision_cyclic}.  A chain map $f \colon
X(\Tanil\MA[E]: \Tanil\MA[Q]) \to X(\Tanil\MA[K])$ and a homotopy $h \colon
X(\Tanil\MA[E]: \Tanil\MA[Q]) \to X(\Tanil\MA[E]: \Tanil\MA[Q])$ are
implicitly constructed in the proof of the Excision Theorem in
Section~\ref{sec:excision}.  We will first define the maps~$f$ and~$h$ that we
are going to use.  Then we will estimate how~$f$ shifts the Hodge filtration.
Finally, we will estimate the behavior of~$h$.  Since the estimates in
Theorem~\ref{the:excision_cyclic} are optimal, we have to do very careful
bookkeeping.

Throughout this section, we adopt the conventions $l, l_j \in \LL$, $q, q_j
\in \MA[Q]$, $x, x_j \in \Tanil\MA[E]$, $\ma[e], \ma[e]_j \in \MA[E]$,
$\ma[k], \ma[k]_j \in \MA[K]$, $\ma[q], \ma[q]_j \in \MA[Q]$, $X, X_j \in
\Tanil\LL$, for all $j \in \Z_+$.

The first step in the proof of excision is Theorem~\ref{the:excision_step_I}
asserting that the natural map $\psi \colon X(\LL) \to X(\Tanil\MA[E]:
\Tanil\MA[Q])$ induced by the inclusion $\LL \subset \Tanil\MA[E]$ is a
homotopy equivalence.  In fact, we have proved more in
Theorem~\ref{the:P_free}.  Namely, the map~$\psi$ is split injective and
$X(\Tanil\MA[E]: \Tanil\MA[Q]) \cong \psi \bigl( X(\LL) \bigr) \oplus
\VS[C]'_\bullet$ with the contractible complex $\VS[C]'_0 \defeq [\LL,
s_L(\Tanil\MA[Q])]$, $\VS[C]'_1 \defeq \LL \,Ds_L(\Tanil\MA[Q]) \bmod [,]$.
In fact, the boundary $\partial_1 \colon \VS[C]'_1 \to \VS[C]'_0$ is a
bornological isomorphism.  Thus we obtain a homotopy inverse $g \colon
X(\Tanil\MA[E]: \Tanil\MA[Q]) \to X(\LL)$ for~$\psi$ by letting $g = 0$ on
$\VS[C]'_\bullet$ and $g = \psi^{-1}$ on $\psi\bigl( X(\LL) \bigr)$.  By
construction, $g \circ \psi = \ID$ on $X(\LL)$ and $\psi \circ g$ is the
projection onto the range of~$\psi$ annihilating~$\VS[C]'_\bullet$.  Let $h_1
\in \Endo \bigl( X(\Tanil\MA[E]: \Tanil\MA[Q]) \bigr)$ be the obvious
contraction of $\VS[C]'_\bullet$ defined by $h_1 = 0$ on $\psi \bigl( X(\LL)
\bigr) + \VS[C]'_1$ and $h_1 = \partial_1 |_{\VS[C]'_1} ^{-1} \colon \VS[C]'_0
\to \VS[C]'_1$ on $\VS[C]'_0$.  Then $\psi \circ g = \ID - [h_1, \partial]$.

The second step in the proof of excision is Theorem~\ref{the:excision_step_II}
asserting that~$\LL$ is analytically quasi-free.  In
Section~\ref{sec:excision_step_II} we have constructed an explicit bounded
splitting homomorphism $\upsilon \colon \LL \to \Tanil\LL$.  Actually, we have
remarked in Section~\ref{sec:pro_algebras} that we do not need
Theorem~\ref{the:excision_step_II} for excision in periodic cyclic cohomology.
It suffices that the pro-algebra $\bigl( \LL/ (\LL \cap \Janil\MA[E]^n)
\bigr)$ is quasi-free, which follows from the pro-algebraic analogue of
Theorem~\ref{the:P_free}.  Nevertheless, since we have constructed the
homomorphism~$\upsilon$ we may just as well use it.  Let $j \colon
\Tanil\MA[K] \to \LL \subset \Tanil\MA[E]$ be the natural map induced by the
inclusion $\MA[K] \subset \MA[E]$.  The analytic quasi-freeness of~$\LL$ and
the Uniqueness Theorem~\ref{the:uniqueness} imply that $X(j) \colon
X(\Tanil\MA[K]) \to X(\LL)$ is a homotopy equivalence.  We make this argument
more explicit.

Let $\tilde\tau \colon \LL \to \MA[K]$ be the natural projection with kernel
$\LL \cap \Janil\MA[E]$.  It induces a bounded homomorphism $\Tanil \tilde\tau
\colon \Tanil\LL \to \Tanil\MA[K]$.  The composition $\Tanil\tilde\tau \circ
\upsilon \colon \LL \to \Tanil\MA[K]$ is a section for~$j$, that is, $\Tanil
\tilde\tau \circ \upsilon \circ j = \ID$ on $\Tanil\MA[K]$.  It suffices to
verify this on $\sigma(\MA[K])$, where it is trivial.  Hence the bounded chain
map
$$
f \defeq X(\Tanil \tilde\tau) \circ X(\upsilon) \circ g \colon
X(\Tanil\MA[E] : \Tanil\MA[Q]) \to X(\LL) \to X(\Tanil\LL) \to X(\Tanil\MA[K])
$$
is a section for $\psi \circ X(j) = \varrho$.  This is the map we take
for~$f$ in Theorem~\ref{the:excision_cyclic}.

Let $\tau_{\LL} \colon \Tanil\LL \to \LL$ be the natural projection
annihilating $\Janil\LL$.  The homomorphisms $j \circ \tilde\tau$
and~$\tau_{\LL}$ are smoothly homotopic.  To construct a homotopy $H \colon
\Tanil\LL \to \CINF([0,1]; \LL)$ between them, it suffices to construct a
\lanilcur $H \colon \LL \to \CINF([0,1]; \LL)$ with $\ev_0 \circ H = j \circ
\Tanil\tilde\tau$ and $\ev_1 \circ H = \tau_{\LL}$.  There are many
possibilities for~$H$, we use the following not so obvious definition.

We define a grading on $\Tanil\MA[E] \supset \LL$ by $\# x = n$ iff $x \in
\Omega^{2n}\MA[E]$.  That is, elements of $\Omega^{2n}\MA[E]$ are homogeneous
of degree~$n$.  We define $H(l) \defeq t^{\# l} \otimes l \in \CINF([0,1])
\hot \LL$ for all homogeneous $l \in \LL$.  This extends to a bounded linear
map $H \colon \LL \to \CINF([0,1]; \LL)$.  Since $\tilde\tau \circ H =
\CINF([0,1]; \tilde\tau)$ is a homomorphism, the curvature of~$H$ has values
in $\CINF([0,1]; \LL \cap \Janil\MA[E])$.  Thus~$H$ is a \lanilcur and can be
extended to a bounded homomorphism $H \colon \Tanil\LL \to \CINF([0,1]; \LL)$.
We have
\begin{multline*}
  H(Dl_1 \,Dl_2) =
  H(l_1 \odot l_2) - H(l_1) \odot H(l_2) =
  H(l_1\cdot l_2) - H(dl_1 dl_2) - H(l_1) \odot H(l_2)
  \\ =
  t^{\# l_1 + \# l_2} l_1 \cdot l_2 -
  t^{\# l_1 + \# l_2 + 1} dl_1 dl_2 -
  t^{\# l_1} l_1 \odot  t^{\# l_2} l_2 =
  t^{\# l_1 + \# l_2}(1-t) dl_1 dl_2
\end{multline*}
and hence $H(\opt{l_0} \,Dl_1 \dots \,Dl_{2n}) = t^{\# \opt{l_0} + \dots + \#
  l_{2n}}(1-t)^n \opt{l_0} dl_1 \dots dl_{2n}$.  We define a bigrading on
$\Tanil\Tanil\MA[E] \supset \Tanil\LL$ as follows.  The \emph{internal
  degree}~$\#_i$ and the \emph{external degree}~$\#_e$ are defined by
$$
\#_i(\opt{x_0} \,Dx_1 \dots \,Dx_{2n}) \defeq \sum_{j=0}^{2n} \# x_j,
\qquad
\#_e(\opt{x_0} \,Dx_1 \dots \,Dx_{2n}) \defeq n.
$$
The \emph{total degree}~$\#_t$ is defined by $\#_t \defeq \#_e + \#_i$.  If
we let $\tilde{H}(\opt{x_0} \,Dx_1 \dots \,Dx_{2n}) \defeq \opt{x_0} dx_1
\dots dx_{2n}$, then we can write
\begin{equation}
  \label{eq:H_explicit}
  H(X) = t^{\#_i X} (1-t)^{\#_e X} \tilde{H}(X)
  \qquad \forall X \in \Tanil\LL.
\end{equation}
If $X \in \Tanil\LL$ is homogeneous with $\#_t X = n$, then $\ev_t \circ H(X)$
is homogeneous of degree $\# H_t(X) = n$ for all $t \in [0,1]$.  We abbreviate
this as $\# H(X) = \#_t X$.  Similarly, $\# \Tanil\tilde\tau (X) = \#_t X$.

We collect the homomorphisms obtained so far in the following commutative
diagram:
\begin{equation}  \label{eq:estimate_homomorphisms}
  \begin{gathered}
    \xymatrix{
      {\:\Tanil\MA[K]\:} \ar@{>->>}[r]^{\ID} \ar[d]_{j} &
        {\:\Tanil\MA[K]\:} \ar@{>->}[r]^{j} &
          {\LL} \\
      {\LL\;} \ar@{>->}[r]_{\upsilon} \ar@{>->>}[dr]_{\ID} &
        {\Tanil\LL} \ar@{->>}[u]_{\Tanil\tilde{\tau}} \ar[r]^-{H}
          \ar@{->>}[d]^{\tau_\LL} &
          {\CINF([0,1]; \LL)} \ar@{->>}[u]_{\ev_0} \ar@{->>}[d]^{\ev_1} \\
      {} &
        {\LL} \ar@{=}[r] &
          {\LL}
      }
  \end{gathered}
\end{equation}

The algebra $\Tanil\LL$ is quasi-free in a natural way.  Since the
homomorphisms $j \circ \Tanil\tilde\tau$ and $\tau_{\LL}$ are smoothly
homotopic, the homotopy invariance of the X-complex~\ref{pro:X_homotopy}
implies that there is a bounded linear map $h_2 \colon X(\Tanil\LL) \to
X(\LL)$ of degree~$1$ such that $X(j \circ \Tanil\tilde\tau) = X(\tau_{\LL}) -
[\partial, h_2]$.  The map~$h_2$ is constructed explicitly during the proof of
Proposition~\ref{pro:X_homotopy} in Appendix~\ref{app:X_homotopy}.  Another
ingredient of the construction is a graded right connection on $\Omega^1
(\Tanil\LL)$ encoding the quasi-freeness of $\Tanil\LL$.  We use the standard
connection $\nabla \colon \Omega^1(\Tanil\LL) \to \Omega^2 (\Tanil\LL)$
defined by
$$
\nabla(\opt{X_0} \delta(l_1) \opt{X_2}) \defeq
-\opt{X_0} \delta l_1 \delta \opt{X_2}
\qquad \forall \opt{X_0}, \opt{X_2} \in \Unse{(\Tanil\LL)},\ l_1 \in \LL.
$$
This well-defines~$\nabla$ on $\Omega^1 (\Tanil\LL)$ because $\Omega^1
(\Tanil\LL) \cong \Unse{(\Tanil\LL)} \hot \LL \hot \Unse{(\Tanil\LL)}$.  It is
easy to verify that~$\nabla$ satisfies~\eqref{eq:def_nabla} and thus is a
graded right connection.  Let~$\dot{H}$ be the derivative of~$H$ and define
$\eta \colon \Omega^{j} (\Tanil\LL) \to \Omega^{j-1} \LL$ by
$$
\eta(\opt{X_0} \delta X_1 \dots \delta X_j) \defeq
\int_0^1 H_t\opt{X_0} \odot \dot{H_t} (X_1) \,DH_t(X_2) \dots \,DH_t(X_j)
\,dt.
$$
Here $\delta \colon \Omega(\Tanil\LL) \to \Omega(\Tanil\LL)$ and $D \colon
\Omega\LL \to \Omega\LL$ denote the standard derivations.  Copying the
formulas in Appendix~\ref{app:X_homotopy}, we put $h_2 \defeq -\eta \circ
\nabla \circ \delta$ on $\Omega^0(\Tanil\LL)$ and $h_2 \defeq \eta \circ (\ID
- b \circ \nabla)$ on $\Omega^1(\Tanil\LL)$.  This map~$h_2$ descends to a map
$h_2 \colon X(\Tanil\LL) \to X(\LL)$ satisfying $[\partial, h_2] = X(\tau_\LL)
- X(j \circ \Tanil\tilde\tau)$.

Finally, we define
$$
h \defeq h_1 + \psi \circ h_2 \circ X(\upsilon) \circ g
$$
and compute
\begin{multline*}
  \varrho \circ f =
  \psi \circ X(j) \circ X(\Tanil\tilde{\tau}) \circ X(\upsilon) \circ g =
  \psi \circ \bigl( X(\tau_\LL) - [\partial, h_2] \bigr) \circ X(\upsilon)
  \circ g
  \\ =
  \psi \circ \bigl( X(\tau_\LL \circ \upsilon) -
  [\partial, h_2\circ X(\upsilon)] \bigr) \circ g =
  \psi \circ g - [\partial, \psi\circ h_2 \circ X(\upsilon) \circ g]
  \\ =
  \ID - [\partial, h_1 + \psi \circ h_2 \circ X(\upsilon) \circ g] =
  \ID - [\partial, h].
\end{multline*}
Thus~$h$ has the required property for Theorem~\ref{the:excision_cyclic}.

We have to estimate how~$f$ and~$h$ shift the Hodge filtration.  The only
problematic ingredients are the maps $X(\upsilon)$ and~$g$.  The dimension
shift mainly occurs in~$\upsilon$.  We will see that it divides degrees
by~$3$, roughly speaking.  The map~$g$ only looses finitely many degrees.
That follows from the way~$g$ is constructed: The isomorphisms in
Section~\ref{sec:excision_isomorphisms} loose only finitely many $d$'s and~$g$
is constructed from them in finitely many steps.  Thus~$g$ can loose at most
finitely many $d$'s.

To estimate the behavior of~$\upsilon$, we need the \emph{modified external
  degree}~$\#'_e$ defined by
$$
\#'_e(x_0 \, Dx_1 \dots \, Dx_{2n}) \defeq 3n,
\qquad
\#'_e(Dx_1 \dots\, Dx_{2n}) \defeq 3n-2
$$
and the \emph{modified total degree} $\#'_t \defeq \#'_e + \#_i$.  The
total degrees $\#_t$ and~$\#'_t$ are related by a factor of~$3$.  More
precisely, we have the inequalities
$$
\#_t \le
\#'_t \le
3\#_e + \#_i \le
3\#_t.
$$
We will see that $\#'_t \upsilon(l) \ge \# l$.  That is, if~$l$ is
homogeneous of degree~$n$, then $\upsilon(l)$ is a sum of homogeneous terms of
degrees $\#'_t \ge n$.  We need a slightly more general assertion in order to
explore some cancellation.  This will be facilitated if we extend the
action~$\triangleright$ to $\II \circledcirc \Unse{(\Tanil\Tanil\MA[E])}
\supset \Tanil\LL$.  The latter space is equal to $\II \, (D\Tanil\MA[E])
^\even \oplus D\II (D\Tanil\MA[E])^\odd$.

The definition of the action~$\triangleright$ in \eqref{eq:triOpen}
and~\eqref{eq:triClosed} uses the map $\alpha \colon \LL \to \LL \,D\LL$
defined by $\alpha(\vs[g] \odot \opt{x}) = \vs[g] \,D\opt{x}$ if $\vs[g] \in
\VS[G]$ and $\opt{x} \in \Unse{\LL}$.  This formula still makes sense for
$\opt{x} \in \Unse{(\Tanil\MA[E])}$ and defines $\alpha \colon \II \to \VS[G]
\,D\Tanil\MA[E] \subset \II \,D\Tanil\MA[E]$ because $\II \cong \VS[G] \hot
\Unse{(\Tanil\MA[E])}$ by~\eqref{eq:mu12}.  Using this extension of~$\alpha$,
we define $\triangleright \colon \MA[E] \times (\II \circledcirc \Unse{(\Tanil
  \Tanil \MA[E])}) \to \II \circledcirc \Unse{(\Tanil \Tanil \MA[E])}$ by the
same formulas \eqref{eq:triOpen} and~\eqref{eq:triClosed}.

The curvature of this extended action still satisfies \eqref{eq:OmtriOpen}
and~\eqref{eq:OmtriClosed}.  And $\triangleright(x) \in \Endo( \II \circledcirc
\Unse{(\Tanil\Tanil\MA[E])})$ is a right \Mpn{\Tanil\Tanil\MA[E]}module
homomorphism for all $x \in \Tanil\MA[E]$.  The proof carries over without
change.

\begin{lemma}  \label{lem:upsilon_estimate}
  We have $\#'_t(x \triangleright X) \ge \# x + \#'_t(X)$ for all homogeneous
  $x \in \Unse{(\Tanil\MA[E])}$ and $X \in \II \circledcirc
  \Unse{(\Tanil\Tanil\MA[E])}$.  Consequently, $\#'_t \upsilon(x) \ge \# x$
  for all $x \in \LL$.
\end{lemma}

\begin{proof}
  We may assume that~$X$ is a homogeneous monomial and that $x \in \MA[E] \cup
  d\MA[E] d\MA[E]$.  More general monomials $\opt{\ma[e]_0} d\ma[e]_1 \dots
  d\ma[e]_{2n}$ can be treated by induction on~$n$.  Thus we have to verify
  $\#'_t (\ma[e] \triangleright X) \ge \#'_t (X)$ for all $ \ma[e] \in \MA[E]$
  and $\#'_t \omega_\triangleright (\ma[e]_1, \ma[e]_2) (X) \ge \#'_t X + 1$
  for all $\ma[e]_1, \ma[e]_2 \in \MA[E]$.  These assertions follow easily by
  inspecting the summands in \eqref{eq:triOpen}, \eqref{eq:triClosed},
  \eqref{eq:OmtriOpen} and~\eqref{eq:OmtriClosed} and estimating their
  degrees.  The crucial observation is that the map $\alpha \colon \II \to \II
  \,D\Tanil\MA[E]$ may decrease the internal degree by at most~$1$.  That is,
  $\alpha(x)$ is a linear combination of homogeneous monomials $x_0 \,Dx_1$
  with $\# x_0 + \# x_1 \ge \# x - 1$ for all $x \in \II$.  This follows from
  the construction of~$\mu_{\ref{eq:mu8}}^{-1}$ in the proof of
  Lemma~\ref{lem:mu68}.
  
  To conclude that $\#'_t \upsilon(x) \ge \# x$ for all $x \in \LL$, it
  suffices to consider the case $x \in d\MA[E] d\MA[E]$.  This implies the
  assertion for monomials of higher degree:
  $$
  \#'_t \upsilon( \opt{\ma[e]_0} d\ma[e]_1 \dots d\ma[e]_{2n-1}
  d\ma[k]_{2n}) =
  \#'_t \bigl(\opt{\ma[e]_0} d\ma[e]_1 \dots d\ma[e]_{2n-2}
  \triangleright \upsilon(d\ma[e]_{2n-1} d\ma[k]_{2n}) \bigr) \ge
  n-1 + \#'_t \upsilon (d\ma[e]_{2n-1} d\ma[k]_{2n}).
  $$
  It is easy to compute $\upsilon( ds\ma[q]_1 d\ma[k]_2) = ds\ma[q]_1
  d\ma[k]_2$ and $\upsilon( d\ma[k]_1 d\ma[k]_2) = d\ma[k]_1 d\ma[k]_2 +
  D\ma[k]_1 \,D\ma[k]_2$.  These terms have modified total degree at
  least~$1$.
\end{proof}

It is now straightforward to prove a ``poor man's version'' of
Theorem~\ref{the:excision_cyclic} asserting only that $f$ maps $F_{3n +
  \const} (\MA[E]: \MA[Q])$ into $F_n(\MA[K])$ and~$h$ maps $F_{3n + \const}
(\MA[E] : \MA[Q])$ into $F_n(\MA[K])$ for some constant $\const$.  The
factor~$3$ occurs when the degrees $\#_t$ and~$\#'_t$ on~$\Tanil\LL$ are
compared.  The additional constant accommodates the degree shift in~$g$ and
the map~$h_2$.

We use the following linear subspaces of $\Omega^1(\LL)$ to describe the image
$g \bigl( F_n( \Tanil\MA[E]: \Tanil\MA[Q]) \bigr) \subset X(\LL)$.
\begin{align*}
  \VS_{p,q} &\defeq
  \cllin {} \{ \opt{l_0} \,Dl_1 \mid \# l_0 \ge p,\ \# l_1 \ge q \};
  \\
  \VS_{p,q}^{\mathrm{i}} &\defeq
  \cllin {} \{ \opt{x_0} \odot (\ma[k]_1 s\ma[q]_2) \,Dl_1
  - \opt{x_0}\odot \ma[k]_1 \,D(s\ma[q]_2 \odot l_1)
  \mid \# \opt{x_0} \ge p,\ \# l_1 \ge q \};
  \\
  \VS_{p,q}^{\mathrm{o}} &\defeq
  \cllin {} \bigl\{ l_0 \,D\bigl( \opt{x_1} \odot (\ma[k]_2 s\ma[q]_3) \bigr)
  - s\ma[q]_3 \odot l_0 \,D( \opt{x_1} \odot \ma[k]_2)
  \bigm| \# l_0 \ge p,\ \# \opt{x_1} \ge q \bigr\};
  \\
  \VS_{p,q}^{\mathrm{io}} &\defeq
  \cllin {} \bigl\{
  \opt{x_0} \odot (\ma[k]_1 s\ma[q]_2)
  \,D\bigl(\opt{x_3} \odot (\ma[k]_4 s\ma[q]_5) \bigr) -
  \opt{x_0} \odot \ma[k]_1
  \,D\bigl(s\ma[q]_2 \odot \opt{x_3} \odot (\ma[k]_4 s\ma[q]_5) \bigr)
  \\ & \qquad {} -
  s\ma[q]_5\odot \opt{x_0} \odot (\ma[k]_1 s\ma[q]_2)
  \,D\bigl(\opt{x_3} \odot \ma[k]_4 \bigr) +
  s\ma[q]_5\odot \opt{x_0} \odot \ma[k]_1
  \,D\bigl(s\ma[q]_2 \odot \opt{x_3} \odot \ma[k]_4 \bigr)
  \\ & \qquad \qquad {}
  \bigm| \# \opt{x_0} \ge p,\ \# \opt{x_1} \ge q \bigr\}.
\end{align*}
We have $\VS_{p,q}^{\mathrm{io}} \subset \VS_{p,q}^{\mathrm{i}} \subset
\VS_{p,q}$, $\VS_{p,q}^{\mathrm{io}} \subset \VS_{p,q}^{\mathrm{o}} \subset
\VS_{p,q}$, and $\VS_{p,q}^{\cdots} \subset \VS_{p',q'}^{\cdots}$ if $p \ge
p'$ and $q \ge q'$.

\begin{lemma}  \label{lem:g_estimate}
  We have $g( \II \cap \Janil \MA[E]^{n}) \subset \LL \cap \Janil
  \MA[E]^{n-1}$ and
  \begin{displaymath}
    g\bigl( F_{2n}^\odd (\MA[E]: \MA[Q]) \bigr) \subset
    \sum_{p+q=n} \VS_{p,q} +
    \sum_{p+q=n-1} \VS_{p,q}^{\mathrm{i}} +
    \sum_{p+q=n-1} \VS_{p,q}^{\mathrm{o}} +
    \sum_{p+q=n-2} \VS_{p,q}^{\mathrm{io}} +
    \VS_{n-1,0} \bmod[,];
  \end{displaymath}
  \begin{multline*}
    g\bigl( F_{2n-1}^\odd (\MA[E]: \MA[Q]) \bigr) \subset
    \sum_{p+q=n} \VS_{p,q} +
    \sum_{p+q=n-1} \VS_{p,q}^{\mathrm{i}} +
    \sum_{p+q=n-1} \VS_{p,q}^{\mathrm{o}}
    \\ +
    \sum_{p+q=n-2} \VS_{p,q}^{\mathrm{io}} +
    \VS_{n-1,0} +
    \VS_{0,n-1} \bmod[,].
  \end{multline*}
\end{lemma}

\begin{proof}
  We have the following easy recipe to compute~$g$.  If $y \in
  X_0(\Tanil\MA[E]: \Tanil\MA[Q]) = \II$, we subtract from~$y$ an element
  of~$\VS[C]'_0 = [\LL, s_L(\Tanil\MA[Q])]$ such that the result is in~$\LL$.
  If $y \in \Omega^1 (\Tanil\MA[E]: \Tanil\MA[Q])$, we subtract from~$y$
  elements of $\VS[C]'_1 = \LL \,Ds_L(\Tanil\MA[Q])$ and $b\bigl( \Omega^2
  (\Tanil\MA[E]: \Tanil \MA[Q]) \bigr)$ until the result is in $\Unse{\LL}
  \,D\LL$.
  
  On the even part $X_0(\Tanil\MA[E]: \Tanil\MA[Q]) = \II$, we have
  $$
  g\bigl( l \odot s_L\opt{q} \bigr) =
  g\bigl( [l, s_L\opt{q}] + s_L\opt{q} \odot l \bigr) =
  s_L\opt{q} \odot l.
  $$
  This suffices to determine the even part of~$g$ by~\eqref{eq:mu11}.
  We decompose the odd part $X_1(\Tanil\MA[E]: \Tanil\MA[Q])$ into the direct
  sum $\II\, Ds(\MA[Q]) \oplus \Unse{( \Tanil\MA[E] )} \,D\MA[K]$.  The
  summand $\Unse{( \Tanil\MA[E] )} \,D\MA[K]$ is easy to handle.  We have
  \begin{equation}  \label{eq:g_odd_K}
    g( \opt{l_0} \odot s_L\opt{q_1} \,D\ma[k] ) =
    g\bigl( \opt{l_0} \,D( s_L\opt{q_1} \odot \ma[k]) \bigr)
    - g( \opt{l_0} \,(Ds_L\opt{q_1})  \ma[k] ) =
    \opt{l_0} \,D(s_L\opt{q_1} \odot \ma[k]).
  \end{equation}
  The isomorphism~\eqref{eq:mu4} shows that this determines~$g$ on
  $\Unse{(\Tanil\MA[E])} \,D\MA[K]$.
  The summand $\II\, Ds(\MA[Q])$ requires two steps.  As above, we have
  $g\bigl( l_0 \odot s_L\opt{q_1} \,Ds(\ma[q]) \bigr) = g\bigl( l_0 \,
  D(s_L\opt{q_1} \odot s\ma[q]) \bigr)$.  Bringing $s_L\opt{q_1} \odot
  s\ma[q]$ into standard form, we get terms in $s_L(\Tanil\MA[Q])$ and terms
  in~$\II$ involving the curvature of~$s$.  Since~$g$ vanishes on $\LL
  \,Ds_L(\Tanil\MA[Q])$, the summands in $s_L(\Tanil\MA[Q])$ do not
  contribute.  We apply~$\mu_{\ref{eq:mu11}}$ to each curvature term to get an
  element of $\LL \otimes s_L\Unse{(\Tanil\MA[Q])}$.  Then
  \begin{displaymath}
    g\bigl( l_0 \,D(l_1 \odot s_L\opt{q_2} ) \bigr) =
    s_L\opt{q_2} \odot l_0 \,Dl_1.
  \end{displaymath}
  This procedure determines~$g$ on $\II \,Ds(\MA[Q])$ because
  of~\eqref{eq:mu11}.
  
  Thus to obtain~$g$, we have to apply the isomorphism
  $\mu_{\ref{eq:mu4}}^{-1}$ at most once on~$\II$ and twice on
  $X_1(\Tanil\MA[E] : \Tanil\MA[Q])$.  Application of
  $\mu_{\ref{eq:mu4}}^{-1}$ either leaves the internal degree unchanged or
  looses one internal degree.  In the latter case, we split a monomial
  $\omega_0 d\ma[k]_1 ds\ma[q]_2 \odot \omega_3$ into $\omega_0 \odot
  (\ma[k]_1 s\ma[q]_2) \otimes \omega_3 - (\omega_0 \odot \ma[k]_1) \otimes
  (s\ma[q]_2 \odot \omega_3)$, producing the kind of expression that occurs in
  the definition of $\VS_{p,q}^{\mathrm{\cdots}}$.
  
  The assertion of the lemma follows in a straightforward way.  To prove the
  assertion about $F_{2n-1}^\odd (\MA[E] : \MA[Q])$, we also have to consider
  expressions of the form $\omega \,D(d\ma[e]_1 d\ma[e]_2)$ with either
  $\omega \in \II$ or $d\ma[e]_1 d\ma[e]_2 \in \II$
  (Lemma~\ref{lem:Fodd_alternative}).  The details are rather boring and
  therefore omitted.
\end{proof}

\begin{claim}  \label{claim:f_even}
  The even part of~$f$ maps $\II \cap \Janil\MA[E]^{3n-1}$ into
  $\Janil\MA[K]^n$.
\end{claim}

\begin{proof}
  The even part of~$g$ maps $\II \cap \Janil\MA[E]^{3n-1}$ to $\LL \cap
  \Janil\MA[E]^{3n-2}$ by Lemma~\ref{lem:g_estimate}.  This space is mapped
  by~$\upsilon$ to $\cllin {} \{ X \in \Tanil\LL \mid \#'_t X \ge 3n-2 \}$ by
  Lemma~\ref{lem:upsilon_estimate}.  Since $3\#_t \ge \#'_t X$ and~$\#_t$ is
  integer valued, we get only terms with $\#_t \ge n$.  These are mapped by
  $\Tanil\tilde\tau$ into $\Janil\MA[K]^n$ because $\# \Tanil\tilde\tau (X) =
  \#_t X$.
\end{proof}

\begin{claim}  \label{claim:f_odd_I}
  The odd part of $X(\Tanil\tilde\tau) \circ X(\upsilon)$ maps~$\VS_{p,q}$
  with $p \ge 3p' - 2$ and $q \ge 3q' - 2$ into $F_{2(p'+q')-1}^\odd
  (\MA[K])$.
\end{claim}

\begin{proof}
  We get as above that $X(\Tanil\tilde\tau \circ \upsilon)$ maps $\opt{l_0}
  \,Dl_1$ with $\# \opt{l_0} \ge 3p'-2$ and $\# l_1 \ge 3q'-2$ to
  $$
  \cllin {} \{ \opt{y_0} \,Dy_1 \in \Omega^1(\Tanil\MA[K]) \mid
  \# \opt{y_0} \ge p',\ \# y_1 \ge q' \} \bmod [,].
  $$
  (We first take the completant linear hull, then map to the commutator
  quotient.)  This is contained in $F_{2(p'+q')-1}^\odd (\MA[K])$ by
  Lemma~\ref{lem:Fodd_alternative}.
\end{proof}

The previous claim applies in particular to $\VS_{p,q}$ with $p+q \ge 3n-2$.
Hence the two claims above combined with Lemma~\ref{lem:g_estimate} show
that~$f$ maps $F_{6n-1} (\MA[E] : \MA[Q])$ into $F_{2n-1} (\MA[E] : \MA[Q])$.
Recall that $\VS_{p,q}^{\cdots} \subset \VS_{p,q}$ for $\cdots = \mathrm{i},
\mathrm{o}, \mathrm{io}$.  When computing $f\bigl( F_{6n+2}(\MA[E]: \MA[Q])
\bigr)$, we can replace the even part by $\II \cap \Janil\MA[E]^{3n+2}$
because of Lemma~\ref{lem:Feven_alternative} and because~$f$ is a chain map.
Claim~\ref{claim:f_even} shows that $f(\II \cap \Janil\MA[E]^{3n+2}) \subset
\Janil\MA[K]^{n+1} \subset F_{2n}^\even (\MA[K])$ as desired.

The odd part of $F_{6n+2} (\MA[E] : \MA[Q])$ requires more attention.  We can
leave out all summands $\VS_{p,q}^{\cdots}$ with $p \ge 3p' - 2$, $q \ge 3q' -
2$ and $p' + q' \ge n+1$ because these are mapped into $F_{2n+1}^\odd
(\MA[K])$ by Claim~\ref{claim:f_odd_I}.  This leaves the summands
$\VS^{\mathrm{i}}_{3p',3(n-p')}$, $\VS^{\mathrm{o}}_{3p',3(n-p')}$,
$\VS^{\mathrm{io}}_{3p'-1,3(n-p')}$, $\VS^{\mathrm{io}}_{3p',3(n-p')-1}$, and
$\VS_{3n,0}$ in Lemma~\ref{lem:g_estimate} to be considered.

The summand $\VS_{3n,0}$ is mapped into $\cllin {} \{ \opt{y_0} \,Dy_1 \mid \#
\opt{y_0} \ge n \} \bmod [,] \subset F_{2n} (\MA[K])$ and therefore
unproblematic.  The remaining summands are more tricky.  First of all, they
are all contained in $\VS_{3p'-2, 3(n-p')-1}$.  If $\opt{l_0} \,Dl_1 \in
\VS_{3p'-2, 3(n-p') - 1}$, then $\upsilon\opt{l_0} \delta \upsilon(l_1)$ can
be written as a sum of homogeneous monomials $\opt{X_0} \delta X_1$ satisfying
$\#'_t \opt{X_0} \ge 3p'-2$ and $\#'_t X_1 \ge 3(n-p') - 1$.  Hence
$\#_t \opt{X_0} \ge p'$ and either $\#_t X_1 \ge n+1-p'$ or $\#_e X_1 = n -
p'$ and $\#_i X_1 = 0$ because $3\#_e + \#_i \ge \#'_t$.

We need an explicit formula for the image of $\opt{X_0} \delta X_1 \bmod [,]$
in $\Omega^\odd \LL$ due to Cuntz and Quillen
\cite[Lemma~5.3]{cuntz95:cyclic}:
\begin{equation}
  \label{eq:xDy_explicit}
  \opt{X_0} \,\delta X_1 \bmod [,] =
  - \sum_{j=0}^{k-1} \kappa^{2j} b (\opt{X_0} \circledcirc X_1)
  + \sum_{j=0}^{2k-1} \kappa^j D(\opt{X_0} \odot X_1)
  + \kappa^{2k}(\opt{X_0} \odot DX_1)
\end{equation}
for all $X_1 \in \Omega^{2k}\LL$, $\opt{X_0} \in \Unse{(\Tanil\LL)}$.  The
operators~$b$ and~$\kappa$ are the usual operators on $\Omega\LL$.

If $\#_t \opt{X_0} + \#_t X_1 \ge n$, then the summands
in~\eqref{eq:xDy_explicit} containing~$D$ are harmless: They are mapped by
$\Omega_\an(\tilde\tau)$ to $\Janil\MA[K]^{n} \,d\MA[K] \subset F_{2n}^\odd
(\MA[K])$.  If even $\#_t \opt{X_0} + \#_t X_1 \ge n+1$, then the first sum
involving~$b$ is harmless as well.  We only get problems if $\#_t \opt{X_0} =
p'$ and $\#_e X_1 = n-p'$, $\#_i X_1 = 0$.  Therefore, we can replace
$X(\upsilon) (\opt{l_0} \,D l_1)$ by
$$
\sum_{j=0}^{n-p'-1} \kappa^{2j} b(\opt{X_0} \circledcirc X_1) =
\sum_{j=0}^{n-p'-1} \kappa^{2j} b\bigl( \upsilon(\opt{l_0} \odot l_1) \bigr).
$$
The map $\opt{l_0} \,Dl_1 \mapsto \opt{l_0} \odot l_1$ maps
$\VS^{\mathrm{i}}_{3p',3(n-p')} + \VS^{\mathrm{o}}_{3p',3(n-p')} +
\VS^{\mathrm{io}}_{3p'-1,3(n-p')} + \VS^{\mathrm{io}}_{3p', 3(n-p')-1}$ into
$\LL \cap \Janil\MA[E]^{3n+1} + \VS[W]_{3n}$ with
$$
\VS[W]_{k} \defeq
\cllin {} \{ \opt{x} \odot (\ma[k] s\ma[q]) -  s\ma[q] \odot \opt{x}
  \odot \ma[k]
\mid \# \opt{x} \ge k \}.
$$
Thus it remains to consider $b \circ \upsilon (\VS[W]_{3n} + \LL \cap
\Janil\MA[E]^{3n+1})$.

\begin{claim}  \label{claim:upsilon_W}
  The homomorphism~$\upsilon$ maps $\VS[W]_{3n} + \LL \cap
  \Janil\MA[E]^{3n+1}$ into
  $$
  \cllin {} \{ X \in \Tanil\Tanil\MA[E] \mid \#_t X \ge n+1 \} +
  b \bigl(\, \cllin {} \{ \opt{X} \,Ds(\ma[q]) \mid \#_i X = 0,\ \#_e X = n \}
  \bigr).
  $$
\end{claim}

\begin{proof}
  The terms in $\LL \cap \Janil\MA[E]^{n+1}$ are treated by
  Lemma~\ref{lem:upsilon_estimate} and the estimate $3\#_t \ge \#'_t$.
  To handle the summand~$\VS[W]_{3n}$ we need the extension
  of~$\triangleright$ to $\II \circledcirc \Unse{(\Tanil\Tanil\MA[E])}$.

  Pick a generator $l \defeq \opt{x} \odot (\ma[k]s\ma[q]) - s\ma[q] \odot
  \opt{x} \odot \ma[k]$ of~$\VS[W]_{3n}$ with $\opt{x} \ge 3n$.  We compute
  \begin{multline*}
    \upsilon(l) =
    \opt{x} \triangleright (\ma[k] s\ma[q])
    - s\ma[q] \triangleright \opt{x} \triangleright \ma[k]
    =
    \opt{x} \triangleright
    (d\ma[k] ds\ma[q] + D\ma[k] \,Ds\ma[q] + \ma[k] \circledcirc s\ma[q])
    - s\ma[q] \triangleright \opt{x} \triangleright \ma[k]
    \\ =
    \opt{x} \triangleright (d\ma[k] ds\ma[q] + D\ma[k] \,Ds\ma[q])
    + (\opt{x} \triangleright \ma[k]) \circledcirc s\ma[q]
    - s\ma[q] \triangleright \opt{x} \triangleright \ma[k]
    \\ =
    \opt{x} \triangleright (d\ma[k] ds\ma[q] + D\ma[k] \,Ds\ma[q])
    - D(\opt{x} \triangleright \ma[k]) \,Ds\ma[q]
    - b\bigl( (\opt{x} \triangleright \ma[k]) \,Ds\ma[q] \bigr)
    + s\ma[q] \odot (\opt{x} \triangleright \ma[k])
    - s\ma[q] \triangleright (\opt{x} \triangleright \ma[k]).
  \end{multline*}
  Lemma~\ref{lem:upsilon_estimate} implies that $\opt{x} \triangleright
  (d\ma[k] ds\ma[q] + D\ma[k] \,Ds\ma[q])$ and $D(\opt{x} \triangleright
  \ma[k]) \,Ds\ma[q]$ have total degree $\ge n+1$.  Hence the first two
  summands above are of the required form.
  
  Write $\opt{x} \triangleright \ma[k]$ as a sum of homogeneous monomials
  $\sum_j X_j$.  Since $\#'_t (\opt{x} \triangleright \ma[k]) \ge \# \opt{x}
  \ge 3n$, we have $\#_t X_j \ge n+1$ or $\#_e X_j = n$ and $\#_i X_j = 0$.
  If $\#_t X_j \ge n+1$, then $b(X_j \,Ds\ma[q])$ has total degree $\ge n+1$.
  The second summand in Claim~\ref{claim:upsilon_W} is there to accommodate
  $b(X_j \,Ds\ma[q])$ with $\#_e X_j = n$ and $\#_i X_j = 0$.
  
  We even have $\#_t X_j \ge n+2$ or $\#_e X_j \ge n+1$ or $X_j = l_j \,DX_j'$
  with $\#_e DX_j' = n$.  These estimates are needed to handle $s\ma[q] \odot
  (\opt{x} \triangleright \ma[k]) - s\ma[q] \triangleright (\opt{x}
  \triangleright \ma[k])$.  The first two cases are harmless because
  $\triangleright (s\ma[q])$ decreases total degrees $\#_t$ by at most~$1$ and
  does not decrease the external degree.  That is, $\#_t X_j \ge n+2$ implies
  $\#_t s\ma[q] \odot X_j \ge n+1$ and $\#_t s\ma[q] \triangleright X_j \ge
  n+1$, and $\#_e X_j \ge n+1$ implies $\#_t s\ma[q] \triangleright X_j \ge
  \#_e s\ma[q] \triangleright X_j \ge n+1$ and $\#_t s\ma[q] \odot X_j \ge
  n+1$.  The definition of~$\triangleright$ in~\eqref{eq:triOpen} implies
  $$
  s\ma[q] \odot l_j \,DX_j' - s\ma[q] \triangleright l_j \,DX_j' =
  D \alpha(s\ma[q] \odot l_j) \,DX_j'.
  $$
  This tackles the third case because $\#_e D\alpha (s\ma[q] \odot l_j)
  \,DX_j' \ge \#_e DX_j' + 1 \ge n+1$.
\end{proof}

Claim~\ref{claim:upsilon_W} implies that $b \circ \upsilon (\VS[W]_{3n} + \LL
\cap \Janil\MA[E]^{3n+1})$ is contained in $\cllin {} \{ \opt{X_0} \delta l_1
\mid \#_t \opt{X_0} + \# l_1 \ge n \}$ because $b \circ b = 0$.  Hence
$X(\Tanil\tilde\tau)$ maps $b \circ \upsilon (\VS[W]_{3n} + \Janil\MA[E]^{n+1}
\cap \LL)$ into $F_{2n}^{\odd} (\MA[K])$ as desired.  This completes the proof
that $f\bigl( F_{3n+2} (\MA[E]: \MA[Q]) \bigr) \subset F_{n}(\MA[K])$ for all
$n \in \Z_+$.

It remains to consider the operator~$h$.  The part~$h_1$ is easy:

\begin{claim}  \label{claim:h_I}
  We have $h_1 \bigl( \Omega^{2n} (\MA[E]: \MA[Q]) \bigr) \subset \Omega^{\ge
    2n-1} (\MA[E]: \MA[Q])$.
\end{claim}

\begin{proof}
  We have
  $$
  h_1(\omega \cdot (d\ma[e]_1 d\ma[k]_2) ds\ma[q]_3 \dots ds\ma[q]_{2n}) =
  \omega \cdot (d\ma[e]_1 d\ma[k]_2) \, D(ds\ma[q]_3 \dots ds\ma[q]_{2n})
  \bmod [,].
  $$
  Rewriting $D(\cdots)$ according to~\eqref{eq:D_long_form}, we only get terms
  with at least $2n-1$ $d$'s.  Moreover,
  \begin{multline*}
    h_1(\omega \cdot (d\ma[k]_1 ds\ma[q]_2) ds\ma[q]_3 \dots ds\ma[q]_{2n}) =
    \omega \odot (\ma[k]_1 s\ma[q]_2) \, D(ds\ma[q]_3 \dots ds\ma[q]_{2n}) -
    \omega \odot \ma[k]_1 \, D(s\ma[q]_2 ds\ma[q]_3 \dots ds\ma[q]_{2n})
    \\ =
    \omega \odot (\ma[k]_1 s\ma[q]_2 - \ma[k]_1 \odot s\ma[q]_2)
    \, D(ds\ma[q]_3 \dots ds\ma[q]_{2n}) -
    \omega \odot \ma[k]_1 \, (Ds\ma[q]_2)\, ds\ma[q]_3 \dots ds\ma[q]_{2n}
    \bmod [,]
    \\ =
    \omega \cdot (d\ma[k]_1 ds\ma[q]_2) \, D(ds\ma[q]_3 \dots ds\ma[q]_{2n}) -
    ds\ma[q]_3 \dots ds\ma[q]_{2n} \odot \omega \odot \ma[k]_1 \, Ds\ma[q]_2
    \bmod [,].
  \end{multline*}
  Both summands on the right hand side are in $\Omega^{\ge 2n-1} (\MA[E]:
  \MA[Q])$ by~\eqref{eq:D_long_form}.
\end{proof}

The part $\psi \circ h_2 \circ X(\upsilon) \circ g$ is quite complicated to
handle.  In a couple of places, we need some trick to conclude that the worst
possible terms drop out or are at least in the range of~$b$.  What happens
here is essentially the following.  Homotopy invariance does not hold in
cyclic cohomology but it holds after stabilizing once with the
\Mpn{S}operator.  Thus we expect $\psi \circ h_2$ to loose at most one
filtration level (if we cared to define a Hodge filtration on $X(\Tanil\LL)$).
This is indeed the case, but at first sight it may appear that $\psi \circ
h_2$ looses two filtration levels.  Therefore, we will only sketch the ideas
and tricks necessary and leave the detailed computations to the interested
reader.

The restriction of~$h_2$ to odd degrees is most easy to write down.  The
connection~$\nabla$ satisfies $\nabla(\opt{X_0} \delta l_1) = 0$ for all $l_1
\in \LL$.  Hence $h_2( \opt{X_0} \delta l_1) = \eta(\opt{X_0} \delta l_1) =
\int_0^1 H_t\opt{X_0} \odot \dot{H}_t(l_1) \,dt$.  Here~$\dot{H}$ denotes the
derivative of~$H$.  Since both~$H$ and~$\dot{H}$ are homogeneous with respect
to $\#_t$ and~$\#$, it follows that $\psi \circ h_2$ maps $\opt{X_0} \delta
l_1$ into $F_{2n-1}^\even (\MA[E] : \MA[Q])$ if $\#_t \opt{X_0} + \# l_1 \ge
n$.  Using~\eqref{eq:xDy_explicit}, it follows that terms of the form
$\opt{X_0} \delta X_1$ with $\#_t \opt{X_0} + \#_t X_1 \ge n+1$ are also
mapped to $F_{2n-1}^\even (\MA[E] : \MA[Q])$.

In addition, if $\# l_1 = 0$, then $h_2(\opt{X_0} \delta l_1)= 0$ because~$H$
is constant on~$\MA[K]$.  Hence terms of the form $\opt{X_0} \delta X_1$ with
$\#_t \opt{X_0} + \#_t X_1 \ge n$ and $\#_i X_1 = 0$ are mapped by $\psi \circ
h_2$ into $F_{2n-1} (\MA[E] : \MA[Q])$.  This is sufficient to show that $\psi
\circ h_2 \circ X(\upsilon)$ maps most of the summands in $g\bigl(
F_{6n-1}^\odd (\MA[E] : \MA[Q]) \bigr)$ occurring in
Lemma~\ref{lem:g_estimate} into $F_{2n-1}^\even (\MA[E] : \MA[Q])$.  The only
problematic summands are $\VS_{3p', 3(n-p')-2}^{\mathrm{io}}$.  They can be
handled using the following symmetry:

\begin{claim}
  $\VS_{p,q}^{\mathrm{io}} \subset \VS_{0,p+q+1} + \VS_{q,p}^{\mathrm{io}}
  \bmod [,]$.
\end{claim}

\begin{proof}
  Up to commutators, we have
  \begin{multline*}
    l_0 \,D\bigl(\opt{x_1} \odot (\ma[k]_1 s\ma[q]_1) \bigr)
    - s\ma[q]_1 \odot l_0 \,D(\opt{x_1} \odot \ma[k]_1)
    \\ =
    D\bigl(\opt{x_1} \odot (\ma[k]_1 s\ma[q]_1) \bigr) \, l_0
    - D(\opt{x_1} \odot \ma[k]_1) \, s\ma[q]_1 \odot l_0
    \\ =
    D( \opt{x_1} \odot d\ma[k]_1 ds\ma[q]_1 \odot l_0)
    - \opt{x_1} \odot (\ma[k]_1 s\ma[q]_1) Dl_0
    + \opt{x_1} \odot \ma[k]_1) \,D( s\ma[q]_1 \odot l_0).
  \end{multline*}
  Thus $\VS_{p,q}^{\mathrm{o}} \subset \VS_{0,p+q+1} + \VS_{q,p}^{\mathrm{i}} +
  [,]$.  A similar computation proves the claim.
\end{proof}

Thus we can replace $\VS_{3p',3(n-p')-2}^{\mathrm{io}}$ by $\VS_{3(n-p')-2,
  3p'} + \VS_{0,3n-1}$.  It follows that $\psi \circ h_2 \circ X(\upsilon)
\circ g$ maps $F_{6n-1}^\odd (\MA[E] : \MA[Q])$ into $F_{2n-1}^\even (\MA[E] :
\MA[Q])$.

Let $b \colon \Omega_\an \LL \to \Omega_\an \LL$ be the \Hochschild boundary.
We have
\begin{multline}  \label{eq:eta_b}
  \eta \circ b( \opt{X} \,Dl_1\,Dl_2) =
  \int_0^1
  [H_t(l_2), H_t\opt{X} \odot \dot{H}_t(l_1)]
  \\
  + H_t(Dl_2 \,DX) \odot \dot{H}_t(l_1)
  + H_t(DX \,Dl_1) \odot \dot{H}_t(l_2)
  - H_t\opt{X} \odot \dot{H}_t(Dl_1 \,Dl_2) \,dt.
\end{multline}
Thus the terms of lowest degree in $\psi \circ \eta \circ b(\opt{X} \,Dl_1
\,Dl_2)$ are \Mpn{b}boundaries.  It follows that all terms $\opt{X_0} \delta
X_1$ with $\#_t \opt{X_0} + \#_t X_1 \ge n+1$ are mapped by $\psi \circ h_2$
into $F_{2n}^\even(\MA[E] : \MA[Q])$.  Using also Claim~\ref{claim:upsilon_W},
we obtain that $\psi \circ h_2 \circ X(\upsilon) \circ g$ maps $F_{6n+2}^\odd
(\MA[E] : \MA[Q])$ into $F_{2n}^\even (\MA[E] : \MA[Q])$.

The formula for~$h_2$ on the even part is slightly more complicated because it
involves the map $\nabla \circ D \colon \Tanil\LL \to \Unse{(\Tanil\LL)} \hot
\LL \hot \Unse{(\Tanil\LL)}$.  This map is described
by~\eqref{eq:D_long_form}.  We have to compose $\nabla \circ D$ with the map
$\opt{X_0} \otimes l_1 \otimes \opt{X_2} \mapsto \eta(\opt{X_0} \delta l_1
\delta \opt{X_2})$.  Finally, we have to bring the result in
$\Omega^1(\Tanil\LL)$ into standard form by~\eqref{eq:xDy_explicit}.  Going
through this computation, we find that $\psi \circ h_2(X) \in F_{2n-1}^\odd
(\MA[E] : \MA[Q])$ if $\#_t X \ge n+1$.  If $\#_e X = n$, $\#_i X = 0$, then
some terms drop out because $\dot{H}|_{\MA[K]} = 0$, so that still $\psi \circ
h_2(X) \in F_{2n-1}^\odd (\MA[E] : \MA[Q])$.  Using
Lemma~\ref{lem:upsilon_estimate} and Lemma~\ref{lem:g_estimate}, we conclude
that $\psi \circ h_2 \circ X(\upsilon) \circ g$ maps $F_{6n-1}^\even (\MA[E] :
\MA[Q])$ into $F_{2n-1}^\odd (\MA[E]: \MA[Q])$.

It is more difficult to show that $\psi \circ h_2 \circ X(\upsilon) \circ g$
maps $F_{6n+2}^\even (\MA[E] : \MA[Q])$ to $F_{2n}^\even (\MA[E] : \MA[Q])$.
It can be shown as above that $\psi \circ h_2 (X) \in F_{2n}^\odd (\MA[E] :
\MA[Q])$ if $\#_t X \ge n+2$ or $\#_e X \ge n+1$ or ($\#_e X \ge n$ and $\#_i
X \ge 1$).  In the first two cases, $X$ is even mapped into $F_{2n+1}^\odd
(\MA[E] : \MA[Q])$.  This suffices to handle $\II \cap \Janil\MA[E]^{3n+2}
\subset F_{6n+2}^\even (\MA[E] : \MA[Q])$.  The remaining summand $\partial_1
\bigl( \Omega^{6n+3} (\MA[E] : \MA[Q]) \bigr)$ has to be considered this time
because~$h_2$ is not a chain map.  If $x \in \partial_1 \bigl( \Omega^{6n+3}
(\MA[E] : \MA[Q]) \bigr)$, then we split $X(\upsilon) \circ g (x)$ into
homogeneous components~$X_j$.  By Lemma~\ref{lem:upsilon_estimate}, we have
$\#_t X_j \ge n+2$ or $\#_e X_j = n+1$ or $\#_e X_j = n$ and $\#_i X_j \in
\{0,1\}$.  The only problematic possibility is $\#_e X_j = n$ and $\#_i X_j =
0$.  These terms may yield a contribution to $\psi \circ h_2(x)$ in $b\bigl(
\Omega^{2n} (\MA[E] : \MA[Q]) \bigr)$.  However, since $f(x) \in F_{2n}^\even
(\MA[K])$, these lowest order components~$X_j$ must be in $b(\Omega^{2n+1}
\MA[K])$ up to terms of higher degree.  This implies that the possible
contribution in $b\bigl( \Omega^{2n} (\MA[E] : \MA[Q]) \bigr)$ cancels.

This proves the first half of Theorem~\ref{the:excision_cyclic}.  We do not
prove the assertion in Theorem~\ref{the:excision_cyclic} about split
extensions.  Suffice it to indicate how to modify the grading~$\#'_t$ to
improve the estimates about~$\upsilon$.  In the split case, the map~$s_L$ is
multiplicative.  Hence $s(\MA[Q]) \odot \VS[G] \subset \VS[G]$ is contained in
the kernel of~$\alpha$.  Thus some of the worst terms in the formulas for
$\triangleright$ and $\omega_\triangleright$ vanish.  We define another
modified total degree by $\#''_t \defeq 2\#_e + \#_i - \epsilon$, where
$\epsilon( l_0 \,DX) = 0$, $\epsilon( Dg \,DX) = 1$ if $g \in \VS[G]$, and
$\epsilon(Dl_1 \,DX) = 2$ if $l_1 = d\ma[e]_1 \dots d\ma[e]_{2n}
\opt{\ma[e]_{2n+1}}$ and not $l_1 \in \VS[G]$.  It is easy to verify that
$\#''_t (\ma[e] \triangleright X) \ge \#''_t X$ and $\#''_t \bigl(
\omega_\triangleright (\ma[e]_1, \ma[e]_2) (X) \bigr) \ge \#''_t X + 1$.  We
use that if $\triangleright$ or $\omega_\triangleright$ produces terms of the
form $Dl_1 \,DX$, then automatically $l_1 \in \VS[G]$.  Since the external
degree only occurs with a factor of~$2$, we get an overall factor of~$2$
instead of~$3$.  The remaining estimations are quite similar to those above.

\section{\Frechet algebras and fine algebras}
\label{sec:permanence}

Let $\FINE$ be the functor from the category of vector spaces to the category
of complete bornological vector spaces that endows a vector space with the
fine bornology.  Then $\FINE(V) \hot \FINE(W)$ is naturally isomorphic to
$\FINE(V \otimes W)$ with $V \otimes W$ the algebraic tensor product.  Thus
$\Omega^j A$ is the usual uncompleted space of differential forms if~$A$ is a
fine algebra.  Furthermore, the functor $\FINE$ is fully faithful.  It follows
that for algebras $A$ and~$B$ without additional structure, $\HP^\ast
\bigl(\FINE(A); \FINE(B) \bigr) = \HP^\ast(A; B)$.  That is, the
``bornological'' bivariant periodic cyclic cohomology of fine algebras is the
usual theory for algebras without additional structure.  The same applies to
cyclic (co)homology and \Hochschild (co)homology.

Let $\COMP$ be the functor from the category of \Frechet spaces (with
continuous linear maps as morphisms) to the category of complete bornological
vector spaces that endows a \Frechet space with the precompact bornology.  The
bornological tensor product $\COMP(V) \hot \COMP(W)$ is naturally isomorphic
to $\COMP(V \prot W)$ with $V \prot W$ the projective tensor product by
Theorem~\ref{the:tensor_Frechet}.  Moreover, the functor $\COMP$ is fully
faithful.  It follows that for \Frechet algebras $A$ and~$B$, $\HP^\ast \bigl(
\COMP(A); \COMP(B) \bigr) = \HP^\ast_{\mathrm{top}} (A; B)$.  That is, the
``bornological'' periodic cyclic cohomology of \Frechet algebras endowed with
the precompact bornology is the usual topological bivariant periodic cyclic
cohomology.  The same applies to cyclic (co)homology and \Hochschild
(co)homology.

Another special case that deserves attention is the following.  Let~$\GR$ be a
separated smooth groupoid.  Let $\CCINF(\GR)$ be the algebra of smooth
compactly supported functions on~$\GR$ with convolution as multiplication and
the topology of locally uniform convergence.  The convolution product is not
jointly continuous unless~$\GR$ is compact or discrete and countable (see
Example~\ref{exa:R_convolution}).  Thus $\CCINF(\GR)$ usually is not a
topological algebra.  However, the convolution is obviously separately
continuous.  Hence $\CCINF(\GR)$ is a complete bornological algebra when
endowed with the precompact bornology (which equals the bounded bornology in
this case).  The bornological tensor products $\CCINF(\GR)^{\hot n}$ are
naturally isomorphic to $\CCINF(\GR^n)$ with the precompact bornology by
Example~\ref{exa:tensor_smooth_manifold}.

This is the completed tensor product used by J.-L.\ Brylinski and Victor
Nistor in~\cite{brylinski94:cyclic} for smooth \'etale groupoids.  Hence they
compute the bornological cyclic and \Hochschild homology of $\CCINF(\GR)$.
The results they obtain are related to the cohomology (in the sense of
algebraic topology) of suitable classifying spaces.

However, neither the purely algebraic theory nor the topological theory are
satisfactory.  The topological theory does not apply at all because the
multiplication is not jointly continuous and thus the differential~$b$ is not
continuous with respect to the projective tensor product topology.  The purely
algebraic theory is applicable but apparently not easy to compute.  Already
the case of the algebra $\CINF([0,1])$ appears to be intractable.  The algebra
$\CINF([0,1])$ is the smooth convolution algebra of the manifold $[0,1]$
viewed as an \'etale groupoid.  The problem is that $\CINF([0,1])$ with the
fine bornology is no longer smoothly homotopy equivalent to~$\C$.

Unfortunately, the computations in~\cite{brylinski94:cyclic} give no
information about the cyclic cohomology of $\CCINF(\GR)$ and about bivariant
groups.  The problem is that in~\cite{brylinski94:cyclic}, the authors
complete $\CCINF(\GR)^{\otimes n}$ to $\CCINF(\GR^n)$ and after that work
purely algebraically.  For the homology of a complex, the bornology does not
matter.  It matters, however, for the cohomology.  It is not obvious how to
rearrange the computations in~\cite{brylinski94:cyclic} so as to respect the
bornologies on the complexes.  Already the very first step, reduction to
loops, does not work because the contractible complex that is used
in~\cite{brylinski94:cyclic} is not contractible in a bounded way, its closure
being no longer contractible.  It is quite plausible that a different argument
can be used to relate the bornological cyclic cohomology to suitable homology
groups of classifying spaces.

Since bornological algebras contain both fine algebras and \Frechet algebras
as special cases, they make it easier to study the interplay between an
algebra and dense subalgebras.  For example, we can compare the algebras
$\CINF(S^1)$ of smooth functions on the circle and its dense, fine subalgebra
$\C[u, u^{-1}]$ of Laurent polynomials and assert that they are
\Mpn{\HP}equivalent.  That is, the class of the inclusion homomorphism is an
invertible element in $\HP^0\bigl(\C[u, u^{-1}]; \CINF(S^1)\bigr)$.  It occurs
quite frequently in examples that a (suitably small) \Frechet algebra has the
same cyclic cohomology as a suitable fine, finitely presented, dense
subalgebra.

\begin{appendix}

\chapter[Appendices]{}

\section{Bornologies and inductive systems}
\label{app:bornologies}

\subsection{Functorial constructions with bornological vector spaces}
\label{app:bornology_construct}

As announced in Section~\ref{sec:bornology_construct}, we give here the
definitions of direct products, direct sums, projective limits, and inductive
limits.  See~\cite{Hogbe-Nlend:Born} for more details.

Let $(\VS_i)_{i\in I}$ be a projective system of bornological vector spaces.
The \emph{projective limit} of~$(\VS_i)_{i\in I}$ is constructed as follows.
As a vector space, $\prolim {(\VS_i)_{i\in I}}$ is equal to the ordinary
projective limit in the category of vector spaces.  There are natural maps
$\pi_i\colon \prolim {(\VS_i)_{i\in I}} \to \VS_i$ for all $i\in I$.  The
bornology on $\prolim {(\VS_i)_{i\in I}}$ is taken to be the \emph{coarsest}
bornology making all maps~$\pi_i$ bounded.  That is, $S\subset \prolim
{(\VS_i)_{i\in I}}$ is small iff $\pi_i(S) \subset \VS_i$ is small for all
$i\in I$.  If all the spaces~$\VS_i$ are separated or complete, so is $\prolim
{(\VS_i)_{i\in I}}$.  The \emph{direct product} $\prod_{i\in I} \VS_i$ of a
set of bornological vector spaces is constructed similarly: Endow the vector
space direct product with the coarsest bornology making all coordinate
projections bounded.

\begin{example}
  The functors $\BOUND$ and $\COMP$ from locally convex topological vector
  spaces to convex bornological vector spaces are compatible with projective
  limits and direct products.  That is, if $(V_i)_{i\in I}$ is a projective
  system of topological vector spaces, then $\prolim V_i \cong \prolim
  {(V_i,\BOUND)} \cong \prolim {(V_i,\COMP)}$ as vector spaces; a subset of
  $\prolim V_i$ is bounded or precompact iff it is small in $\prolim
  {(V_i,\BOUND)}$ or $\prolim {(V_i,\COMP)}$, respectively.
\end{example}

The \emph{separated quotient} $\sep(\VS)$ of a complete bornological vector
space~$\VS$ is the quotient of~$\VS$ by the closure of~$\{0\}$.  It is always
separated.  It is characterized by the universal property that any map $\VS
\to \VS[W]$ with separated range~$\VS[W]$ factors uniquely
through~$\sep(\VS)$.

Let $(\VS_i)_{i\in I}$ be an inductive system of bornological vector spaces.
The \emph{inductive limit} of~$(\VS_i)_{i\in I}$ is constructed as follows.
As a vector space, $\indlim {(\VS_i)_{i\in I}}$ is equal to the ordinary
inductive limit in the category of vector spaces.  There are natural maps
$\iota_i\colon \VS_i \to \indlim {(\VS_i)_{i\in I}}$ for all $i\in I$.  The
bornology on $\indlim {(\VS_i)_{i\in I}}$ is taken to be the \emph{finest}
bornology making all maps~$\iota_i$ bounded.  That is, $S\subset \indlim
{(\VS_i)_{i\in I}}$ is small iff there are $i\in I$, $T\in \CBS(\VS_i)$ such
that $S \subset \iota_i(T)$.

The \emph{direct sum} $\sum_{i\in I} \VS_i$ of a set of bornological vector
spaces is constructed similarly: Endow the vector space direct sum with the
finest bornology making all coordinate injections bounded.  The direct sum is
complete (separated) iff all summands~$\VS_i$ are complete (separated).

The inductive limit need not be separated even if all~$\VS_i$ are separated.
However, if the structure maps $\VS_i \to \VS_j$ of the inductive system are
all \emph{injective}, then $\indlim {(\VS_i)_{i\in I}}$ is separated.  An
inductive system is \emph{injective} iff all the structure maps are injective.
This notion is not invariant under equivalence of inductive systems.  That is,
if $(\VS_i)_{i \in I}$ is an injective inductive system and isomorphic
to~$(\VS[W]_j)_{j\in J}$, then the structure maps in~$(\VS[W]_j)_{j\in J}$
need not be injective.  Therefore, we call an inductive system
\emph{essentially injective} iff it is equivalent to an injective inductive
system.

To obtain inductive limits in the category of separated convex bornological
vector spaces, we combine the functors $\indlim$ and $\sep$.  We call $\sep (
\indlim {(\VS_i)_{i\in I}} )$ the \emph{separated inductive limit}
of $(\VS_i)_{i\in I}$.  It has the right universal property for an inductive
limit in the category of separated bornological vector spaces.  We also write
$\indlim$ for this separated inductive limit if it is clear that the separated
inductive limit is intended.  If all~$\VS_i$ are complete, so is the separated
inductive limit.

We write $\Lin_n(\VS_1, \dots, \VS_n; \VS[W])$ for the vector space of bounded
\Mpn{n}linear maps $\VS_1 \times \dots \VS_n \to \VS[W]$, endowed with the
equibounded bornology.  The following proposition formulates a strengthening
of the universal property of inductive limits.  The proof is elementary.

\begin{proposition}  \label{pro:indlim_universal}
  Let~$(\VS_i)_{i\in I}$ be an inductive system of convex bornological vector
  spaces and let~$\VS[W]$ be a convex bornological vector space.  There is a
  natural bornological isomorphism
  \begin{equation}  \label{eq:indlim_universal}
    \Lin\bigl(\indlim {(\VS_i)_{i\in I}}; \VS[W]\bigr) \cong
    \prolim_{i\in I} \Lin(\VS_i; \VS[W]).
  \end{equation}
  More generally, if $(\VS_i^j)_{i \in I_j}$, $j=1,\dots,n$, are inductive
  systems of convex bornological vector spaces, then there is a natural
  bornological isomorphism
  \begin{equation}  \label{eq:indlim_universal_multi}
    \Lin_n\bigl(\indlim {(\VS_{i_1}^1)_{i_1\in I_1}}, \dots,
    \indlim {(\VS_{i_n}^n)_{i_n\in I_n}}; \VS[W]\bigr) \cong
    \prolim_{i_1\in I_1,\dots, i_n \in I_n}
    \Lin_n(\VS_{i_1}^n, \dots, \VS_{i_n}^n; \VS[W]).
  \end{equation}
  The inductive limits above are non-separated.  If we assume that all
  $\VS_i^j$ and all~$\VS[W]$ are separated convex bornological vector spaces,
  then analogous statements hold for the separated inductive limit.
\end{proposition}

In particular, an inductive limit of a family of bornological algebras is
again a bornological algebra.  That is, the natural multiplication is bounded.
This usually fails for topological algebras.

\begin{example}  \label{exa:LF_indlim}
  Let~$V$ be an LF-space.  That is, $V = \bigcup V_n$ for an increasing
  sequence of \Frechet subspaces~$V_n$ and the topology on~$V$ is the finest
  one that makes the inclusions $V_n \to V$ continuous.  A bounded
  subspace~$S$ of~$V$ is already contained in~$V_n$ for some $n \in \N$.  Of
  course, if~$S$ is bounded or precompact in~$V$, then it is bounded or
  precompact as a subset of~$V_n$, respectively.  Thus $(V, \BOUND) \cong
  \indlim {(V_n, \BOUND)}$ and $(V, \COMP) \cong \indlim {(V_n, \COMP)}$.
\end{example}

\subsection{Bornological vector spaces and inductive systems}
\label{app:bornologies_inductive}

Recall that $\CBS_c(\VS)$ denotes the directed set of all completant small
disks in~$\VS$.  If $S_1\subset S_2$, there is a canonical inclusion
$\VS_{S_1}\to \VS_{S_2}$.  Thus $(\VS_S)_{S \in \CBS_c(\VS)}$ is an injective
inductive system of Banach spaces.  The map $\VS \mapsto (\VS_S)_{S \in
  \CBS_c(\VS)}$ can be extended to a functor from the category of complete
bornological vector spaces to the category of inductive systems of Banach
spaces.  If $l \colon \VS \to \VS[W]$ is a bounded map, then $l(S) \in
\CBS_c(\VS[W])$ for all $S \in \CBS_c(\VS)$ and~$l$ restricts to a bounded map
$l_S\colon \VS_S \to \VS[W]_{l(S)}$.  These maps piece together to a morphism
of inductive systems from $(\VS_S)$ to $(\VS[W]_T)$.  We call this functor the
\emph{dissection functor}~$\dissect$.

Let $V$ and~$W$ be inductive systems of Banach spaces.  We write $\Mor(V; W)$
for the space of morphisms of inductive systems $V \to W$.  If $V = (V_i)_{i
\in I}$ and $W = (W_j)_{j \in J}$, then
\begin{equation}  \label{eq:inductive_morphisms}
  \Mor(V; W) =
  \prolim_{i\in I} \indlim_{j\in J} \Lin( V_i; W_j).
\end{equation}
Here the limits are taken in the category of (bornological) vector spaces and
are \emph{non-separated}.  By~\eqref{eq:inductive_morphisms}, we have a
natural bornology on the space $\Mor(V; W)$.

The \emph{inductive limit functor} $\indlim$ associates to an inductive system
of Banach spaces the corresponding \emph{separated} inductive limit in the
category of complete bornological vector spaces.

\begin{theorem}  \label{the:bornologies_inductive}
  The functors $\indlim$ and $\dissect$ are adjoint, that is,
  \begin{equation}  \label{eq:indlim_dissect_adjoint}
    \Lin( \indlim {(V_i)_{i\in I}}; \VS[W]) \cong
    \prolim_{i\in I} \indlim_{S\in \CBS_c(\VS[W])} \Lin( V_i; \VS[W]_S) =
    \Mor\bigl( (V_i)_{i\in I}; \dissect \VS[W] \bigr)
  \end{equation}
  if~$(V_i)_{i\in I}$ is an inductive system of Banach spaces and~$\VS[W]$ is
  a complete bornological vector space.

  The composition $\indlim \circ \dissect$ is equivalent to the identity.
  That is, there is a natural isomorphism
  \begin{equation}  \label{eq:indlim_dissect_inverse}
    \VS \cong \indlim_{S\in\CBS_c(\VS)} \VS_S
  \end{equation}
  for all complete bornological vector spaces~$\VS$.  We have a natural
  bornological isomorphism
  \begin{equation}
    \label{eq:dissect_fully_faithful}
    \Lin(\VS; \VS[W]) \cong
    \prolim_{S\in\CBS_c(\VS)} \indlim_{T\in\CBS_c(\VS[W])}
    \Lin(\VS_S; \VS[W]_T) =
    \Mor( \dissect \VS; \dissect \VS[W])
  \end{equation}
  for all complete bornological vector spaces $\VS,\VS[W]$.  In particular,
  the functor $\dissect$ is fully faithful.

  More generally, if\/ $\VS_1,\dots,\VS_n,\VS[W]$ are complete bornological
  vector spaces, we have a natural bornological isomorphism
  \begin{equation}
    \label{eq:dissect_fully_faithful_multi}
    \Lin_n(\VS_1, \dots,\VS_n; \VS[W]) \cong
    \prolim_{S_j\in\CBS_c(\VS_j)} \indlim_{T \in \CBS_c(\VS[W])}
    \Lin_n\bigl( (\VS_1)_{S_1}, \dots, (\VS_n)_{S_n}; \VS[W]_T \bigr).
  \end{equation}
  
  Let $(V_i)_{i\in I}$ be an injective inductive system of Banach spaces.
  Then $\dissect \circ \indlim {(V_i)_{i\in I}}$ is naturally isomorphic
  to~$(V_i)_{i\in I}$.  Hence the inductive limit functor induces a
  bornological isomorphism
  \begin{equation}  \label{eq:indlim_iso_essin}
    \Mor( W; V) \cong \Lin( \indlim W; \indlim V)
  \end{equation}
  if $W$ and~$V$ are inductive systems and~$V$ is essentially injective.

  The functor $\dissect$ is an equivalence between the category of complete
  bornological vector spaces and the category of injective inductive systems
  of Banach spaces.
  
  These assertions remain true if we replace complete bornological vector
  spaces, $\CBS_c(\VS)$, separated inductive limits, and Banach spaces by
  convex bornological vector spaces, $\CBS_d(\VS)$, non-separated inductive
  limits, and semi-normed spaces, respectively; or by separated convex
  bornological vector spaces, $\CBS_d(\VS)$, separated inductive limits, and
  normed spaces, respectively.
\end{theorem}

\begin{proof}
  Since the~$V_i$ in~\eqref{eq:indlim_dissect_adjoint} are Banach spaces, we
  have bornological isomorphisms
  \begin{displaymath}
    \Lin( V_i; \VS[W]) \cong
    \indlim_{S \in \CBS_c(\VS[W])} \Lin( V_i; \VS[W]_S)
    \qquad \forall i \in I.
  \end{displaymath}
  This follows directly from the definitions of bounded linear maps.
  Equation~\eqref{eq:indlim_dissect_adjoint} follows if we plug this into the
  universal property of inductive limits~\eqref{eq:indlim_universal}.
  
  Let~$\VS$ be a complete bornological vector space.  The universal property
  of the inductive limit gives rise to a natural bounded linear map $\indlim
  \VS_S \to \VS$ extending the natural inclusions $\VS_S \subset \VS$.  The
  concrete construction of inductive limits shows that this map is a
  bornological isomorphism.
  
  If we apply~\eqref{eq:indlim_dissect_adjoint} to $\dissect \VS$
  and~$\VS[W]$, we get~\eqref{eq:dissect_fully_faithful}.  Thus the functor
  $\dissect$ is fully faithful.
  Equation~\eqref{eq:dissect_fully_faithful_multi} can be proved similarly.
  By~\eqref{eq:indlim_universal_multi}, we can reduce to the case where
  $\VS_1,\dots,\VS_n$ are all primitive spaces.  If $\VS_1,\dots,\VS_n$ are
  primitive, the definition of an equibounded family of \Mpn{n}linear maps
  implies immediately that
  \begin{displaymath}
    \Lin_n( \VS_1,\dots \VS_n; \VS[W]) \cong
    \indlim_{S \in \CBS_c(\VS[W])} \Lin_n( \VS_1, \dots, \VS_n; \VS[W]_S).
  \end{displaymath}
  
  Let $(V_i)_{i\in I}$ be an injective system of Banach spaces and let
  $\VS\defeq \indlim {(V_i)_{i\in I}}$ be its inductive limit.  Let $B_i
  \subset \VS$ be the image of the unit ball of~$V_i$ under the structure map
  $V_i \to \VS$.  By the concrete definition of the inductive limit bornology,
  a subset of~$\VS$ is small iff it is contained in $c\cdot B_i$ for some
  $i\in I$ for some constant $c\in\R$.  It follows that $\dissect \VS$ is
  equivalent to the inductive system $(\VS_{B_i})_{i\in I}$.  Since~$(V_i)$ is
  injective, the structure maps $V_i \to \VS$ are all injective.  Hence we
  have natural isomorphisms $V_i \cong \VS_{B_i}$ for all $i\in I$.  It
  follows that $(V_i)_{i\in I}$ is equivalent to $\dissect \VS$.  Conversely,
  inductive systems of the form $\dissect \VS$ are injective by construction.
  Thus $\dissect$ is an equivalence between the category of complete
  bornological vector spaces and the category of injective inductive systems
  of Banach spaces.  Its inverse is the inductive limit functor $\indlim$.
  
  Of course, we still have $\dissect \circ \indlim V \cong V$ if~$V$ is only
  essentially injective.  Therefore, if~$V$ is essentially injective,
  then~\eqref{eq:indlim_dissect_adjoint} implies
  $$
  \Mor( W; V) \cong
  \Mor(W; \dissect \circ \indlim V) \cong
  \Lin(\indlim W; \indlim V).
  $$
  
  If we consider separated convex bornological vector spaces instead, we have
  to replace $\CBS_c(\VS)$ by $\CBS_d(\VS)$.  The spaces~$\VS_S$ for $S\in
  \CBS_d(\VS)$ are then only normed spaces.  Otherwise, nothing changes.  If
  we drop the separation assumption, then the spaces~$\VS_S$ are only
  semi-normed.  We have to take the non-separated inductive limit instead of
  the separated inductive limit to use~\eqref{eq:indlim_universal}.
\end{proof}

\subsection{Construction of completions}
\label{app:completions}

Let~$\VS$ be a convex bornological vector space.  Write $\VS = \indlim
{(\VS_S)_{S \in \CBS_d(\VS)}}$.  Let~$\Cpl{\VS_S}$ be the Hausdorff completion
of~$\VS_S$.  Since taking Hausdorff completions is functorial,
$(\Cpl{\VS_S})_{S \in \CBS_d(\VS)}$ is an inductive system of Banach spaces.
Let~$\Cpl{\VS}$ be its separated inductive limit.  The universal property of
Hausdorff completions implies that $\Mor \bigl( (\Cpl{\VS_S}); W\bigr) \cong
\Mor( \dissect \VS; W)$ for all inductive systems of Banach spaces~$W$.  The
same holds for multi-linear morphisms between inductive systems.  These are
defined by the right hand side of~\eqref{eq:dissect_fully_faithful_multi}.
Applying~\eqref{eq:indlim_dissect_adjoint}
and~\eqref{eq:dissect_fully_faithful}, we get bornological isomorphisms
\begin{displaymath}
  \Lin(\VS; \VS[W])
  \cong
  \Mor( \dissect \VS; \dissect \VS[W])
  \cong
  \Mor\bigl( (\Cpl{\VS_S})_{S \in \CBS_d(\VS)};  \dissect \VS[W] \bigr)
  \cong
  \Lin\bigl( \indlim (\Cpl{\VS_S})_{S \in \CBS_d(\VS)};  \VS[W] \bigr)
\end{displaymath}
for all complete bornological vector spaces~$\VS[W]$.  Thus the
space~$\Cpl{\VS}$ constructed above has the right universal property for a
completion.  Using a multi-linear version
of~\eqref{eq:indlim_dissect_adjoint}, we obtain a bornological isomorphism
$$
\Lin_n( \VS_1,\dots,\VS_n; \VS[W]) \cong
\Lin_n( \Cpl{\VS_1},\dots,\Cpl{\VS_n}; \VS[W])
$$
for any complete bornological vector space.  This is the meaning of
Lemma~\ref{lem:Cpl_multi}.

The good incomplete spaces are those for which the map $\natural \colon \VS
\to \Cpl{\VS}$ is a bornological isomorphism onto its range, endowed with the
subspace bornology from~$\Cpl{\VS}$.  That is, a set $S \subset \VS$ is small
iff $\natural(S) \subset \Cpl{\VS}$ is small.  This implies that~$\VS$ is a
bornological subspace of a complete space.  Conversely, a bornological
subspace of a complete space is a subspace of its completion.  This follows
easily from the following lemma.

\begin{lemma}  \label{lem:completion_injective}
  Let~$\VS$ be a convex bornological vector space, $\VS[W]$ a complete
  bornological vector space, let $l \colon \VS \to \VS[W]$ be a bounded linear
  map.  Assume that $\lim l(x_n) = 0$ in~$\VS[W]$ implies that $\lim x_n = 0$
  for all Cauchy sequences $(x_n)_{n \in \N}$ in~$\VS$.  In particular, $l$ is
  injective.
  
  Let $\CBS' \subset \CBS(\VS[W])$ be the collection of all subsets of sets of
  the form $\coco{\bigl( l(S) \bigr)}$ with $S \in \CBS(\VS)$.  Let $\VS'
  \subset \VS[W]$ be the linear span of\/ $\bigcup \CBS'$.

  Then~$\CBS'$ is a completant bornology on~$\VS'$.  The map~$l$ can be viewed
  as an injective bounded linear map $l' \colon \VS \to \VS'$.  The complete
  bornological vector space $(\VS', \CBS')$ has the universal property of the
  completion of\/~$\VS$ (with respect to the map $l' \colon \VS \to \VS'$) and
  therefore is equal to~$\Cpl{\VS}$.
\end{lemma}

\begin{proof}
  The sets $\coco{\bigl( l(S) \bigr)}$ are completant disks by construction.
  Each set of~$\CBS'$ is contained in a completant small disk and $\{ x \} \in
  \CBS'$ for all $x \in \VS'$.  Thus~$\CBS'$ is a completant bornology
  on~$\VS'$.  For example, condition~(iii) for a bornology follows because
  $\coco{l(S_1 + S_2)} = \coco{l(S_1)} + \coco{l(S_2)}$.  By construction,
  $l(S) \in \CBS'$ for all $S \in \CBS(\VS)$.  Hence we can view~$l$ as a
  bounded linear map $l' \colon \VS \to \VS'$.
  
  It remains to verify that any bounded linear map $f \colon \VS \to \VS[X]$
  with complete range~$\VS[X]$ can be factored as $f = f' \circ l'$ with a
  bounded linear map $f' \colon (\VS', \CBS') \to \VS[X]$.  Let $S \in
  \CBS_d(\VS)$, then~$f(S)$ is contained in a small completant disk $T \in
  \CBS_c(\VS[X])$ because~$\VS[X]$ is complete.  We claim that $f_S \defeq
  f|_{\VS_S} \colon \VS_S \to \VS[X]_T$ can be extended to a bounded linear
  map $f_S' \colon \VS'_{\coco{S}} \to \VS[X]_T$.
  
  It is clear that~$f_S$ can be extended to a bounded linear map on the
  completion~$\Cpl{\VS_S}$ of~$\VS_S$.  The problem is that the natural map
  $\Cpl{\VS_S} \to \VS'_{\coco{S}}$ may fail to be injective.  An element in
  the kernel is the limit~$x$ of a (bornological) Cauchy sequence $(x_n)$
  in~$\VS_S$ for which $\lim l(x_n) = 0 \in \VS[W]$.  By assumption, since
  $l(x_n)$ is a bornological null sequence, $(x_n)$ is a bornological null
  sequence in~$\VS$ and hence $\lim f(x_n) = 0$ because~$f$ is bounded.  Thus
  the kernel of the map $\Cpl{\VS_S} \to \VS'_{\coco{S}}$ is annihilated by
  the extension $\Cpl{f_S}$.  Consequently, $\Cpl{f_S}$ descends to a bounded
  linear map $f_S' \colon \VS'_{\coco{S}} \to \VS[X]_T$.
  
  Finally, we observe that the extensions~$f_S'$, $S \in \CBS(\VS)$ piece
  together to a bounded linear map $(\VS', \CBS') \to \VS[X]$.  Thus $(\VS',
  \CBS')$ is isomorphic to the completion of~$\VS$.
\end{proof}

\subsection{Tensor products of \Frechet spaces}
\label{app:tensor_Frechet}

We prove Theorem~\ref{the:tensor_Frechet} and Corollary~\ref{cor:tensor_LF}.

{
\theoremstyle{plain}
\newtheorem*{thmfrechet}{Theorem~\ref{the:tensor_Frechet}}

\begin{thmfrechet}
  Let\/ $\VS_1$ and~$\VS_2$ be \Frechet spaces and let\/ $\VS_1 \prot \VS_2$ be
  their completed projective tensor product~\cite{grothendieck55:produits}.
  The natural bilinear map $\natural \colon \VS_1 \times \VS_2 \to \VS_1 \prot
  \VS_2$ induces bornological isomorphisms
  \begin{gather*}
    (\VS_1, \COMP ) \hot (\VS_2, \COMP ) \cong (\VS_1 \prot \VS_2, \COMP);
    \\
    (\VS_1, \BOUND) \hot (\VS_2, \BOUND) \cong (\VS_1 \prot \VS_2, \CBS ).
  \end{gather*}
  Here~$\CBS$ denotes the bornology of all $S \subset \VS_1 \prot \VS_2$ that
  are contained in a set of the form $\coco{(B_1 \otimes B_2)}$ with bounded
  sets $B_1 \in \BOUND(\VS_1)$ and $B_2 \in \BOUND(\VS_2)$.

  The bornology $\CBS$ is equal to $\BOUND( \VS_1 \prot \VS_2 )$ in the
  following cases.  If both\/ $\VS_1$ and\/~$\VS_2$ are Banach spaces;
  if\/~$\VS_2$ is arbitrary and\/ $\VS_1 = L^1(M, \mu)$ is the space of
  integrable functions on a locally compact space with respect to some Borel
  measure; or if\/~$\VS_2$ is arbitrary and\/~$\VS_1$ is nuclear.
\end{thmfrechet}

}
\begin{proof}
  Recall \Grothendieck's fundamental theorem about compact subsets of the
  projective tensor product of two \Frechet spaces.

  \begin{theorem}[{\cite[p.~51]{grothendieck55:produits}}]
    \label{the:Grothendieck}
    Let $\VS_1$ and~$\VS_2$ be two \Frechet spaces.  Let $K \subset \VS_1
    \prot \VS_2$ be a compact subset.  Then there are null-sequences
    $(\vs_{1,n})$ and~$(\vs_{2,n})$ in $\VS_1$ and~$\VS_2$, respectively, and
    a compact subset~$K_0$ of the unit ball of the space~$\ell^1(\N)$ of
    absolutely summable sequences such that
    $$
    K \subset \biggl\{ \sum_{n=1}^\infty \lambda_n
    \vs_{1,n} \otimes \vs_{2,n} \biggm| (\lambda_n) \in K_0 \biggr\}.
    $$
  \end{theorem}

  In particular, it follows that each point of $\VS_1 \prot \VS_2$ is contained
  in a set of the form $\coco{(B_1 \otimes B_2)}$ with bounded $B_1$ and~$B_2$.
  Thus~$\CBS$ is a bornology on $\VS_1 \prot \VS_2$.  Since $\VS_1 \prot \VS_2$
  is topologically complete, the bornologies $\COMP(\VS_1 \prot \VS_2)$ and
  $\CBS$ are completant.
  
  The natural bilinear map $(\VS_1, \COMP) \times (\VS_2, \COMP) \to (\VS_1
  \prot \VS_2, \COMP)$ is bounded.  Hence there is a bounded linear map $f_C
  \colon (\VS_1, \COMP) \hot (\VS_2, \COMP) \to (\VS_1 \prot \VS_2, \COMP)$.
  Theorem~\ref{the:Grothendieck} implies that each precompact subset of $\VS_1
  \prot \VS_2$ is contained in a set of the form $\coco{(K_1 \otimes K_2)}$
  with $K_j \defeq \{\vs_{j,n} \mid n \in \N\}$, $j = 1,2$, and null-sequences
  $(\vs_{j,n})$.  Since the points of a null-sequence form a precompact set,
  $K_1$ and~$K_2$ are precompact.  Thus~$f_C$ is a quotient map.  Similarly,
  the bilinear map $(\VS_1, \BOUND) \times (\VS_2, \BOUND) \to (\VS_1 \prot
  \VS_2, \CBS)$ is bounded and induces a bounded linear map $f_B \colon
  (\VS_1, \BOUND) \hot (\VS_2, \BOUND) \to (\VS_1 \prot \VS_2, \CBS)$.  By
  definition of~$\CBS$, each $S \in \CBS$ is of the form $f_B(\hat{S})$ for a
  small set $\hat{S} \subset (\VS_1, \BOUND) \hot (\VS_2, \BOUND)$.  That is,
  $f_B$ is a quotient map.
  
  It remains to show that $f_B$ and~$f_C$ are injective.  For an injective
  quotient map is automatically a bornological isomorphism.  By the concrete
  definition of the completion, $(\VS_1, \COMP) \hot (\VS_2, \COMP)$ is the
  separated inductive limit of the spaces $(\VS_1)_{K_1} \prot (\VS_2)_{K_2}$
  with~$K_j$ running through $\COMP(\VS_j)$, $j = 1,2$.  Pick any $x \in
  (\VS_1, \COMP) \hot (\VS_2, \COMP)$ with $f_C(x) = 0$.  We have to show that
  $x = 0$.  There are completant precompact disks $K_j \subset \VS_j$,
  $j = 1,2$, such that $x \in (\VS_1)_{K_1} \prot (\VS_2)_{K_2}$.
  
  We view~$x$ as an element in the kernel of the natural map $(\VS_1)_{K_1}
  \prot (\VS_2)_{K_2} \to \VS_1 \prot \VS_2$.  This natural map can indeed
  fail to be injective because there are Banach spaces that do not have
  \Grothendieck's approximation property.  However, for~$x$ to describe the
  zero element of $(\VS_1, \COMP) \hot (\VS_2, \COMP)$, it suffices that there
  are completant precompact disks $L_j \subset \VS_j$, $j = 1,2$, such
  that~$x$ is annihilated by the natural map $(\VS_1)_{K_1} \prot
  (\VS_2)_{K_2} \to (\VS_1)_{L_1} \prot (\VS_2)_{L_2}$.  This follows from the
  following corollary of Theorem~\ref{the:Grothendieck} also due to
  \Grothendieck.

  \begin{proposition}[{\cite[Remarque~4, p.~57]{grothendieck55:produits}}]
    Let $\VS_1$ and~$\VS_2$ be \Frechet spaces.  Let~$(\vs[w]_N)$ be a
    null-sequence in $\VS_1 \prot \VS_2$.  Then there are null-sequences
    $(\vs'_{j,n})_{n \in \N}$ in~$\VS_j$, $j = 1,2$, and a null-sequence
    $(\Lambda_N)_{N\in\N}$ in~$\ell^1(\N)$ such that $\vs[w]_N =
    \sum_{i=1}^\infty \Lambda_N(i) \vs'_{1,i} \otimes \vs'_{2,i}$ for all $n
    \in \N$.
  \end{proposition}
  
  There is a Cauchy sequence $(x_N)$ in the uncompleted tensor product
  $(\VS_1)_{K_1} \protu (\VS_2)_{K_2}$ with $x = \lim x_N$.  The sequence
  $\vs[w]_N \defeq f_C(x_N)$ converges bornologically towards $f_C(x) = 0$.
  We apply the proposition to this null-sequence.  Let~$L_j$ be the disked
  hull of $\{\vs'_{j,n} \mid n \in \N\} \cup K_j$.  Thus~$L_j$ is a precompact
  subset of~$\VS_j$ and~$(\vs[w]_N)$ is a null-sequence in $(\VS_1)_{L_1} \hot
  (\VS_2)_{L_2}$.  Hence the image of~$x$ in $(\VS_1)_{L_1} \hot
  (\VS_2)_{L_2}$ is zero.  Thus $x = 0$ in $(\VS_1, \COMP) \hot (\VS_2,
  \COMP)$ as desired.

  The proof of Theorem~\ref{the:tensor_Frechet} is almost finished.  It
  remains to compare $\CBS$ and $\BOUND(\VS_1 \prot \VS_2)$.  It is trivial
  that $\CBS = \BOUND(\VS_1 \prot \VS_2)$ if both $\VS_1$ and~$\VS_2$ are
  Banach spaces.  For $L^1(M,\mu)$ and nuclear spaces, see Grothendieck
  \cite[p.~68--69]{grothendieck55:produits}
  \cite[p.~73--74]{grothendieck55:espaces}.
\end{proof}

{
\theoremstyle{plain}
\newtheorem*{corfrechet}{Corrollary~\ref{cor:tensor_LF}}

\begin{corfrechet}
  Let\/ $\VS_1$ and\/~$\VS_2$ be nuclear LF-spaces.  Then $(\VS_1, \BOUND) \hot
  (\VS_2, \BOUND)$ is isomorphic to \Grothendieck's inductive tensor product\/
  $\VS_1 \indt \VS_2$, endowed with the bounded bornology.  The inductive
  tensor product\/ $\VS_1 \indt \VS_2$ is again a nuclear LF-space.
\end{corfrechet}

}
\begin{proof}
  Write $\VS_1$ and~$\VS_2$ as inductive limits of nuclear \Frechet spaces,
  $\VS_j = \indlim_{n} \VS_{j,n}$, $j=1,2$.  The inductive tensor product
  $\VS_1 \indt \VS_2$ is isomorphic to $\indlim \VS_{1,n} \prot \VS_{2,n}$
  (see~\cite{grothendieck55:produits}).  Since all~$\VS_{j,n}$ are nuclear,
  the maps $\VS_{1,n} \prot \VS_{2,n} \to \VS_{1,n+1} \prot \VS_{2,n+1}$ are
  topological embeddings for all $n \in \N$.  Hence $\VS_1 \indt \VS_2$ is a
  nuclear LF-space.  We claim that
  \begin{multline*}
    (\VS_1 \indt \VS_2, \BOUND) \cong
    \indlim {(\VS_{1,n} \prot \VS_{2,n}, \BOUND)} \cong
    \indlim {\bigl( (\VS_{1,n}, \BOUND) \hot (\VS_{2,n}, \BOUND) \bigr)}
    \\ \cong
    \indlim {(\VS_{1,n}, \BOUND)} \hot \indlim {(\VS_{2,n}, \BOUND)} \cong
    (\VS_1, \BOUND) \hot (\VS_2, \BOUND).
  \end{multline*}
  The first and the last isomorphism follow from
  Example~\ref{exa:LF_indlim}.  The second isomorphism follows from
  Theorem~\ref{the:tensor_Frechet}.  The third isomorphism follows because
  completed bornological tensor products commute with inductive limits.
\end{proof}

\subsection{Admissible \Frechet algebras}
\label{app:proof_admissible}

\Frechet algebras satisfying one of the following equivalent conditions are
called \emph{admissible} by Puschnigg~\cite{puschnigg96:asymptotic}.

{
\theoremstyle{plain}
\newtheorem*{thadmiss}{Theorem~\ref{the:tensoring_Frechet}}

\begin{thadmiss}
  Let~$\MA[C]$ be a \Frechet algebra.  Then the following are equivalent:
  \begin{enumerate}[(i)]%

  \item for all precompact subsets $S \subset \MA[C]$, there is $\lambda > 0$
    such that $(\lambda S)^\infty$ is precompact;
    
  \item for each null-sequence $(x_n)_{n \in \N}$ in~$\MA[C]$, there is $N \in
    \N$ such that $\{x_n \mid n \ge N \}^\infty$ is precompact;
    
  \item there is a neighborhood~$U$ of the origin in~$\MA[C]$ such
    that~$S^\infty$ is precompact for all precompact sets $S \subset U$.

  \end{enumerate}
\end{thadmiss}

}

\begin{proof}
  \textbf{(i) implies~(ii)}.  Let $(x_n)_{n \in \N}$ be a null-sequence
  in~$\MA[C]$.  Since~$\MA[C]$ is \Frechet, any topological null-sequence is a
  bornological null-sequence.  That is, there is a null-sequence of positive
  real numbers~$(\epsilon_n)$ such that $\{ \epsilon_n^{-1} \cdot x_n \mid n
  \in \N\}$ is precompact.  By~(i), there is $\lambda > 0$ such that $T =
  \lambda \cdot \{\epsilon_n^{-1} \cdot x_n\}$ is bornopotent.  The circled
  hull of~$T$ is still bornopotent.  There is $N \in \N$ such that $\lambda >
  \epsilon_n$ for all $n \ge N$.  Thus for $n \ge N$, $x_n \in (\epsilon_n /
  \lambda) \cdot T$ is contained in the circled hull of~$T$.  It follows that
  the set $\{x_n \mid n \ge N\}$ is bornopotent.

  \textbf{(ii) implies~(iii)}.  Assume that~(iii) is false.  Choose a
  countable basis $(U_n)_{n \in \N}$ of neighborhoods of the origin
  in~$\MA[C]$.  Since~(iii) is false, there are precompact sets $S_n \subset
  U_n$ for all $n \in \N$ such that~$S_n^\infty$ not precompact.
  Since~$\MA[C]$ is a \Frechet space, there is a null-sequence $(x_{n,j})_{j
    \in \N}$ such that~$S_n$ is contained in the completant disked hull of $\{
  x_{n,j} \mid j \in \N\}$ and $x_{n,j}\in S_n$ for all $n,j$.  Rearrange the
  numbers~$(x_{n,j})$ in a sequence.  We claim that this sequence is a
  null-sequence.  Fix a neighborhood of the origin~$U$.  Since
  $(x_{n,j})_{j}$ is a null-sequence for each $n \in \N$, there is $j(n) \in
  \N$ such that $x_{n,j} \in U$ for all $j\ge j(n)$.  Since $U_N \subset U$
  for some $N \in \N$, $x_{n,j} \in S_n \subset U$ for all $n \ge N$.  Thus
  all but finitely many of the~$x_{n,j}$ are in~$U$.  That is, $(x_{n,j})_{n,
    j \in \N}$ is a null-sequence for any ordering of the indices.
  
  Let $F \subset \N \times \N$ be finite and let $S \defeq \{ x_{n,j} \mid
  (n,j) \in \N^2 \setminus F\}$.  There is $N \in \N$ for which $(\{N\} \times
  \N) \cap F = \emptyset$.  Thus~$S_N$ is contained in the disked hull of~$S$.
  Therefore, $S^\infty$ is not precompact.  Consequently, the null-sequence
  $(x_{n,j})_{n, j \in \N^2}$ violates~(ii).  Thus if~(iii) is false,
  then~(ii) is false.
  
  \textbf{(iii) implies~(i)} is evident because neighborhoods of the origin
  absorb all bounded sets.
\end{proof}

\section{Homological algebra and universal algebra}
\label{app:homological_algebra}

We carry over some algebra from algebras without additional structure to
complete bornological algebras.  We obtain the \emph{long exact homology
sequence} for allowable extensions of complexes.  We define \emph{modules}
over complete bornological algebras and introduces \emph{free} and
\emph{projective} modules.  Modules are always assumed to be complete.  We
state the \emph{comparison theorem} for allowable projective resolutions.  We
define the \emph{bar resolution} and \emph{\Hochschild homology}.  These
things work out as in \emph{relative} homological
algebra~\cite{MacLane:Homology}.  Relative homological algebra treats an
algebra~$\Lambda$ over a commutative ground ring~$\mathsf{K}$ as if the ground
ring were a field.  In order to exclude non-trivial homology contributions
from~$\mathsf{K}$, only extensions that split as extensions of
\Mpn{\mathsf{K}}modules are allowed.

We identify the \emph{tensor algebra} $\Tens\MA$ of a complete bornological
algebra~$\MA$ with the even part of the algebra of non-commutative
differential forms $\Omega\MA$ endowed with the Fedosov product.  The reader
is assumed to be familiar with the definition of $\Omega\MA$ over a
\emph{non-unital} algebra~$\MA$ in~\cite{connes94:ncg}.  The identification
$\Tens\MA \cong (\Omega\MA, \odot)$ is due to Cuntz and
Quillen~\cite{cuntz95:algebra}.  We prove the equivalence of several
definitions for \emph{quasi-free} algebras.  The notion of quasi-freeness was
introduced (for \emph{unital} algebras without additional structure) by Cuntz
and Quillen~\cite{cuntz95:algebra}.

Most results carry over to objects in an arbitrary additive category with
finite direct sums, tensor products, and a ground ring~$\C$ satisfying $\C
\hot X \cong X \cong X \hot \C$.  In particular, they continue to hold in
categories of projective or inductive systems.  However, the algebras
$\Omega\MA$ and $\Tens\VS$ of~\ref{app:Omega_Tens} are infinite direct sums.
In more general categories, a direct sum of algebras need not be an algebra
any more.  This happens, for example, for topological algebras or
pro-algebras.  However, the application to quasi-freeness only needs the
quotients $\Tens\VS / (\Jens\VS)^k$ that do not cause any problems because
they are finite direct sums.

\subsection{Complexes of complete bornological vector spaces}
\label{app:complexes}

A \emph{\Mpn{\Z}graded complex of complete bornological vector spaces} is a
collection $(K_n)_{n \in \Z}$ of complete bornological vector spaces with
bounded linear maps $\partial_n \colon K_n \to K_{n-1}$ satisfying $\partial_n
\circ \partial_{n+1} = 0$ for all $n \in \Z$.  A \emph{\Mpn{\Ztwo}graded
  complex of complete bornological vector spaces} consists of complete
bornological vector spaces $K_0$ and~$K_1$ and bounded linear maps $\partial_0
\colon K_0 \to K_1$ and $\partial_1 \colon K_1 \to K_0$ satisfying $\partial_0
\circ \partial_1 = 0$ and $\partial_1 \circ \partial_0 = 0$.

Let $J_\bullet$ and~$K_\bullet$ be complexes of complete bornological vector
spaces.  We define the complex $\Lin(J_\bullet; K_\bullet)$ as follows.  Let
$\Lin(J_\bullet; K_\bullet)_k \defeq \prod_{n} \Lin(J_n; K_{k+n})$ be the
space of bounded linear maps of degree~$k$.  If the complexes $J_\bullet$
and~$K_\bullet$ are both \Mpn{\Ztwo}graded, then $n, k \in \Ztwo$.  The
boundary in $\Lin(J_\bullet; K_\bullet)$ is $\delta(f) \defeq [\delta, f]
\defeq \delta_{K}\circ f - (-1)^{\deg f} f \circ \delta_{J}$.  Thus the cycles
in $\Lin(J_\bullet; K_\bullet)_0$ are the bounded chain maps $J_\bullet \to
K_\bullet$ and the boundaries are those chain maps homotopic to zero.  We
write
$$
H^n(J_\bullet; K_\bullet) \defeq H^n \bigl( \Lin(J_\bullet; K_\bullet) \bigr).
$$

The composition of linear maps descends to a well-defined bilinear map on
homology
\begin{equation}  \label{eq:Hn_product}
  \circ \colon H^n(J_\bullet; K_\bullet) \times H^m(I_\bullet; J_\bullet) \to
  H^{n+m} (I_\bullet; K_\bullet).
\end{equation}
This product is associative in the usual sense.  If $f \in H^n(J_\bullet;
K_\bullet)$ and $g \in H^m(I_\bullet; J_\bullet)$, we write $f_\ast(g) \defeq
f\circ g$ and $g^\ast(f) \defeq (-1)^{\deg g \cdot \deg f} f \circ g$.

\begin{theorem}[Long Exact Homology Sequence]  \label{the:long_exact_homology}
  Let $(i,p) \colon K_\bullet \injto E_\bullet \prto Q_\bullet$ be an
  allowable extension of\/ \Mpn{\Ztwo}graded complexes of complete bornological
  vector spaces.  Let~$L_\bullet$ be a \Mpn{\Ztwo}graded complex of complete
  bornological vector spaces.  Then there are natural six-term exact sequences

  \begin{gather}
    \label{eq:long_exact_homology}
    \begin{gathered}
      \xymatrix{
        {H^0(L_\bullet; K_\bullet)} \ar[r]^{i_\ast} &
          {H^0(L_\bullet; E_\bullet)} \ar[r]^{p_\ast} &
            {H^0(L_\bullet; Q_\bullet)} \ar[d]^{\partial_\ast} \\
        {H^1(L_\bullet; Q_\bullet)} \ar[u]_{\partial_\ast} &
          {H^1(L_\bullet; E_\bullet)} \ar[l]_{p_\ast} &
            {H^1(L_\bullet; K_\bullet)} \ar[l]_{i_\ast}
        }
    \end{gathered}
    \displaybreak[0] \\
    \label{eq:long_exact_cohomology}
    \begin{gathered}
      \xymatrix{
        {H_0(Q_\bullet; L_\bullet)} \ar[r]^{p^\ast} &
          {H_0(E_\bullet; L_\bullet)} \ar[r]^{i^\ast} &
            {H_0(K_\bullet; L_\bullet)} \ar[d]^{\partial^\ast} \\
        {H_1(K_\bullet; L_\bullet)} \ar[u]_{\partial^\ast} &
          {H_1(E_\bullet; L_\bullet)} \ar[l]_{i^\ast} &
            {H_1(Q_\bullet; L_\bullet)} \ar[l]_{p^\ast}
        }
    \end{gathered}
  \end{gather}
  
  $i_\ast$, $i^\ast$, $p_\ast$, $p^\ast$, $\partial_\ast$, $\partial^\ast$ are
  the (signed) composition products with the homology classes of the chain
  maps $i$ and~$p$ and with a certain natural element $\partial \in
  H_1(Q_\bullet; K_\bullet)$.
  
  If $s \colon Q_\bullet \to E_\bullet$ is a bounded linear section,
  $\partial$ can be described as follows.  The linear map $[\partial, s] =
  \partial_E \circ s - s \circ \partial_Q$ maps $Q_\bullet$ into $K_\bullet
  \subset E_\bullet$ and satisfies $[\partial, [\partial, s]] = 0$.  Thus
  $[\partial, s]$ yields an element $\partial \in H_1(Q_\bullet; K_\bullet)$.
  The homology class of~$\partial$ does not depend on the choice of~$s$.
\end{theorem}

\begin{proof}
  Since the extension is allowable, $E_\bullet \cong K_\bullet \oplus
  Q_\bullet$ as graded bornological vector spaces.  Thus $\Lin(
  L_\bullet; E_\bullet) \cong \Lin( L_\bullet; K_\bullet) \oplus \Lin(
  L_\bullet; Q_\bullet)$ and $\Lin( E_\bullet; L_\bullet) \cong \Lin(
  K_\bullet; L_\bullet) \oplus \Lin( Q_\bullet; L_\bullet)$.  Consequently,
  \begin{displaymath} 
    \begin{split}
      & \Lin(L_\bullet; K_\bullet) \overset{i_\ast}{\injto}
        \Lin(L_\bullet; E_\bullet) \overset{p_\ast}{\prto}
        \Lin(L_\bullet; Q_\bullet),\\
      & \Lin(Q_\bullet; L_\bullet) \overset{p^\ast}{\injto}
        \Lin(E_\bullet; L_\bullet) \overset{i^\ast}{\prto}
        \Lin(K_\bullet; L_\bullet)
    \end{split}
  \end{displaymath}
  are short exact sequences of complexes of vector spaces.  These induce the
  long exact homology sequences \eqref{eq:long_exact_homology}
  and~\eqref{eq:long_exact_cohomology} as usual.
\end{proof}

The naturality of the long exact sequence means that chain maps $L_\bullet \to
L_\bullet'$ and morphisms between allowable extensions of complexes from
$K_\bullet \injto E_\bullet \prto Q_\bullet$ to $K'_\bullet \injto E'_\bullet
\prto Q'_\bullet$ give rise to commuting diagrams of the resulting six-term
exact sequences.

\subsection{Modules over a complete bornological algebra~$\MA$}
\label{app:modules}

A \emph{(left) \Mpn{\MA}module} is a \emph{complete} bornological vector
space~$\VS$ with a bounded bilinear map $\MA \times \VS \to \VS$ satisfying
the usual associativity condition.  Since~$\MA$ may be non-unital, there is no
unitarity condition for modules.  Let $\Lin_{\MA} (\VS; \VS[W])$ be the space
of \Mpn{\MA}module homomorphisms from $\VS \to \VS[W]$.

Let~$\Unse{\MA}$ be the \emph{unitarization} of~$\MA$.  That is, $\Unse{\MA}
\defeq \MA \oplus \C$ as a bornological vector space with the usual
multiplication for which $1 \in \C$ acts as the identity of~$\Unse{\MA}$.  A
\emph{free left \Mpn{\MA}module} is defined as $\Unse{\MA} \hot \VS$ with left
\Mpn{\MA}action $\ma_1 \cdot (\opt{\ma_2} \otimes \vs) \defeq (\ma_1 \cdot
\opt{\ma_2}) \otimes \vs$.  Free modules have the usual universal property:
Composition with the natural inclusion map $\VS \to \Unse{\MA} \hot \VS$ gives
rise to a natural isomorphism $\Lin_{\MA}(\Unse{\MA} \hot \VS; \VS[W]) \cong
\Lin(\VS; \VS[W])$. A \emph{free \Mpn{\MA}bimodule} is defined as $\Unse{\MA}
\hot \VS \hot \Unse{\MA}$ with the obvious left and right \Mpn{\MA}action.  A
module is \emph{projective} iff it is a direct summand of a free module.

Projective modules~$\VS[P]$ have the usual projectivity property for allowable
extensions.  That is, the functor $\Lin_{\MA}(\VS[P]; \blank)$ is exact in the
sense that it maps allowable extensions of \Mpn{\MA}modules to allowable
extensions of complete bornological vector spaces.  This is evident for free
modules because $\Lin_{\MA}(\VS[P]; \blank) \cong \Lin(\VS; \blank)$ naturally
if~$\VS[P]$ is free on~$\VS$.  The exactness of the functor $\Lin_{\MA}
(\VS[P]; \blank)$ carries over to direct summands of free modules in the usual
way.

Let~$\VS$ be an \Mpn{\MA}bimodule.  A complex $(\VS[P]_n, \delta)_{n \ge 0}$
of \Mpn{\MA}bimodules with an augmentation $\VS[P]_0 \prto \VS[P]_{-1} \defeq
\VS$ is called an \emph{allowable resolution} of~$\VS$ if the augmented
complex $(\VS[P]_n, \delta)_{n\ge-1}$ is contractible as a complex of
bornological vector spaces.  That is, there are bounded linear maps $h_n
\colon \VS[P]_n \to \VS[P]_{n+1}$, $n \ge -1$, such that $[(h_n)_{n \ge -1},
\delta] = \ID$.  If all the modules $\VS[P]_n$, $n \ge 0$, in an allowable
resolution are free (projective) we have an \emph{allowable free (projective)
resolution} of~$\VS$.

\begin{theorem}[Comparison Theorem]  \label{the:comparison}
  Let $\VS[P]_\bullet \to \VS$ and $\VS[Q]_\bullet \to \VS'$ be complexes of
  \Mpn{\MA}modules over $\VS$ and\/~$\VS'$, respectively.  Let $f \colon \VS
  \to \VS'$ be an \Mpn{\MA}module homomorphism.  Assume that~$\VS[P]_\bullet$
  is projective and that~$\VS[Q]_\bullet$ is an allowable resolution
  of~$\VS'$.  Then~$f$ can be lifted to a chain map $\VS[P]_\bullet \to
  \VS[Q]_\bullet$.  Any two liftings of~$f$ are chain homotopic.

  In particular, any two allowable projective resolutions of~$\VS$ are
  homotopy equivalent.
\end{theorem}

The usual proof for modules over a ring~\cite{MacLane:Homology} carries over
easily.

If~$\VS$ is a right and~$\VS[W]$ is a left \Mpn{\MA}module, let $\VS
\hot_{\MA} \VS[W]$ be the quotient of $\VS \hot \VS[W]$ by the range of the
bounded linear map $\VS \hot \MA \hot \VS[W] \to \VS \hot \VS[W]$ sending $\vs
\otimes \ma \otimes \vs[w] \mapsto \vs \cdot \ma \otimes \vs[w] - \vs \otimes
\ma \cdot \vs[w]$.  The space $\VS \hot_{\MA} \VS[W]$ may fail to be
separated.  If $\VS$ is a free right \Mpn{\MA}module, that is, $\VS \cong \VS'
\hot \Unse{\MA}$ for some complete bornological vector space~$\VS'$, then we
have a bornological isomorphism $\VS \hot_{\MA} \VS[W] \cong \VS' \hot
\VS[W]$.  In particular, $\VS \hot_{\MA} \VS[W]$ is separated.  The same
applies if~$\VS[W]$ is a free left \Mpn{\MA}module.

If~$\VS$ is an \Mpn{\MA}bimodule, let $[\VS, \MA]$ be the range of the bounded
homomorphism $\VS \hot \MA \to \VS$ sending $\vs \otimes \ma \mapsto [\vs,
\ma] \defeq \vs \cdot \ma - \ma \cdot \vs$.  Write $\VS / [,] \defeq \VS /
[\VS, \MA]$ for the \emph{commutator quotient} of~$\VS$.  Although $[\VS,
\MA]$ is bigger than the purely algebraic linear span of the commutators,
$[\VS, \MA]$ need not be bornologically closed.  Thus $\VS / [,]$ need not be
separated.  If $\VS \cong \Unse{\MA} \hot \VS'$ is free as a left
\Mpn{\MA}module, then $\VS / [,] \cong \VS'$.  The same applies if~$\VS$ is
free as a right \Mpn{\MA}module.

\subsection{The bar resolution and \Hochschild homology}
\label{app:bar}

Let~$\MA$ be a complete bornological algebra.  Define $\Barn_n(\MA) \defeq
\Unse{\MA} \hot \MA^{\hot n} \hot \Unse{\MA}$ for $n \ge 0$, $\Barn_{-1} (\MA)
\defeq \Unse{\MA}$.  For $n \ge 0$, $\Barn_n(\MA)$ is a free
\Mpn{\MA}bimodule; $\Barn_{-1}(\MA)$ is an \Mpn{\MA}bimodule with respect to
the usual module structure.  For $n \ge 0$, define $b' = b'_n \colon
\Barn_n(\MA) \to \Barn_{n-1}(\MA)$ by
\begin{equation}  \label{eq:def_bPrime}
  b'( \ma_0 \otimes \dots \otimes \ma_{n+1}) \defeq
  \sum_{j=0}^n (-1)^j \ma_0 \otimes \dots \otimes \ma_{j-1} \otimes
  (\ma_j \cdot \ma_{j+1}) \otimes \ma_{j+2} \otimes \dots \otimes
  \ma_{n+1}.
\end{equation}
The linear map~$b'_n$ exists and is bounded because of the boundedness of the
multiplication in~$\MA$ and the universal property of the completed
bornological tensor product.  By definition, $b'_n$ is an \Mpn{\MA}bimodule
homomorphism for all $n\ge0$.  One verifies as usual that $b'_n \circ
b'_{n+1}=0$.  That is, $\Barn_\bullet(\MA) \to \Unse{\MA}$ is a complex of
\Mpn{\MA}bimodules over~$\Unse{\MA}$.  In fact, $\Barn_\bullet(\MA)$ is a
resolution, called the \emph{bar resolution of~$\MA$}.  A natural contracting
homotopy $h_n^L \colon \Barn_n(\MA) \to \Barn_{n+1}(\MA)$ is defined by
$$
h_n^L(\ma_0 \otimes \dots \otimes \ma_{n+1}) \defeq
1 \otimes [\ma_0] \otimes \dots \otimes \ma_{n+1}.
$$
Here $[\ma_0] \defeq 0$ if $\ma_0 \in \C \subset \Unse{\MA}$ and $[\ma_0]
\defeq \ma_0$ if $\ma_0 \in \MA$.  The boundedness of~$h_n^L$ is trivial.  A
computation shows that $h^L \circ h^L = 0$ and $h^L \circ b' + b' \circ h^L =
\ID$.  That is, $h^L$ is a contracting homotopy.  By definition, $h^L$
consists of right \Mpn{\MA}module homomorphisms.  Another contracting
homotopy~$h_\bullet^R$ consisting of left \Mpn{\MA}module homomorphisms is
defined symmetrically by
$$
h_n^R(\ma_0 \otimes \dots \otimes \ma_{n+1}) \defeq
(-1)^{n+1} \ma_0 \otimes \dots \otimes [\ma_{n+1}] \otimes 1.
$$

\begin{proposition}  \label{pro:bar}
  $\Barn_\bullet(\MA)$ is an allowable free \Mpn{\MA}bimodule resolution
  of~$\Unse{\MA}$.  It is contractible as a resolution of left or right
  \Mpn{\MA}bimodules.
\end{proposition}

Since $\Barn_n(\MA)$ is free, we have $\Barn_n(\MA) / [,] \cong \Unse{\MA}
\hot \MA^{\hot n}$.  The differential~$b'$ induces on the commutator quotients
the differential~$b$ defined by
\begin{multline}  \label{eq:def_b}
  b(\ma_0 \otimes \dots \otimes \ma_n) \defeq
  \sum_{j=0}^{n-1} (-1)^j \ma_0 \otimes \dots \otimes \ma_{j-1} \otimes
  (\ma_j\cdot\ma_{j+1}) \otimes \ma_{j+2} \otimes \dots \otimes \ma_{n}
  \\ {} +
  (-1)^n (\ma_n\cdot\ma_0) \otimes \ma_1 \otimes \dots \otimes \ma_{n-1}.
\end{multline}
Define $\tilde\Omega^n\MA \defeq \Unse{\MA} \hot \MA^{\hot n}$, $\Omega^n\MA
\defeq \Unse{\MA} \hot \MA^{\hot n}$ for $n \ge 1$ and $\Omega^0\MA \defeq
\MA$.  Thus $\tilde{\Omega}^n\MA = \Omega^n \MA$ for $n\ge 1$ and
$$
\Barn_\bullet(\MA) / [,] \cong
(\tilde{\Omega}^n\MA, b) \cong
(\Omega^n\MA, b) \oplus \C[0].
$$
The homology of the complex $(\Omega^n\MA, b)$ is the \emph{\Hochschild
  homology} $\HH_\ast(\MA)$ of~$\MA$, its cohomology is the \emph{\Hochschild
  cohomology} $\HH^\ast(\MA)$ of~$\MA$.  If $\VS[P]_\bullet \to \Unse{\MA}$ is
another allowable projective \Mpn{\MA}bimodule resolution of~$\Unse{\MA}$,
then $\Barn_n(\MA)$ is chain homotopic to~$\VS[P]_\bullet$ as a complex of
\Mpn{\MA}bimodules by the Comparison Theorem~\ref{the:comparison}.  Therefore,
the associated commutator quotients $\VS[P]_\bullet/[,]$ and $(\Omega^n\MA, b)
\oplus \C[0]$ are chain homotopic and have the same homology and cohomology.

Since $h^R \circ h^R = 0$ and $[h^R, b'] = 0$, it follows that
$\Barn_\bullet(\MA) \cong \Ker (h^R) \oplus \Ker b'$ as vector spaces.  The
restriction of~$b'$ to $\Ker h^R$ is an isomorphism onto $\Ker b'$.  Of
course, the kernel of~$h_n^R$ is equal to $\Unse{\MA} \hot \MA^{\hot n} \hot
1$ for all $n \ge 0$ and thus naturally isomorphic to $\tilde\Omega^n\MA$.
Thus $\tilde\Omega^n\MA \cong \Ker b'$ carries a natural \Mpn{\MA}bimodule
structure.  The left multiplication is the obvious one, $\ma \cdot
(\opt{\ma_0} \otimes \dots \otimes \ma_{n}) = (\ma \cdot \opt{\ma_0}) \otimes
\dots \otimes \ma_n$.  The right multiplication is
$$
(\opt{\ma_0} \otimes \dots \otimes \ma_n) \cdot \ma =
(h_n^R)^{-1} (\ID - h_n^R\circ b'_{n+1})
(\opt{\ma_0} \otimes \dots \otimes \ma_n \otimes \ma) =
(-1)^n b'_{n+1}(\opt{\ma_0} \otimes \dots \otimes \ma_n \otimes \ma).
$$
Here~$b'_{n+1}$ is viewed as a map $b'_{n+1} \colon \tilde{\Omega}^{n+1}\MA
\to \tilde{\Omega}^{n}\MA$ using the same formula~\eqref{eq:def_bPrime}.  We
can replace $\tilde{\Omega}^n\MA$ by $\Omega^n\MA$.  If we write
$\opt{\ma_0} d\ma_1 \dots d\ma_n$ for $\opt{\ma_0} \otimes \dots \otimes \ma_n
\in \Omega^n\MA$, then we get
$$
b'(\omega d\ma) = (-1)^{\deg \omega} \omega \cdot \ma, \qquad
b(\omega d\ma) = (-1)^{\deg \omega} [\omega, \ma]
$$
for all homogeneous $\omega \in \Unse{(\Omega\MA)}$ and $\ma \in \MA$.

Let $(\VS[P]_\bullet, \partial_\bullet)$ be an allowable resolution.  Then
$\Ker \partial_n \injto \VS[P]_n \prto \Ker \partial_{n-1}$ are allowable
extensions for all $n \in \N$ and~$\VS[P]_\bullet$ can be obtained by
concatenating these extensions.  Especially, the bar resolution yields
allowable extensions
\begin{equation}  \label{eq:bar_split}
  \begin{gathered}
    \xymatrix{
      {\tilde\Omega^{n+1}\MA\;} \ar@{>->}[r]^-{\iota_n} &
        {\Barn_{n}(\MA)} \ar@{->>}[r]^-{\pi_n} &
          {\tilde\Omega^{n}\MA}
      }
    \\
    \iota_n(\opt{\ma_0} \otimes \ma_1 \otimes \dots \otimes \ma_{n+1}) =
    b'_{n+2}( \opt{\ma_0} \otimes \dots \otimes \ma_{n+1} \otimes 1)
    \\
    \pi_n(\opt{\ma_0} \otimes \ma_1 \otimes \dots \otimes \ma_n \otimes
    \opt{ \ma_{n+1}}) =
    (\opt{\ma_0} \otimes \ma_1 \otimes \dots \otimes \ma_n)
    \cdot \opt{\ma_{n+1}}
  \end{gathered}
\end{equation}
for all $n \ge 0$.

\begin{deflemma}  \label{deflem:n_dimensional}
  A complete bornological algebra~$\MA$ is \emph{\Mpn{n}dimensional}, $n \in
  \Z_+$, iff it satisfies one of the following equivalent conditions:
  \begin{enumerate}[(i)]%

  \item $\tilde\Omega^n\MA$ is a projective \Mpn{\MA}bimodule;

  \item the extension $\tilde\Omega^{n+1}\MA \injto \Barn_n(\MA) \prto
    \tilde\Omega^n\MA$ in~\eqref{eq:bar_split} splits as an extension of
    \Mpn{\MA}bimodules;

  \item there is a bounded linear map $\varphi \colon \MA^{\hot n} \to
    \Omega^{n+1}\MA$ satisfying
    \begin{multline}
      \label{eq:def_varphi_n}
      \ma_0 \cdot \varphi(\ma_1 \otimes \dots \otimes \ma_n) -
      \sum_{j=0}^{n-1} (-1)^j \varphi(\ma_0 \otimes \dots \otimes
      (\ma_j\ma_{j+1}) \otimes \dots \otimes \ma_n)
      \\ \quad {} +
      (-1)^{n+1} \varphi(\ma_0 \otimes \dots \otimes \ma_{n-1}) \cdot \ma_n =
      d\ma_0 \dots d\ma_{n}
      \qquad \forall \ma_0, \dots, \ma_n \in \MA.
    \end{multline}

  \item there is a bounded linear map $\nabla \colon \tilde\Omega^n\MA \to
    \Omega^{n+1}\MA$ satisfying
    \begin{equation}  \label{eq:def_nabla_n}
      \nabla(\ma \cdot \omega) = \ma \cdot \nabla(\omega),\quad
      \nabla(\omega \cdot \ma) = \nabla(\omega) \cdot \ma + (-1)^n \omega d\ma
      \qquad
      \forall \ma \in \MA,\ \omega \in \tilde\Omega^n\MA;
    \end{equation}
    
  \item $\Unse{\MA}$ has an allowable projective \Mpn{\MA}bimodule resolution
    $0\to \VS[P]_n \to \VS[P]_{n-1} \to \cdots \to \VS[P]_0 \to \Unse{\MA}$ of
    length~$n$.

  \end{enumerate}
\end{deflemma}

\begin{proof}
  If $\tilde\Omega^n\MA$ is projective, then~\eqref{eq:bar_split} splits by an
  \Mpn{\MA}bimodule homomorphism.  Conversely, if~\eqref{eq:bar_split} splits
  by an \Mpn{\MA}bimodule homomorphism, then $\tilde\Omega^n\MA$ is projective
  as a direct summand of the free \Mpn{\MA}bimodule $\Barn_n(\MA)$.  Thus (i)
  and~(ii) are equivalent.
  
  The extension~\eqref{eq:bar_split} splits iff there is an \Mpn{\MA}bimodule
  homomorphism $f \colon \Unse{\MA} \hot \MA^{\hot n} \hot \Unse{\MA} =
  \Barn_n(\MA) \to \tilde\Omega^{n+1}\MA = \Omega^{n+1}\MA$ satisfying $f
  \circ \iota_n = \ID$.  The bimodule homomorphisms $f \colon \Unse{\MA} \hot
  \MA^{\hot n} \hot \Unse{\MA} \to \Omega^{n+1}\MA$ correspond to bounded
  linear maps $\varphi \colon \MA^{\hot n} \to \Omega^{n+1}\MA$ via
  $f(\opt{\ma_0} \otimes \ma_1 \otimes \dots \otimes \ma_n \otimes
  \opt{\ma_{n+1}}) \defeq \opt{\ma_0} \cdot \phi(\ma_1 \otimes \dots \otimes
  \ma_n) \cdot \opt{\ma_{n+1}}$.  One verifies easily that $f\circ \iota_n =
  \ID$ is equivalent to~\eqref{eq:def_varphi_n}.  Thus~(ii) is equivalent
  to~(iii).
  
  Since $\tilde\Omega^n\MA$ is the free left \Mpn{\MA}module on $\MA^{\hot
    n}$, left \Mpn{\MA}module homomorphisms $\nabla \colon \tilde\Omega^n\MA
  \to \Omega^{n+1}\MA$ correspond to bounded linear maps $\varphi \colon
  \MA^{\hot n} \to \Omega^{n+1}\MA$ via $\nabla (\opt{\ma_0} d\ma_1 \dots
  d\ma_n) \defeq \opt{\ma_0} \cdot \varphi (\ma_1 \otimes \dots \otimes
  \ma_n)$.  One verifies easily that~$\varphi$
  satisfies~\eqref{eq:def_varphi_n} iff~$\nabla$
  satisfies~\eqref{eq:def_nabla_n}.  Thus~(iii) and~(iv) are equivalent.
  
  If~$\tilde\Omega^n\MA$ is a projective \Mpn{\MA}bimodule, then~$\Unse{\MA}$
  has an allowable projective \Mpn{\MA}bimodule resolution of length~$n$,
  namely $0 \to \tilde\Omega^n\MA \to \Barn_{n-1}(\MA) \to \cdots \to
  \Barn_0(\MA) \to \Unse{\MA}$.  Thus (i) implies~(v).  The converse direction
  is usually proved using the Ext functors.  If we cared to define and study
  these functors, we could follow that line.  Instead, we use the following
  Theorem~\ref{the:Schanuel} of Schanuel.
\end{proof}

\begin{theorem}  \label{the:Schanuel}
  Let $\VS[K] \injto \VS[P] \overset{\partial}{\prto} \VS$ and $\VS[K]' \injto
  \VS[P]' \overset{\partial'}{\prto} \VS$ be allowable resolutions of
  length~$1$ with projective $\VS[P]$ and~$\VS[P]'$.  Then $\VS[K] \oplus
  \VS[P]' \cong \VS[K]' \oplus \VS[P]$.  Thus~$\VS[K]$ is projective
  iff\/~$\VS[K]'$ is projective.
  
  More generally, let $\VS[K]_n \injto \VS[P]_{n-1} \to \dots \to \VS[P]_0
  \prto \VS$ and $\VS[K]'_n \injto \VS[P]'_{n-1} \to \dots \to \VS[P]'_0 \prto
  \VS'$ be allowable resolutions with all $\VS[P]_j$ and all~$\VS[P]'_j$
  projective.  Assume that there are projective modules $\VS[Q]$ and~$\VS[Q]'$
  such that $\VS \oplus \VS[Q] \cong \VS' \oplus \VS[Q]'$.  Then there are
  projective modules $\VS[R]$ and~$\VS[R]'$ such that $\VS[K]_n \oplus \VS[R]
  \cong \VS[K]'_n \oplus \VS[R]'$.  Thus~$\VS[K]_n$ is projective
  iff\/~$\VS[K]'_n$ is projective.
\end{theorem}

\begin{proof}
  Let $\VS[P]'' \defeq \VS[P] \oplus \VS[P]'$ and let $\partial'' \colon
  \VS[P]'' \to \VS$ be the map $\partial \oplus \partial'$.  The claim is that
  both $\VS[K] \oplus \VS[P]'$ and $\VS[K]' \oplus \VS[P]$ are isomorphic as
  modules to the kernel of $\partial''$.  It suffices to prove this for
  $\VS[K]' \oplus \VS[P]$.  As vector spaces, $\VS[P] \cong \VS[K] \oplus \VS$
  and $\VS[P]' \cong \VS[K]' \oplus \VS$ and therefore $\VS[P]'' \cong \VS
  \oplus \VS[K]' \oplus \VS[P]$ and $\Ker \partial'' \cong \VS[K]' \oplus
  \VS[P]$.  Thus we have an allowable extension of modules $\VS[K]' \injto
  \Ker \partial'' \prto \VS[P]$.  Since~$\VS[P]$ is projective, this extension
  splits.  That is, $\Ker \partial'' \cong \VS[K]' \oplus \VS[P]$ as modules.
  Since a module is projective iff it is a direct summand of a projective
  module, it follows that~$\VS[K]$ is projective iff~$\VS[K]'$ is projective.
  
  Since we can get long allowable resolutions by concatenating allowable
  extensions, it suffices to prove the second assertion for the case $n=1$.
  Adding the allowable resolutions $\ID \colon \VS[Q] \to \VS[Q]$ and $\ID
  \colon \VS[Q]' \to \VS[Q]'$, we can reduce to the special case $\VS \cong
  \VS'$.  In that special case we have $\VS[K] \oplus \VS[P]' \cong \VS[P]
  \oplus \VS[K]'$.  The second assertion follows.
\end{proof}

The map~$\nabla$ in~(iv) is a \emph{graded right connection} in the following
sense:

\begin{definition}  \label{def:right_connection}
  Let~$\MA$ be a complete bornological algebra and let~$\VS$ be a graded
  \Mpn{\MA}bimodule.  A \emph{graded right connection} on~$\VS$ is a bounded
  linear map map $\nabla\colon \VS \to \VS \hot_{\MA} \Omega^1\MA$ satisfying
  \begin{equation}  \label{eq:def_right_connection}
    \nabla(\ma \cdot \vs) =
    \ma \cdot \nabla(\vs), \quad
    \nabla(\vs \cdot \ma) =
    \nabla(\vs) \cdot \ma + (-1)^{\deg\vs} \vs \otimes d\ma \qquad
    \forall \ma \in \MA,\ \vs \in \VS;
  \end{equation}
\end{definition}

Of course, $\tilde\Omega^n\MA$ is graded by $\deg(\omega) = n$ for all $\omega
\in \tilde\Omega^n\MA$.  Condition~\eqref{eq:def_nabla} asserts that~$\nabla$
is a graded right connection on $\tilde\Omega^n\MA$ with respect to this
grading because $\Omega^{n+1} \MA \cong \tilde\Omega^n \MA \hot_{\MA}
\Omega^1\MA$.  In \cite{cuntz95:algebra} and~\cite{connes94:ncg}, no grading
is taken into account in the definition of a connection.  A map $\nabla \colon
\tilde\Omega^n\MA \to \Omega^{n+1} \MA$ is a graded connection on
$\tilde\Omega^n\MA$ iff $(-1)^n \nabla$ is an ungraded connection.

\subsection{The tensor algebra and differential forms}
\label{app:Omega_Tens}

Let~$\MA$ be a complete bornological algebra.  The \emph{differential
  envelope} $\Omega\MA$ is the universal complete bornological differential
algebra generated by~$\MA$.  We can describe $\Omega\MA$ in a concrete way as
follows.  Let $\Omega^n \MA \defeq \Unse{\MA} \hot \MA^{\hot n}$.  Then
$\Omega\MA \defeq \sum_{j=0}^\infty \Omega^j\MA$ with the direct sum
bornology.  The usual differential and multiplication of non-commutative
differential forms are bounded with respect to the direct sum bornology and
thus make $\Omega\MA$ into a complete bornological differential algebra.  The
multiplication is bounded because direct sums commute with bornological tensor
products.

The \emph{Fedosov product} on $\Omega\MA$ is defined by
\begin{equation}  \label{eq:def_Fedosov}
  x \odot y \defeq x \cdot y - (-1)^{\deg x} dx dy
  \qquad \forall x,y \in \Omega\MA.
\end{equation}
It is bounded because the multiplication, differential, and grading on
$\Omega\MA$ are bounded.  It is trivial to verify that $\odot$ is associative.
Decompose $\Omega\MA$ into the differential forms of even degree
$\Omega^\even\MA$ and the differential forms of odd degree $\Omega^\odd\MA$.
The Fedosov product of two even forms is again even, so that
$(\Omega^\even\MA, \odot)$ is a complete bornological algebra.  Let
$\sigma_{\MA} \colon \MA\to (\Omega^\even\MA, \odot)$ be the natural bounded
linear map coming from the isomorphism $\MA \cong \Omega^0\MA$.

\begin{proposition}  \label{pro:Tens_universal}
  Let $l \colon \MA \to \MA[B]$ be a bounded linear map between complete
  bornological algebras.  There is a unique bounded homomorphism $\LLH{l}
  \colon (\Omega^\even\MA, \odot) \to \MA[B]$ satisfying $\LLH{l} \circ
  \sigma_{\MA} = l$.  The homomorphism~$\LLH{l}$ is uniquely determined by
  \begin{equation}
    \label{eq:LLH_Tens}
    \LLH{l}(\opt{\ma_0} d\ma_1 \dots d\ma_{2n}) \defeq
    l\opt{\ma_0} \omega_l(\ma_1, \ma_2) \cdots \omega_l(\ma_{2n-1}, \ma_{2n}).
  \end{equation}
\end{proposition}

\begin{proof}
  It is clear that~$\LLH{l}$ above defines a bounded linear map
  $\Omega^\even\MA \to \MA[B]$ with $\LLH{l} \circ \sigma_{\MA} = l$.
  Furthermore, it is the only map that has a chance to be multiplicative: If
  $\psi \colon (\Omega^\even\MA, \odot) \to \MA[B]$ is a homomorphism with
  $\psi \circ \sigma_{\MA} = l$, then necessarily $\psi \circ
  \omega_\sigma (\ma_1, \ma_2) = \omega_l (\ma_1, \ma_2)$ and thus
  $$
  \psi(\opt{\ma_0} d\ma_1 \dots d\ma_{2n}) = \psi\bigl(\opt{\ma_0} \odot
  \omega_\sigma (\ma_1, \ma_2) \odot \cdots \odot \omega_\sigma (\ma_{2n-1},
  \ma_{2n}) \bigr) = l\opt{\ma_0} \cdot \omega_l (\ma_1, \ma_2) \cdots
  \omega_l (\ma_{2n-1}, \ma_{2n}).
  $$
  It remains to verify that $\hat{l} \defeq \LLH{l}$ is indeed
  multiplicative.  Let $X \subset \Omega^\even\MA$ be the set of all $\omega'
  \in \Omega^\even\MA$ satisfying $\hat{l} (\omega' \odot \omega) = \hat{l}
  (\omega') \cdot \hat{l} (\omega)$ for all $\omega \in \Omega^\even\MA$.
  Evidently, $X$ is a bornologically closed subalgebra of $\Omega^\even\MA$.
  It is easy to verify that $\hat{l}(\ma \odot \omega) = \hat{l}(\ma) \cdot
  \hat{l}(\omega)$ for all $\ma \in \MA$ and~$\omega$ of the form $\opt{\ma_0}
  d\ma_1 \dots d\ma_{2n}$.  It follows that $\hat{l} (\ma \odot \omega) =
  \hat{l}(\ma) \cdot \hat{l} (\omega)$ for all $\ma \in \MA$, $\omega \in
  \Omega^\even\MA$, that is, $\sigma(\MA) \subset X$.  Consequently,
  $\omega_\sigma (\MA, \MA) \subset \MA \odot \MA + \MA \subset X$ and
  therefore $\Unse{\MA} \odot \omega_\sigma (\MA, \MA)^{\odot n} \subset X$
  for all $n \in \N$.  Hence $X = \Omega^\even\MA$.  That is, $\hat{l}$ is
  multiplicative.
\end{proof}

The universal property described in Proposition~\ref{pro:Tens_universal} is
precisely the universal property of the \emph{tensor algebra} $\Tens\MA$
of~$\MA$.  Hence we get a natural isomorphism $\Tens\MA \cong (\Omega^\even
\MA, \odot)$ as in~\cite{cuntz95:algebra}.  Let $(\Jens\MA)^k \defeq
\sum_{n=k}^\infty \Omega^{2n}\MA$.  The natural projection $\tau_{\MA} \colon
\Tens\MA \to \MA$ with kernel $\Jens\MA$ is the unique bounded homomorphism
satisfying $\tau_{\MA} \circ \sigma_{\MA} = \ID$.

\subsection{Quasi-free algebras}
\label{app:quasi_free}

A complete bornological algebra~$\MA[N]$ is \emph{\Mpn{k}nilpotent} iff
$\MA[N]^k = \{0\}$.  Thus~$\MA[N]$ is \Mpn{2}nilpotent iff the multiplication
in~$\MA$ is the zero map.  We call~$\MA[N]$ \emph{nilpotent} iff it is
\Mpn{k}nilpotent for some $k \in \N$.  An extension of complete bornological
algebras $\MA[K] \injto \MA[E] \prto \MA[Q]$ is called \emph{\Mpn{k}nilpotent}
iff~$\MA[K]$ is \Mpn{k}nilpotent and \emph{nilpotent} iff~$\MA[K]$ is
nilpotent.

\begin{deflemma}  \label{deflem:quasi_free}
  A complete bornological algebra~$\MA[R]$ is \emph{quasi-free} iff it
  satisfies one of the equivalent conditions:
  \begin{enumerate}[(i)]%

  \item \label{qf_i}
    There is a bounded linear map $\varphi \colon \MA[R] \to \Omega^2\MA[R]$
    satisfying
    \begin{equation}
      \label{eq:def_varphi}
      \varphi(x_1 x_2) =
      x_1 \varphi(x_2) + \varphi(x_1) x_2 - dx_1 dx_2
      \qquad \forall x_1, x_2 \in \MA[R].
    \end{equation}
    
  \item \label{qf_viii}
    There is a bounded linear map $\nabla \colon \Omega^1\MA[R] \to
    \Omega^2\MA[R]$ satisfying
    \begin{equation}  \label{eq:def_nabla}
      \nabla(x \cdot \omega) = x \cdot \nabla(\omega),\quad
      \nabla(\omega \cdot x) = \nabla(\omega) \cdot x - \omega dx \qquad
      \forall x \in \MA[R],\ \omega \in \Omega^1\MA[R].
    \end{equation}
    
  \item \label{qf_vii}
    The \Mpn{\MA[R]}bimodule $\Omega^1 \MA[R]$ is projective.

  \item \label{qf_ix}
    There is an allowable projective \Mpn{\MA[R]}bimodule resolution of
    length~$1$
    $$
    0 \longrightarrow
    P_1 \overset{\partial_1}{\longrightarrow}
    P_0 \overset{\partial_0}{\longrightarrow}
    \Unse{\MA[R]}
    $$
    of\/~$\Unse{\MA[R]}$ with the standard bimodule structure.

  \item \label{qf_ii}
    There is a bounded splitting homomorphism $\upsilon \colon \MA[R] \to
    \Tens\MA[R]/ (\Jens\MA[R])^2$ for $\tau \colon \Tens\MA[R] /
    (\Jens\MA[R])^2 \to \MA[R]$.
    
  \item \label{qf_iii}
    Let $\MA[K] \injto \MA[E] \prto \MA[Q]$ be an allowable \Mpn{2}nilpotent
    extension and $f \colon \MA[R] \to \MA[Q]$ a bounded homomorphism; then
    there is a bounded homomorphism $\hat{f} \colon \MA[R] \to \MA[E]$
    lifting~$f$, that is, making the following diagram commute:
    \begin{equation}  \label{eq:lift_nilpotent}
      \begin{gathered}
        \xymatrix{
          & & {\MA[R]} \ar[d]^{f} \ar@{-->}[dl]_{\hat{f}} \\
          {\MA[K]\;} \ar@{>->}[r]^{i} &
            {\MA[E]} \ar@{->>}[r]^{p} &
              {\MA[Q].} \ar@/^/@{.>}[l]^{s}
          }
      \end{gathered}
    \end{equation}

  \item \label{qf_iv}
    Each allowable \Mpn{2}nilpotent extension $\MA[K] \injto \MA[E] \prto
    \MA[R]$ has a bounded splitting homomorphism.

  \item \label{qf_v}
    There is a family of bounded homomorphisms $\upsilon_n \colon \MA[R] \to
    \Tens\MA[R] / (\Jens\MA[R])^n$, $n \ge 1$, such that $\upsilon_1 =
    \ID[{\MA[R]}]$ and each~$\upsilon_{n+1}$ lifts~$\upsilon_n$.  That is, the
    following diagram commutes:
    \begin{displaymath}
      \xymatrix{
        {\MA[R]} \ar[d]_{\upsilon_1}^{\cong} \ar[dr]_{\upsilon_2}
          \ar[drrr]_{\upsilon_n} \ar[drrrr]^{\upsilon_{n+1}} \\
        {\Tens\MA[R]/\Jens\MA[R]} &
          {\Tens\MA[R]/(\Jens\MA[R])^2} \ar[l] &
            {\cdots} \ar[l] &
              {\Tens\MA[R]/(\Jens\MA[R])^n} \ar[l] &
                {\Tens\MA[R]/(\Jens\MA[R])^{n+1}.} \ar[l] &
                  {\cdots} \ar[l]
        }
    \end{displaymath}
    
  \item \label{qf_vi}
    Let $\MA[K] \injto \MA[E] \prto \MA[Q]$ be an allowable nilpotent
    extension and $f \colon \MA[R] \to \MA[Q]$ a bounded homomorphism.  Then
    there is a bounded homomorphism $\hat{f} \colon \MA[R] \to \MA[E]$
    lifting~$f$ as in~\eqref{qf_iii}.

  \end{enumerate}
\end{deflemma}

\begin{proof}
  The equivalence of (i)--(iv) is the special case $n = 1$ of
  Lemma~\ref{deflem:n_dimensional}.
  
  \textbf{(i) is equivalent to~(v)}.  Identify $\Tens\MA[R] / (\Jens\MA[R])^2
  \cong \MA[R] \oplus \Omega^2\MA[R]$ as a bornological vector space.  The
  multiplication is the Fedosov product, where terms of degree~$4$ or higher
  are omitted.  The bounded linear sections for the natural projection
  $\Tens\MA[R] / (\Jens\MA[R])^2\to \MA[R]$ are of the form $\sigma_{\MA[R]} +
  \varphi$ with bounded linear maps $\varphi \colon \MA[R] \to
  \Omega^2\MA[R]$.  The section $\sigma_{\MA[R]} + \varphi$ is multiplicative
  iff~$\varphi$ satisfies~\eqref{eq:def_varphi} because
  $$
  (\sigma + \varphi)(x_1 x_2) - (\sigma + \varphi)(x_1) \odot
  (\sigma + \varphi)(x_2) =
  dx_1 dx_2 + \varphi(x_1 x_2) - \varphi(x_1) \cdot x_2
  - x_1 \cdot \varphi(x_2)
  $$
  up to terms of degree~$4$ and higher.  Thus \eqref{qf_i}
  and~\eqref{qf_ii} are equivalent.
  
  \textbf{(v) implies~(vi)}.  Let $\upsilon \colon \MA[R] \to \Tens\MA[R] /
  (\Jens\MA[R])^2$ be a bounded splitting homomorphism.  Let $\MA[K] \injto
  \MA[E] \prto \MA[Q]$, $f,i,p,s$ be as in~\eqref{eq:lift_nilpotent}.
  Consider the bounded linear map $s \circ f \colon \MA[R] \to \MA[E]$.  Its
  curvature is $\omega_{s \circ f} (x_1, x_2) = \omega_s \bigl( f(x_1), f(x_2)
  \bigr)$ because~$f$ is a homomorphism.  Since $p \circ s = \ID[{\MA[Q]}]$ is
  a homomorphism, it follows that $\omega_{s \circ f}$ maps~$\MA[R]$
  into~$\MA[K]$.  Thus $\omega_{s \circ f} (x_1, x_2) \cdot \omega_{s \circ f}
  (x_3, x_4) = 0$ for all $x_1, \dots, x_4 \in \MA[R]$.  Consequently, the
  bounded homomorphism $\LLH{s \circ f} \colon \Tens\MA[R] \to \MA[E]$ of
  Proposition~\ref{pro:Tens_universal} annihilates $(\Jens\MA[R])^2$.
  Therefore, it descends to a bounded homomorphism $g \colon \Tens\MA[R] /
  (\Jens\MA[R])^2 \to \MA[E]$.  Furthermore, $p \circ g = f\circ \tau$ with
  $\tau \colon \Tens\MA[R] / (\Jens\MA[R])^2 \to \MA[R]$ the natural
  projection.  The composition $\hat{f} \defeq g \circ \upsilon \colon \MA[R]
  \to \MA[E]$ is a bounded homomorphism satisfying $p \circ \hat{f} = (p \circ
  g) \circ \upsilon = f \circ (\tau \circ \upsilon) = f$.  Thus~\eqref{qf_ii}
  implies~\eqref{qf_iii}.
  
  Of course, \eqref{qf_iv} is the special case $\MA[Q] = \MA[R]$ and $f =
  \ID[{\MA[R]}]$ of~\eqref{qf_iii}.  Thus~\eqref{qf_iii}
  implies~\eqref{qf_iv}.  Since $\Jens\MA[R] / (\Jens\MA[R])^2 \injto
  \Tens\MA[R] / (\Jens\MA[R])^2 \prto \MA[R]$ is an allowable \Mpn{2}nilpotent
  extension, \eqref{qf_ii} is a special case of~\eqref{qf_iv}.
  
  Therefore, the conditions \eqref{qf_ii}--\eqref{qf_iv} are equivalent.
  
  \textbf{(vi) implies~(viii)}.  The homomorphisms~$\upsilon_n$ are
  constructed by induction, starting with $\upsilon_1 = \ID$.  In each
  induction step we apply~(vi) to the \Mpn{2}nilpotent extension
  $(\Jens\MA[R])^n / (\Jens\MA[R])^{n+1} \injto \Tens\MA[R] /
  (\Jens\MA[R])^{n+1} \prto \Tens\MA[R] / (\Jens\MA[R])^n$ to lift $\upsilon_n
  \colon \MA[R] \to \Tens\MA[R] / (\Jens\MA[R])^n$ to a bounded homomorphism
  $\upsilon_{n+1} \colon \MA[R] \to \Tens\MA[R] / (\Jens\MA[R])^{n+1}$.
  Thus~\eqref{qf_iii} implies~\eqref{qf_v}.

  The implication \eqref{qf_v}$\Longrightarrow$\eqref{qf_vi} is proved as the
  implication \eqref{qf_ii}$\Longrightarrow$\eqref{qf_iii}.
  Since~\eqref{qf_iii} is a special case of~\eqref{qf_vi}, the conditions
  \eqref{qf_ii}--\eqref{qf_vi} are equivalent.
\end{proof}

\begin{digression}
  The sign in front of $-dx_1 dx_2$ in~\eqref{eq:def_varphi} is a matter of
  convention.  In \cite{cuntz95:algebra} and~\cite{cuntz95:cyclic}, the other
  sign is used, but in \cite{cuntz97:excision}
  and~\cite{cuntz97:excision_top}, Cuntz and Quillen switch to the sign
  convention in~\eqref{eq:def_varphi}.
\end{digression}

\section{Cyclic homology and cohomology}
\label{app:cyclic}

We define cyclic homology and cohomology using the Hodge filtration on
$\Omega\MA$.  This definition is equivalent to the more standard definition
using bicomplexes.  First, we introduce the Karoubi operator and the
boundary~$B$ on $\Omega\MA$.

\subsection{The Karoubi operator}
\label{app:kappa}

In order to write down the boundary~$\partial$ of the X-complex explicitly, we
need the \emph{Karoubi operator} $\kappa\colon \Omega\MA\to \Omega\MA$ defined
by
$$
\kappa(\omega d\ma) \defeq (-1)^{\deg \omega} d\ma \cdot \omega \qquad
\forall \omega \in \Omega\MA,\ \ma \in \MA.
$$
It is clear that~$\kappa$ is a bounded operator of degree zero on
$\Omega\MA$.  It is also bounded with respect to the bornology~$\CBS_\an$ and
thus extends to a bounded operator on $\Omega_\an\MA$.  We recall some
computations of Cuntz and Quillen~\cite{cuntz95:operators}.  We have
\begin{multline*}
  (bd +db) (\omega d\ma) =
  (-1)^{\deg \omega + 1} [d\omega, \ma] + (-1)^{\deg \omega} d [\omega, \ma] =
  (-1)^{\deg \omega} \bigl(
    \ma d\omega - (d\omega)\ma + d(\omega \cdot \ma) - d(\ma \cdot \omega)
  \bigr)
  \\ =
  (-1)^{\deg \omega} \bigl(
    -(d\ma) \omega + (-1)^{\deg \omega} \omega d\ma
  \bigr) =
  \omega d\ma - \kappa(\omega d\ma)
\end{multline*}
because~$d$ is a graded derivation.  That is,
\begin{equation}  \label{eq:kappa_i}
  bd+db = \ID - \kappa.
\end{equation}
It follows that~$\kappa$ is homotopic to the identity with respect to both
differentials $b$ and~$d$.  In particular, $\kappa$ commutes with $b$ and~$d$
(because $b^2 = 0 = d^2$).  Iterating~$\kappa$, we get
\begin{equation}  \label{eq:kappa_ii}
  \kappa^j (\opt{\ma_0} d\ma_1 \dots d\ma_n) =
  (-1)^{j(n-1)} d\ma_{n-j+1} \dots d\ma_n \cdot \opt{\ma_0}
  d\ma_1 \dots d\ma_{n-j} \qquad
  \forall 0 \le j \le n.
\end{equation}
In particular, for $j = n$, this yields $\kappa^n (\opt{\ma_0} d\ma_1 \dots
d\ma_n) = d\ma_1 \dots d\ma_n \cdot \opt{\ma_0} = \opt{\ma_0} d\ma_1 \dots
d\ma_n + [d\ma_1 \dots d\ma_n, \opt{\ma_0}] = \opt{\ma_0} d\ma_1 \dots d\ma_n
+ (-1)^n b(d\ma_1 \dots d\ma_n d\opt{\ma_0})$.  That is,
\begin{equation}  \label{eq:kappa_iii}
  \kappa^n = \ID + b \kappa^n d \qquad \text{on $\Omega^n\MA$.}
\end{equation}
Especially, $d^2 = 0$ implies $\kappa^n = \ID$ on $d\Omega^{n-1}\MA \subset
\Omega^n\MA$ and thus
\begin{equation}  \label{eq:kappa_iv}
  \kappa^{n+1} d = d \qquad \text{on $\Omega^n\MA$.}
\end{equation}
Combining this with~\eqref{eq:kappa_iii}, it follows that
\begin{equation}  \label{eq:kappa_vii}
  (\kappa^n - \ID) (\kappa^{n+1} - \ID) =
  b\kappa^n d (\kappa^{n+1} - \ID) =
  b\kappa^n (\kappa^{n+1}d - d) = 0
  \qquad \text{on $\Omega^n\MA$.}
\end{equation}
Moreover, \eqref{eq:kappa_iii}, \eqref{eq:kappa_iv}, and~\eqref{eq:kappa_i}
yield
\begin{equation}  \label{eq:kappa_v}
  \kappa^{n+1} = \kappa + b \kappa^{n+1} d = \kappa + bd = \ID - db
  \qquad \text{on $\Omega^n\MA$.}
\end{equation}
Since $b^2 = 0$, this implies $\kappa^{n+1} = \ID$ on $b(\Omega^{n+1}\MA)
\subset \Omega^n\MA$ and thus
\begin{equation}  \label{eq:kappa_vi}
  \kappa^{n} b = b \qquad \text{on $\Omega^n\MA$.}
\end{equation}
Define Connes's operator $B \colon \Omega\MA \to \Omega\MA$ by
$$
B(\opt{\ma_0} d\ma_1 \dots d\ma_n) \defeq
\sum_{j=0}^n (-1)^{nj}
d\ma_j \dots d\ma_n d\opt{\ma_0} d\ma_1 \dots d\ma_{j-1}
\qquad \forall \ma_0, \dots, \ma_n \in \MA.
$$
That is,
\begin{equation}
  \label{eq:kappa_ix}
  B = \sum_{j=0}^n \kappa^j d \qquad \text{on $\Omega^n\MA$.}
\end{equation}
Thus $B \kappa = \kappa B$ and $B^2 = Bd = dB = 0$.  Using geometric series,
we get
$$
(\kappa^{n})^{n+1} - \ID =
\sum_{j=0}^n \kappa^{nj} (\kappa^n - \ID) =
\sum_{j=0}^n \kappa^{nj} b\kappa^n d =
b\sum_{j=0}^n \kappa^{n(j+1)} d =
b\sum_{j=0}^n \kappa^{n-j} d =
bB
$$
on $\Omega^n\MA$ because $n(j+1) \equiv -(j+1) \equiv n-j \bmod n+1$ and
$\kappa^{n+1}d = d$.  Similarly,
$$
(\kappa^{n+1})^{n} - \ID =
\sum_{j=0}^{n-1} \kappa^{(n+1)j} (\kappa^{n+1} - \ID) =
- \sum_{j=0}^{n-1} \kappa^{(n+1)j} db =
- \sum_{j=0}^{n-1} \kappa^{j} db =
-B b
$$
on $\Omega^n\MA$ because~$b$ lowers degree by~$1$ and $\kappa^n d = d$ on
$\Omega^{n-1}\MA$.  Thus
\begin{equation}  \label{eq:kappa_viii}
  \kappa^{n(n+1)} - 1 = bB = - Bb \qquad \text{on $\Omega^n\MA$.}
\end{equation}
Thus we get a bicomplex if we endow $\Omega\MA$ with the boundaries $(b,B)$.
It follows that
\begin{equation}  \label{eq:Bb_bicomplex}
  (B+b)^2 = B^2 + (Bb+bB) + b^2 = 0.
\end{equation}

\subsection{The Hodge filtration and cyclic cohomology}
\label{app:Hodge}

Let~$\MA$ be a complete bornological algebra.  The \emph{Hodge filtration} is
a certain filtration $F_n(\MA)$ on the bicomplex $(\Omega\MA, b, B)$ such that
the quotients $(\Omega\MA / F_n(\MA), B+b)$ compute the cyclic cohomology
of~$\MA$.  It was introduced by Cuntz and Quillen~\cite{cuntz95:cyclic}.
The Hodge filtration on $\Omega\MA$ is defined by
\begin{equation}  \label{eq:def_Hodge}
  F_n(\MA) \defeq b(\Omega^{n+1}\MA) \oplus \sum_{j=n+1}^\infty \Omega^j\MA.
\end{equation}
Clearly, $F_n(\MA)$ is invariant under $b$ and~$d$.  Hence $F_n(\MA)$ is
invariant under the operator~$\kappa$ constructed from $b$ and~$d$
by~\eqref{eq:kappa_i}.  Thus $F_n(\MA)$ is invariant under the boundary~$B$ and
the spectral operators $P$ and~$H$ constructed from~$\kappa$.  We are mainly
interested in the corresponding quotients
\begin{equation}  \label{eq:def_Xtower}
  \Xtower_n(\MA) \defeq
  \Omega\MA/ F_n(\MA) \cong
  \sum_{j=0}^{n-1} \Omega^j\MA \oplus \Omega^n\MA / b(\Omega^{n+1}\MA).
\end{equation}
The space $\Xtower_n(\MA)$ may fail to be separated.  The operators $b$, $d$,
$\kappa$, $B$, $P$, and~$H$ descend to bounded operators on $\Xtower_n(\MA)$
because they map $F_n(\MA)$ into itself.

We work in $X(\Tanil\MA) \cong \Omega_\an\MA$, not in $\Omega\MA$ throughout
this thesis.  This only makes life more difficult, but since we have
introduced most notation in $\Omega_\an\MA$, it is convenient to define the
Hodge filtration in $\Omega_\an\MA$ as well.  We define $F_n(\MA) \subset
\Omega_\an\MA$ so that $\Omega_\an\MA / F_n(\MA)$ is still the space
$\Xtower_n(\MA)$ above.  Thus we complete the sum $\sum_{j=n+1} \Omega^j\MA$
inside $\Omega_\an\MA$ but do not complete $b(\Omega^{n+1}\MA)$.

We make $\Xtower_n(\MA)$ a \Mpn{\Ztwo}graded complex using the boundary $b+B$
and the even/odd grading.  The complexes $(\Xtower_n(\MA), B+b)$ and
$(\Xtower_n(\MA), \partial)$ are homotopy equivalent.  The argument in
Section~\ref{sec:X_Tanil} carries over easily (and with some simplifications).
Since we divide out all but finitely many $\Omega^n\MA$, the rescaling by
$[j/2]!$ on $\Omega^j\MA$ is bounded on $\Xtower_n(\MA)$.  Since $F_n(\MA)
\subset F_m(\MA)$ for $m\le n$, there are natural projection maps $\Xmap_{nm}
\colon \Xtower_n(\MA) \to \Xtower_m(\MA)$ that make $\bigl( \Xtower_n(\MA)
\bigr)_{n \in \N}$ a \emph{tower of complexes} in the sense
of~\cite{cuntz95:cyclic}.  Especially, for $n = 1$, the complex
$\Xtower_1(\MA)$ is equal to the X-complex of~$\MA$.

If we are interested in properties of the projective system $\bigl(
\Xtower_n(\MA) \bigr)_{n \in \N}$ of complexes, then it is appropriate to view
it as the X-complex of the pro-algebra $\Tpnil \MA$ defined in
Section~\ref{sec:pro_algebras}.  However, this gives little information about
the \emph{individual} complexes $\Xtower_n(\MA)$.  These complexes can be
studied by homological algebra, using their relation to \Hochschild homology.
Therefore, the boundary $B+b$ is preferable over~$\partial$ on
$\Xtower_n(\MA)$ because it is closer to the \Hochschild boundary~$b$.

As in~\cite{cuntz95:cyclic}, we define cyclic (co)homology groups using
the Hodge tower.  Let
\begin{equation}  \label{eq:def_cyclic}
  \begin{aligned}
    \HC_n(\MA) &\defeq H_n \bigl( \Xtower_ n   (\MA) \bigr), \\
    \HC^n(\MA) &\defeq H^n \bigl( \Xtower_ n   (\MA) \bigr),
  \end{aligned} \qquad
  \begin{aligned}
    \HD_n(\MA) &\defeq H_n \bigl( \Xtower_{n+1}(\MA) \bigr), \\
    \HD^n(\MA) &\defeq H^n \bigl( \Xtower_{n+1}(\MA) \bigr).
  \end{aligned}
\end{equation}
By definition, $H^n\bigl( \Xtower_n(\MA) \bigr)$ is the homology of the
dual complex $\Xtower_k'(\MA) \defeq \Lin(\Xtower_k(\MA); \C[0])$ of bounded
linear functionals on $\Xtower_k(\MA)$.  The groups $\HC_n(\MA)$ and
$\HC^n(\MA)$ are the \emph{cyclic (co)homology groups} of~$\MA$, whereas
$\HD_n(\MA)$ and~$\HD^n(\MA)$ are called the \emph{de Rham homology groups}
of~$\MA$.

Consider the complex $G_n \defeq (F_{n-1} (\MA) / F_n (\MA), B+b)$.  This
complex has $b(\Omega^n\MA)$ in degree $n-1 \bmod 2$ and $\Omega^n\MA /
b(\Omega^{n+1}\MA)$ in degree $n \bmod 2$.  The homology of~$G_n$ is equal to
$\HH_n(\MA) = \Ker b|_{\Omega^n\MA} / b(\Omega^{n+1}\MA)$ in degree~$n$
and~$0$ in degree $n-1$.  For each $n \in \Z_+$, we have a (usually not
allowable) extension $G_n \injto \Xtower_n(\MA) \prto \Xtower_{n-1}(\MA)$ of
complexes and get an associated exact sequence
$$
\HD_{n-1}(\MA) \injto
\HC_{n-1}(\MA) \to
\HH_n(\MA) \to
\HC_n(\MA) \prto
\HD_{n-2}(\MA)
$$
for the homology of these complexes.  Splicing these extensions for
different~$n$, we get the homology version of Connes's exact sequence:
$$
\cdots \to 
\HC_{n+1}(\MA) \overset{S}{\longrightarrow}
\HC_{n-1}(\MA) \overset{I}{\longrightarrow}
\HH_n(\MA) \overset{B}{\longrightarrow}
\HC_n(\MA) \overset{S}{\longrightarrow}
\HC_{n-2}(\MA) \to \cdots.
$$
Thus $\HD_n(\MA)$ can be interpreted as the range of the operator $S \colon
\HC_{n+2}(\MA) \to \HC_n(\MA)$.  A dual exact sequence exists for the
cohomology groups.  It allows to interpret $\HD^n(\MA)$ as the range of the
operator $S \colon \HC^n(\MA) \to \HC^{n+2} (\MA)$.  Therefore, the cyclic
cohomology groups $\bigl( \HC^{2n+\ast}(\MA) \bigr)_{n \in \N}$ form an
inductive system of vector spaces for $\ast = 0,1$.  The inductive limit is
equal to $\HP^\ast(\MA)$, the cohomology of the complex $\prod_{j=0}^\infty
\Omega^j\MA$ with boundary $B+b$.  Dually, the homology groups $\bigl(
\HC_{2n+\ast}(\MA) \bigr)_{n \in \N}$ form a projective system.  The limit is
not quite equal to $\HP_\ast(\MA)$, but there is a \Mp{\prolim^1}exact
sequence that computes $\HP_\ast(\MA)$ from this projective system of cyclic
homology groups.

The above definition of cyclic (co)homology yields the same result as the more
standard definition using bicomplexes (see~\cite{cuntz95:cyclic}).  The
bicomplex approach avoids the use of non-separated spaces and non-allowable
extensions and readily gives Connes's exact sequence in homology and
cohomology.  Moreover, all cyclic cohomology groups are computed by the same
complex.

However, the definition in~\eqref{eq:def_cyclic} is good when relating
operations in \emph{periodic} cyclic cohomology back to cyclic cohomology.
This is the situation in Theorem~\ref{the:excision_cyclic}.  We already have
excision in the periodic theory and want to estimate how the dimensions in
cyclic (co)homology are shifted.

Let $(i,p) \colon \MA[K] \injto \MA[E] \prto \MA[Q]$ be an allowable
extension.  Let $\Omega^n(\MA[E]: \MA[Q]) \defeq \Ker p_\ast \colon
\Omega^n\MA[E] \to \Omega^n\MA[Q]$.  We define the \emph{relative Hodge
  filtration} $F_n(\MA[E]: \MA[Q])$, the \emph{relative Hodge tower}
$\Xtower_n(\MA[E]: \MA[Q])$, and the relative cyclic (co)homology groups
$\HC_n(\MA[E]: \MA[Q])$, etc., as in \eqref{eq:def_Hodge},
\eqref{eq:def_Xtower}, and~\eqref{eq:def_cyclic}.  We only replace $\MA$ by
$\MA[E]: \MA[Q]$ everywhere.  For example,
\begin{equation}  \label{eq:def_Hodge_relative}
  F_n (\MA[E]: \MA[Q]) \defeq
  \sum_{j\ge n+1} \Omega^j(\MA[E]: \MA[Q]) \oplus
  b \bigl( \Omega^{n+1} (\MA[E]: \MA[Q]) \bigr).
\end{equation}
Everything said above for the absolute cyclic (co)homology groups carries over
to the relative case.  Furthermore, there is a long exact sequence
\begin{equation}  \label{eq:cyclic_relative_absolute}
  \cdots \to
  \HC_{n+1} (\MA[Q]) \to
  \HC_{n} (\MA[E]: \MA[Q]) \to
  \HC_{n} (\MA[E]) \to
  \HC_{n} (\MA[Q]) \to
  \HC_{n-1} (\MA[E]: \MA[Q]) \to
  \cdots
\end{equation}
and a dual sequence in cohomology (with all indexes raised and all arrows
reversed).  The exact sequence~\eqref{eq:cyclic_relative_absolute} follows
from the alternative definition of cyclic homology in terms of bicomplexes.

The Hodge filtration lives on $X(\Tanil\MA)$ but is defined in terms of
operators on $\Omega_\an\MA$.  It is important that we can describe it using
the algebra structure on $\Tanil\MA$ and powers of the ideal $\Janil\MA$ (see
also~\cite{cuntz95:cyclic}).  The $n$th power $(\Janil\MA)^n$ is the subspace
of $\Omega_\an^\even\MA$ spanned by monomials $\opt{\ma_0} d\ma_1 \dots
d\ma_{2j}$ with $j \ge n$.  We interpret $(\Janil\MA)^0 \defeq
\Unse{(\Tanil\MA)}$.

\begin{lemma}  \label{lem:Feven_alternative}
  Let $\hat{F}_{2n}^\even(\MA) \defeq (\Janil\MA)^{n+1} \subset
  \Omega_\an^\even$ and let $\hat{F}_{2n}^\odd(\MA) \subset
  \Omega^1(\Tanil\MA)$ be the completant linear hull of $(\Janil\MA)^n
  D(\Tanil\MA)$.  The odd part $F_{2n}^\odd(\MA)$ of $F_{2n}(\MA)$ is equal to
  the range of $\hat{F}_{2n}^\odd(\MA)$ under the quotient map
  $\Omega^1(\Tanil\MA) \to \Omega_\an\MA$.  The even part $F_{2n}^\odd(\MA)$
  is equal to $\hat{F}_{2n}^\even(\MA) + \partial_1 \bigl(F_{2n}^\odd (\MA)
  \bigr)$.
\end{lemma}

\begin{proof}
  Since $\partial_1 = b - (1+\kappa)d$, the operators $b$ and~$\partial_1$
  agree up to the perturbation $(1+\kappa)d$, which produces terms of higher
  order.  Hence it makes no difference whether we add $b(\Omega^{2n+1}\MA)$ or
  $\partial_1(\Omega^{2n+1}\MA)$ to $(\Janil\MA)^{n+1}$.
\end{proof}

A similar description of $F_{2n}(\MA[E]: \MA[Q])$ can be given in the relative
case.

\begin{lemma}  \label{lem:Fodd_alternative}
  Let $\hat{F}_{2n-1}^\odd (\MA) \subset \Omega^1 (\Tanil\MA)$ be the
  completant linear hull of $\sum_{k + l = n} (\Janil\MA)^k
  \,D(\Janil\MA)^l$.  Then $F_{2n-1}^\odd(\MA)$ is the image of
  $\hat{F}_{2n-1}^\odd(\MA)$ under the quotient map $\Omega^1(\Tanil\MA) \to
  \Omega_\an^\odd\MA$.  Furthermore, $F_{2n-1}^\even(\MA) = (\Janil\MA)^n$.
  In the relative case, we get an analogous statement with
  $$
  \hat{F}_{2n-1}^\odd(\MA[E]: \MA[Q]) \defeq
  \cllin {} \sum_{k + l = n}
  \bigl((\Janil\MA[E])^k \cap \II\bigr) \,D(\Janil\MA[E])^l +
  (\Janil\MA[E])^k \,D \bigr( \Janil\MA[E]^l  \cap \II \bigr).
  $$
\end{lemma}

\begin{proof}
  We can replace $\hat{F}_{2n-1}^\odd(\MA)$ by $(\Janil\MA)^{n-1} \,D(d\MA
  d\MA) + (\Janil\MA)^n \,D\MA$ without changing the image under the map to
  $\Omega_\an^\odd\MA$.  This follows from the derivation rule for~$D$ as
  in~\eqref{eq:D_long_form}.  We compute
  \begin{multline}  \label{eq:Fodd_alternative}
    \opt{x} \,D(d\ma_1 d\ma_2) \bmod [,] =
    \opt{x} \,D(\ma_1 \ma_2) -
    \opt{x} \odot \ma_1 \,D\ma_2 -
    \ma_2\odot \opt{x} \,D\ma_1 \bmod [,]
    \\ =
    d\opt{x} d\ma_1 d\ma_2 +
    d\ma_2 d\opt{x} d\ma_1 -
    b(\opt{x} d\ma_1 d\ma_2).
  \end{multline}
  Hence the map $\Omega_\an^\even\MA \,d\MA d\MA \cong \Unse{(\Tanil\MA)}
  \,D(d\MA d\MA) \to \Omega_\an^\odd\MA$ sending $\opt{x} d\ma_1d\ma_2$ to
  $\opt{x} \,D(d\ma_1 d\ma_2) + [,]$ is equal to~$-b$ up to perturbations of
  higher degree.  Thus it makes no difference whether we add
  $b(\Omega^{2n}\MA)$ or $\hat{F}_{2n-1}^\odd(\MA) \bmod [,]$ to
  $\Omega_\an^{\ge 2n}\MA$.  This proves the assertion about the odd part of
  $F_{2n-1}(\MA)$ in the absolute case.  The assertion about the even part is
  trivial.  The relative case is handled similarly.
\end{proof}

\section{The homotopy invariance of the X-complex}
\label{app:X_homotopy}

We prove the homotopy invariance of the X-complex:

{
\theoremstyle{plain}
\newtheorem*{proX}{Proposition~\ref{pro:X_homotopy}}

\begin{proX}
  Let $\MA$ and~$\MA[B]$ be complete bornological algebras.  Let $\phi_0,
  \phi_1 \colon \MA \to \MA[B]$ be AC-homotopic bounded homomorphisms.  Assume
  that~$\MA$ is quasi-free.

  Then the induced maps $X(\phi_t)$, $t = 0,1$, are chain homotopic.  That is,
  there is a bounded linear map $h \colon X(\MA) \to X(\MA[B])$ of degree~$1$
  such that $X(\phi_0) = X(\phi_1) - [\partial, h]$ with $[\partial, h] \defeq
  \partial_{X(\MA[B])} \circ h + h\circ \partial_{X(\MA)}$.
\end{proX}

}
We will derive the formula for~$h$ that is written down
in~\cite{cuntz95:cyclic}.  The ingredients are the quasi-freeness of the
source and an absolutely continuous homotopy $\Phi \colon \MA \to \ABC([0,1];
\MA[B])$.  The quasi-freeness of~$\MA$ is encoded in a graded
connection~$\nabla$ on $\Omega^1\MA$.  Since~$\MA$ is quasi-free, the natural
map $\Xmap_{n} \colon \Xtower_n(\MA) \to X(\MA)$ is a homotopy equivalence of
complexes.  Using the graded connection~$\nabla$, we can explicitly construct
a chain map section $\nabla_\ast \colon X(\MA) \to \Xtower_2(\MA)$.

The absolutely continuous homotopy $\Phi \colon \MA \to \ABC([0,1]; \MA[B])$
is used to construct a chain homotopy $\eta \colon \Xtower_2(\MA) \to
X(\MA[B])$ between the chain maps $X(\phi_0) \circ \Xmap_2$ and $X(\phi_1)
\circ \Xmap_2$.  Finally, $\eta \circ \nabla_\ast \colon X(\MA) \to X(\MA[B])$
is a chain homotopy between $X(\phi_t) \circ \Xmap_2 \circ \nabla_\ast =
X(\phi_t)$ for $t = 0,1$.

The proof of the homotopy invariance is quite explicit and therefore works in
many different categories.  However, the map~$\eta$ involves integration of
functions and hence division by integers.  As a result, we only get homotopy
invariance for the X-complex in a \Mpn{\Q}linear category.

\subsection{Contracting the \Hochschild complex}
\label{app:bar_contract}

If~$\MA$ is \Mpn{n}dimensional, then the bar resolution is contractible as a
bimodule resolution above degree~$n$.  Hence the \Hochschild complex is
contractible above degree~$n$.  Connections can be used to write down an
explicit contraction (see also~\cite{cuntz95:cyclic}).

Let~$\VS$ be a graded \Mpn{\MA}bimodule and let $\nabla \colon \VS \to \VS
\hot_{\MA} \Omega^1\MA$ be a graded right connection on~$\VS$ in the sense of
Definition~\ref{def:right_connection}.  The complete bornological vector space
$\VS \hot_{\MA} \Unse{(\Omega\MA)} \cong \sum_{n=0}^\infty \VS \hot \MA^{\hot
  n}$ is a graded \Mp{\MA}\Mp{\Omega\MA}bimodule in a natural way.  The graded
connection $\nabla \colon \VS \to \VS \hot_{\MA} \Omega^1\MA$ can be extended
uniquely to an endomorphism $\nabla \in \Endo(\VS \hot_{\MA}
\Unse{(\Omega\MA)})$ of degree~$1$ satisfying
\begin{equation}  \label{eq:nabla_extend}
  \nabla(\ma \cdot x) = \ma \cdot \nabla(x),
  \quad
  \nabla(x \cdot \omega) =
  \nabla(x) \cdot \omega + (-1)^{\deg x} x \cdot d\omega
\end{equation}
for all $\ma \in \MA$, $x \in \VS \hot_{\MA} \Unse{(\Omega\MA)}$, and $\omega
\in \Omega\MA$.  Equation~\eqref{eq:nabla_extend} implies that $\nabla(\vs
\hot d\ma_1 \dots d\ma_n) = \nabla(\vs) \cdot d\ma_1 \dots d\ma_n$.  This
uniquely determines a bounded linear map $\nabla \in \Endo(\VS \hot_{\MA}
\Unse{(\Omega\MA)})$.  This extension~$\nabla$
satisfies~\eqref{eq:nabla_extend} because of~\eqref{eq:def_right_connection}.

Especially, for a graded connection on $\Omega^n\MA$ the extension is a map
$\sum_{j \ge n} \Omega^j\MA \to \sum_{j \ge n+1} \Omega^j\MA$ of
degree~$1$ satisfying~\eqref{eq:nabla_extend}.  It is defined by
\begin{equation}  \label{eq:def_nabla_extend}
  \nabla( \opt{\ma_0} d\ma_1 \dots d\ma_j) \defeq
  \nabla( \opt{\ma_0} d\ma_1 \dots d\ma_n) d\ma_{n+1} \dots d\ma_j
  \qquad \forall \opt{\ma_0} \in \Unse{\MA},\ \ma_1, \dots, \ma_n \in \MA.
\end{equation}
Equation~\eqref{eq:nabla_extend} implies that $\nabla|_{\Omega^j\MA}$ is a
graded right connection for all $j \ge n$.

\begin{proposition}  \label{pro:Hochschild_contract}
  Let~$\MA$ be an \Mpn{n}dimensional complete bornological algebra.  Extend a
  graded connection $\nabla \colon \Omega^n\MA \to \Omega^{n+1}\MA$ as
  in~\eqref{eq:def_nabla_extend} and let $\nabla \equiv 0$ on $\Omega^j\MA$
  for $j < n$.
  
  Then $[b, \nabla] \colon \Omega\MA \to \Omega\MA$ is an idempotent whose
  range is equal to $F_n(\MA)$.  That is, $[b, \nabla] = \ID$ on $F_n(\MA)$
  and $\Ker [b, \nabla]$ is a \Mpn{b}invariant direct complement for
  $F_n(\MA)$.
\end{proposition}

\begin{proof}
  Let $x \in \Omega^j\MA$, $j \ge n+1$.  For $\omega = \ma \in \Omega^0\MA$,
  equation~\eqref{eq:nabla_extend} implies
  \begin{multline*}
    [b, \nabla] (xd\ma) =
    b \bigl(\nabla(x) d\ma\bigr) + \nabla \circ b(xd\ma) =
    (-1)^{j+1} [\nabla(x), \ma]  + (-1)^j \nabla([x,\ma])
    \\ =
    (-1)^j \bigl( \ma \cdot \nabla(x) - \nabla(x) \cdot \ma
      + \nabla(x \cdot \ma) - \nabla(\ma \cdot x) \bigr)
    \\ =
    (-)^j \bigl( \ma \cdot \nabla(x) - \nabla(x) \cdot \ma
     + \nabla(x) \cdot \ma + (-1)^j x d\ma - \ma \cdot \nabla(x) \bigr)
    =
    x d\ma.
  \end{multline*}
  Thus $[b, \nabla] = \ID$ on $\sum_{j \ge n+1} \Omega^j\MA$.  Since $[b,
  \nabla]$ commutes with~$b$, this also holds on $b(\Omega^{n+1}\MA) \subset
  \Omega^n\MA$.
  
  If $j < n$, then $[b, \nabla] = 0$ on $\Omega^j\MA$ because $\nabla = 0$ on
  $\Omega^j\MA$ and $\Omega^{j-1}\MA$.  We have $[b, \nabla] = b \circ \nabla$
  on $\Omega^n\MA$ because~$\nabla$ vanishes on $\Omega^{n-1}\MA$.  Hence $[b,
  \nabla] \circ [b,\nabla] = b \nabla b \nabla = b(\ID - b\nabla) \nabla =
  b\nabla = [b,\nabla]$ on $\Omega^n\MA$.  Thus $[b, \nabla]$ is idempotent.
  In addition, the range of $[b, \nabla]|_{\Omega^n\MA} = b
  \nabla|_{\Omega^n\MA}$ is contained in $F_n(\MA)$.  The range of
  $[b,\nabla]$ is equal to $F_n(\MA)$ because $[b,\nabla] = \ID$ on
  $F_n(\MA)$.  Thus $\Ker [b, \nabla]$ is a direct complement for $F_n(\MA)$.
  It is invariant under~$b$ because $[b, \nabla]$ commutes with~$b$.
\end{proof}

For cyclic cohomology, we must work with the boundary $B+b$ instead of~$b$.
We have $[\nabla, B+b] = [\nabla, B] + [\nabla, b] = [\nabla, B] + \pi_n$,
where $\pi_n = [\nabla, b]$ is an idempotent projection onto $F_n(\MA)$.  The
perturbation $[\nabla, B]$ maps $F_n(\MA)$ into itself and commutes there
with~$b$ because $[\nabla,b] = \ID$ on $F_n(\MA)$ commutes with~$B$ and $[B,b]
= 0$.  Evidently, $[\nabla, B]$ has degree~$+2$ and thus maps $\Omega^j\MA$ to
$\Omega^{j+2}\MA$ for all $j \in \N$.  Thus $\ID - [\nabla, B+b]$ maps
$F_j(\MA)$ into $F_{j+2}(\MA)$ for all $j \ge n$.  The map $[\nabla, B+b]$ is
a chain map homotopic to the identity by construction.  Iterating it $k$
times, we obtain chain maps $\Omega\MA \to \Omega\MA$ that map $F_n(\MA) \to
F_{n+2k}(\MA)$.  Hence we get induced chain maps $\Xtower_n(\MA) \to
\Xtower_{n+2k}(\MA)$ that are inverse up to homotopy to the natural projection
$\Xtower_{n+2k}(\MA) \to \Xtower_n(\MA)$.

This technique to contract the complex $(\Omega\MA, B+b)$ using a connection
is important for many applications.  Let me mention here only Khalkhali's
result that entire and periodic cyclic cohomology coincide for \emph{Banach}
algebras of finite homological dimension~\cite{khalkhali94:connections}.  This
is proved by observing that $\lim_{n \to \infty} (\ID - [\nabla, B+b])^n$
makes sense and is a bounded operator $\Omega_\an\MA \to \Omega_\an\MA$
homotopic to the identity.  However, the proof of boundedness depends on the
assumption that~$\MA$ is a Banach algebra.  The assertion becomes false
already for \Frechet algebras.

Here, we only use Proposition~\ref{pro:Hochschild_contract} to construct a
chain map $\nabla_\ast \colon X(\MA) \to \Xtower_2(\MA)$ if~$\MA$ is
quasi-free in the sense of Definition~\ref{deflem:quasi_free}.  Let $\nabla
\colon \Omega^1\MA \to \Omega^2\MA$ be a graded connection extended
by~\eqref{eq:def_nabla_extend}.  We have seen above that $\ID - [\nabla, B+b]$
maps $F_1(\MA)$ into $F_3(\MA) \subset F_2(\MA)$ and thus descends to a
bounded chain map $\nabla_\ast \colon X(\MA) = \Xtower_1(\MA) \to
\Xtower_2(\MA)$.  In addition, $\Xmap_2 \circ \nabla_\ast = \ID$.  Explicitly,
\begin{equation}  \label{eq:nabla_ast}
  \begin{alignedat}{2}{2}
    \nabla_\ast & = \ID - \nabla \circ d &
    \qquad & \text{on $\Omega^0\MA$,}
    \\
    \nabla_\ast & = \ID -  [\nabla, b] = \ID - b \circ \nabla &
    \qquad & \text{on $\Omega^1\MA / b(\Omega^2\MA)$.}
  \end{alignedat}{}
\end{equation}

\subsection{Absolutely continuous homotopies and chain homotopies}
\label{sec:X_tensor}

Let $\Phi \colon \MA \to \ABC([0,1]; \MA[B])$ be an absolutely continuous
homotopy.  Its derivative is a bounded linear map $\Phi' \colon \MA \to
L^1([0,1]; \MA[B])$.  If we view $L^1([0,1]; \MA[B])$ as an
\Mp{\ABC([0,1];\MA[B])}bimodule, then $\Phi'$ is a derivation relative
to~$\Phi$ in the sense that $\Phi'(xy) = \Phi(x)\Phi'(y) + \Phi'(x)\Phi(y)$
for all $x, y \in \MA$.  For $j \in \N$, define $\eta \colon \Omega^j\MA \to
\Omega^{j-1}\MA[B]$ by
$$
\eta( \opt{\ma_0} d\ma_1 \dots d\ma_j) \defeq
\int_0^1 \Phi_t\opt{\ma_0} \cdot \Phi'_t(\ma_1) \,d\Phi_t(\ma_2) \dots
\,d\Phi_t(\ma_n) \,dt.
$$
This is a well-defined bounded linear map because integration $\int_0^1
\,dt \colon L^1([0,1]; \Omega^{j-1}\MA[B]) \to \Omega^{j-1}\MA[B]$ is bounded.
We let $\eta|_{\Omega^0\MA} = 0$ and view~$\eta$ as an operator $\Omega\MA \to
\Omega\MA[B]$ of degree~$-1$.  Using that~$\Phi'$ is a derivation relative
to~$\Phi$, we compute that $[\eta,b] = \eta b + b \eta = 0$:
\begin{multline*}
  \eta \circ b(\opt{\ma_0} d\ma_1 \dots d\ma_n)
  =
  \int_0^1 \Phi_t(\opt{\ma_0} \ma_1) \Phi_t'(\ma_2) \,d\Phi_t(\ma_3) \dots
  \,d\Phi_t(\ma_n)
  \displaybreak[0] \\
  - \Phi_t\opt{\ma_0} \Phi_t'(\ma_1 \ma_2) \,d\Phi_t(\ma_3) \dots
  \,d\Phi_t(\ma_n)
  + \Phi_t\opt{\ma_0} \Phi_t'(\ma_1) \Phi_t(\ma_2) \,d\Phi_t(\ma_3) \dots
  \,d\Phi_t(\ma_n)
  \displaybreak[0] \\
  - (-1)^n \Phi_t\opt{\ma_0} \Phi_t'(\ma_1) \,\bigl( d\Phi_t(\ma_2) \dots
  \,d\Phi_t(\ma_{n-1}) \bigr) \cdot \Phi_t(\ma_n)
  \displaybreak[0] \\
  + (-1)^n \Phi_t(\ma_n\opt{\ma_0}) \Phi_t'(\ma_1) \,d\Phi_t(\ma_2) \dots
  \,d\Phi_t(\ma_{n-1})
  \,dt
  \displaybreak[1] \\ =
  \int_0^1 (-1)^{n-1}
  [\Phi_t\opt{\ma_0} \Phi_t'(\ma_1) \,d\Phi_t(\ma_2) \dots
  \,d\Phi_t(\ma_{n-1}), \Phi_t(\ma_n)] \,dt
  =
  -b\circ \eta(\opt{\ma_0} d\ma_1 \dots d\ma_n).
\end{multline*}

Thus~$\eta$ maps $b(\Omega^3\MA)$ into $b(\Omega^2\MA[B])$ and descends to a
well-defined bounded linear map $\eta \colon \Xtower_2(\MA) \to X(\MA[B])$.
Recall that $\Xmap_2 \colon \Xtower_2(\MA) \to X(\MA)$ is the natural
projection.  Using the map~$\eta$, we can finish the proof of
Proposition~\ref{pro:X_homotopy}:

\begin{lemma}  \label{lem:X_II_homotopy}
  We have $X(\phi_0) \circ \Xmap_2 = X(\phi_1) \circ \Xmap_2 - [\partial,
  \eta]$.  Thus the chain maps $X(\phi_t) \circ \Xmap_2 \colon \Xtower_2(\MA)
  \to X(\MA[B])$ induced by AC-homotopic homomorphisms $\phi_t \colon \MA \to
  \MA[B]$ are chain homotopic without any quasi-freeness assumptions on~$\MA$.
\end{lemma}

\begin{proof}
  We compute $[\partial, \eta] = \partial \eta + \eta \partial = X(\phi_1)
  \circ \Xmap_2 - X(\phi_0) \circ \Xmap_2$ on $\Omega^j\MA$ for $j = 0,1,2$:
  \begin{displaymath}
    [\partial, \eta](\ma) =
    \eta(d\ma) =
    \int_0^1 \Phi'_t(\ma) \,dt =
    \Phi_1(\ma) - \Phi_0(\ma).
  \end{displaymath}
  \begin{multline*}
    [\partial, \eta](\opt{\ma_0} d\ma_1) =
    d\circ \eta(\opt{\ma_0} d\ma_1) + \eta \circ B(\opt{\ma_0} d\ma_1)
    \\ =
    \int_0^1
      d\bigl( \Phi_t\opt{\ma_0} \Phi'_t(\ma_1) \bigr)
    + \Phi'_t\opt{\ma_0} d\Phi_t(\ma_1)
    - \Phi'_t(\ma_1) d\Phi_t\opt{\ma_0} \,dt
    \\ =
    \int_0^1
      [d\Phi_t\opt{\ma_0}, \Phi'_t(\ma_1)]
    + \frac{\partial}{\partial t}
      \bigl( \Phi_t\opt{\ma_0} d\Phi_t(\ma_1) \bigr) \,dt
    =
    \Phi_1\opt{\ma_0} d\Phi_1(\ma_1) - \Phi_0\opt{\ma_0} d\Phi_0(\ma_1).
  \end{multline*}
  Since the range of~$\eta$ is $X(\MA[B])$, we can omit $\int_0^1
  [d\Phi_t\opt{\ma_0}, \Phi'_t(\ma_1)] \,dt \in b(\Omega^2\MA[B])$.  On
  $\Omega^2\MA/ b(\Omega^3\MA)$ we have $\eta \partial + \partial \eta = \eta
  b + b \eta = 0$ as desired.
\end{proof}

The operator $h \colon X(\MA) \to X(\MA[B])$ in
Proposition~\ref{pro:X_homotopy} is the composition $\eta \circ \nabla_\ast$
with~$\eta$ as in Lemma~\ref{lem:X_II_homotopy} and~$\nabla_\ast$ as
in~\eqref{eq:nabla_ast}.  We can replace the even part of $\nabla_\ast$ by $-
\nabla \circ d$ because~$\eta$ vanishes on $\MA \subset \Xtower_2(\MA)$.  We
get the same formula as Cuntz and Quillen \cite[p.~408]{cuntz95:cyclic}.

\section{Some computations in $X(\Tanil\MA)$}
\label{app:X_Tanil}

We do some computations in $X(\Tanil\MA)$ that were left out in
Chapter~\ref{cha:HA}.  We derive the formula for the X-complex
boundary~$\partial$ on $\Omega_\an\MA$ and provide some details of the proof
of Proposition~\ref{pro:X_equals_C} that were left out in
Section~\ref{sec:X_Tanil}.

\subsection{The boundary of $X(\Tanil\MA)$}
\label{app:X_boundary}

Firstly, we write down the inverse of the map $\mu_{\ref{eq:mu3}} \colon
\Unse{(\Tanil\MA)} \hot \MA \hot \Unse{(\Tanil\MA)} \to \Omega^1(\Tanil\MA)$.
Since we already know from Proposition~\ref{pro:Omega_I_Tanil} that
$\mu_{\ref{eq:mu3}}$ is a bornological isomorphism, it suffices to write down
$\mu_{\ref{eq:mu3}}^{-1} \colon \Omega^1(\Tens\MA) \to \Unse{(\Tens\MA)} \hot
\MA \hot \Unse{(\Tens\MA)}$ such that $\mu_{\ref{eq:mu3}} \circ
\mu_{\ref{eq:mu3}}^{-1} = \ID$.  Then the map $\mu_{\ref{eq:mu3}}^{-1}$
extends to a bounded linear map $\mu_{\ref{eq:mu3}}^{-1} \colon
\Omega^1(\Tanil\MA) \to \Unse{(\Tanil\MA)} \hot \MA \hot \Unse{(\Tanil\MA)}$
inverse to $\mu_{\ref{eq:mu3}}$.

The derivation rule for~$D$ implies $D(dxdy) = D(xy) - x\,Dy - (Dx)\,y$ and
hence
\begin{multline}  \label{eq:D_long_form}
  D( \opt{\ma_0} d\ma_1 \dots d\ma_{2n}) =
  D\bigl( \opt{\ma_0} \odot (d\ma_1 d\ma_2) \odot (d\ma_3 d\ma_4)
  \odot \dots \odot (d\ma_{2n-1} d\ma_{2n}) \bigr)
  \displaybreak[0] \\ {} =
  (D\opt{\ma_0})\; d\ma_1 \dots d\ma_{2n} +
  \sum_{j=1}^n \opt{\ma_0} d\ma_1 \dots d\ma_{2j-2}
  \bigl( D(d\ma_{2j-1} d\ma_{2j}) \bigr) d\ma_{2j+1} \dots d\ma_{2n}
  \displaybreak[0] \\ {} =
  (D\opt{\ma_0})\; d\ma_1 \dots d\ma_{2n} +
  \sum_{j=1}^n \opt{\ma_0} d\ma_1 \dots d\ma_{2j-2}
  \bigl( D(\ma_{2j-1} \ma_{2j}) \bigr)\; d\ma_{2j+1} \dots d\ma_{2n}
  \\ \quad {} -
  \sum_{j=1}^n \opt{\ma_0} d\ma_1 \dots d\ma_{2j-2} \odot \ma_{2j-1}
  ( D\ma_{2j} )\; d\ma_{2j+1} \dots d\ma_{2n}
  \\ \quad {} -
  \sum_{j=1}^n \opt{\ma_0} d\ma_1\dots d\ma_{2j-2}
  ( D\ma_{2j-1} )\; \ma_{2j} d\ma_{2j+1} \dots d\ma_{2n}.
\end{multline}
This yields an explicit \Mp{\mu_{\ref{eq:mu3}}}pre-image of $D( \opt{\ma_0}
d\ma_1 \dots d\ma_n)$.  The resulting map $f \colon D(\Tens\MA) \to
\Unse{(\Tens\MA)} \hot \MA \hot \Unse{(\Tens\MA)}$ is extended to a left
\Mp{\Tens\MA}module homomorphism on $\Omega^1(\Tens\MA)$ by
$\mu_{\ref{eq:mu3}}^{-1} (\omega_1 \,D\omega_2) \defeq \omega_1 f(D\omega_2)$
for $\omega_1 \in \Unse{(\Tens\MA)}$, $\omega_2 \in \Tens\MA$.  By
construction, $\mu_{\ref{eq:mu3}} \circ \mu_{\ref{eq:mu3}}^{-1} = \ID$ on
$D(\Tens\MA)$ and hence on $\Omega^1(\Tens\MA)$.  It is not hard to show
directly that $\mu_{\ref{eq:mu3}}^{-1}$ extends to a bounded linear map
$\Omega^1(\Tanil\MA) \to \Unse{(\Tanil\MA)} \hot \MA \hot \Unse{(\Tanil\MA)}$
that is inverse to $\mu_{\ref{eq:mu3}}$.

Let $\partial \colon \Omega_\an\MA \to \Omega_\an\MA$ be the boundary of the
X-complex transported to $\Omega_\an\MA$ by the natural vector space
isomorphism $X(\Tanil\MA) \cong \Omega_\an\MA$.  We repeat the proof of the
explicit formulas for~$\partial$ in~\cite{cuntz95:cyclic}.  For $\omega \in
\Omega_\an^\even\MA$, $\ma \in \MA$, we have
$$
\partial_1(\omega d\ma) =
b( \omega \,D\ma) =
\omega \odot \ma - \ma \odot \omega =
\omega \cdot \ma - \ma \cdot \omega - d\omega d\ma + d\ma d\omega =
(b - d - \kappa \circ d)(\omega d\ma).
$$
That is, we obtain $\partial_1 = b - (1+\kappa) d$ as asserted
in~\eqref{eq:partial_odd}.

If $x \in \Omega_\an^\even\MA$, then $\partial_0(x) \in \Omega_\an^\odd\MA$ is
computed as follows.  We map $Dx \in \Omega^1(\Tanil\MA)$ to
$\Unse{(\Tanil\MA)} \hot \MA \hot \Unse{(\Tanil\MA)}$ via
$\mu_{\ref{eq:mu3}}^{-1}$ and compose with the bounded linear map
$\Unse{(\Tanil\MA)} \hot \MA \hot \Unse{(\Tanil\MA)} \to \Omega_\an^\odd\MA$
sending $\opt{x_0} \otimes \ma \otimes \opt{x_1}$ to $\opt{x_1} \odot
\opt{x_0} d\ma \in \Omega_\an^\odd\MA$.  Thus~\eqref{eq:D_long_form}
implies~\eqref{eq:partial_even}:
\begin{multline*}
  \partial_0( \opt{\ma_0} d\ma_1 \dots d\ma_{2n}) = 
  D( \opt{\ma_0} d\ma_1 \dots d\ma_{2n}) \bmod [,]
  =
  d\ma_1 \dots d\ma_{2n} D\ma_0
  \displaybreak[0] \\ {} +
  \sum_{j=1}^n d\ma_{2j+1} \dots d\ma_{2n} \odot
    \opt{\ma_0} d\ma_1 \dots d\ma_{2j-2} D(\ma_{2j-1} \ma_{2j})
  \displaybreak[0] \\ {} -
  \sum_{j=1}^n d\ma_{2j+1} \dots d\ma_{2n} \odot
    \opt{\ma_0} d\ma_1 \dots d\ma_{2j-2} \odot \ma_{2j-1} D\ma_{2j}
  \displaybreak[0] \\ {} -
  \sum_{j=1}^n \ma_{2j} \odot d\ma_{2j+1} \dots d\ma_{2n} \odot
    \opt{\ma_0} d\ma_1 \dots d\ma_{2j-2} D\ma_{2j-1}  \bmod [,]
  \displaybreak[1] \\ =
  \sum_{j=0}^{2n}
    d\ma_{j} \dots d\ma_{2n} d\opt{\ma_{0}} d\ma_1 \dots d\ma_{j-1}
  - \sum_{j=1}^n b( d\ma_{2j+1} \dots d\ma_{2n} \cdot
    \opt{\ma_0} d\ma_1 \dots d\ma_{2j-1} d\ma_{2j})
  \displaybreak[1] \\ =
  \sum_{j=0}^{2n} \kappa^j d(\opt{\ma_0} d\ma_1 \dots d\ma_{2n}) -
  \sum_{j=0}^{n-1} b\kappa^{2j} (\opt{\ma_0} d\ma_1 \dots d\ma_{2n}).
\end{multline*}

\subsection{Boundedness of the spectral operators}
\label{app:X_boundedness}

We verify some details for the proof of Proposition~\ref{pro:X_equals_C}.  We
have to show that the complexes $(\Omega_\an\MA, \partial)$ and
$(\Omega_\an\MA, \delta)$ with~$\delta$ as defined in~\eqref{eq:def_delta} are
chain homotopic.  The proof uses the spectral projection~$P$ of~$\kappa^2$
associated to the eigenvalue~$1$ and the operator~$H$ that annihilates the
range of~$P$ and is inverse to $1 - \kappa^2$ on the range of $P^\bot \defeq
\ID - P$.  We need explicit formulas for $P$ and~$H$ to verify that they are
bounded.

By~\eqref{eq:kappa_vii}, $\kappa$ satisfies the identities $(\kappa^n -
1)(\kappa^{n+1} - 1) = 0$ on $\Omega^n\MA$.  This carries over to~$\kappa^2$
because $(1-\kappa^{2n}) = (1-\kappa^n) (1+\kappa^n)$.  The \Mpn{1}eigenspace
of~$\kappa^2$ is the sum of the eigenspaces of~$\kappa$ for the
eigenvalues~$\pm 1$.  The following linear algebra applies to any
operator~$\lambda$ on a complex vector space~$V$ that satisfies the identity
$(\lambda^n-1) (\lambda^{n+1} - 1) = 0$.

Since~$\lambda$ satisfies a polynomial identity and we are working over the
complex numbers, it follows that~$V$ decomposes into eigenspaces $V =
\bigoplus_{\omega \in N} V_\omega$ with~$N$ the set of roots of the polynomial
$(z^n-1) (z^{n+1}-1)$.  Thus $\omega^n = 1$ or $\omega^{n+1} = 1$.  In
particular, $1 \in N$ is an eigenvalue with multiplicity at most~$2$.  All
other eigenvalues $\omega \in N \setminus \{1\}$ have no multiplicity, so that
$\lambda|_{V_\omega}$ is scalar multiplication by~$\omega$.  On~$V_1$, we have
$\lambda|_{V_1} = 1 + \epsilon$ with $\epsilon \circ \epsilon = 0$.  If the
eigenvalue~$1$ has multiplicity~$1$, then $\epsilon = 0$.  Linear algebra
yields a polynomial~$f_n$ such that $f_n(\lambda)$ is the
projection~$P_\lambda$ onto~$V_1$ that annihilates all other
eigenspaces~$V_\omega$ with $\omega \neq 1$.  Let $P_\lambda^\bot = \ID -
P_\lambda$ be the complementary projection.  There is a polynomial~$g_n$ such
that $P_\lambda \circ g_n(\lambda) = 0$ and $g_n(\lambda) (1 - \lambda) =
P_\lambda^\bot$.  That is, $H_\lambda \defeq g_n(\lambda)$ lives on the range
of $P_\lambda^\bot$ and is the inverse of $1-\lambda$ there.  The polynomials
$f_n$ and~$g_n$ are not uniquely determined.  However, the operators
$P_\lambda$ and~$H_\lambda$ are uniquely determined by~$\lambda$.

We need explicit formulas for~$f_n$ and~$g_n$.  Define $N_n \defeq \frac{1}{n}
(1 + z + \dots + z^{n-1})$.  Thus $N_n(1) = 1$ and $(z-1) \cdot N_n =
\frac{1}{n} (z^n-1)$.  It follows that $N_n(\omega) = 0$ for all $\omega \neq
1$ with $\omega^n = 1$.  Let
$$
f_n(z) \defeq
N_n(z) \cdot N_{n+1}(z)
\bigl( 1 - \bigl (n- \tfrac{1}{2} \bigr) (z-1) \bigr).
$$
By construction, $f_n(\omega) = 0$ for all $\omega \in N \setminus \{1\}$.
One verifies easily that $f_n(1) = 1$ and $f_n'(1) = 0$.  Thus $f_n(\lambda)$
annihilates~$V_\omega$ for $\omega \neq 1$.  On~$V_1$ it is the identity
operator because $\lambda|_{V_1} = 1 + \epsilon$ with $\epsilon^2 = 0$ and
hence $f_n(1 + \epsilon) = f_n(1) + f'_n(1) \epsilon = 1$.  Thus $f_n(\lambda)
P_\lambda$.  Since $f_n(1) = 1$, the function $\tilde{g}_n \defeq (f_n-1) /
(z-1)$ is a polynomial.  By definition, $\tilde{g}_n (\lambda) \circ (1 -
\lambda) = 1 - f_n(\lambda) = 1 - P_\lambda$.  Evidently,
$$
\tilde{g}_n(z) =
-\bigl( n- \tfrac{1}{2} \bigr) N_n(z) \cdot N_{n+1}(z) +
N_n(z) \cdot \frac{N_{n+1}(z) - 1}{z - 1} + \frac{N_n(z) - 1}{z - 1}.
$$
An easy computation shows that
$$
\frac{N_{n}(z) - 1}{z - 1} =
\tfrac{1}{n} (z^{n-2} + 2z^{n-1} + 3z^{n-2} + \dots + (n-1)z^0) =
\sum_{j=0}^{n-2} \frac{n - 1 - j}{n} z^j.
$$
The relation $g_n(\lambda) \cdot P_\lambda = 0$ can be enforced by defining
$g_n \defeq (1-f_n) \cdot \tilde{g}_n$.  This does not alter the relation
$g_n(\lambda) \circ (1 - \lambda) = 1 - P_\lambda$ and achieves that
$g_n(\lambda) \circ P_\lambda = 0$ because~$P_\lambda$ is idempotent.

We now apply this linear algebra to~$\kappa^2$.  We define $P \colon \Omega\MA
\to \Omega\MA$ and $H \colon \Omega\MA \to \Omega\MA$ such that
$P|_{\Omega^n\MA}$ is the projection onto the \Mpn{1}eigenspace of~$\kappa^2$
and~$H$ satisfies $PH = HP = 0$ and $H\cdot (1- \kappa^2) = (1-\kappa^2) \cdot
H = P^\bot$.  We claim that these operators $P$ and~$H$ are bounded with
respect to the bornology~$\CBS_\an$ on $\Omega\MA$.

For each $n \in \Z_+$, we have constructed above polynomials $f_n$ and~$g_n$
such that $P|_{\Omega^n\MA} = f_n(\kappa^2)$ and $H|_{\Omega^n\MA} =
g_n(\kappa^2)$.  Let $Q_n \colon \Omega\MA \to \Omega^n\MA \subset \Omega\MA$
be the projection annihilating $\Omega^k\MA$ with $k \neq n$.  The family of
operators $\{ C^n \kappa^j Q_n \mid 0 \le j \le n,\ n \in \Z_+\}$ is
equibounded with respect to the bornology~$\CBS_\an$ for any $C \in \R$.  This
follows easily from~\eqref{eq:kappa_ii}.  Furthermore, \eqref{eq:kappa_ii}
implies that $\{ \kappa^n Q_n \mid n \in \Z_+\}$ is equibounded.  Hence the
set of operators $\{C^n \kappa^j Q_n \mid 0\le j \le kn,\ n \in \Z_+\}$ is
equibounded for all $k \in \N$.  By the way, the set of operators $\{ \kappa^j
Q_n \mid n, j \in \Z_+\}$ is \emph{not} equibounded
because~\eqref{eq:kappa_viii} and $bB bB = 0$ imply $\kappa^{kn(n+1)} = 1 +
kbB$.

Hence an operator of the form $\sum_{j=0}^\infty Q_n h_n(\kappa)$ is bounded
with respect to the bornology~$\CBS_\an$ if $(h_n)_{n \in \Z_+}$ is a sequence
of polynomials whose degrees increase at most linearly and whose
\emph{absolute coefficient sum} grows at most exponentially.  The absolute
coefficient sum of a polynomial $a_0 + a_1x + \dots + a_nx^n$ is simply $|a_0|
+ \dots + |a_n|$.  Under these assumptions, the sum $\sum Q_n h_n(\kappa)$ is
in the completant disked hull of the equibounded set $\{ C^n \kappa^j Q_n \mid
0 \le j \le kn \}$ for suitable $C \in \R$, $k \in \N$.  The degrees of the
polynomials $f_n$ and~$g_n$ grow linearly in~$n$ and their absolute
coefficient sums grows polynomially.  Hence $P$ and~$H$ are bounded with
respect to the bornology $\CBS_\an$ and thus can be extended to bounded
operators on $\Omega_\an\MA$.
  
The X-complex boundary $\partial$ and $b+B$ commute with~$\kappa$ because~$d$
and~$b$ commute with~$\kappa$ and~$\kappa$ preserves the degree.  The range
of~$P$ is equal to the kernel of $(\kappa^2-1)^2$, whereas the kernel of~$P$
is equal to the range of $(\kappa^2-1)^2$.  These spaces are invariant under
any operator commuting with~$\kappa$.  Thus~$P$ commutes with $\partial$ and
$b+B$.  The operator~$H$ can be characterized as the inverse of $P^\bot
(1-\kappa^2) P^\bot$ in $\Endo(P^\bot \Omega\MA)$.  Thus any operator
commuting with~$\kappa$ commutes with~$H$.

It remains to verify that the restriction of $\partial$ to the range of~$P$ is
equal to~$\delta$.  Equation~\eqref{eq:kappa_iv} implies that~$1$ is a single
eigenvalue of $\kappa|_{d\Omega_\an\MA}$.  Thus $P \circ (\kappa^2 - \ID)$
vanishes on $d(\Omega_\an\MA)$.  Using~\eqref{eq:kappa_vi} we get similarly
that $P \circ (\kappa^2 - \ID)$ vanishes on $b(\Omega_\an\MA)$.  Consequently,
we have $\kappa^2 d = d$ and $\kappa^2 b = b$ on the subspace $P
\Omega_\an\MA$.  Hence the restriction of~$\partial_0$ to $P \Omega^{2n}\MA$
is equal to
$$
\sum_{j=0}^{2n} \kappa^j d - \sum_{j=0}^{n-1} \kappa^{2j} b =
B - n b =
\delta.
$$
The restriction of~$\partial_1$ to $P \Omega^{2n+1}\MA$ is equal to
$$
b - (1 + \kappa) d =
b - \frac{1}{n+1} \sum_{j=0}^n (1 + \kappa) \kappa^{2j} d =
b - \frac{1}{n+1} B =
\delta.
$$
These are the details needed to fill the gaps in the proof of
Proposition~\ref{pro:X_equals_C} in Section~\ref{sec:X_Tanil}.

\end{appendix}

\bibliographystyle{plain}
\bibliography{cyclic}

\end{document}